\documentclass[letterpaper,12pt]{article}
\usepackage{amsmath,amssymb,amsfonts,epsf,epsfig,graphicx,enumerate,placeins,bm}
\usepackage{relsize}
\newtheorem{theo}{Theorem}[section]
\newtheorem{prop}[theo]{Proposition}

\newtheorem{defn}[theo]{Definition}

\newtheorem{conj}[theo]{Conjecture}

\newenvironment{proof}{\noindent {\sc Proof}.}
                {\phantom{a} \hfill \framebox[2.2mm]{ } \bigskip}

\newenvironment{itemizenew}
{\begin{list}{$\bullet$}{
\leftmargin=5mm
\labelwidth=3mm
\labelsep=3mm
}}{\end{list}}

\newcommand{\ZZ}{\mathbb{Z}}

\def\int{{\rm int}}

\newcommand{\D}{{\mathcal{D}}}
\renewcommand{\phi}{\varphi}
\newcommand{\n}{{\bullet}}
\renewcommand{\t}{{\circ}}

\newcommand{\Circ}{{\rm Circ}}

\title{The Honeymoon Oberwolfach Problem: small cases}
\author{{\large \sc Marie Rose Jerade}  and  {\large \sc Mateja \v{S}ajna}\footnote{Corresponding author. Email: msajna$@$uottawa.ca. Mailing address: Department of Mathematics and Statistics, University of Ottawa, Ottawa, ON, Canada.}  \medskip \\
University of Ottawa }

\begin{document}
\maketitle \baselineskip 17pt

\begin{abstract}
The Honeymoon Oberwolfach Problem HOP$(2m_1,2m_2,\ldots,2m_t)$ asks the following question. Given $n=m_1+m_2+\ldots +m_t$ newlywed couples at a conference and $t$ round tables of sizes $2m_1,2m_2,\ldots,2m_t$, is it possible to arrange the $2n$ participants at these tables for $2n-2$ meals so that each participant sits next to their spouse at every meal, and sits next to every other participant exactly once? A solution to HOP$(2m_1,2m_2,\ldots,2m_t)$ is a decomposition of $K_{2n}+(2n-3)I$, the complete graph $K_{2n}$ with $2n-3$ additional copies of a fixed 1-factor $I$, into 2-factors, each consisting of disjoint $I$-alternating cycles of lengths $2m_1,2m_2,\ldots,2m_t$.

The Honeymoon Oberwolfach Problem was introduced in a 2019 paper by Lepine and \v{S}ajna. The authors conjectured that HOP$(2m_1,2m_2,\ldots,$ $2m_t)$ has a solution whenever the obvious necessary conditions are satisfied, and proved the conjecture for several large cases, including the uniform cycle length case $m_1=\ldots=m_t$, and the small cases with $n \le 9$. In the present paper, we extend the latter result to all cases with $n \le 20$ using a computer search.

\end{abstract}

\section{Introduction}

The well-known Oberwolfach Problem, denoted OP$(m_1,\ldots,m_t)$, asks whether $n=m_1+\ldots+m_t$  participants can be seated at $t$ tables of sizes $m_1,\ldots,m_t$ for several nights in a row so that each participant gets to sit next to every other participant exactly once. Thus, we are asking whether $K_n$, the complete graph on $n$ vertices, admits a 2-factorization such that each 2-factor is a disjoint union of $t$ cycles of lengths $m_1,\ldots,m_t$. The problem has been solved in many special cases --- see \cite{AdaBry,AlsHag,AlsSch,BryDan,Dez,FraHol,FraRos,GloJoo,HofSch,SalDra,Tra}
--- but is in general still open.

In the 2019 paper \cite{LepSaj}, Lepine and the second author introduced a new variant of the Oberwolfach Problem, called the Honeymoon Oberwolfach Problem. This problem, denoted HOP$(m_1,\ldots,m_t)$, can be described as follows. We have $n=\frac{1}{2}(m_1+\ldots+m_t)$ newlywed couples attending a conference and $t$ tables of sizes $m_1,\ldots,m_t$ (where each $m_i  \ge 3$). Is it possible to arrange the participants  at these $t$ round tables on $2n-2$ consecutive nights so that each couple sit together every night, and every participant sits next to every other participant exactly once?

In graph-theoretic terms we are asking whether $K_{2n}+(2n-3)I$, the multigraph obtained from the complete graph $K_{2n}$ by adjoining $2n-3$ additional copies of a chosen 1-factor $I$, admits a decomposition into 2-factors, each a vertex-disjoint union of cycles of lengths $m_1,\ldots,m_t$, so that in each of these cycles, every other edge is a copy of an edge of $I$.
A solution to HOP$(m_1,m_2,\ldots,m_t)$ is equivalent to a {\em semi-uniform 1-factorization of $K_{2n}$ of type $(m_1,m_2,\ldots,m_t)$}; that is, a 1-factorization $\{ F_1,F_2,\ldots,F_{2n-1} \}$ such that for all $i \ne 1$, the 2-factor $F_1 \cup F_i$ consists of disjoint cycles of lengths $m_1,m_2,\ldots,m_t$.

If $m_1=\ldots=m_t=m$ and $tm=2n$, then the symbol HOP$(m_1,\ldots,m_t)$ is abbreviated as HOP$(2n;m)$. Note that if HOP$(m_1,\ldots,m_t)$ has a solution, then the $m_i$ are all even and at least 4; these are the obvious necessary conditions.

In \cite{LepSaj}, the authors proposed the following conjecture.

\begin{conj}\cite{LepSaj}\label{conj:main}
The obvious necessary conditions for HOP$(m_1,\ldots,m_t)$ to have a solution are also sufficient.
\end{conj}

They also proved the conjecture in the following cases.

\begin{theo}\cite{LepSaj} \label{the:intro2}
Let $m$ and $n$ be positive integers, $2 \le m \le n$. Then HOP$(2n;2m)$ has a solution if and only if $n \equiv 0 \pmod{m}$.
\end{theo}

\begin{theo}\cite{LepSaj} \label{the:intro1}
Let $2 \le m_1 \le \ldots \le m_t$ be integers, and $n=m_1+\ldots +m_t$. Then HOP$(2m_1,\ldots,2m_t)$ has a solution in each of the following cases.
\begin{enumerate}[(i)]
\item $m_i \equiv 0 \pmod{4}$ for all $i$.
\item $n$ is odd and OP$(m_1,\ldots,m_t)$ has a solution.
\item $n$ is odd and $t=2$.
\item $n$ is odd, $n < 40$, and $m_1 \ge 3$.
\item $n \le 9$.
\end{enumerate}
\end{theo}

\noindent We remark that case (iv) of Theorem~\ref{the:intro1} is extended to $n \le 60$ by the result of \cite{SalDra}.
In this paper, we extend case (v) of Theorem~\ref{the:intro1}, thus proving the following result.

\begin{theo}\label{the:main}
Let $m_1, \ldots, m_t$ be integers with $m_i \ge 2$ for all $i$, and $n=m_1+\ldots +m_t$ such that $n \le 20$. Then HOP$(2m_1,\ldots,2m_t)$ has a solution.
\end{theo}

To prove Theorem~\ref{the:main}, we use the approach described in \cite{LepSaj}, combined with a computer-assisted search. In Sections~\ref{sec:prereqs} and \ref{sec:tools}, we present the relevant terminology and tools from \cite{LepSaj}, in Section~\ref{sec:main}, we give the framework of the proof of Theorem~\ref{the:main}, while in the appendices, we list the computational results supporting the proof.

\section{Terminology}\label{sec:prereqs}

In this paper, graphs may contain parallel edges, but not loops.
As usual, $K_n$ and $\lambda K_n$ denote the complete graph and the $\lambda$-fold complete graph, respectively, on $n$ vertices.   For $m \ge 2$, the symbol $C_m$ denotes the cycle of length $m$, or $m$-cycle.

Let $G$ be a graph, and let $H_1,\ldots,H_t$ be subgraphs of $G$. The collection $\{ H_1,\ldots,H_t\}$ is called a {\em decomposition} of $G$ if $\{ E(H_1),\ldots,E(H_t)\}$ is a partition of $E(G)$.

An {\em $r$-factor} in a graph $G$ is an $r$-regular spanning subgraph of $G$, and an {\em $r$-factorization} of $G$ is a decomposition of $G$ into $r$-factors. A 2-factor of $G$ consisting of disjoint  cycles of lengths $m_1,\ldots, m_t$, respectively, is called a {\em $(C_{m_1},\ldots,C_{m_t})$-factor} of $G$, and a decomposition into $(C_{m_1},\ldots,C_{m_t})$-factors  is called a {\em $(C_{m_1},\ldots,C_{m_t})$-factorization} of $G$.

For a positive integer $n$ and $S \subseteq \ZZ_n^\ast$ such that $S=-S$, we define a {\em circulant} $\Circ(n;S)$ as the graph with vertex set $\{ x_i: i \in \ZZ_n\}$ and edge set $\{ x_i x_{i+d}: i \in \ZZ_n, d \in S\}$. An edge of the form $x_i x_{i+d}$ is said to be of {\em difference} $d$. Note that an edge of difference $d$ is also of difference $n-d$, so we may assume that each difference is in $\{ 1,2,\ldots, \lfloor \frac{n}{2} \rfloor \}$.

In this paper, the complete graph $K_n$ will be viewed as the join of the circulant $\Circ(n-1;\ZZ_{n-1}^\ast)$ and the complete graph $K_1$ with vertex $x_{\infty}$. Thus, $V(K_n)=\{ x_i: i \in \ZZ_{n-1}\} \cup \{ x_{\infty} \}$ and $E(K_n)=\{ x_i x_j: i,j \in \ZZ_{n-1}, i \ne j \} \cup \{ x_i x_{\infty}: i \in \ZZ_{n-1} \}$. If this is the case, then an edge of the form $x_i x_{\infty}$ will be called of {\em difference infinity}.

Let $I$ be a chosen 1-factor in the graph $K_{2n}$. An edge of $K_{2n}$ is said to be an {\em $I$-edge} if it belongs to $E(I)$, and a {\em non-$I$-edge} otherwise. The symbol  $K_{2n}+\lambda I$ denotes the graph $K_{2n}$ with $\lambda$ additional copies of each $I$-edge, for a total of $\lambda+1$ copies of each $I$-edge.  (Note that these additional copies of $I$-edges of $K_{2n}$ are then also considered to be $I$-edges of $K_{2n}+\lambda I$.)
A cycle $C$ of $K_{2n}+\lambda I$, necessarily of even length, is said to be {\em $I$-alternating} if the $I$-edges and non-$I$-edges along $C$ alternate. A 2-factor (or 2-factorization) of $K_{2n}+\lambda I$ is said to be {\em $I$-alternating} if each of its cycles is $I$-alternating.

Thus, a solution to the Honeymoon Oberwolfach Problem HOP$(m_1,\ldots,m_t)$ is an $I$-alternating $(C_{m_1},\ldots,C_{m_t})$-factorization of $K_{2n}+(2n-3)I$ for $2n=m_1+m_2+\ldots +m_t$.

\section{The Tools}\label{sec:tools}

As in \cite{LepSaj}, we use the symbol  $4K_n^\n$ to denote the 4-fold complete graph with $n$ vertices whose edges are coloured pink, blue, and black, and black edges are oriented so that each 4-set of parallel edges contains one pink edge, one blue edge, and two opposite black arcs.

\begin{defn} \cite{LepSaj}\label{def:HOP} {\rm
A 2-factorization ${\cal D}$  of  $4K_n^\n$  is said to be {\em HOP} if each cycle of ${\cal D}$ satisfies the following condition:
\begin{description}
\item{\bf (C)} any two adjacent (that is, consecutive) edges satisfy one of the following:
    \begin{itemize}
    \item one is blue and the other pink; or
    \item both are black and directed in the same way;
    \item one is blue and the other black, directed towards the blue edge; or
    \item one is pink and the other black, directed away from the pink edge.
    \end{itemize}
\end{description}}
\end{defn}

\begin{theo}\cite{LepSaj} \label{the:NASC}
Let $m_1,\ldots,m_t$ be integers greater than 1, and let $n=m_1+\ldots+m_t$. Then HOP$(2m_1,\ldots,2m_t)$ has a solution if and only if $4K_n^\n$ admits an HOP $(C_{m_1},\ldots,C_{m_t})$-factorization.
\end{theo}

In the next proposition, the symbol $2K_n^\t$ denotes the multigraph $2K_n$ whose edges are coloured pink and black so that each 2-set of parallel edges contains one pink edge and one black edge.

\begin{prop}\cite{LepSaj} \label{pro:main1}
Assume $n$ is even, and let the vertex set of $2K_n^\t$ be $\{ x_i: i \in \ZZ_{n-1}\} \cup \{ x_{\infty} \}$. Let $\rho$ be the permutation $\rho=(x_{\infty})(x_0 \; x_1 \; x_2 \; \ldots \; x_{n-2})$, and let $\rho_\t$ denote the permutation on the edge set of $2K_n^\t$ that is induced by $\rho$ and that preserves the colour of the edges.

Suppose $2K_n^\t$ admits a $(C_{m_1},\ldots,C_{m_t})$-factor $F$  such that
\begin{description}
\item[{\bf (A1)}] each cycle in $F$ of length at least 3 contains an even number of pink edges, and
\item[{\bf (A2)}] $F$  contains exactly one edge from each of the orbits of $\langle \rho_\t \rangle$.
\end{description}
 Then $4K_n^{\n}$ admits  an HOP $(C_{m_1},\ldots,C_{m_t})$-factorization.
\end{prop}

\begin{prop}\cite{LepSaj} \label{pro:main}
Assume $n$ is even, and let the vertex set of $4K_n^{\n}$ be $\{ x_i: i \in \ZZ_{n-1}\} \cup \{ x_{\infty} \}$. Let $\rho$ be the permutation $\rho=(x_{\infty})(x_0 \; x_1 \; x_2 \; \ldots \; x_{n-2})$, and let $\rho_{\n}$ denote the permutation on the edge set of $4K_n^{\n}$ that is induced by $\rho$ and that preserves the colour (and orientation) of the edges.

Suppose $4K_n^{\n}$ admits edge-disjoint $(C_{m_1},\ldots,C_{m_t})$-factors $F_1$ and $F_2$ such that
\begin{description}
\item[{\bf (D1)}] each cycle in $F_1$ and $F_2$ satisfies Condition (C) in Definition~\ref{def:HOP}, and
\item[{\bf (D2)}] $F_1$ and $F_2$ jointly contain exactly one edge from each of the orbits of $\langle \rho_{\n} \rangle$.
\end{description}
 Then  $\D=\{ \rho_{\n}^i(F_1), \rho_{\n}^i(F_2): i \in \ZZ_{n-1} \}$ is an  HOP $(C_{m_1},\ldots,C_{m_t})$-factorization of $4K_n^{\n}$.
\end{prop}

\begin{prop}\cite{LepSaj} \label{pro:main2}
Assume $n$ is odd, and let the vertex set of $4K_n^{\n}$ be $\{ x_i: i \in \ZZ_{n-1} \} \cup \{ x_{\infty} \}$. Let $\rho$ be the permutation $\rho=(x_{\infty})(x_0 \; x_1 \; x_2 \; \ldots \; x_{n-2})$, and let $\rho_{\n}$ denote the permutation on the edge set of $4K_n^{\n}$ that is induced by $\rho$ and that preserves the colour (and orientation) of the edges.

Suppose $4K_n^{\n}$ admits pairwise edge-disjoint $(C_{m_1},\ldots,C_{m_t})$-factors $F_1$, $F_2$,  and $F_3$ such that
\begin{description}
\item[{\bf (E1)}] each cycle in $F_1$, $F_2$, and $F_3$ satisfies Condition (C) in Definition~\ref{def:HOP};
\item[{\bf (E2)}] each orbit of $\langle \rho_{\n} \rangle$ has edges  either in $F_1 \cup F_2$ or in $F_3$;
\item[{\bf (E3)}] if $e \in E(F_1 \cup F_2)$, then $\rho_{\n}^{\frac{n-1}{2}}(e)\in E(F_1 \cup F_2)$; and
\item[{\bf (E4)}] $F_1 \cup F_2$ contains a pink and a blue edge of difference $\frac{n-1}{2}$.
\end{description}
 Then  $\D=\{ \rho_{\n}^i(F_1), \rho_{\n}^i(F_2): i=0,1,\ldots,\frac{n-3}{2} \} \cup \{ \rho_{\n}^i(F_3): i \in \ZZ_{n-1} \}$ is an  HOP $(C_{m_1},\ldots,C_{m_t})$-factorization of $4K_n^{\n}$.
\end{prop}

The required 2-factors $F$ from Proposition~\ref{pro:main1}, $F_1$ and $F_2$ from Proposition~\ref{pro:main}, and $F_1$, $F_2$, and $F_3$ from Proposition~\ref{pro:main2} will are called the {\em starter 2-factors} (or {\em starters}) of the resulting 2-factorizations. Thus, we are referring to Propositions~\ref{pro:main1}, \ref{pro:main}, and \ref{pro:main2} as the one-starter, two-starter, and three-starter approach, respectively.

\section{Proof of Theorem~\ref{the:main}}\label{sec:main}

By the {\em type} of a $(C_{m_1},\ldots,C_{m_t})$-factor we mean the multiset $[m_1,\ldots,m_t]$.

\medskip

\begin{proof}
By Theorem~\ref{the:intro1}(v), it suffices to consider $10 \le n \le 20$. For each value of $n$, we list all possible 2-factor types $[m_1,\ldots,m_t]$, and refer to the result that guarantees existence of a solution to HOP$(2m_1,\ldots,2m_t)$. In cases where we are referring to Proposition~\ref{pro:main1}, \ref{pro:main}, or \ref{pro:main2}, appropriate starter 2-factors are given in the appendices.

\bigskip

\noindent {\sc Case} $n=10$: see Appendix~\ref{app:10}.

\bigskip

\begin{tabular}{|l|l|}
\hline
2-factor type $[m_1,\ldots,m_t]$ & HOP$(2m_1,\ldots,2m_t)$ has a solution by... \\ \hline
$[2, 2, 2, 2, 2]$ & Theorem~\ref{the:intro2} \\
$[4, 2, 2, 2]$ & Proposition~\ref{pro:main1} (one starter) \\
$[3, 3, 2, 2]$ & Proposition~\ref{pro:main} (two starters) \\
$[6, 2, 2]$ & Proposition~\ref{pro:main} (two starters) \\
$[5, 3, 2]$ & Proposition~\ref{pro:main1} (one starter)\\
$[4, 4, 2]$ & Proposition~\ref{pro:main1} (one starter)\\
$[4, 3, 3]$ & Proposition~\ref{pro:main} (two starters) \\
$[8,2]$ & Proposition~\ref{pro:main1} (one starter)\\
$[7,3]$ & Proposition~\ref{pro:main} (two starters) \\
$[6,4]$ & Proposition~\ref{pro:main} (two starters) \\
$[5,5]$ & Theorem~\ref{the:intro2} \\
$[10]$ & Theorem~\ref{the:intro2} \\
\hline
\end{tabular}

\newpage

\noindent {\sc Case} $n=11$: see Appendix~\ref{app:11}.

\bigskip

\begin{tabular}{|l|l|}
\hline
2-factor type $[m_1,\ldots,m_t]$ & HOP$(2m_1,\ldots,2m_t)$ has a solution by... \\ \hline
$[3, 2, 2, 2, 2]$ & Proposition~\ref{pro:main2} (three starters) \\
$[5, 2, 2, 2]$ & Proposition~\ref{pro:main2} (three starters) \\
$[4, 3, 2, 2]$ & Proposition~\ref{pro:main2} (three starters) \\
$[3, 3, 3, 2]$ & Proposition~\ref{pro:main2} (three starters) \\
$[7, 2, 2]$ & Proposition~\ref{pro:main2} (three starters) \\
$[6, 3, 2]$ & Proposition~\ref{pro:main2} (three starters) \\
$[5, 4, 2]$ & Proposition~\ref{pro:main2} (three starters) \\
$[5, 3, 3]$ & Theorem~\ref{the:intro2} \\
$[4, 4, 3]$ & Theorem~\ref{the:intro2} \\
$[9,2]$ & Theorem~\ref{the:intro2} \\
$[8,3]$ & Theorem~\ref{the:intro2} \\
$[7,4]$ & Theorem~\ref{the:intro2} \\
$[6,5]$ & Theorem~\ref{the:intro2} \\
$[11]$ & Theorem~\ref{the:intro2} \\
\hline
\end{tabular}

\bigskip\bigskip

\noindent {\sc Case} $n=12$: see Appendix~\ref{app:12}.

\bigskip

\begin{tabular}{|l|l|}
\hline
2-factor type $[m_1,\ldots,m_t]$ & HOP$(2m_1,\ldots,2m_t)$ has a solution by... \\ \hline
$[2, 2, 2, 2, 2, 2]$ & Theorem~\ref{the:intro2} \\
$[4, 2, 2, 2, 2]$ & Proposition~\ref{pro:main1} (one starter) \\
$[3, 3, 2, 2, 2]$ & Proposition~\ref{pro:main} (two starters) \\
$[6, 2, 2, 2]$ & Proposition~\ref{pro:main} (two starters) \\
$[5, 3, 2, 2]$ & Proposition~\ref{pro:main1} (one starter) \\
$[4, 4, 2, 2]$ & Proposition~\ref{pro:main1} (one starter) \\
$[4, 3, 3, 2]$ & Proposition~\ref{pro:main} (two starters) \\
$[3, 3, 3, 3]$ & Theorem~\ref{the:intro2} \\
$[8, 2, 2]$ & Proposition~\ref{pro:main1} (one starter) \\
$[7, 3, 2]$ & Proposition~\ref{pro:main} (two starters) \\
$[6, 4, 2]$ & Proposition~\ref{pro:main} (two starters) \\
$[5, 5, 2]$ & Proposition~\ref{pro:main} (two starters) \\
$[6, 3, 3]$ & Proposition~\ref{pro:main1} (one starter) \\
$[5, 4, 3]$ & Proposition~\ref{pro:main1} (one starter) \\
$[4, 4, 4]$ & Theorem~\ref{the:intro2} \\
$[10,2]$ & Proposition~\ref{pro:main} (two starters) \\
$[9,3]$ & Proposition~\ref{pro:main1} (one starter) \\
$[8,4]$ & Theorem~\ref{the:intro2} \\
$[7,5]$ & Proposition~\ref{pro:main1} (one starter) \\
$[6,6]$ & Theorem~\ref{the:intro2} \\
$[12]$ & Theorem~\ref{the:intro2} \\
\hline
\end{tabular}

\newpage

\noindent {\sc Case} $n=13$: see Appendix~\ref{app:13}.

\bigskip

\begin{tabular}{|l|l|}
\hline
2-factor type $[m_1,\ldots,m_t]$ & HOP$(2m_1,\ldots,2m_t)$ has a solution by... \\ \hline
$[3, 2, 2, 2, 2, 2]$ & Proposition~\ref{pro:main2} (three starters) \\
$[5, 2, 2, 2, 2]$ & Proposition~\ref{pro:main2} (three starters) \\
$[4, 3, 2, 2, 2]$ & Proposition~\ref{pro:main2} (three starters) \\
$[3, 3, 3, 2, 2]$ & Proposition~\ref{pro:main2} (three starters) \\
$[7, 2, 2, 2]$ & Proposition~\ref{pro:main2} (three starters) \\
$[6, 3, 2, 2]$ & Proposition~\ref{pro:main2} (three starters) \\
$[5, 4, 2, 2]$ & Proposition~\ref{pro:main2} (three starters) \\
$[5, 3, 3, 2]$ & Proposition~\ref{pro:main2} (three starters) \\
$[4, 4, 3, 2]$ & Proposition~\ref{pro:main2} (three starters) \\
$[4, 3, 3, 3]$ & Theorem~\ref{the:intro2} \\
$[9, 2, 2]$ & Proposition~\ref{pro:main2} (three starters) \\
$[8, 3, 2]$ & Proposition~\ref{pro:main2} (three starters) \\
$[7, 4, 2]$ & Proposition~\ref{pro:main2} (three starters) \\
$[6, 5, 2]$ & Proposition~\ref{pro:main2} (three starters) \\
$[7, 3, 3]$ & Theorem~\ref{the:intro2} \\
$[6, 4, 3]$ & Theorem~\ref{the:intro2} \\
$[5, 5, 3]$ & Theorem~\ref{the:intro2} \\
$[5, 4, 4]$ & Theorem~\ref{the:intro2} \\
$[11,2]$ & Theorem~\ref{the:intro2} \\
$[10,3]$ & Theorem~\ref{the:intro2} \\
$[9,4]$ & Theorem~\ref{the:intro2} \\
$[8,5]$ & Theorem~\ref{the:intro2} \\
$[7,6]$ & Theorem~\ref{the:intro2} \\
$[13]$ & Theorem~\ref{the:intro2} \\
\hline
\end{tabular}

\bigskip\bigskip

\noindent {\sc Case} $n=14$: see Appendix~\ref{app:14}.

\bigskip

\begin{tabular}{|l|l|}
\hline
2-factor type $[m_1,\ldots,m_t]$ & HOP$(2m_1,\ldots,2m_t)$ has a solution by... \\ \hline
$[2, 2, 2, 2, 2, 2, 2]$ & Theorem~\ref{the:intro2} \\
$[4, 2, 2, 2, 2, 2]$ & Proposition~\ref{pro:main1} (one starter) \\
$[3, 3, 2, 2, 2, 2]$ & Proposition~\ref{pro:main} (two starters) \\
$[6, 2, 2, 2, 2]$ & Proposition~\ref{pro:main} (two starters) \\
$[5, 3, 2, 2, 2]$ & Proposition~\ref{pro:main1} (one starter) \\
$[4, 4, 2, 2, 2]$ & Proposition~\ref{pro:main1} (one starter) \\
$[4, 3, 3, 2, 2]$ & Proposition~\ref{pro:main} (two starters) \\
$[3, 3, 3, 3, 2]$ & Proposition~\ref{pro:main1} (one starter) \\
$[8, 2, 2, 2]$ & Proposition~\ref{pro:main1} (one starter) \\
$[7, 3, 2, 2]$ & Proposition~\ref{pro:main} (two starters) \\
$[6, 4, 2, 2]$ & Proposition~\ref{pro:main} (two starters) \\
$[5, 5, 2, 2]$ & Proposition~\ref{pro:main} (two starters) \\
$[6, 3, 3, 2]$ & Proposition~\ref{pro:main1} (one starter) \\
$[5, 4, 3, 2]$ & Proposition~\ref{pro:main1} (one starter) \\
$[4, 4, 4, 2]$ & Proposition~\ref{pro:main1} (one starter) \\
$[5, 3, 3, 3]$ & Proposition~\ref{pro:main} (two starters) \\
$[4, 4, 3, 3]$ & Proposition~\ref{pro:main} (two starters) \\
\hline
\end{tabular}

\begin{tabular}{|l|l|}
\hline
2-factor type $[m_1,\ldots,m_t]$ & HOP$(2m_1,\ldots,2m_t)$ has a solution by... \\ \hline
$[10, 2, 2]$ & Proposition~\ref{pro:main} (two starters) \\
$[9, 3, 2]$ & Proposition~\ref{pro:main1} (one starter) \\
$[8, 4, 2]$ & Proposition~\ref{pro:main1} (one starter) \\
$[7, 5, 2]$ & Proposition~\ref{pro:main1} (one starter) \\
$[6, 6, 2]$ & Proposition~\ref{pro:main1} (one starter) \\
$[8, 3, 3]$ & Proposition~\ref{pro:main} (two starters) \\
$[7, 4, 3]$ & Proposition~\ref{pro:main} (two starters) \\
$[6, 5, 3]$ & Proposition~\ref{pro:main} (two starters) \\
$[6, 4, 4]$ & Proposition~\ref{pro:main} (two starters) \\
$[5, 5, 4]$ & Proposition~\ref{pro:main} (two starters) \\
$[12,2]$ & Proposition~\ref{pro:main1} (one starter) \\
$[11,3]$ & Proposition~\ref{pro:main} (two starters) \\
$[10,4]$ & Proposition~\ref{pro:main} (two starters) \\
$[9,5]$ & Proposition~\ref{pro:main} (two starters) \\
$[8,6]$ & Proposition~\ref{pro:main} (two starters) \\
$[7,7]$ & Theorem~\ref{the:intro2} \\
$[14]$ & Theorem~\ref{the:intro2} \\
\hline
\end{tabular}

\bigskip\bigskip

\noindent {\sc Case} $n=15$: see Appendix~\ref{app:15}.

\bigskip

\begin{tabular}{|l|l|}
\hline
2-factor type $[m_1,\ldots,m_t]$ & HOP$(2m_1,\ldots,2m_t)$ has a solution by... \\ \hline
$[3, 2, 2, 2, 2, 2, 2]$ & Proposition~\ref{pro:main2} (three starters)  \\
$[5, 2, 2, 2, 2, 2]$ & Proposition~\ref{pro:main2} (three starters)  \\
$[4, 3, 2, 2, 2, 2]$ & Proposition~\ref{pro:main2} (three starters)  \\
$[3, 3, 3, 2, 2, 2]$ & Proposition~\ref{pro:main2} (three starters)  \\
$[7, 2, 2, 2, 2]$ & Proposition~\ref{pro:main2} (three starters)  \\
$[6, 3, 2, 2, 2]$ & Proposition~\ref{pro:main2} (three starters)  \\
$[5, 4, 2, 2, 2]$ & Proposition~\ref{pro:main2} (three starters)  \\
$[5, 3, 3, 2, 2]$ & Proposition~\ref{pro:main2} (three starters)  \\
$[4, 4, 3, 2, 2]$ & Proposition~\ref{pro:main2} (three starters)  \\
$[4, 3, 3, 3, 2]$ & Proposition~\ref{pro:main2} (three starters)  \\
$[3, 3, 3, 3, 3]$ & Theorem~\ref{the:intro2} \\
$[9, 2, 2, 2]$ & Proposition~\ref{pro:main2} (three starters)  \\
$[8, 3, 2, 2]$ & Proposition~\ref{pro:main2} (three starters)  \\
$[7, 4, 2, 2]$ & Proposition~\ref{pro:main2} (three starters)  \\
$[6, 5, 2, 2]$ & Proposition~\ref{pro:main2} (three starters)  \\
$[7, 3, 3, 2]$ & Proposition~\ref{pro:main2} (three starters)  \\
$[6, 4, 3, 2]$ & Proposition~\ref{pro:main2} (three starters)  \\
$[5, 5, 3, 2]$ & Proposition~\ref{pro:main2} (three starters)  \\
$[5, 4, 4, 2]$ & Proposition~\ref{pro:main2} (three starters)  \\
$[6, 3, 3, 3]$ & Theorem~\ref{the:intro2} \\
$[5, 4, 3, 3]$ & Theorem~\ref{the:intro2} \\
$[4, 4, 4, 3]$ & Theorem~\ref{the:intro2} \\
\hline
\end{tabular}

\begin{tabular}{|l|l|}
\hline
2-factor type $[m_1,\ldots,m_t]$ & HOP$(2m_1,\ldots,2m_t)$ has a solution by... \\ \hline
$[11, 2, 2]$ & Proposition~\ref{pro:main2} (three starters)  \\
$[10, 3, 2]$ & Proposition~\ref{pro:main2} (three starters)  \\
$[9, 4, 2]$ & Proposition~\ref{pro:main2} (three starters)  \\
$[8, 5, 2]$ & Proposition~\ref{pro:main2} (three starters)  \\
$[7, 6, 2]$ & Proposition~\ref{pro:main2} (three starters)  \\
$[9, 3, 3]$ & Theorem~\ref{the:intro2} \\
$[8, 4, 3]$ & Theorem~\ref{the:intro2} \\
$[7, 5, 3]$ & Theorem~\ref{the:intro2} \\
$[6, 6, 3]$ & Theorem~\ref{the:intro2} \\
$[7, 4, 4]$ & Theorem~\ref{the:intro2} \\
$[6, 5, 4]$ & Theorem~\ref{the:intro2} \\
$[5, 5, 5]$ & Theorem~\ref{the:intro2} \\
$[13,2]$ & Theorem~\ref{the:intro2} \\
$[12,3]$ & Theorem~\ref{the:intro2} \\
$[11,4]$ & Theorem~\ref{the:intro2} \\
$[10,5]$ & Theorem~\ref{the:intro2} \\
$[9,6]$ & Theorem~\ref{the:intro2} \\
$[8,7]$ & Theorem~\ref{the:intro2} \\
$[15]$ & Theorem~\ref{the:intro2} \\
\hline
\end{tabular}

\bigskip\bigskip

\noindent {\sc Case} $n=16$: see Appendix~\ref{app:16}.

\bigskip

\begin{tabular}{|l|l|}
\hline
2-factor type $[m_1,\ldots,m_t]$ & HOP$(2m_1,\ldots,2m_t)$ has a solution by... \\ \hline
$[2, 2, 2, 2, 2, 2, 2, 2]$ & Theorem~\ref{the:intro2} \\
$[4, 2, 2, 2, 2, 2, 2]$ & Proposition~\ref{pro:main1} (one starter) \\
$[3, 3, 2, 2, 2, 2, 2]$ & Proposition~\ref{pro:main} (two starters) \\
$[6, 2, 2, 2, 2, 2]$ & Proposition~\ref{pro:main} (two starters) \\
$[5, 3, 2, 2, 2, 2]$ & Proposition~\ref{pro:main1} (one starter) \\
$[4, 4, 2, 2, 2, 2]$ & Proposition~\ref{pro:main1} (one starter) \\
$[4, 3, 3, 2, 2, 2]$ & Proposition~\ref{pro:main} (two starters) \\
$[3, 3, 3, 3, 2, 2]$ & Proposition~\ref{pro:main1} (one starter) \\
$[8, 2, 2, 2, 2]$ & Proposition~\ref{pro:main1} (one starter) \\
$[7, 3, 2, 2, 2]$ & Proposition~\ref{pro:main} (two starters) \\
$[6, 4, 2, 2, 2]$ & Proposition~\ref{pro:main} (two starters) \\
$[5, 5, 2, 2, 2]$ & Proposition~\ref{pro:main} (two starters) \\
$[6, 3, 3, 2, 2]$ & Proposition~\ref{pro:main1} (one starter) \\
$[5, 4, 3, 2, 2]$ & Proposition~\ref{pro:main1} (one starter) \\
$[4, 4, 4, 2, 2]$ & Proposition~\ref{pro:main1} (one starter) \\
$[5, 3, 3, 3, 2]$ & Proposition~\ref{pro:main} (two starters) \\
$[4, 4, 3, 3, 2]$ & Proposition~\ref{pro:main} (two starters) \\
$[4, 3, 3, 3, 3]$ & Proposition~\ref{pro:main1} (one starter) \\
$[10, 2, 2, 2]$ & Proposition~\ref{pro:main} (two starters) \\
$[9, 3, 2, 2]$ & Proposition~\ref{pro:main1} (one starter) \\
$[8, 4, 2, 2]$ & Proposition~\ref{pro:main1} (one starter) \\
$[7, 5, 2, 2]$ & Proposition~\ref{pro:main1} (one starter) \\
$[6, 6, 2, 2]$ & Proposition~\ref{pro:main1} (one starter) \\
\hline
\end{tabular}

\begin{tabular}{|l|l|}
\hline
2-factor type $[m_1,\ldots,m_t]$ & HOP$(2m_1,\ldots,2m_t)$ has a solution by... \\ \hline
$[8, 3, 3, 2]$ & Proposition~\ref{pro:main} (two starters) \\
$[7, 4, 3, 2]$ & Proposition~\ref{pro:main} (two starters) \\
$[6, 5, 3, 2]$ & Proposition~\ref{pro:main} (two starters) \\
$[6, 4, 4, 2]$ & Proposition~\ref{pro:main} (two starters) \\
$[5, 5, 4, 2]$ & Proposition~\ref{pro:main} (two starters) \\
$[7, 3, 3, 3]$ & Proposition~\ref{pro:main1} (one starter) \\
$[6, 4, 3, 3]$ & Proposition~\ref{pro:main1} (one starter) \\
$[5, 5, 3, 3]$ & Proposition~\ref{pro:main1} (one starter) \\
$[5, 4, 4, 3]$ & Proposition~\ref{pro:main1} (one starter) \\
$[4, 4, 4, 4]$ & Theorem~\ref{the:intro2} \\
$[12, 2, 2]$ & Proposition~\ref{pro:main1} (one starter) \\
$[11, 3, 2]$ & Proposition~\ref{pro:main} (two starters) \\
$[10, 4, 2]$ & Proposition~\ref{pro:main} (two starters) \\
$[9, 5, 2]$ & Proposition~\ref{pro:main} (two starters) \\
$[8, 6, 2]$ & Proposition~\ref{pro:main} (two starters) \\
$[7, 7, 2]$ & Proposition~\ref{pro:main} (two starters) \\
$[10, 3, 3]$ & Proposition~\ref{pro:main1} (one starter) \\
$[9, 4, 3]$ & Proposition~\ref{pro:main1} (one starter) \\
$[8, 5, 3]$ & Proposition~\ref{pro:main1} (one starter) \\
$[7, 6, 3]$ & Proposition~\ref{pro:main1} (one starter) \\
$[8, 4, 4]$ & Theorem~\ref{the:intro2} \\
$[7, 5, 4]$ & Proposition~\ref{pro:main1} (one starter) \\
$[6, 6, 4]$ & Proposition~\ref{pro:main1} (one starter) \\
$[6, 5, 5]$ & Proposition~\ref{pro:main1} (one starter) \\
$[14, 2]$ & Proposition~\ref{pro:main} (two starters) \\
$[13, 3]$ & Proposition~\ref{pro:main1} (one starter) \\
$[12, 4]$ & Theorem~\ref{the:intro2} \\
$[11, 5]$ & Proposition~\ref{pro:main1} (one starter) \\
$[10, 6]$ & Proposition~\ref{pro:main1} (one starter) \\
$[9, 7]$ & Proposition~\ref{pro:main1} (one starter) \\
$[8, 8]$ & Theorem~\ref{the:intro2} \\
$[16]$ & Theorem~\ref{the:intro2} \\
\hline
\end{tabular}

\bigskip\bigskip

\noindent {\sc Case} $n=17$: see Appendix~\ref{app:17}.

\bigskip

\begin{tabular}{|l|l|}
\hline
2-factor type $[m_1,\ldots,m_t]$ & HOP$(2m_1,\ldots,2m_t)$ has a solution by... \\ \hline
$[3, 2, 2, 2, 2, 2, 2, 2]$ & Proposition~\ref{pro:main2} (three starters)  \\
$[5, 2, 2, 2, 2, 2, 2]$ & Proposition~\ref{pro:main2} (three starters)  \\
$[4, 3, 2, 2, 2, 2, 2]$ & Proposition~\ref{pro:main2} (three starters)  \\
$[3, 3, 3, 2, 2, 2, 2]$ & Proposition~\ref{pro:main2} (three starters)  \\
$[7, 2, 2, 2, 2, 2]$ & Proposition~\ref{pro:main2} (three starters)  \\
$[6, 3, 2, 2, 2, 2]$ & Proposition~\ref{pro:main2} (three starters)  \\
$[5, 4, 2, 2, 2, 2]$ & Proposition~\ref{pro:main2} (three starters)  \\
$[5, 3, 3, 2, 2, 2]$ & Proposition~\ref{pro:main2} (three starters)  \\
$[4, 4, 3, 2, 2, 2]$ & Proposition~\ref{pro:main2} (three starters)  \\
$[4, 3, 3, 3, 2, 2]$ & Proposition~\ref{pro:main2} (three starters)  \\
$[3, 3, 3, 3, 3, 2]$ & Proposition~\ref{pro:main2} (three starters)  \\
\hline
\end{tabular}

\begin{tabular}{|l|l|}
\hline
2-factor type $[m_1,\ldots,m_t]$ & HOP$(2m_1,\ldots,2m_t)$ has a solution by... \\ \hline
$[9, 2, 2, 2, 2]$ & Proposition~\ref{pro:main2} (three starters)  \\
$[8, 3, 2, 2, 2]$ & Proposition~\ref{pro:main2} (three starters)  \\
$[7, 4, 2, 2, 2]$ & Proposition~\ref{pro:main2} (three starters)  \\
$[6, 5, 2, 2, 2]$ & Proposition~\ref{pro:main2} (three starters)  \\
$[7, 3, 3, 2, 2]$ & Proposition~\ref{pro:main2} (three starters)  \\
$[6, 4, 3, 2, 2]$ & Proposition~\ref{pro:main2} (three starters)  \\
$[5, 5, 3, 2, 2]$ & Proposition~\ref{pro:main2} (three starters)  \\
$[5, 4, 4, 2, 2]$ & Proposition~\ref{pro:main2} (three starters)  \\
$[6, 3, 3, 3, 2]$ & Proposition~\ref{pro:main2} (three starters)  \\
$[5, 4, 3, 3, 2]$ & Proposition~\ref{pro:main2} (three starters)  \\
$[4, 4, 4, 3, 2]$ & Proposition~\ref{pro:main2} (three starters)  \\
$[5, 3, 3, 3, 3]$ & Theorem~\ref{the:intro2} \\
$[4, 4, 3, 3, 3]$ & Theorem~\ref{the:intro2} \\
$[11, 2, 2, 2]$ & Proposition~\ref{pro:main2} (three starters)  \\
$[10, 3, 2, 2]$ & Proposition~\ref{pro:main2} (three starters)  \\
$[9, 4, 2, 2]$ & Proposition~\ref{pro:main2} (three starters)  \\
$[8, 5, 2, 2]$ & Proposition~\ref{pro:main2} (three starters)  \\
$[7, 6, 2, 2]$ & Proposition~\ref{pro:main2} (three starters)  \\
$[9, 3, 3, 2]$ & Proposition~\ref{pro:main2} (three starters)  \\
$[8, 4, 3, 2]$ & Proposition~\ref{pro:main2} (three starters)  \\
$[7, 5, 3, 2]$ & Proposition~\ref{pro:main2} (three starters)  \\
$[6, 6, 3, 2]$ & Proposition~\ref{pro:main2} (three starters)  \\
$[7, 4, 4, 2]$ & Proposition~\ref{pro:main2} (three starters)  \\
$[6, 5, 4, 2]$ & Proposition~\ref{pro:main2} (three starters)  \\
$[5, 5, 5, 2]$ & Proposition~\ref{pro:main2} (three starters)  \\
$[8, 3, 3, 3]$ & Theorem~\ref{the:intro2} \\
$[7, 4, 3, 3]$ & Theorem~\ref{the:intro2} \\
$[6, 5, 3, 3]$ & Theorem~\ref{the:intro2} \\
$[6, 4, 4, 3]$ & Theorem~\ref{the:intro2} \\
$[5, 5, 4, 3]$ & Theorem~\ref{the:intro2} \\
$[5, 4, 4, 4]$ & Theorem~\ref{the:intro2} \\
$[13, 2, 2]$ & Proposition~\ref{pro:main2} (three starters)  \\
$[12, 3, 2]$ & Proposition~\ref{pro:main2} (three starters)  \\
$[11, 4, 2]$ & Proposition~\ref{pro:main2} (three starters)  \\
$[10, 5, 2]$ & Proposition~\ref{pro:main2} (three starters)  \\
$[9, 6, 2]$ & Proposition~\ref{pro:main2} (three starters)  \\
$[8, 7, 2]$ & Proposition~\ref{pro:main2} (three starters)  \\
$[11, 3, 3]$ & Theorem~\ref{the:intro2} \\
$[10, 4, 3]$ & Theorem~\ref{the:intro2} \\
$[9, 5, 3]$ & Theorem~\ref{the:intro2} \\
$[8, 6, 3]$ & Theorem~\ref{the:intro2} \\
$[7, 7, 3]$ & Theorem~\ref{the:intro2} \\
$[9, 4, 4]$ & Theorem~\ref{the:intro2} \\
$[8, 5, 4]$ & Theorem~\ref{the:intro2} \\
$[7, 6, 4]$ & Theorem~\ref{the:intro2} \\
$[7, 5, 5]$ & Theorem~\ref{the:intro2} \\
$[6, 6, 5]$ & Theorem~\ref{the:intro2} \\%
\hline
\end{tabular}

\begin{tabular}{|l|l|}
\hline
2-factor type $[m_1,\ldots,m_t]$ & HOP$(2m_1,\ldots,2m_t)$ has a solution by... \\ \hline
$[15, 2]$ & Theorem~\ref{the:intro2} \\
$[14, 3]$ & Theorem~\ref{the:intro2} \\
$[13, 4]$ & Theorem~\ref{the:intro2} \\
$[12, 5]$ & Theorem~\ref{the:intro2} \\
$[11, 6]$ & Theorem~\ref{the:intro2} \\
$[10, 7]$ & Theorem~\ref{the:intro2} \\
$[9, 8]$ & Theorem~\ref{the:intro2} \\
$[17]$ & Theorem~\ref{the:intro2} \\
\hline
\end{tabular}

\bigskip\bigskip

\noindent {\sc Case} $n=18$: see Appendix~\ref{app:18}.

\bigskip

\begin{tabular}{|l|l|}
\hline
2-factor type $[m_1,\ldots,m_t]$ & HOP$(2m_1,\ldots,2m_t)$ has a solution by... \\ \hline
$[2, 2, 2, 2, 2, 2, 2, 2, 2]$ & Theorem~\ref{the:intro2} \\
$[4, 2, 2, 2, 2, 2, 2, 2]$ & Proposition~\ref{pro:main1} (one starter) \\
$[3, 3, 2, 2, 2, 2, 2, 2]$ & Proposition~\ref{pro:main} (two starters) \\
$[6, 2, 2, 2, 2, 2, 2]$ & Proposition~\ref{pro:main} (two starters) \\
$[5, 3, 2, 2, 2, 2, 2]$ & Proposition~\ref{pro:main1} (one starter) \\
$[4, 4, 2, 2, 2, 2, 2]$ & Proposition~\ref{pro:main1} (one starter) \\
$[4, 3, 3, 2, 2, 2, 2]$ & Proposition~\ref{pro:main} (two starters) \\
$[3, 3, 3, 3, 2, 2, 2]$ & Proposition~\ref{pro:main1} (one starter) \\
$[8, 2, 2, 2, 2, 2]$ & Proposition~\ref{pro:main1} (one starter) \\
$[7, 3, 2, 2, 2, 2]$ & Proposition~\ref{pro:main} (two starters) \\
$[6, 4, 2, 2, 2, 2]$ & Proposition~\ref{pro:main} (two starters) \\
$[5, 5, 2, 2, 2, 2]$ & Proposition~\ref{pro:main} (two starters) \\
$[6, 3, 3, 2, 2, 2]$ & Proposition~\ref{pro:main1} (one starter) \\
$[5, 4, 3, 2, 2, 2]$ & Proposition~\ref{pro:main1} (one starter) \\
$[4, 4, 4, 2, 2, 2]$ & Proposition~\ref{pro:main1} (one starter) \\
$[5, 3, 3, 3, 2, 2]$ & Proposition~\ref{pro:main} (two starters) \\
$[4, 4, 3, 3, 2, 2]$ & Proposition~\ref{pro:main} (two starters) \\
$[4, 3, 3, 3, 3, 2]$ & Proposition~\ref{pro:main1} (one starter) \\
$[3, 3, 3, 3, 3, 3]$ & Theorem~\ref{the:intro2} \\
$[10, 2, 2, 2, 2]$ & Proposition~\ref{pro:main} (two starters) \\
$[9, 3, 2, 2, 2]$ & Proposition~\ref{pro:main1} (one starter) \\
$[8, 4, 2, 2, 2]$ & Proposition~\ref{pro:main1} (one starter) \\
$[7, 5, 2, 2, 2]$ & Proposition~\ref{pro:main1} (one starter) \\
$[6, 6, 2, 2, 2]$ & Proposition~\ref{pro:main1} (one starter) \\
$[8, 3, 3, 2, 2]$ & Proposition~\ref{pro:main} (two starters) \\
$[7, 4, 3, 2, 2]$ & Proposition~\ref{pro:main} (two starters) \\
$[6, 5, 3, 2, 2]$ & Proposition~\ref{pro:main} (two starters) \\
$[6, 4, 4, 2, 2]$ & Proposition~\ref{pro:main} (two starters) \\
$[5, 5, 4, 2, 2]$ & Proposition~\ref{pro:main} (two starters) \\
$[7, 3, 3, 3, 2]$ & Proposition~\ref{pro:main1} (one starter) \\
$[6, 4, 3, 3, 2]$ & Proposition~\ref{pro:main1} (one starter) \\
$[5, 5, 3, 3, 2]$ & Proposition~\ref{pro:main1} (one starter) \\
$[5, 4, 4, 3, 2]$ & Proposition~\ref{pro:main1} (one starter) \\
$[4, 4, 4, 4, 2]$ & Proposition~\ref{pro:main1} (one starter) \\
\hline
\end{tabular}

\begin{tabular}{|l|l|}
\hline
2-factor type $[m_1,\ldots,m_t]$ & HOP$(2m_1,\ldots,2m_t)$ has a solution by... \\ \hline
$[6, 3, 3, 3, 3]$ & Proposition~\ref{pro:main} (two starters) \\
$[5, 4, 3, 3, 3]$ & Proposition~\ref{pro:main} (two starters) \\
$[4, 4, 4, 3, 3]$ & Proposition~\ref{pro:main} (two starters) \\
$[12, 2, 2, 2]$ & Proposition~\ref{pro:main1} (one starter) \\
$[11, 3, 2, 2]$ & Proposition~\ref{pro:main} (two starters) \\
$[10, 4, 2, 2]$ & Proposition~\ref{pro:main} (two starters) \\
$[9, 5, 2, 2]$ & Proposition~\ref{pro:main} (two starters) \\
$[8, 6, 2, 2]$ & Proposition~\ref{pro:main} (two starters) \\
$[7, 7, 2, 2]$ & Proposition~\ref{pro:main} (two starters) \\
$[10, 3, 3, 2]$ & Proposition~\ref{pro:main1} (one starter) \\
$[9, 4, 3, 2]$ & Proposition~\ref{pro:main1} (one starter) \\
$[8, 5, 3, 2]$ & Proposition~\ref{pro:main1} (one starter) \\
$[7, 6, 3, 2]$ & Proposition~\ref{pro:main1} (one starter) \\
$[8, 4, 4, 2]$ & Proposition~\ref{pro:main1} (one starter) \\
$[7, 5, 4, 2]$ & Proposition~\ref{pro:main1} (one starter) \\
$[6, 6, 4, 2]$ & Proposition~\ref{pro:main1} (one starter) \\
$[6, 5, 5, 2]$ & Proposition~\ref{pro:main1} (one starter) \\
$[9, 3, 3, 3]$ & Proposition~\ref{pro:main} (two starters) \\
$[8, 4, 3, 3]$ & Proposition~\ref{pro:main} (two starters) \\
$[7, 5, 3, 3]$ & Proposition~\ref{pro:main} (two starters) \\
$[6, 6, 3, 3]$ & Proposition~\ref{pro:main} (two starters) \\
$[7, 4, 4, 3]$ & Proposition~\ref{pro:main} (two starters) \\
$[6, 5, 4, 3]$ & Proposition~\ref{pro:main} (two starters) \\
$[5, 5, 5, 3]$ & Proposition~\ref{pro:main} (two starters) \\
$[6, 4, 4, 4]$ & Proposition~\ref{pro:main} (two starters) \\
$[5, 5, 4, 4]$ & Proposition~\ref{pro:main} (two starters) \\
$[14, 2, 2]$ & Proposition~\ref{pro:main} (two starters) \\
$[13, 3, 2]$ & Proposition~\ref{pro:main1} (one starter) \\
$[12, 4, 2]$ & Proposition~\ref{pro:main1} (one starter) \\
$[11, 5, 2]$ & Proposition~\ref{pro:main1} (one starter) \\
$[10, 6, 2]$ & Proposition~\ref{pro:main1} (one starter) \\
$[9, 7, 2]$ & Proposition~\ref{pro:main1} (one starter) \\
$[8, 8, 2]$ & Proposition~\ref{pro:main1} (one starter) \\
$[12, 3, 3]$ & Proposition~\ref{pro:main} (two starters) \\
$[11, 4, 3]$ & Proposition~\ref{pro:main} (two starters) \\
$[10, 5, 3]$ & Proposition~\ref{pro:main} (two starters) \\
$[9, 6, 3]$ & Proposition~\ref{pro:main} (two starters) \\
$[8, 7, 3]$ & Proposition~\ref{pro:main} (two starters) \\
$[10, 4, 4]$ & Proposition~\ref{pro:main} (two starters) \\
$[9, 5, 4]$ & Proposition~\ref{pro:main} (two starters) \\
$[8, 6, 4]$ & Proposition~\ref{pro:main} (two starters) \\
$[7, 7, 4]$ & Proposition~\ref{pro:main} (two starters) \\
$[8, 5, 5]$ & Proposition~\ref{pro:main} (two starters) \\
$[7, 6, 5]$ & Proposition~\ref{pro:main} (two starters) \\
$[6, 6, 6]$ & Theorem~\ref{the:intro2} \\
\hline
\end{tabular}

\begin{tabular}{|l|l|}
\hline
2-factor type $[m_1,\ldots,m_t]$ & HOP$(2m_1,\ldots,2m_t)$ has a solution by... \\ \hline
$[16, 2]$ & Proposition~\ref{pro:main1} (one starter) \\
$[15, 3]$ & Proposition~\ref{pro:main} (two starters) \\
$[14, 4]$ & Proposition~\ref{pro:main} (two starters) \\
$[13, 5]$ & Proposition~\ref{pro:main} (two starters) \\
$[12, 6]$ & Proposition~\ref{pro:main} (two starters) \\
$[11, 7]$ & Proposition~\ref{pro:main} (two starters) \\
$[10, 8]$ & Proposition~\ref{pro:main} (two starters) \\
$[9, 9]$ & Theorem~\ref{the:intro2} \\
$[18]$ & Theorem~\ref{the:intro2} \\
\hline
\end{tabular}

\bigskip\bigskip

\noindent {\sc Case} $n=19$: see Appendix~\ref{app:19}.

\bigskip

\begin{tabular}{|l|l|}
\hline
2-factor type $[m_1,\ldots,m_t]$ & HOP$(2m_1,\ldots,2m_t)$ has a solution by... \\ \hline
$[3, 2, 2, 2, 2, 2, 2, 2, 2]$ & Proposition~\ref{pro:main2} (three starters)  \\
$[5, 2, 2, 2, 2, 2, 2, 2]$ & Proposition~\ref{pro:main2} (three starters)  \\
$[4, 3, 2, 2, 2, 2, 2, 2]$ & Proposition~\ref{pro:main2} (three starters)  \\
$[3, 3, 3, 2, 2, 2, 2, 2]$ & Proposition~\ref{pro:main2} (three starters)  \\
$[7, 2, 2, 2, 2, 2, 2]$ & Proposition~\ref{pro:main2} (three starters)  \\
$[6, 3, 2, 2, 2, 2, 2]$ & Proposition~\ref{pro:main2} (three starters)  \\
$[5, 4, 2, 2, 2, 2, 2]$ & Proposition~\ref{pro:main2} (three starters)  \\
$[5, 3, 3, 2, 2, 2, 2]$ & Proposition~\ref{pro:main2} (three starters)  \\
$[4, 4, 3, 2, 2, 2, 2]$ & Proposition~\ref{pro:main2} (three starters)  \\
$[4, 3, 3, 3, 2, 2, 2]$ & Proposition~\ref{pro:main2} (three starters)  \\
$[3, 3, 3, 3, 3, 2, 2]$ & Proposition~\ref{pro:main2} (three starters)  \\
$[9, 2, 2, 2, 2, 2]$ & Proposition~\ref{pro:main2} (three starters)  \\
$[8, 3, 2, 2, 2, 2]$ & Proposition~\ref{pro:main2} (three starters)  \\
$[7, 4, 2, 2, 2, 2]$ & Proposition~\ref{pro:main2} (three starters)  \\
$[6, 5, 2, 2, 2, 2]$ & Proposition~\ref{pro:main2} (three starters)  \\
$[7, 3, 3, 2, 2, 2]$ & Proposition~\ref{pro:main2} (three starters)  \\
$[6, 4, 3, 2, 2, 2]$ & Proposition~\ref{pro:main2} (three starters)  \\
$[5, 5, 3, 2, 2, 2]$ & Proposition~\ref{pro:main2} (three starters)  \\
$[5, 4, 4, 2, 2, 2]$ & Proposition~\ref{pro:main2} (three starters)  \\
$[6, 3, 3, 3, 2, 2]$ & Proposition~\ref{pro:main2} (three starters)  \\
$[5, 4, 3, 3, 2, 2]$ & Proposition~\ref{pro:main2} (three starters)  \\
$[4, 4, 4, 3, 2, 2]$ & Proposition~\ref{pro:main2} (three starters)  \\
$[5, 3, 3, 3, 3, 2]$ & Proposition~\ref{pro:main2} (three starters)  \\
$[4, 4, 3, 3, 3, 2]$ & Proposition~\ref{pro:main2} (three starters)  \\
$[4, 3, 3, 3, 3, 3]$ & Theorem~\ref{the:intro2} \\
$[11, 2, 2, 2, 2]$ & Proposition~\ref{pro:main2} (three starters)  \\
$[10, 3, 2, 2, 2]$ & Proposition~\ref{pro:main2} (three starters)  \\
$[9, 4, 2, 2, 2]$ & Proposition~\ref{pro:main2} (three starters)  \\
$[8, 5, 2, 2, 2]$ & Proposition~\ref{pro:main2} (three starters)  \\
$[7, 6, 2, 2, 2]$ & Proposition~\ref{pro:main2} (three starters)  \\
$[9, 3, 3, 2, 2]$ & Proposition~\ref{pro:main2} (three starters)  \\
$[8, 4, 3, 2, 2]$ & Proposition~\ref{pro:main2} (three starters)  \\
$[7, 5, 3, 2, 2]$ & Proposition~\ref{pro:main2} (three starters)  \\
\hline
\end{tabular}

\begin{tabular}{|l|l|}
\hline
2-factor type $[m_1,\ldots,m_t]$ & HOP$(2m_1,\ldots,2m_t)$ has a solution by... \\ \hline
$[6, 6, 3, 2, 2]$ & Proposition~\ref{pro:main2} (three starters)  \\
$[7, 4, 4, 2, 2]$ & Proposition~\ref{pro:main2} (three starters)  \\
$[6, 5, 4, 2, 2]$ & Proposition~\ref{pro:main2} (three starters)  \\
$[5, 5, 5, 2, 2]$ & Proposition~\ref{pro:main2} (three starters)  \\
$[8, 3, 3, 3, 2]$ & Proposition~\ref{pro:main2} (three starters)  \\
$[7, 4, 3, 3, 2]$ & Proposition~\ref{pro:main2} (three starters)  \\
$[6, 5, 3, 3, 2]$ & Proposition~\ref{pro:main2} (three starters)  \\
$[6, 4, 4, 3, 2]$ & Proposition~\ref{pro:main2} (three starters)  \\
$[5, 5, 4, 3, 2]$ & Proposition~\ref{pro:main2} (three starters)  \\
$[5, 4, 4, 4, 2]$ & Proposition~\ref{pro:main2} (three starters)  \\
$[7, 3, 3, 3, 3]$ & Theorem~\ref{the:intro2} \\
$[6, 4, 3, 3, 3]$ & Theorem~\ref{the:intro2} \\
$[5, 5, 3, 3, 3]$ & Theorem~\ref{the:intro2} \\
$[5, 4, 4, 3, 3]$ & Theorem~\ref{the:intro2} \\
$[4, 4, 4, 4, 3]$ & Theorem~\ref{the:intro2} \\
$[13, 2, 2, 2]$ & Proposition~\ref{pro:main2} (three starters)  \\
$[12, 3, 2, 2]$ & Proposition~\ref{pro:main2} (three starters)  \\
$[11, 4, 2, 2]$ & Proposition~\ref{pro:main2} (three starters)  \\
$[10, 5, 2, 2]$ & Proposition~\ref{pro:main2} (three starters)  \\
$[9, 6, 2, 2]$ & Proposition~\ref{pro:main2} (three starters)  \\
$[8, 7, 2, 2]$ & Proposition~\ref{pro:main2} (three starters)  \\
$[11, 3, 3, 2]$ & Proposition~\ref{pro:main2} (three starters)  \\
$[10, 4, 3, 2]$ & Proposition~\ref{pro:main2} (three starters)  \\
$[9, 5, 3, 2]$ & Proposition~\ref{pro:main2} (three starters)  \\
$[8, 6, 3, 2]$ & Proposition~\ref{pro:main2} (three starters)  \\
$[7, 7, 3, 2]$ & Proposition~\ref{pro:main2} (three starters)  \\
$[9, 4, 4, 2]$ & Proposition~\ref{pro:main2} (three starters)  \\
$[8, 5, 4, 2]$ & Proposition~\ref{pro:main2} (three starters)  \\
$[7, 6, 4, 2]$ & Proposition~\ref{pro:main2} (three starters)  \\
$[7, 5, 5, 2]$ & Proposition~\ref{pro:main2} (three starters)  \\
$[6, 6, 5, 2]$ & Proposition~\ref{pro:main2} (three starters)  \\
$[10, 3, 3, 3]$ & Theorem~\ref{the:intro2} \\
$[9, 4, 3, 3]$ & Theorem~\ref{the:intro2} \\
$[8, 5, 3, 3]$ & Theorem~\ref{the:intro2} \\
$[7, 6, 3, 3]$ & Theorem~\ref{the:intro2} \\
$[8, 4, 4, 3]$ & Theorem~\ref{the:intro2} \\
$[7, 5, 4, 3]$ & Theorem~\ref{the:intro2} \\
$[6, 6, 4, 3]$ & Theorem~\ref{the:intro2} \\
$[6, 5, 5, 3]$ & Theorem~\ref{the:intro2} \\
$[7, 4, 4, 4]$ & Theorem~\ref{the:intro2} \\
$[6, 5, 4, 4]$ & Theorem~\ref{the:intro2} \\
$[5, 5, 5, 4]$ & Theorem~\ref{the:intro2} \\
$[15, 2, 2]$ & Proposition~\ref{pro:main2} (three starters)  \\
$[14, 3, 2]$ & Proposition~\ref{pro:main2} (three starters)  \\
$[13, 4, 2]$ & Proposition~\ref{pro:main2} (three starters)  \\
$[12, 5, 2]$ & Proposition~\ref{pro:main2} (three starters)  \\
\hline
\end{tabular}

\begin{tabular}{|l|l|}
\hline
2-factor type $[m_1,\ldots,m_t]$ & HOP$(2m_1,\ldots,2m_t)$ has a solution by... \\ \hline
$[11, 6, 2]$ & Proposition~\ref{pro:main2} (three starters)  \\
$[10, 7, 2]$ & Proposition~\ref{pro:main2} (three starters)  \\
$[9, 8, 2]$ & Proposition~\ref{pro:main2} (three starters)  \\
$[13, 3, 3]$ & Theorem~\ref{the:intro2} \\
$[12, 4, 3]$ & Theorem~\ref{the:intro2} \\
$[11, 5, 3]$ & Theorem~\ref{the:intro2} \\
$[10, 6, 3]$ & Theorem~\ref{the:intro2} \\
$[9, 7, 3]$ & Theorem~\ref{the:intro2} \\
$[8, 8, 3]$ & Theorem~\ref{the:intro2} \\
$[11, 4, 4]$ & Theorem~\ref{the:intro2} \\
$[10, 5, 4]$ & Theorem~\ref{the:intro2} \\
$[9, 6, 4]$ & Theorem~\ref{the:intro2} \\
$[8, 7, 4]$ & Theorem~\ref{the:intro2} \\
$[9, 5, 5]$ & Theorem~\ref{the:intro2} \\
$[8, 6, 5]$ & Theorem~\ref{the:intro2} \\
$[7, 7, 5]$ & Theorem~\ref{the:intro2} \\
$[7, 6, 6]$ & Theorem~\ref{the:intro2} \\
$[17, 2]$ & Theorem~\ref{the:intro2} \\
$[16, 3]$ & Theorem~\ref{the:intro2} \\
$[15, 4]$ & Theorem~\ref{the:intro2} \\
$[14, 5]$ & Theorem~\ref{the:intro2} \\
$[13, 6]$ & Theorem~\ref{the:intro2} \\
$[12, 7]$ & Theorem~\ref{the:intro2} \\
$[11, 8]$ & Theorem~\ref{the:intro2} \\
$[10, 9]$ & Theorem~\ref{the:intro2} \\
$[19]$ & Theorem~\ref{the:intro2} \\
\hline
\end{tabular}

\bigskip\bigskip

\noindent {\sc Case} $n=20$: see Appendix~\ref{app:20}.

\bigskip

\begin{tabular}{|l|l|}
\hline
2-factor type $[m_1,\ldots,m_t]$ & HOP$(2m_1,\ldots,2m_t)$ has a solution by... \\ \hline
$[2, 2, 2, 2, 2, 2, 2, 2, 2, 2]$ & Theorem~\ref{the:intro2} \\
$[4, 2, 2, 2, 2, 2, 2, 2, 2]$ & Proposition~\ref{pro:main1} (one starter) \\
$[3, 3, 2, 2, 2, 2, 2, 2, 2]$ & Proposition~\ref{pro:main} (two starters) \\
$[6, 2, 2, 2, 2, 2, 2, 2]$ & Proposition~\ref{pro:main} (two starters) \\
$[5, 3, 2, 2, 2, 2, 2, 2]$ & Proposition~\ref{pro:main1} (one starter) \\
$[4, 4, 2, 2, 2, 2, 2, 2]$ & Proposition~\ref{pro:main1} (one starter) \\
$[4, 3, 3, 2, 2, 2, 2, 2]$ & Proposition~\ref{pro:main} (two starters) \\
$[3, 3, 3, 3, 2, 2, 2, 2]$ & Proposition~\ref{pro:main1} (one starter) \\
$[8, 2, 2, 2, 2, 2, 2]$ & Proposition~\ref{pro:main1} (one starter) \\
$[7, 3, 2, 2, 2, 2, 2]$ & Proposition~\ref{pro:main} (two starters) \\
$[6, 4, 2, 2, 2, 2, 2]$ & Proposition~\ref{pro:main} (two starters) \\
$[5, 5, 2, 2, 2, 2, 2]$ & Proposition~\ref{pro:main} (two starters) \\
$[6, 3, 3, 2, 2, 2, 2]$ & Proposition~\ref{pro:main1} (one starter) \\
$[5, 4, 3, 2, 2, 2, 2]$ & Proposition~\ref{pro:main1} (one starter) \\
$[4, 4, 4, 2, 2, 2, 2]$ & Proposition~\ref{pro:main1} (one starter) \\
$[5, 3, 3, 3, 2, 2, 2]$ & Proposition~\ref{pro:main} (two starters) \\
\hline
\end{tabular}

\begin{tabular}{|l|l|}
\hline
2-factor type $[m_1,\ldots,m_t]$ & HOP$(2m_1,\ldots,2m_t)$ has a solution by... \\ \hline
$[4, 4, 3, 3, 2, 2, 2]$ & Proposition~\ref{pro:main} (two starters) \\
$[4, 3, 3, 3, 3, 2, 2]$ & Proposition~\ref{pro:main1} (one starter) \\
$[3, 3, 3, 3, 3, 3, 2]$ & Proposition~\ref{pro:main} (two starters) \\
$[10, 2, 2, 2, 2, 2]$ & Proposition~\ref{pro:main} (two starters) \\
$[9, 3, 2, 2, 2, 2]$ & Proposition~\ref{pro:main1} (one starter) \\
$[8, 4, 2, 2, 2, 2]$ & Proposition~\ref{pro:main1} (one starter) \\
$[7, 5, 2, 2, 2, 2]$ & Proposition~\ref{pro:main1} (one starter) \\
$[6, 6, 2, 2, 2, 2]$ & Proposition~\ref{pro:main1} (one starter) \\
$[8, 3, 3, 2, 2, 2]$ & Proposition~\ref{pro:main} (two starters) \\
$[7, 4, 3, 2, 2, 2]$ & Proposition~\ref{pro:main} (two starters) \\
$[6, 5, 3, 2, 2, 2]$ & Proposition~\ref{pro:main} (two starters) \\
$[6, 4, 4, 2, 2, 2]$ & Proposition~\ref{pro:main} (two starters) \\
$[5, 5, 4, 2, 2, 2]$ & Proposition~\ref{pro:main} (two starters) \\
$[7, 3, 3, 3, 2, 2]$ & Proposition~\ref{pro:main1} (one starter) \\
$[6, 4, 3, 3, 2, 2]$ & Proposition~\ref{pro:main1} (one starter) \\
$[5, 5, 3, 3, 2, 2]$ & Proposition~\ref{pro:main1} (one starter) \\
$[5, 4, 4, 3, 2, 2]$ & Proposition~\ref{pro:main1} (one starter) \\
$[4, 4, 4, 4, 2, 2]$ & Proposition~\ref{pro:main1} (one starter) \\
$[6, 3, 3, 3, 3, 2]$ & Proposition~\ref{pro:main} (two starters) \\
$[5, 4, 3, 3, 3, 2]$ & Proposition~\ref{pro:main} (two starters) \\
$[4, 4, 4, 3, 3, 2]$ & Proposition~\ref{pro:main} (two starters) \\
$[5, 3, 3, 3, 3, 3]$ & Proposition~\ref{pro:main1} (one starter) \\
$[4, 4, 3, 3, 3, 3]$ & Proposition~\ref{pro:main1} (one starter) \\
$[12, 2, 2, 2, 2]$ & Proposition~\ref{pro:main1} (one starter) \\
$[11, 3, 2, 2, 2]$ & Proposition~\ref{pro:main} (two starters) \\
$[10, 4, 2, 2, 2]$ & Proposition~\ref{pro:main} (two starters) \\
$[9, 5, 2, 2, 2]$ & Proposition~\ref{pro:main} (two starters) \\
$[8, 6, 2, 2, 2]$ & Proposition~\ref{pro:main} (two starters) \\
$[7, 7, 2, 2, 2]$ & Proposition~\ref{pro:main} (two starters) \\
$[10, 3, 3, 2, 2]$ & Proposition~\ref{pro:main1} (one starter) \\
$[9, 4, 3, 2, 2]$ & Proposition~\ref{pro:main1} (one starter) \\
$[8, 5, 3, 2, 2]$ & Proposition~\ref{pro:main1} (one starter) \\
$[7, 6, 3, 2, 2]$ & Proposition~\ref{pro:main1} (one starter) \\
$[8, 4, 4, 2, 2]$ & Proposition~\ref{pro:main1} (one starter) \\
$[7, 5, 4, 2, 2]$ & Proposition~\ref{pro:main1} (one starter) \\
$[6, 6, 4, 2, 2]$ & Proposition~\ref{pro:main1} (one starter) \\
$[6, 5, 5, 2, 2]$ & Proposition~\ref{pro:main1} (one starter) \\
$[9, 3, 3, 3, 2]$ & Proposition~\ref{pro:main} (two starters) \\
$[8, 4, 3, 3, 2]$ & Proposition~\ref{pro:main} (two starters) \\
$[7, 5, 3, 3, 2]$ & Proposition~\ref{pro:main} (two starters) \\
$[6, 6, 3, 3, 2]$ & Proposition~\ref{pro:main} (two starters) \\
$[7, 4, 4, 3, 2]$ & Proposition~\ref{pro:main} (two starters) \\
$[6, 5, 4, 3, 2]$ & Proposition~\ref{pro:main} (two starters) \\
$[5, 5, 5, 3, 2]$ & Proposition~\ref{pro:main} (two starters) \\
$[6, 4, 4, 4, 2]$ & Proposition~\ref{pro:main} (two starters) \\
$[5, 5, 4, 4, 2]$ & Proposition~\ref{pro:main} (two starters) \\
\hline
\end{tabular}

\begin{tabular}{|l|l|}
\hline
2-factor type $[m_1,\ldots,m_t]$ & HOP$(2m_1,\ldots,2m_t)$ has a solution by... \\ \hline
$[8, 3, 3, 3, 3]$ & Proposition~\ref{pro:main1} (one starter) \\
$[7, 4, 3, 3, 3]$ & Proposition~\ref{pro:main1} (one starter) \\
$[6, 5, 3, 3, 3]$ & Proposition~\ref{pro:main1} (one starter) \\
$[6, 4, 4, 3, 3]$ & Proposition~\ref{pro:main1} (one starter) \\
$[5, 5, 4, 3, 3]$ & Proposition~\ref{pro:main1} (one starter) \\
$[5, 4, 4, 4, 3]$ & Proposition~\ref{pro:main1} (one starter) \\
$[4, 4, 4, 4, 4]$ & Theorem~\ref{the:intro2} \\
$[14, 2, 2, 2]$ & Proposition~\ref{pro:main} (two starters) \\
$[13, 3, 2, 2]$ & Proposition~\ref{pro:main1} (one starter) \\
$[12, 4, 2, 2]$ & Proposition~\ref{pro:main1} (one starter) \\
$[11, 5, 2, 2]$ & Proposition~\ref{pro:main1} (one starter) \\
$[10, 6, 2, 2]$ & Proposition~\ref{pro:main1} (one starter) \\
$[9, 7, 2, 2]$ & Proposition~\ref{pro:main1} (one starter) \\
$[8, 8, 2, 2]$ & Proposition~\ref{pro:main1} (one starter) \\
$[12, 3, 3, 2]$ & Proposition~\ref{pro:main} (two starters) \\
$[11, 4, 3, 2]$ & Proposition~\ref{pro:main} (two starters) \\
$[10, 5, 3, 2]$ & Proposition~\ref{pro:main} (two starters) \\
$[9, 6, 3, 2]$ & Proposition~\ref{pro:main} (two starters) \\
$[8, 7, 3, 2]$ & Proposition~\ref{pro:main} (two starters) \\
$[10, 4, 4, 2]$ & Proposition~\ref{pro:main} (two starters) \\
$[9, 5, 4, 2]$ & Proposition~\ref{pro:main} (two starters) \\
$[8, 6, 4, 2]$ & Proposition~\ref{pro:main} (two starters) \\
$[7, 7, 4, 2]$ & Proposition~\ref{pro:main} (two starters) \\
$[8, 5, 5, 2]$ & Proposition~\ref{pro:main} (two starters) \\
$[7, 6, 5, 2]$ & Proposition~\ref{pro:main} (two starters) \\
$[6, 6, 6, 2]$ & Proposition~\ref{pro:main} (two starters) \\
$[11, 3, 3, 3]$ & Proposition~\ref{pro:main1} (one starter) \\
$[10, 4, 3, 3]$ & Proposition~\ref{pro:main1} (one starter) \\
$[9, 5, 3, 3]$ & Proposition~\ref{pro:main1} (one starter) \\
$[8, 6, 3, 3]$ & Proposition~\ref{pro:main1} (one starter) \\
$[7, 7, 3, 3]$ & Proposition~\ref{pro:main1} (one starter) \\
$[9, 4, 4, 3]$ & Proposition~\ref{pro:main1} (one starter) \\
$[8, 5, 4, 3]$ & Proposition~\ref{pro:main1} (one starter) \\
$[7, 6, 4, 3]$ & Proposition~\ref{pro:main1} (one starter) \\
$[7, 5, 5, 3]$ & Proposition~\ref{pro:main1} (one starter) \\
$[6, 6, 5, 3]$ & Proposition~\ref{pro:main1} (one starter) \\
$[8, 4, 4, 4]$ & Theorem~\ref{the:intro2} \\
$[7, 5, 4, 4]$ & Proposition~\ref{pro:main1} (one starter) \\
$[6, 6, 4, 4]$ & Proposition~\ref{pro:main1} (one starter) \\
$[6, 5, 5, 4]$ & Proposition~\ref{pro:main1} (one starter) \\
$[5, 5, 5, 5]$ & Theorem~\ref{the:intro2} \\
$[16, 2, 2]$ & Proposition~\ref{pro:main1} (one starter) \\
$[15, 3, 2]$ & Proposition~\ref{pro:main} (two starters) \\
$[14, 4, 2]$ & Proposition~\ref{pro:main} (two starters) \\
$[13, 5, 2]$ & Proposition~\ref{pro:main} (two starters) \\
$[12, 6, 2]$ & Proposition~\ref{pro:main} (two starters) \\
\hline
\end{tabular}

\begin{tabular}{|l|l|}
\hline
2-factor type $[m_1,\ldots,m_t]$ & HOP$(2m_1,\ldots,2m_t)$ has a solution by... \\ \hline
$[11, 7, 2]$ & Proposition~\ref{pro:main} (two starters) \\
$[10, 8, 2]$ & Proposition~\ref{pro:main} (two starters) \\
$[9, 9, 2]$ & Proposition~\ref{pro:main} (two starters) \\
$[14, 3, 3]$ & Proposition~\ref{pro:main1} (one starter) \\
$[13, 4, 3]$ & Proposition~\ref{pro:main1} (one starter) \\
$[12, 5, 3]$ & Proposition~\ref{pro:main1} (one starter) \\
$[11, 6, 3]$ & Proposition~\ref{pro:main1} (one starter) \\
$[10, 7, 3]$ & Proposition~\ref{pro:main1} (one starter) \\
$[9, 8, 3]$ & Proposition~\ref{pro:main1} (one starter) \\
$[12, 4, 4]$ & Theorem~\ref{the:intro2} \\
$[11, 5, 4]$ & Proposition~\ref{pro:main1} (one starter) \\
$[10, 6, 4]$ & Proposition~\ref{pro:main1} (one starter) \\
$[9, 7, 4]$ & Proposition~\ref{pro:main1} (one starter) \\
$[8, 8, 4]$ & Theorem~\ref{the:intro2} \\
$[10, 5, 5]$ & Proposition~\ref{pro:main1} (one starter) \\
$[9, 6, 5]$ & Proposition~\ref{pro:main1} (one starter) \\
$[8, 7, 5]$ & Proposition~\ref{pro:main1} (one starter) \\
$[8, 6, 6]$ & Proposition~\ref{pro:main1} (one starter) \\
$[7, 7, 6]$ & Proposition~\ref{pro:main1} (one starter) \\
$[18, 2]$ & Proposition~\ref{pro:main} (two starters) \\
$[17, 3]$ & Proposition~\ref{pro:main1} (one starter) \\
$[16, 4]$ & Theorem~\ref{the:intro2} \\
$[15, 5]$ & Proposition~\ref{pro:main1} (one starter) \\
$[14, 6]$ & Proposition~\ref{pro:main1} (one starter) \\
$[13, 7]$ & Proposition~\ref{pro:main1} (one starter) \\
$[12, 8]$ & Theorem~\ref{the:intro2} \\
$[11, 9]$ & Proposition~\ref{pro:main1} (one starter) \\
$[10, 10]$ & Theorem~\ref{the:intro2} \\
$[20]$ & Theorem~\ref{the:intro2} \\
\hline
\end{tabular}

\end{proof}

\appendix

\section{Guide to computational results}

In cases where Proposition~\ref{pro:main1}, \ref{pro:main}, or \ref{pro:main2} is used to prove existence of a solution to HOP$(2m_1,\ldots,2m_t)$, the remaining sections provide appropriate starter 2-factors.
\begin{itemize}
\item Where Proposition~\ref{pro:main1} is used: one starter, namely $F=\{ C_0,C_1,\ldots,C_{t-1} \}$.
\item Where Proposition~\ref{pro:main} is used: two starters, namely $F_1=\{ C_0,C_1,\ldots,C_{t-1} \}$ and $F_2=\{ C_0',C_1',\ldots,C_{t-1}' \}$.
\item Where Proposition~\ref{pro:main2} is used: two starters, namely $F_1=\{ C_0,C_1,\ldots,C_{t-1} \}$ and $F_3=\{ C_0',C_1',\ldots,C_{t-1}' \}$. Note that in all cases, $F_1$ contains a 2-cycle consisting of a pink and blue edge of difference $\frac{n-1}{2}$.
    The starter 2-factor $F_2$ (as denoted in Proposition~\ref{pro:main2}) is obtained from $F_1$ by applying $\rho_{\n}^{\frac{n-1}{2}}$ and then replacing the 2-cycle containing a pink and blue edge of difference $\frac{n-1}{2}$ with a directed 2-cycle (with the same vertices) containing black arcs of length $\frac{n-1}{2}$.
\end{itemize}
The vertices of $2K_n^\t$ or $4K_n^\n$ are labelled $0,1,\ldots,n-1$; label $n-1$ represents vertex $x_{\infty}$. Edge orbits are given in the form $[d,c]$ where $d \in \{1,2,\ldots, \lfloor \frac{n-1}{2} \rfloor \} \cup \{ n-1 \}$ is the difference (with $n-1$ representing the infinity difference), and $c$ is the colour-orientation: 0 for pink, 2 for blue, 1 for black (single starter cases) or black orientated forward (two- and three-starter cases), and -1 for black oriented backward. We consider an arc $(u,u+d)$, for $u,d\in \ZZ_{n-1}$,  oriented forward if $d < \frac{n-1}{2}$, and oriented backward if  $d > \frac{n-1}{2}$. The case $d=\frac{n-1}{2}$ can occur only in the three-starter approach (that is, for odd $n$); however, in this case, black arcs of difference $\frac{n-1}{2}$ appear only in the starter 2-factor $F_2$, which is not specified below.

Each $t$-cycle in a starter 2-factor is presented in the form
$$C=[ v_0, [d_0,c_0], v_1, [d_1,c_1], v_2,\ldots, v_{t-2}, [d_{t-2},c_{t-2}],v_{t-1}, [d_{t-1},c_{t-1}]],$$
where $v_0,v_1,\ldots,v_{t-1}$ are the consecutive vertices of the cycle, and $[d_i,c_i]$ represents the orbit of the edge $v_iv_{i+1}$, for $i \in \ZZ_t$, as explained above.


\small

\section{Computational results for $n=10$}\label{app:10}

\begin{itemizenew}
\item 2-factor type $[4, 2, 2, 2]$: one starter \vspace*{-3mm} \\
\hspace*{-10mm}
\begin{minipage}{8cm}
\begin{align*}
C_0 =& \; [6, [4, 1], 2, [4, 0], 7, [9, 1], 9, [9, 0]] \\
C_1 =& \; [5, [2, 0], 3, [2, 1]] \\
C_2 =& \; [4, [3, 0], 1, [3, 1]] \\
C_3 =& \; [0, [1, 0], 8, [1, 1]]
\end{align*}
\end{minipage}

\item 2-factor type $[3, 3, 2, 2]$: two starters \vspace*{-3mm} \\
\hspace*{-17mm}
\begin{minipage}{8cm}
\begin{align*}
C_0 =& \; [1, [9, 2], 9, [9, 0], 0, [1, 1]] \\
C_1 =& \; [8, [4, 1], 3, [4, 2], 7, [1, 0]] \\
C_2 =& \; [6, [2, -1], 4, [2, 1]] \\
C_3 =& \; [5, [3, 0], 2, [3, 2]]
\end{align*}
\end{minipage}
\begin{minipage}{8cm}
\begin{align*}
C_0' =& \; [9, [9, 1], 5, [1, -1], 6, [9, -1]] \\
C_1' =& \; [4, [4, -1], 8, [4, 0], 3, [1, 2]] \\
C_2' =& \; [2, [2, 0], 0, [2, 2]] \\
C_3' =& \; [7, [3, -1], 1, [3, 1]]
\end{align*}
\end{minipage}

\item 2-factor type $[6, 2, 2]$: two starters \vspace*{-3mm} \\
\begin{minipage}{8cm}
\begin{align*}
C_0 =& \; [9, [9, 2], 3, [3, 0], 0, [2, -1], 7, [4, 1], 2, [2, 1], 4, [9, 1]] \\
C_1 =& \; [6, [4, 0], 1, [4, 2]] \\
C_2 =& \; [8, [3, -1], 5, [3, 1]]
\end{align*}
\end{minipage}

\vspace*{-4mm}
\begin{minipage}{8cm}
\begin{align*}
C_0' =& \; [5, [4, -1], 0, [9, -1], 9, [9, 0], 4, [1, 2], 3, [1, 0], 2, [3, 2]] \\
C_1' =& \; [7, [1, -1], 6, [1, 1]] \\
C_2' =& \; [1, [2, 0], 8, [2, 2]]
\end{align*}
\end{minipage}

\item 2-factor type $[5, 3, 2]$: one starter \vspace*{-3mm} \\
\hspace*{-4mm}
\begin{minipage}{8cm}
\begin{align*}
C_0 =& \; [3, [2, 0], 5, [1, 0], 4, [4, 0], 8, [2, 1], 6, [3, 0]] \\
C_1 =& \; [1, [3, 1], 7, [4, 1], 2, [1, 1]] \\
C_2 =& \; [0, [9, 0], 9, [9, 1]]
\end{align*}
\end{minipage}

\item 2-factor type $[4, 4, 2]$: one starter \vspace*{-3mm} \\
\hspace*{-10mm}
\begin{minipage}{8cm}
\begin{align*}
C_0 =& \; [5, [2, 1], 7, [3, 0], 1, [1, 1], 0, [4, 0]] \\
C_1 =& \; [3, [1, 0], 4, [2, 0], 2, [3, 1], 8, [4, 1]] \\
C_2 =& \; [9, [9, 0], 6, [9, 1]]
\end{align*}
\end{minipage}

\item 2-factor type $[4, 3, 3]$: two starters \vspace*{-3mm} \\
\hspace*{-10mm}
\begin{minipage}{8cm}
\begin{align*}
C_0 =& \; [0, [1, 2], 1, [2, 1], 8, [2, 0], 6, [3, 1]] \\
C_1 =& \; [3, [9, 1], 9, [9, -1], 7, [4, -1]] \\
C_2 =& \; [5, [1, -1], 4, [2, 2], 2, [3, 0]]
\end{align*}
\end{minipage}
\begin{minipage}{8cm}
\begin{align*}
C_0' =& \; [6, [2, -1], 4, [9, 2], 9, [9, 0], 0, [3, -1]] \\
C_1' =& \; [8, [4, 1], 3, [4, 2], 7, [1, 0]] \\
C_2' =& \; [1, [4, 0], 5, [3, 2], 2, [1, 1]]
\end{align*}
\end{minipage}

\item 2-factor type $[8, 2]$: one starter \vspace*{-3mm} \\
\begin{minipage}{8cm}
\begin{align*}
C_0 =& \; [0, [3, 0], 3, [4, 0], 7, [9, 1], 9, [9, 0], 1, [4, 1], 6, [2, 1], 8, [3, 1], 2, [2, 0]] \\
C_1 =& \; [4, [1, 0], 5, [1, 1]]
\end{align*}
\end{minipage}

\item 2-factor type $[7, 3]$: two starters \vspace*{-3mm} \\
\begin{minipage}{8cm}
\begin{align*}
C_0 =& \; [1, [4, -1], 5, [4, 0], 0, [9, 1], 9, [9, 2], 2, [3, 1], 8, [2, 0], 6, [4, 2]] \\
C_1 =& \; [3, [4, 1], 7, [3, -1], 4, [1, -1]]
\end{align*}
\end{minipage}

\vspace*{-4mm}
\begin{minipage}{8cm}
\begin{align*}
C_0' =& \; [2, [2, -1], 0, [1, 1], 1, [3, 2], 7, [3, 0], 4, [1, 2], 5, [2, 1], 3, [1, 0]] \\
C_1' =& \; [8, [9, -1], 9, [9, 0], 6, [2, 2]]
\end{align*}
\end{minipage}

\item 2-factor type $[6, 4]$: two starters \vspace*{-3mm} \\
\begin{minipage}{8cm}
\begin{align*}
C_0 =& \; [0, [1, 0], 1, [4, 1], 5, [3, 2], 2, [3, 0], 8, [4, 2], 3, [3, 1]] \\
C_1 =& \; [6, [9, -1], 9, [9, 0], 7, [3, -1], 4, [2, 2]]
\end{align*}
\end{minipage}

\vspace*{-4mm}
\begin{minipage}{8cm}
\begin{align*}
C_0' =& \; [0, [4, -1], 4, [1, 1], 3, [2, 1], 1, [2, 0], 8, [1, 2], 7, [2, -1]] \\
C_1' =& \; [2, [4, 0], 6, [1, -1], 5, [9, 2], 9, [9, 1]]
\end{align*}
\end{minipage}

\end{itemizenew}

\section{Computational results for $n=11$}\label{app:11}

\begin{itemizenew}
\item 2-factor type $[3, 2, 2, 2, 2]$: three starters \vspace*{-3mm} \\
\hspace*{-12mm}
\begin{minipage}{8cm}
\begin{align*}
C_0 =& \; [7, [10, 1], 10, [10, -1], 6, [1, 1]] \\
C_1 =& \; [3, [1, 0], 4, [1, 2]] \\
C_2 =& \; [5, [3, 0], 8, [3, 2]] \\
C_3 =& \; [9, [2, -1], 1, [2, 1]] \\
C_4 =& \; [0, [2, 0], 2, [2, 2]]
\end{align*}
\end{minipage}
\begin{minipage}{8cm}
\begin{align*}
C_0' =& \; [1, [1, -1], 2, [10, 0], 10, [10, 2]] \\
C_1' =& \; [0, [5, 0], 5, [5, 2]] \\
C_2' =& \; [3, [4, -1], 7, [4, 1]] \\
C_3' =& \; [4, [4, 0], 8, [4, 2]] \\
C_4' =& \; [9, [3, -1], 6, [3, 1]]
\end{align*}
\end{minipage}

\newpage

\item 2-factor type $[5, 2, 2, 2]$: three starters \vspace*{-3mm} \\
\begin{minipage}{8cm}
\begin{align*}
C_0 =& \; [9, [4, 0], 5, [1, -1], 4, [4, 2], 8, [1, 1], 7, [2, -1]] \\
C_1 =& \; [6, [3, -1], 3, [3, 1]] \\
C_2 =& \; [0, [2, 0], 2, [2, 2]] \\
C_3 =& \; [10, [10, -1], 1, [10, 1]]
\end{align*}
\end{minipage}

\vspace*{-4mm}
\hspace*{-2mm}
\begin{minipage}{8cm}
\begin{align*}
C_0' =& \; [8, [4, 1], 4, [2, 1], 2, [1, 0], 1, [4, -1], 7, [1, 2]] \\
C_1' =& \; [0, [5, 0], 5, [5, 2]] \\
C_2' =& \; [9, [3, 0], 6, [3, 2]] \\
C_3' =& \; [3, [10, 0], 10, [10, 2]]
\end{align*}
\end{minipage}

\item 2-factor type $[4, 3, 2, 2]$: three starters \vspace*{-3mm} \\
\hspace*{-8mm}
\begin{minipage}{8cm}
\begin{align*}
C_0 =& \; [7, [2, 1], 5, [2, 0], 3, [4, -1], 9, [2, 2]] \\
C_1 =& \; [0, [2, -1], 8, [4, 2], 4, [4, 0]] \\
C_2 =& \; [1, [1, -1], 2, [1, 1]] \\
C_3 =& \; [10, [10, -1], 6, [10, 1]]
\end{align*}
\end{minipage}
\begin{minipage}{8cm}
\begin{align*}
C_0' =& \; [6, [3, 0], 9, [3, 1], 2, [1, 2], 3, [3, -1]] \\
C_1' =& \; [1, [4, 1], 7, [1, 0], 8, [3, 2]] \\
C_2' =& \; [0, [5, 0], 5, [5, 2]] \\
C_3' =& \; [10, [10, 0], 4, [10, 2]]
\end{align*}
\end{minipage}

\item 2-factor type $[3, 3, 3, 2]$: three starters \vspace*{-3mm} \\
\hspace*{-15mm}
\begin{minipage}{8cm}
\begin{align*}
C_0 =& \; [4, [2, -1], 6, [4, 0], 0, [4, 2]] \\
C_1 =& \; [5, [2, 2], 7, [1, 0], 8, [3, -1]] \\
C_2 =& \; [3, [2, 0], 1, [1, 2], 2, [1, -1]] \\
C_3 =& \; [10, [10, -1], 9, [10, 1]]
\end{align*}
\end{minipage}
\begin{minipage}{8cm}
\begin{align*}
C_0' =& \; [8, [1, 1], 7, [3, 1], 4, [4, -1]] \\
C_1' =& \; [3, [2, 1], 1, [10, 0], 10, [10, 2]] \\
C_2' =& \; [9, [3, 2], 6, [4, 1], 2, [3, 0]] \\
C_3' =& \; [0, [5, 0], 5, [5, 2]]
\end{align*}
\end{minipage}

\item 2-factor type $[7, 2, 2]$: three starters \vspace*{-3mm} \\
\begin{minipage}{8cm}
\begin{align*}
C_0 =& \; [4, [3, 0], 1, [1, 2], 0, [3, 1], 7, [2, -1], 9, [4, 0], 3, [2, 2], 5, [1, 1]] \\
C_1 =& \; [2, [4, -1], 8, [4, 1]] \\
C_2 =& \; [10, [10, 0], 6, [10, 2]]
\end{align*}
\end{minipage}

\vspace*{-4mm}
\begin{minipage}{8cm}
\begin{align*}
C_0' =& \; [8, [2, 0], 6, [3, 2], 3, [1, 0], 2, [3, -1], 9, [2, 1], 1, [4, 2], 7, [1, -1]] \\
C_1' =& \; [0, [5, 0], 5, [5, 2]] \\
C_2' =& \; [4, [10, -1], 10, [10, 1]]
\end{align*}
\end{minipage}

\item 2-factor type $[6, 3, 2]$: three starters \vspace*{-3mm} \\
\begin{minipage}{8cm}
\begin{align*}
C_0 =& \; [9, [4, -1], 5, [10, 1], 10, [10, 2], 7, [1, 1], 6, [3, 1], 3, [4, 0]] \\
C_1 =& \; [2, [2, 1], 4, [3, 2], 1, [1, 0]] \\
C_2 =& \; [8, [2, 0], 0, [2, 2]]
\end{align*}
\end{minipage}

\vspace*{-4mm}
\begin{minipage}{8cm}
\begin{align*}
C_0' =& \; [6, [2, -1], 4, [3, -1], 1, [1, 2], 2, [4, 1], 8, [1, -1], 9, [3, 0]] \\
C_1' =& \; [7, [10, 0], 10, [10, -1], 3, [4, 2]] \\
C_2' =& \; [0, [5, 0], 5, [5, 2]]
\end{align*}
\end{minipage}

\item 2-factor type $[5, 4, 2]$: three starters \vspace*{-3mm} \\
\begin{minipage}{8cm}
\begin{align*}
C_0 =& \; [0, [1, -1], 1, [3, 1], 8, [3, 0], 5, [1, 2], 6, [4, -1]] \\
C_1 =& \; [10, [10, -1], 7, [2, 2], 9, [4, 0], 3, [10, 1]] \\
C_2 =& \; [2, [2, -1], 4, [2, 1]]
\end{align*}
\end{minipage}
\hspace{2mm}
\begin{minipage}{8cm}
\begin{align*}
C_0' =& \; [9, [3, 2], 2, [10, 0], 10, [10, 2], 6, [2, 0], 8, [1, 1]] \\
C_1' =& \; [4, [1, 0], 3, [4, 1], 7, [4, 2], 1, [3, -1]] \\
C_2' =& \; [0, [5, 0], 5, [5, 2]]
\end{align*}
\end{minipage}

\end{itemizenew}


\section{Computational results for $n=12$}\label{app:12}

\begin{itemizenew}
\item 2-factor type $[4, 2, 2, 2, 2]$: one starter \vspace*{-3mm} \\
\hspace*{-10mm}
\begin{minipage}{8cm}
\begin{align*}
C_0 =& \; [6, [1, 0], 5, [3, 0], 2, [1, 1], 3, [3, 1]] \\
C_1 =& \; [11, [11, 0], 0, [11, 1]] \\
C_2 =& \; [1, [4, 0], 8, [4, 1]] \\
C_3 =& \; [7, [2, 0], 9, [2, 1]] \\
C_4 =& \; [10, [5, 0], 4, [5, 1]]
\end{align*}
\end{minipage}

\item 2-factor type $[3, 3, 2, 2, 2]$: two starters \vspace*{-3mm} \\
\hspace*{-12mm}
\begin{minipage}{8cm}
\begin{align*}
C_0 =& \; [3, [1, 1], 4, [11, 1], 11, [11, -1]] \\
C_1 =& \; [0, [5, -1], 5, [4, -1], 9, [2, -1]] \\
C_2 =& \; [6, [2, 0], 8, [2, 2]] \\
C_3 =& \; [1, [1, 0], 2, [1, 2]] \\
C_4 =& \; [10, [3, -1], 7, [3, 1]]
\end{align*}
\end{minipage}
\begin{minipage}{8cm}
\begin{align*}
C_0' =& \; [6, [1, -1], 5, [5, 2], 0, [5, 0]] \\
C_1' =& \; [10, [5, 1], 4, [4, 1], 8, [2, 1]] \\
C_2' =& \; [9, [3, 0], 1, [3, 2]] \\
C_3' =& \; [11, [11, 0], 2, [11, 2]] \\
C_4' =& \; [7, [4, 0], 3, [4, 2]]
\end{align*}
\end{minipage}

\item 2-factor type $[6, 2, 2, 2]$: two starters \vspace*{-3mm} \\
\begin{minipage}{8cm}
\begin{align*}
C_0 =& \; [11, [11, 2], 4, [1, 0], 5, [5, 1], 10, [2, 2], 1, [4, 1], 8, [11, 0]] \\
C_1 =& \; [3, [3, 0], 6, [3, 2]] \\
C_2 =& \; [2, [5, 0], 7, [5, 2]] \\
C_3 =& \; [9, [2, -1], 0, [2, 1]]
\end{align*}
\end{minipage}

\vspace*{-4mm}
\begin{minipage}{8cm}
\begin{align*}
C_0' =& \; [10, [1, 2], 0, [5, -1], 5, [11, -1], 11, [11, 1], 1, [2, 0], 3, [4, -1]] \\
C_1' =& \; [9, [1, -1], 8, [1, 1]] \\
C_2' =& \; [6, [4, 0], 2, [4, 2]] \\
C_3' =& \; [7, [3, -1], 4, [3, 1]]
\end{align*}
\end{minipage}

\item 2-factor type $[5, 3, 2, 2]$: one starter \vspace*{-3mm} \\
\hspace*{-2mm}
\begin{minipage}{8cm}
\begin{align*}
C_0 =& \; [4, [2, 1], 2, [5, 1], 7, [2, 0], 9, [1, 1], 10, [5, 0]] \\
C_1 =& \; [6, [11, 1], 11, [11, 0], 5, [1, 0]] \\
C_2 =& \; [3, [3, 0], 0, [3, 1]] \\
C_3 =& \; [1, [4, 0], 8, [4, 1]]
\end{align*}
\end{minipage}

\item 2-factor type $[4, 4, 2, 2]$: one starter \vspace*{-3mm} \\
\hspace*{-6mm}
\begin{minipage}{8cm}
\begin{align*}
C_0 =& \; [11, [11, 0], 2, [3, 0], 5, [1, 1], 4, [11, 1]] \\
C_1 =& \; [0, [3, 1], 8, [2, 1], 10, [1, 0], 9, [2, 0]] \\
C_2 =& \; [3, [4, 0], 7, [4, 1]] \\
C_3 =& \; [1, [5, 0], 6, [5, 1]]
\end{align*}
\end{minipage}

\item 2-factor type $[4, 3, 3, 2]$: two starters \vspace*{-3mm} \\
\hspace*{-8mm}
\begin{minipage}{8cm}
\begin{align*}
C_0 =& \; [2, [5, 0], 8, [1, 1], 9, [3, -1], 6, [4, 2]] \\
C_1 =& \; [11, [11, 2], 7, [3, 1], 4, [11, 0]] \\
C_2 =& \; [10, [4, -1], 3, [2, 0], 1, [2, 2]] \\
C_3 =& \; [5, [5, -1], 0, [5, 1]]
\end{align*}
\end{minipage}
\begin{minipage}{8cm}
\begin{align*}
C_0' =& \; [9, [4, 0], 2, [3, 2], 10, [2, 1], 8, [1, -1]] \\
C_1' =& \; [3, [3, 0], 6, [5, 2], 1, [2, -1]] \\
C_2' =& \; [7, [4, 1], 0, [11, 1], 11, [11, -1]] \\
C_3' =& \; [5, [1, 0], 4, [1, 2]]
\end{align*}
\end{minipage}

\item 2-factor type $[8, 2, 2]$: one starter \vspace*{-3mm} \\
\begin{minipage}{8cm}
\begin{align*}
C_0 =& \; [1, [1, 1], 0, [4, 1], 4, [2, 1], 6, [4, 0], 10, [3, 1], 7, [2, 0], 5, [3, 0], 2, [1, 0]] \\
C_1 =& \; [8, [5, 0], 3, [5, 1]] \\
C_2 =& \; [9, [11, 0], 11, [11, 1]]
\end{align*}
\end{minipage}

\item 2-factor type $[7, 3, 2]$: two starters \vspace*{-3mm} \\
\begin{minipage}{8cm}
\begin{align*}
C_0 =& \; [5, [11, 1], 11, [11, -1], 3, [5, 2], 9, [1, 1], 8, [2, 1], 6, [5, 0], 1, [4, 1]] \\
C_1 =& \; [4, [3, 2], 7, [5, 1], 2, [2, 0]] \\
C_2 =& \; [0, [1, 0], 10, [1, 2]]
\end{align*}
\end{minipage}

\vspace*{-4mm}
\begin{minipage}{8cm}
\begin{align*}
C_0' =& \; [3, [2, -1], 1, [4, 2], 5, [1, -1], 6, [4, -1], 10, [3, 0], 7, [11, 2], 11, [11, 0]] \\
C_1' =& \; [4, [5, -1], 9, [4, 0], 2, [2, 2]] \\
C_2' =& \; [0, [3, -1], 8, [3, 1]]
\end{align*}
\end{minipage}

\item 2-factor type $[6, 4, 2]$: two starters \vspace*{-3mm} \\
\begin{minipage}{8cm}
\begin{align*}
C_0 =& \; [0, [1, -1], 1, [3, -1], 4, [2, 1], 2, [5, 0], 7, [4, 2], 3, [3, 1]] \\
C_1 =& \; [9, [4, 0], 5, [5, 1], 10, [4, -1], 6, [3, 2]] \\
C_2 =& \; [8, [11, 0], 11, [11, 2]]
\end{align*}
\end{minipage}

\vspace*{-4mm}
\begin{minipage}{8cm}
\begin{align*}
C_0' =& \; [10, [1, 1], 0, [1, 2], 1, [11, -1], 11, [11, 1], 4, [5, -1], 9, [1, 0]] \\
C_1' =& \; [3, [4, 1], 7, [5, 2], 2, [3, 0], 5, [2, -1]] \\
C_2' =& \; [6, [2, 0], 8, [2, 2]]
\end{align*}
\end{minipage}

\item 2-factor type $[5, 5, 2]$: two starters \vspace*{-3mm} \\
\begin{minipage}{8cm}
\begin{align*}
C_0 =& \; [11, [11, 2], 9, [1, 0], 8, [4, 2], 1, [1, 1], 0, [11, 0]] \\
C_1 =& \; [5, [1, -1], 6, [4, 0], 10, [3, 2], 7, [4, 1], 3, [2, -1]] \\
C_2 =& \; [4, [2, 0], 2, [2, 2]]
\end{align*}
\end{minipage}

\vspace*{-4mm}
\begin{minipage}{8cm}
\begin{align*}
C_0' =& \; [8, [4, -1], 1, [2, 1], 10, [3, -1], 2, [11, -1], 11, [11, 1]] \\
C_1' =& \; [4, [3, 1], 7, [1, 2], 6, [3, 0], 3, [5, 2], 9, [5, 0]] \\
C_2' =& \; [5, [5, -1], 0, [5, 1]]
\end{align*}
\end{minipage}

\item 2-factor type $[6, 3, 3]$: one starter \vspace*{-3mm} \\
\begin{minipage}{8cm}
\begin{align*}
C_0 =& \; [9, [3, 0], 1, [1, 1], 2, [1, 0], 3, [4, 1], 7, [3, 1], 4, [5, 1]] \\
C_1 =& \; [0, [5, 0], 5, [11, 0], 11, [11, 1]] \\
C_2 =& \; [6, [2, 0], 8, [2, 1], 10, [4, 0]]
\end{align*}
\end{minipage}

\item 2-factor type $[5, 4, 3]$: one starter \vspace*{-3mm} \\
\hspace*{-4mm}
\begin{minipage}{8cm}
\begin{align*}
C_0 =& \; [1, [2, 0], 3, [1, 1], 2, [4, 0], 9, [5, 1], 4, [3, 1]] \\
C_1 =& \; [10, [3, 0], 7, [1, 0], 8, [2, 1], 6, [4, 1]] \\
C_2 =& \; [11, [11, 0], 0, [5, 0], 5, [11, 1]]
\end{align*}
\end{minipage}

\newpage

\item 2-factor type $[10, 2]$: two starters \vspace*{-3mm} \\
\begin{minipage}{8cm}
\begin{align*}
C_0 =& \; [2, [1, 0], 1, [2, 1], 3, [1, 2], 4, [4, 0], 8, [11, 1], 11, [11, -1], 9, [3, -1], 6, [1, 1], 7, [4, 1], 0, [2, 2]] \\
C_1 =& \; [5, [5, -1], 10, [5, 1]]
\end{align*}
\end{minipage}

\vspace*{-4mm}
\begin{minipage}{8cm}
\begin{align*}
C_0' =& \; [9, [1, -1], 8, [5, 2], 3, [3, 0], 6, [4, 2], 2, [3, 1], 10, [2, -1], 1, [4, -1], 5, [2, 0], 7, [3, 2], 4, [5, 0]] \\
C_1' =& \; [0, [11, 0], 11, [11, 2]]
\end{align*}
\end{minipage}

\item 2-factor type $[9, 3]$: one starter \vspace*{-3mm} \\
\begin{minipage}{8cm}
\begin{align*}
C_0 =& \; [0, [3, 1], 8, [4, 0], 1, [4, 1], 5, [2, 1], 3, [3, 0], 6, [2, 0], 4, [5, 1], 9, [1, 1], 10, [1, 0]] \\
C_1 =& \; [2, [5, 0], 7, [11, 1], 11, [11, 0]]
\end{align*}
\end{minipage}

\item 2-factor type $[7, 5]$: one starter \vspace*{-3mm} \\
\begin{minipage}{8cm}
\begin{align*}
C_0 =& \; [10, [4, 0], 3, [1, 0], 2, [1, 1], 1, [5, 0], 6, [2, 1], 8, [3, 0], 5, [5, 1]] \\
C_1 =& \; [9, [2, 0], 7, [3, 1], 4, [4, 1], 0, [11, 1], 11, [11, 0]]
\end{align*}
\end{minipage}

\end{itemizenew}


\section{Computational results for $n=13$}\label{app:13}

\begin{itemizenew}
\item 2-factor type $[3, 2, 2, 2, 2, 2]$: three starters \vspace*{-3mm} \\
\hspace*{-14mm}
\begin{minipage}{8cm}
\begin{align*}
C_0 =& \; [10, [5, 2], 3, [3, -1], 6, [4, 0]] \\
C_1 =& \; [0, [5, -1], 5, [5, 1]] \\
C_2 =& \; [12, [12, 0], 11, [12, 2]] \\
C_3 =& \; [7, [2, 0], 9, [2, 2]] \\
C_4 =& \; [2, [1, 0], 1, [1, 2]] \\
C_5 =& \; [8, [4, -1], 4, [4, 1]]
\end{align*}
\end{minipage}
\begin{minipage}{8cm}
\begin{align*}
C_0' =& \; [1, [5, 0], 8, [4, 2], 4, [3, 1]] \\
C_1' =& \; [0, [6, 0], 6, [6, 2]] \\
C_2' =& \; [7, [12, -1], 12, [12, 1]] \\
C_3' =& \; [10, [1, -1], 9, [1, 1]] \\
C_4' =& \; [3, [2, -1], 5, [2, 1]] \\
C_5' =& \; [11, [3, 0], 2, [3, 2]]
\end{align*}
\end{minipage}

\item 2-factor type $[5, 2, 2, 2, 2]$: three starters \vspace*{-3mm} \\
\hspace*{-2mm}
\begin{minipage}{8cm}
\begin{align*}
C_0 =& \; [9, [5, 2], 2, [4, 1], 10, [2, 1], 8, [3, 1], 5, [4, 0]] \\
C_1 =& \; [1, [2, 0], 3, [2, 2]] \\
C_2 =& \; [0, [12, -1], 12, [12, 1]] \\
C_3 =& \; [7, [1, -1], 6, [1, 1]] \\
C_4 =& \; [4, [5, -1], 11, [5, 1]]
\end{align*}
\end{minipage}

\vspace*{-4mm}
\begin{minipage}{8cm}
\begin{align*}
C_0' =& \; [1, [2, -1], 3, [4, -1], 7, [5, 0], 2, [4, 2], 10, [3, -1]] \\
C_1' =& \; [0, [6, 0], 6, [6, 2]] \\
C_2' =& \; [11, [3, 0], 8, [3, 2]] \\
C_3' =& \; [9, [12, 0], 12, [12, 2]] \\
C_4' =& \; [5, [1, 0], 4, [1, 2]]
\end{align*}
\end{minipage}

\item 2-factor type $ [4, 3, 2, 2, 2]$: three starters \vspace*{-3mm} \\
\hspace*{-6mm}
\begin{minipage}{8cm}
\begin{align*}
C_0 =& \; [10, [4, 0], 6, [5, 2], 11, [2, -1], 1, [3, 1]] \\
C_1 =& \; [7, [4, 2], 3, [5, 0], 8, [1, -1]] \\
C_2 =& \; [9, [12, -1], 12, [12, 1]] \\
C_3 =& \; [2, [2, 0], 4, [2, 2]] \\
C_4 =& \; [0, [5, -1], 5, [5, 1]]
\end{align*}
\end{minipage}
\begin{minipage}{8cm}
\begin{align*}
C_0' =& \; [12, [12, 0], 10, [1, 1], 11, [3, -1], 8, [12, 2]] \\
C_1' =& \; [3, [1, 2], 2, [1, 0], 1, [2, 1]] \\
C_2' =& \; [0, [6, 0], 6, [6, 2]] \\
C_3' =& \; [7, [3, 0], 4, [3, 2]] \\
C_4' =& \; [9, [4, -1], 5, [4, 1]]
\end{align*}
\end{minipage}

\item 2-factor type $[3, 3, 3, 2, 2]$: three starters \vspace*{-3mm} \\
\hspace*{-11mm}
\begin{minipage}{8cm}
\begin{align*}
C_0 =& \; [8, [12, -1], 12, [12, 1], 7, [1, -1]] \\
C_1 =& \; [11, [2, 0], 9, [1, 1], 10, [1, 2]] \\
C_2 =& \; [2, [2, 2], 4, [4, 0], 0, [2, 1]] \\
C_3 =& \; [6, [3, 0], 3, [3, 2]] \\
C_4 =& \; [1, [4, -1], 5, [4, 1]]
\end{align*}
\end{minipage}
\begin{minipage}{8cm}
\begin{align*}
C_0' =& \; [9, [4, 2], 1, [3, 1], 10, [1, 0]] \\
C_1' =& \; [5, [2, -1], 3, [5, 1], 8, [3, -1]] \\
C_2' =& \; [4, [12, 0], 12, [12, 2], 11, [5, -1]] \\
C_3' =& \; [0, [6, 0], 6, [6, 2]] \\
C_4' =& \; [2, [5, 0], 7, [5, 2]]
\end{align*}
\end{minipage}

\item 2-factor type $[7, 2, 2, 2]$: three starters \vspace*{-3mm} \\
\begin{minipage}{8cm}
\begin{align*}
C_0 =& \; [5, [3, 2], 2, [12, 0], 12, [12, 2], 4, [5, -1], 9, [3, 1], 6, [5, 1], 1, [4, 0]] \\
C_1 =& \; [11, [4, -1], 3, [4, 1]] \\
C_2 =& \; [0, [2, 0], 10, [2, 2]] \\
C_3 =& \; [8, [1, 0], 7, [1, 2]]
\end{align*}
\end{minipage}

\vspace*{-4mm}
\begin{minipage}{8cm}
\begin{align*}
C_0' =& \; [2, [1, 1], 1, [5, 0], 8, [4, 2], 4, [3, -1], 7, [3, 0], 10, [1, -1], 9, [5, 2]] \\
C_1' =& \; [0, [6, 0], 6, [6, 2]] \\
C_2' =& \; [11, [12, -1], 12, [12, 1]] \\
C_3' =& \; [3, [2, -1], 5, [2, 1]]
\end{align*}
\end{minipage}

\item 2-factor type $[6, 3, 2, 2]$: three starters \vspace*{-3mm} \\
\begin{minipage}{8cm}
\begin{align*}
C_0 =& \; [1, [5, 2], 8, [1, 1], 7, [2, 0], 9, [2, 2], 11, [12, 0], 12, [12, -1]] \\
C_1 =& \; [4, [1, -1], 3, [3, 2], 0, [4, 0]] \\
C_2 =& \; [5, [5, -1], 10, [5, 1]] \\
C_3 =& \; [2, [4, -1], 6, [4, 1]]
\end{align*}
\end{minipage}

\vspace*{-4mm}
\begin{minipage}{8cm}
\begin{align*}
C_0' =& \; [1, [12, 2], 12, [12, 1], 9, [2, 1], 7, [5, 0], 2, [4, 2], 10, [3, 0]] \\
C_1' =& \; [3, [2, -1], 5, [1, 0], 4, [1, 2]] \\
C_2' =& \; [0, [6, 0], 6, [6, 2]] \\
C_3' =& \; [8, [3, -1], 11, [3, 1]]
\end{align*}
\end{minipage}

\item 2-factor type $[5, 4, 2, 2]$: three starters \vspace*{-3mm} \\
\begin{minipage}{8cm}
\begin{align*}
C_0 =& \; [3, [1, -1], 4, [5, -1], 9, [2, 0], 7, [1, 2], 6, [3, 1]] \\
C_1 =& \; [1, [1, 1], 0, [2, -1], 2, [12, 0], 12, [12, 2]] \\
C_2 =& \; [5, [5, 0], 10, [5, 2]] \\
C_3 =& \; [8, [3, 0], 11, [3, 2]]
\end{align*}
\end{minipage}

\vspace*{-4mm}
\hspace*{-2mm}
\begin{minipage}{8cm}
\begin{align*}
C_0' =& \; [5, [3, -1], 2, [5, 1], 7, [4, 2], 3, [2, 1], 1, [4, 0]] \\
C_1' =& \; [12, [12, -1], 11, [2, 2], 9, [1, 0], 10, [12, 1]] \\
C_2' =& \; [0, [6, 0], 6, [6, 2]] \\
C_3' =& \; [4, [4, -1], 8, [4, 1]]
\end{align*}
\end{minipage}

\item 2-factor type $[5, 3, 3, 2]$: three starters \vspace*{-3mm} \\
\begin{minipage}{8cm}
\begin{align*}
C_0 =& \; [12, [12, -1], 4, [3, 2], 1, [2, 1], 11, [2, 0], 9, [12, 1]] \\
C_1 =& \; [10, [2, 2], 8, [3, 0], 5, [5, 1]] \\
C_2 =& \; [6, [4, 1], 2, [1, -1], 3, [3, -1]] \\
C_3 =& \; [0, [5, 0], 7, [5, 2]]
\end{align*}
\end{minipage}

\newpage

\begin{minipage}{8cm}
\begin{align*}
C_0' =& \; [5, [5, -1], 10, [12, 0], 12, [12, 2], 1, [4, 0], 9, [4, 2]] \\
C_1' =& \; [11, [3, 1], 8, [1, 1], 7, [4, -1]] \\
C_2' =& \; [2, [2, -1], 4, [1, 0], 3, [1, 2]] \\
C_3' =& \; [0, [6, 0], 6, [6, 2]]
\end{align*}
\end{minipage}

\item 2-factor type $[4, 4, 3, 2]$: three starters \vspace*{-3mm} \\
\begin{minipage}{8cm}
\begin{align*}
C_0 =& \; [8, [4, 1], 0, [1, -1], 11, [5, 2], 4, [4, 0]] \\
C_1 =& \; [10, [12, -1], 12, [12, 1], 1, [4, -1], 5, [5, -1]] \\
C_2 =& \; [3, [4, 2], 7, [5, 0], 2, [1, 1]] \\
C_3 =& \; [6, [3, -1], 9, [3, 1]]
\end{align*}
\end{minipage}
\begin{minipage}{8cm}
\begin{align*}
C_0' =& \; [5, [1, 2], 4, [2, 0], 2, [5, 1], 7, [2, -1]] \\
C_1' =& \; [10, [12, 2], 12, [12, 0], 3, [2, 2], 1, [3, 0]] \\
C_2' =& \; [11, [3, 2], 8, [1, 0], 9, [2, 1]] \\
C_3' =& \; [0, [6, 0], 6, [6, 2]]
\end{align*}
\end{minipage}

\item 2-factor type $[9, 2, 2]$: three starters \vspace*{-3mm} \\
\begin{minipage}{8cm}
\begin{align*}
C_0 =& \; [0, [2, 2], 10, [5, -1], 3, [4, -1], 7, [5, 1], 2, [1, 0], 1, [4, 2], 5, [4, 0], 9, [12, 2], 12, [12, 0]] \\
C_1 =& \; [4, [2, -1], 6, [2, 1]] \\
C_2 =& \; [8, [3, 0], 11, [3, 2]]
\end{align*}
\end{minipage}

\vspace*{-4mm}
\begin{minipage}{8cm}
\begin{align*}
C_0' =& \; [7, [2, 0], 9, [5, 2], 2, [3, -1], 5, [4, 1], 1, [5, 0], 8, [12, 1], 12, [12, -1], 11, [1, 2], 10, [3, 1]] \\
C_1' =& \; [0, [6, 0], 6, [6, 2]] \\
C_2' =& \; [4, [1, -1], 3, [1, 1]]
\end{align*}
\end{minipage}

\item 2-factor type $[8, 3, 2]$: three starters \vspace*{-3mm} \\
\begin{minipage}{8cm}
\begin{align*}
C_0 =& \; [7, [5, 1], 2, [3, -1], 5, [4, -1], 9, [5, 0], 4, [4, 2], 0, [2, 1], 10, [4, 0], 6, [1, 2]] \\
C_1 =& \; [8, [5, -1], 1, [2, 0], 11, [3, 2]] \\
C_2 =& \; [12, [12, -1], 3, [12, 1]]
\end{align*}
\end{minipage}

\vspace*{-4mm}
\begin{minipage}{8cm}
\begin{align*}
C_0' =& \; [12, [12, 2], 5, [2, -1], 7, [4, 1], 3, [1, 1], 2, [1, 0], 1, [3, 1], 4, [5, 2], 9, [12, 0]] \\
C_1' =& \; [11, [3, 0], 8, [2, 2], 10, [1, -1]] \\
C_2' =& \; [0, [6, 0], 6, [6, 2]]
\end{align*}
\end{minipage}

\item 2-factor type $[7, 4, 2]$: three starters \vspace*{-3mm} \\
\begin{minipage}{8cm}
\begin{align*}
C_0 =& \; [8, [1, 0], 9, [5, 2], 2, [2, 1], 0, [4, 0], 4, [3, -1], 1, [4, 1], 5, [3, 2]] \\
C_1 =& \; [7, [4, -1], 3, [3, 1], 6, [12, 1], 12, [12, -1]] \\
C_2 =& \; [11, [1, -1], 10, [1, 1]]
\end{align*}
\end{minipage}

\vspace*{-4mm}
\begin{minipage}{8cm}
\begin{align*}
C_0' =& \; [8, [4, 2], 4, [3, 0], 7, [2, 2], 5, [5, -1], 10, [5, 0], 3, [12, 2], 12, [12, 0]] \\
C_1' =& \; [9, [5, 1], 2, [1, 2], 1, [2, 0], 11, [2, -1]] \\
C_2' =& \; [0, [6, 0], 6, [6, 2]]
\end{align*}
\end{minipage}

\item 2-factor type $[6, 5, 2]$: three starters \vspace*{-3mm} \\
\begin{minipage}{8cm}
\begin{align*}
C_0 =& \; [5, [2, 1], 3, [12, 0], 12, [12, -1], 7, [5, 2], 2, [1, 1], 1, [4, -1]] \\
C_1 =& \; [10, [4, 1], 6, [2, -1], 8, [4, 0], 0, [1, -1], 11, [1, 2]] \\
C_2 =& \; [9, [5, -1], 4, [5, 1]]
\end{align*}
\end{minipage}

\vspace*{-4mm}
\begin{minipage}{8cm}
\begin{align*}
C_0' =& \; [11, [12, 1], 12, [12, 2], 9, [5, 0], 4, [3, 2], 7, [1, 0], 8, [3, 1]] \\
C_1' =& \; [3, [2, 0], 5, [3, -1], 2, [4, 2], 10, [3, 0], 1, [2, 2]] \\
C_2' =& \; [0, [6, 0], 6, [6, 2]]
\end{align*}
\end{minipage}

\end{itemizenew}


\section{Computational results for $n=14$}\label{app:14}

\begin{itemizenew}
\item 2-factor type $[4, 2, 2, 2, 2, 2]$: one starter \vspace*{-3mm} \\
\hspace*{-7mm}
\begin{minipage}{8cm}
\begin{align*}
C_0 =& \;  [6, [3, 0], 9, [13, 1], 13, [13, 0], 3, [3, 1]] \\
C_1 =& \;  [0, [4, 1], 4, [4, 0]] \\
C_2 =& \;  [12, [5, 1], 7, [5, 0]] \\
C_3 =& \;  [10, [2, 0], 8, [2, 1]] \\
C_4 =& \;  [1, [1, 0], 2, [1, 1]] \\
C_5 =& \;  [5, [6, 1], 11, [6, 0]]
\end{align*}
\end{minipage}

\item 2-factor type $[3, 3, 2, 2, 2, 2]$: two starters \vspace*{-3mm} \\
\hspace*{-14mm}
\begin{minipage}{8cm}
\begin{align*}
C_0 =& \; [5, [3, 2], 2, [5, 0], 10, [5, -1]] \\
C_1 =& \; [12, [2, 0], 1, [6, 1], 7, [5, 2]] \\
C_2 =& \; [4, [2, -1], 6, [2, 1]] \\
C_3 =& \; [3, [6, 0], 9, [6, 2]] \\
C_4 =& \; [11, [3, -1], 8, [3, 1]] \\
C_5 =& \; [13, [13, -1], 0, [13, 1]]
\end{align*}
\end{minipage}
\begin{minipage}{8cm}
\begin{align*}
C_0' =& \; [8, [6, -1], 2, [13, 2], 13, [13, 0]] \\
C_1' =& \; [4, [3, 0], 1, [5, 1], 6, [2, 2]] \\
C_2' =& \; [7, [4, -1], 3, [4, 1]] \\
C_3' =& \; [11, [1, -1], 10, [1, 1]] \\
C_4' =& \; [9, [4, 0], 5, [4, 2]] \\
C_5' =& \; [12, [1, 0], 0, [1, 2]]
\end{align*}
\end{minipage}

\item 2-factor type $[6, 2, 2, 2, 2]$: two starters \vspace*{-3mm} \\
\begin{minipage}{8cm}
\begin{align*}
C_0 =& \; [6, [3, 2], 3, [13, -1], 13, [13, 1], 0, [5, 1], 8, [6, 0], 2, [4, 1]] \\
C_1 =& \; [11, [1, 0], 10, [1, 2]] \\
C_2 =& \; [5, [6, -1], 12, [6, 1]] \\
C_3 =& \; [1, [5, 0], 9, [5, 2]] \\
C_4 =& \; [4, [3, -1], 7, [3, 1]]
\end{align*}
\end{minipage}

\vspace*{-4mm}
\begin{minipage}{8cm}
\begin{align*}
C_0' =& \; [5, [3, 0], 8, [5, -1], 3, [6, 2], 9, [2, 0], 11, [4, -1], 7, [2, 2]] \\
C_1' =& \; [6, [13, 0], 13, [13, 2]] \\
C_2' =& \; [10, [4, 0], 1, [4, 2]] \\
C_3' =& \; [12, [1, -1], 0, [1, 1]] \\
C_4' =& \; [4, [2, -1], 2, [2, 1]]
\end{align*}
\end{minipage}

\item 2-factor type $[5, 3, 2, 2, 2]$: one starter \vspace*{-3mm} \\
\begin{minipage}{8cm}
\begin{align*}
C_0 =& \;  [6, [2, 1], 8, [3, 0], 11, [13, 0], 13, [13, 1], 3, [3, 1]] \\
C_1 =& \;  [2, [2, 0], 0, [1, 0], 1, [1, 1]] \\
C_2 =& \;  [4, [6, 1], 10, [6, 0]] \\
C_3 =& \;  [5, [4, 0], 9, [4, 1]] \\
C_4 =& \;  [12, [5, 0], 7, [5, 1]]
\end{align*}
\end{minipage}

\item 2-factor type $[4, 4, 2, 2, 2]$: one starter \vspace*{-3mm} \\
\hspace*{-8mm}
\begin{minipage}{8cm}
\begin{align*}
C_0 =& \;  [5, [2, 1], 3, [5, 1], 8, [1, 0], 7, [2, 0]] \\
C_1 =& \;  [12, [3, 1], 9, [5, 0], 1, [3, 0], 11, [1, 1]] \\
C_2 =& \;  [13, [13, 0], 0, [13, 1]] \\
C_3 =& \;  [4, [6, 0], 10, [6, 1]] \\
C_4 =& \;  [2, [4, 1], 6, [4, 0]]
\end{align*}
\end{minipage}

\newpage

\item 2-factor type $[4, 3, 3, 2, 2]$: two starters \vspace*{-3mm} \\
\hspace*{-6mm}
\begin{minipage}{8cm}
\begin{align*}
C_0 =& \; [12, [4, 0], 3, [6, -1], 10, [4, 1], 1, [2, 2]] \\
C_1 =& \; [4, [5, -1], 9, [2, 0], 11, [6, 2]] \\
C_2 =& \; [2, [6, 0], 8, [13, 2], 13, [13, 1]] \\
C_3 =& \; [0, [5, 0], 5, [5, 2]] \\
C_4 =& \; [7, [1, -1], 6, [1, 1]]
\end{align*}
\end{minipage}
\begin{minipage}{8cm}
\begin{align*}
C_0' =& \; [0, [4, -1], 9, [2, -1], 7, [3, 2], 10, [3, 0]] \\
C_1' =& \; [4, [2, 1], 6, [5, 1], 11, [6, 1]] \\
C_2' =& \; [3, [13, -1], 13, [13, 0], 12, [4, 2]] \\
C_3' =& \; [1, [1, 0], 2, [1, 2]] \\
C_4' =& \; [8, [3, -1], 5, [3, 1]]
\end{align*}
\end{minipage}

\item 2-factor type $[3, 3, 3, 3, 2]$: one starter \vspace*{-3mm} \\
\hspace*{-14mm}
\begin{minipage}{8cm}
\begin{align*}
C_0 =& \;  [0, [5, 0], 8, [13, 1], 13, [13, 0]] \\
C_1 =& \;  [11, [6, 1], 4, [5, 1], 12, [1, 1]] \\
C_2 =& \;  [9, [6, 0], 2, [3, 0], 5, [4, 1]] \\
C_3 =& \;  [6, [1, 0], 7, [3, 1], 10, [4, 0]] \\
C_4 =& \;  [1, [2, 1], 3, [2, 0]]
\end{align*}
\end{minipage}

\item 2-factor type $[8, 2, 2, 2]$: one starter \vspace*{-3mm} \\
\begin{minipage}{8cm}
\begin{align*}
C_0 =& \;  [4, [1, 1], 5, [3, 1], 2, [6, 0], 8, [13, 0], 13, [13, 1], 12, [1, 0], 0, [3, 0], 10, [6, 1]] \\
C_1 =& \;  [11, [2, 0], 9, [2, 1]] \\
C_2 =& \;  [6, [5, 1], 1, [5, 0]] \\
C_3 =& \;  [7, [4, 1], 3, [4, 0]]
\end{align*}
\end{minipage}

\item 2-factor type $[7, 3, 2, 2]$: two starters \vspace*{-3mm} \\
\begin{minipage}{8cm}
\begin{align*}
C_0 =& \; [4, [2, -1], 6, [13, -1], 13, [13, 0], 8, [6, 2], 2, [6, 0], 9, [6, -1], 3, [1, 2]] \\
C_1 =& \; [12, [2, 1], 10, [4, 0], 1, [2, 2]] \\
C_2 =& \; [7, [4, -1], 11, [4, 1]] \\
C_3 =& \; [5, [5, -1], 0, [5, 1]]
\end{align*}
\end{minipage}

\vspace*{-4mm}
\begin{minipage}{8cm}
\begin{align*}
C_0' =& \; [5, [1, 0], 6, [4, 2], 2, [6, 1], 9, [2, 0], 7, [3, 2], 10, [3, 0], 0, [5, 2]] \\
C_1' =& \; [8, [13, 1], 13, [13, 2], 3, [5, 0]] \\
C_2' =& \; [11, [1, -1], 12, [1, 1]] \\
C_3' =& \; [1, [3, -1], 4, [3, 1]]
\end{align*}
\end{minipage}

\item 2-factor type $[6, 4, 2, 2]$: two starters \vspace*{-3mm} \\
\begin{minipage}{8cm}
\begin{align*}
C_0 =& \; [2, [1, 2], 1, [3, 0], 11, [3, -1], 8, [13, 2], 13, [13, 0], 12, [3, 1]] \\
C_1 =& \; [0, [3, 2], 10, [6, 1], 4, [5, -1], 9, [4, 0]] \\
C_2 =& \; [7, [4, -1], 3, [4, 1]] \\
C_3 =& \; [5, [1, -1], 6, [1, 1]]
\end{align*}
\end{minipage}

\vspace*{-4mm}
\begin{minipage}{8cm}
\begin{align*}
C_0' =& \; [1, [6, -1], 8, [5, 2], 3, [6, 0], 10, [6, 2], 4, [5, 0], 12, [2, 1]] \\
C_1' =& \; [11, [5, 1], 6, [1, 0], 5, [4, 2], 9, [2, -1]] \\
C_2' =& \; [13, [13, -1], 7, [13, 1]] \\
C_3' =& \; [2, [2, 0], 0, [2, 2]]
\end{align*}
\end{minipage}

\newpage

\item 2-factor type $[5, 5, 2, 2]$: two starters \vspace*{-3mm} \\
\begin{minipage}{8cm}
\begin{align*}
C_0 =& \; [12, [1, 1], 0, [4, 1], 4, [5, 2], 9, [1, 0], 10, [2, 1]] \\
C_1 =& \; [1, [4, -1], 5, [2, -1], 7, [4, 0], 11, [13, 1], 13, [13, 2]] \\
C_2 =& \; [2, [6, 0], 8, [6, 2]] \\
C_3 =& \; [3, [3, -1], 6, [3, 1]]
\end{align*}
\end{minipage}

\vspace*{-4mm}
\begin{minipage}{8cm}
\begin{align*}
C_0' =& \; [13, [13, 0], 4, [1, 2], 5, [3, 0], 8, [4, 2], 12, [13, -1]] \\
C_1' =& \; [11, [1, -1], 10, [3, 2], 7, [5, 0], 2, [2, 2], 0, [2, 0]] \\
C_2' =& \; [6, [5, -1], 1, [5, 1]] \\
C_3' =& \; [9, [6, -1], 3, [6, 1]]
\end{align*}
\end{minipage}

\item 2-factor type $[6, 3, 3, 2]$: one starter \vspace*{-3mm} \\
\begin{minipage}{8cm}
\begin{align*}
C_0 =& \;  [8, [6, 1], 1, [5, 1], 9, [13, 1], 13, [13, 0], 11, [6, 0], 5, [3, 1]] \\
C_1 =& \;  [4, [5, 0], 12, [4, 1], 3, [1, 0]] \\
C_2 =& \;  [7, [1, 1], 6, [4, 0], 10, [3, 0]] \\
C_3 =& \;  [0, [2, 0], 2, [2, 1]]
\end{align*}
\end{minipage}

\item 2-factor type $[5, 4, 3, 2]$: one starter \vspace*{-3mm} \\
\begin{minipage}{8cm}
\begin{align*}
C_0 =& \;  [6, [4, 1], 10, [5, 0], 5, [6, 1], 12, [5, 1], 7, [1, 0]] \\
C_1 =& \;  [1, [13, 1], 13, [13, 0], 4, [1, 1], 3, [2, 0]] \\
C_2 =& \;  [2, [6, 0], 9, [4, 0], 0, [2, 1]] \\
C_3 =& \;  [8, [3, 0], 11, [3, 1]]
\end{align*}
\end{minipage}

\item 2-factor type $[4, 4, 4, 2]$: one starter \vspace*{-3mm} \\
\hspace*{-7mm}
\begin{minipage}{8cm}
\begin{align*}
C_0 =& \; [4, [4, 0], 0, [13, 0], 13, [13, 1], 6, [2, 1]] \\
C_1 =& \; [8, [3, 1], 11, [3, 0], 1, [1, 1], 2, [6, 0]] \\
C_2 =& \; [9, [1, 0], 10, [6, 1], 3, [2, 0], 5, [4, 1]] \\
C_3 =& \; [12, [5, 0], 7, [5, 1]]
\end{align*}
\end{minipage}

\item 2-factor type $[5, 3, 3, 3]$: two starters \vspace*{-3mm} \\
\begin{minipage}{8cm}
\begin{align*}
C_0 =& \; [7, [2, 0], 9, [3, -1], 6, [6, 1], 12, [4, 1], 3, [4, 2]] \\
C_1 =& \; [2, [2, 2], 0, [3, 0], 10, [5, 1]] \\
C_2 =& \; [4, [4, 0], 8, [3, 2], 5, [1, 1]] \\
C_3 =& \; [13, [13, 0], 11, [3, 1], 1, [13, 2]]
\end{align*}
\end{minipage}

\vspace*{-4mm}
\begin{minipage}{8cm}
\begin{align*}
C_0' =& \; [0, [2, -1], 2, [4, -1], 6, [5, 0], 1, [6, 2], 7, [6, -1]] \\
C_1' =& \; [8, [2, 1], 10, [1, 2], 9, [1, 0]] \\
C_2' =& \; [5, [6, 0], 12, [5, 2], 4, [1, -1]] \\
C_3' =& \; [3, [5, -1], 11, [13, 1], 13, [13, -1]]
\end{align*}
\end{minipage}

\item 2-factor type $[4, 4, 3, 3]$: two starters \vspace*{-3mm} \\
\hspace*{-5mm}
\begin{minipage}{8cm}
\begin{align*}
C_0 =& \; [11, [1, 2], 10, [6, -1], 3, [4, 1], 12, [1, 0]] \\
C_1 =& \; [5, [13, 0], 13, [13, 2], 9, [1, 1], 8, [3, 1]] \\
C_2 =& \; [4, [3, -1], 1, [5, 2], 6, [2, 0]] \\
C_3 =& \; [7, [6, 0], 0, [2, 2], 2, [5, -1]]
\end{align*}
\end{minipage}
\begin{minipage}{8cm}
\begin{align*}
C_0' =& \; [13, [13, -1], 1, [5, 1], 6, [2, 1], 8, [13, 1]] \\
C_1' =& \; [9, [6, 2], 3, [1, -1], 4, [4, 0], 0, [4, -1]] \\
C_2' =& \; [11, [4, 2], 2, [3, 0], 5, [6, 1]] \\
C_3' =& \; [7, [3, 2], 10, [2, -1], 12, [5, 0]]
\end{align*}
\end{minipage}

\item 2-factor type $[10, 2, 2]$: two starters \vspace*{-3mm} \\
\begin{minipage}{8cm}
\begin{align*}
C_0 =& \; [5, [3, 2], 2, [6, 1], 9, [3, 0], 12, [2, 2], 1, [6, -1], 7, [1, 0], 6, [4, 2], 10, [3, -1], 0, [4, 0], 4, [1, 1]] \\
C_1 =& \; [11, [5, 0], 3, [5, 2]] \\
C_2 =& \; [13, [13, 0], 8, [13, 2]]
\end{align*}
\end{minipage}

\vspace*{-4mm}
\begin{minipage}{8cm}
\begin{align*}
C_0' =& \; [1, [2, 0], 3, [5, -1], 11, [6, 2], 5, [6, 0], 12, [5, 1], 4, [3, 1], 7, [1, -1], 6, [13, 1], 13, [13, -1], 2, [1, 2]] \\
C_1' =& \; [9, [4, -1], 0, [4, 1]] \\
C_2' =& \; [10, [2, -1], 8, [2, 1]]
\end{align*}
\end{minipage}

\item 2-factor type $[9, 3, 2]$: one starter \vspace*{-3mm} \\
\begin{minipage}{8cm}
\begin{align*}
C_0 =& \;  [7, [6, 0], 0, [2, 1], 2, [6, 1], 9, [3, 0], 12, [2, 0], 10, [4, 1], 1, [4, 0], 5, [3, 1], 8, [1, 1]] \\
C_1 =& \;  [3, [13, 1], 13, [13, 0], 4, [1, 0]] \\
C_2 =& \;  [6, [5, 1], 11, [5, 0]]
\end{align*}
\end{minipage}

\item 2-factor type $[8, 4, 2]$: one starter \vspace*{-3mm} \\
\begin{minipage}{8cm}
\begin{align*}
C_0 =& \;  [3, [13, 0], 13, [13, 1], 7, [3, 1], 10, [6, 1], 4, [5, 0], 12, [1, 1], 0, [2, 0], 2, [1, 0]] \\
C_1 =& \;  [6, [5, 1], 11, [3, 0], 1, [6, 0], 8, [2, 1]] \\
C_2 =& \;  [5, [4, 0], 9, [4, 1]]
\end{align*}
\end{minipage}

\item 2-factor type $[7, 5, 2]$: one starter \vspace*{-3mm} \\
\begin{minipage}{8cm}
\begin{align*}
C_0 =& \;  [8, [6, 1], 2, [6, 0], 9, [4, 1], 5, [1, 1], 4, [3, 1], 7, [3, 0], 10, [2, 1]] \\
C_1 =& \;  [1, [13, 1], 13, [13, 0], 0, [1, 0], 12, [4, 0], 3, [2, 0]] \\
C_2 =& \;  [11, [5, 1], 6, [5, 0]]
\end{align*}
\end{minipage}

\item 2-factor type $[6, 6, 2]$: one starter \vspace*{-3mm} \\
\begin{minipage}{8cm}
\begin{align*}
C_0 =& \;  [0, [6, 1], 6, [3, 0], 9, [6, 0], 3, [5, 0], 8, [13, 0], 13, [13, 1]] \\
C_1 =& \;  [4, [3, 1], 7, [2, 0], 5, [4, 1], 1, [4, 0], 10, [5, 1], 2, [2, 1]] \\
C_2 =& \;  [12, [1, 1], 11, [1, 0]]
\end{align*}
\end{minipage}

\item 2-factor type $[8, 3, 3]$: two starters \vspace*{-3mm} \\
\begin{minipage}{8cm}
\begin{align*}
C_0 =& \; [8, [6, -1], 1, [4, -1], 5, [6, 1], 12, [5, -1], 4, [6, 0], 10, [1, 2], 9, [2, -1], 11, [3, 1]] \\
C_1 =& \; [0, [13, 0], 13, [13, -1], 6, [6, 2]] \\
C_2 =& \; [3, [1, 1], 2, [5, 0], 7, [4, 2]]
\end{align*}
\end{minipage}

\vspace*{-4mm}
\begin{minipage}{8cm}
\begin{align*}
C_0' =& \; [2, [5, 2], 7, [3, 0], 4, [2, 2], 6, [1, 0], 5, [5, 1], 10, [3, 2], 0, [2, 1], 11, [4, 0]] \\
C_1' =& \; [3, [13, 2], 13, [13, 1], 1, [2, 0]] \\
C_2' =& \; [12, [3, -1], 9, [1, -1], 8, [4, 1]]
\end{align*}
\end{minipage}

\item 2-factor type $[7, 4, 3]$: two starters \vspace*{-3mm} \\
\begin{minipage}{8cm}
\begin{align*}
C_0 =& \; [6, [1, 2], 5, [6, 1], 12, [13, -1], 13, [13, 0], 9, [1, -1], 8, [4, 2], 4, [2, 0]] \\
C_1 =& \; [1, [4, 1], 10, [5, 0], 2, [4, -1], 11, [3, 2]] \\
C_2 =& \; [3, [4, 0], 7, [6, 2], 0, [3, -1]]
\end{align*}
\end{minipage}

\vspace*{-4mm}
\begin{minipage}{8cm}
\begin{align*}
C_0' =& \; [10, [5, 2], 2, [2, 1], 0, [1, 0], 12, [13, 1], 13, [13, 2], 11, [5, -1], 3, [6, 0]] \\
C_1' =& \; [7, [6, -1], 1, [3, 1], 4, [5, 1], 9, [2, -1]] \\
C_2' =& \; [8, [3, 0], 5, [1, 1], 6, [2, 2]]
\end{align*}
\end{minipage}

\item 2-factor type $[6, 5, 3]$: two starters \vspace*{-3mm} \\
\begin{minipage}{8cm}
\begin{align*}
C_0 =& \; [6, [3, 1], 3, [2, 0], 5, [1, 2], 4, [6, 0], 11, [13, 2], 13, [13, 1]] \\
C_1 =& \; [7, [5, 2], 12, [4, 1], 8, [6, -1], 1, [4, 0], 10, [3, -1]] \\
C_2 =& \; [0, [4, -1], 9, [6, 1], 2, [2, -1]]
\end{align*}
\end{minipage}

\vspace*{-4mm}
\begin{minipage}{8cm}
\begin{align*}
C_0' =& \; [2, [2, 2], 0, [5, 1], 8, [1, -1], 9, [5, 0], 1, [4, 2], 5, [3, 0]] \\
C_1' =& \; [12, [1, 1], 11, [5, -1], 3, [1, 0], 4, [2, 1], 6, [6, 2]] \\
C_2' =& \; [10, [3, 2], 7, [13, 0], 13, [13, -1]]
\end{align*}
\end{minipage}

\item 2-factor type $[6, 4, 4]$: two starters \vspace*{-3mm} \\
\begin{minipage}{8cm}
\begin{align*}
C_0 =& \; [13, [13, 0], 1, [1, 2], 2, [5, 0], 7, [4, -1], 3, [2, 1], 5, [13, 2]] \\
C_1 =& \; [0, [6, 2], 6, [3, 0], 9, [2, 2], 11, [2, 0]] \\
C_2 =& \; [8, [4, 2], 12, [5, -1], 4, [6, 0], 10, [2, -1]]
\end{align*}
\end{minipage}

\vspace*{-4mm}
\begin{minipage}{8cm}
\begin{align*}
C_0' =& \; [2, [3, 1], 12, [4, 1], 8, [1, 1], 7, [13, -1], 13, [13, 1], 1, [1, -1]] \\
C_1' =& \; [4, [6, -1], 10, [3, -1], 0, [4, 0], 9, [5, 2]] \\
C_2' =& \; [6, [1, 0], 5, [6, 1], 11, [5, 1], 3, [3, 2]]
\end{align*}
\end{minipage}

\item 2-factor type $[5, 5, 4]$: two starters \vspace*{-3mm} \\
\begin{minipage}{8cm}
\begin{align*}
C_0 =& \; [13, [13, 2], 10, [6, 0], 4, [3, 1], 7, [4, 1], 11, [13, 1]] \\
C_1 =& \; [1, [5, 0], 9, [4, 2], 5, [5, 1], 0, [1, 0], 12, [2, 2]] \\
C_2 =& \; [8, [2, 1], 6, [4, 0], 2, [1, 1], 3, [5, 2]]
\end{align*}
\end{minipage}

\vspace*{-4mm}
\begin{minipage}{8cm}
\begin{align*}
C_0' =& \; [12, [3, 0], 9, [5, -1], 4, [1, 2], 3, [3, -1], 6, [6, -1]] \\
C_1' =& \; [2, [4, -1], 11, [6, 2], 5, [2, -1], 7, [13, -1], 13, [13, 0]] \\
C_2' =& \; [10, [2, 0], 8, [6, 1], 1, [1, -1], 0, [3, 2]]
\end{align*}
\end{minipage}

\item 2-factor type $[12, 2]$: one starter \vspace*{-3mm} \\
\begin{minipage}{8cm}
\begin{align*}
C_0 =& \;  [0, [13, 1], 13, [13, 0], 9, [1, 1], 10, [2, 0], 12, [1, 0], 11, [5, 1], 3, [3, 1], 6, [5, 0], 1, [4, 1], 5, [2, 1], 7, [3, 0],  \\ & \; 4, [4, 0]] \\
C_1 =& \;  [8, [6, 0], 2, [6, 1]]
\end{align*}
\end{minipage}

\item 2-factor type $[11, 3]$: two starters \vspace*{-3mm} \\
\begin{minipage}{8cm}
\begin{align*}
C_0 =& \; [7, [4, 2], 3, [5, 0], 8, [2, 2], 10, [6, 0], 4, [1, 2], 5, [6, 1], 12, [3, 0], 2, [2, -1], 0, [1, 1], 1, [5, 2], 9, [2, 0]] \\
C_1 =& \; [6, [5, 1], 11, [13, 2], 13, [13, 0]]
\end{align*}
\end{minipage}

\vspace*{-4mm}
\begin{minipage}{8cm}
\begin{align*}
C_0' =& \; [9, [2, 1], 7, [4, -1], 11, [1, -1], 12, [6, -1], 5, [3, 1], 2, [1, 0], 1, [3, 2], 4, [13, -1], 13, [13, 1], 8, [5, -1], \\ & \; 0, [4, 1]] \\
C_1' =& \; [10, [6, 2], 3, [3, -1], 6, [4, 0]]
\end{align*}
\end{minipage}

\item 2-factor type $[10, 4]$: two starters \vspace*{-3mm} \\
\begin{minipage}{8cm}
\begin{align*}
C_0 =& \; [0, [4, 1], 4, [5, 1], 9, [6, 1], 2, [5, -1], 10, [2, 2], 8, [13, 0], 13, [13, -1], 11, [5, 2], 3, [4, 0], 12, [1, 1]] \\
C_1 =& \; [7, [2, 0], 5, [1, 2], 6, [5, 0], 1, [6, 2]]
\end{align*}
\end{minipage}

\vspace*{-4mm}
\begin{minipage}{8cm}
\begin{align*}
C_0' =& \; [1, [2, -1], 3, [3, 0], 6, [4, -1], 2, [3, -1], 12, [4, 2], 8, [1, 0], 9, [2, 1], 11, [6, -1], 5, [1, -1], 4, [3, 2]] \\
C_1' =& \; [13, [13, 1], 0, [6, 0], 7, [3, 1], 10, [13, 2]]
\end{align*}
\end{minipage}

\item 2-factor type $[9, 5]$: two starters \vspace*{-3mm} \\
\begin{minipage}{8cm}
\begin{align*}
C_0 =& \; [6, [1, 1], 5, [2, 0], 3, [6, 1], 9, [13, 1], 13, [13, -1], 7, [3, 2], 4, [6, 0], 10, [2, 1], 12, [6, 2]] \\
C_1 =& \; [1, [1, 0], 2, [4, -1], 11, [3, -1], 8, [5, 1], 0, [1, 2]]
\end{align*}
\end{minipage}

\vspace*{-4mm}
\begin{minipage}{8cm}
\begin{align*}
C_0' =& \; [2, [6, -1], 9, [5, 2], 4, [3, 0], 7, [1, -1], 6, [2, 2], 8, [5, 0], 0, [3, 1], 3, [5, -1], 11, [4, 1]] \\
C_1' =& \; [13, [13, 0], 12, [2, -1], 10, [4, 2], 1, [4, 0], 5, [13, 2]]
\end{align*}
\end{minipage}

\item 2-factor type $[8, 6]$: two starters \vspace*{-3mm} \\
\begin{minipage}{8cm}
\begin{align*}
C_0 =& \; [13, [13, -1], 10, [5, 2], 5, [3, -1], 8, [4, 0], 12, [1, 1], 0, [3, 1], 3, [2, -1], 1, [13, 1]] \\
C_1 =& \; [9, [3, 2], 6, [4, 1], 2, [5, -1], 7, [3, 0], 4, [6, 2], 11, [2, 0]]
\end{align*}
\end{minipage}

\vspace*{-4mm}
\begin{minipage}{8cm}
\begin{align*}
C_0' =& \; [1, [1, -1], 0, [2, 1], 2, [2, 2], 4, [5, 0], 9, [6, -1], 3, [4, 2], 7, [1, 0], 8, [6, 1]] \\
C_1' =& \; [13, [13, 2], 5, [6, 0], 12, [1, 2], 11, [5, 1], 6, [4, -1], 10, [13, 0]]
\end{align*}
\end{minipage}

\end{itemizenew}

\section{Computational results for $n=15$}\label{app:15}

\begin{itemizenew}
\item 2-factor type $[3, 2, 2, 2, 2, 2, 2]$: three starters \vspace*{-3mm} \\
\hspace*{-15mm}
\begin{minipage}{8cm}
\begin{align*}
C_0 =& \; [9, [5, 0], 0, [6, 2], 6, [3, -1]] \\
C_1 =& \; [11, [3, 0], 8, [3, 2]] \\
C_2 =& \; [14, [14, -1], 10, [14, 1]] \\
C_3 =& \; [12, [1, -1], 13, [1, 1]] \\
C_4 =& \; [4, [2, 0], 2, [2, 2]] \\
C_5 =& \; [1, [4, -1], 5, [4, 1]] \\
C_6 =& \; [7, [4, 0], 3, [4, 2]]
\end{align*}
\end{minipage}
\begin{minipage}{8cm}
\begin{align*}
C_0' =& \; [4, [5, 2], 9, [6, 0], 1, [3, 1]] \\
C_1' =& \; [6, [5, -1], 11, [5, 1]] \\
C_2' =& \; [3, [1, 0], 2, [1, 2]] \\
C_3' =& \; [8, [14, 0], 14, [14, 2]] \\
C_4' =& \; [12, [2, -1], 10, [2, 1]] \\
C_5' =& \; [13, [6, -1], 5, [6, 1]] \\
C_6' =& \; [0, [7, 0], 7, [7, 2]]
\end{align*}
\end{minipage}

\item 2-factor type $[5, 2, 2, 2, 2, 2]$: three starters \vspace*{-3mm} \\
\begin{minipage}{8cm}
\begin{align*}
C_0 =& \; [11, [6, 2], 3, [6, -1], 9, [6, 0], 1, [1, 2], 0, [3, 0]] \\
C_1 =& \; [2, [14, -1], 14, [14, 1]] \\
C_2 =& \; [13, [5, -1], 4, [5, 1]] \\
C_3 =& \; [8, [4, 0], 12, [4, 2]] \\
C_4 =& \; [7, [2, 0], 5, [2, 2]] \\
C_5 =& \; [6, [4, -1], 10, [4, 1]]
\end{align*}
\end{minipage}

\vspace*{-4mm}
\begin{minipage}{8cm}
\begin{align*}
C_0' =& \; [6, [3, 2], 9, [6, 1], 3, [14, 0], 14, [14, 2], 5, [1, 0]] \\
C_1' =& \; [12, [1, -1], 13, [1, 1]] \\
C_2' =& \; [2, [5, 0], 11, [5, 2]] \\
C_3' =& \; [10, [2, -1], 8, [2, 1]] \\
C_4' =& \; [4, [3, -1], 1, [3, 1]] \\
C_5' =& \; [0, [7, 0], 7, [7, 2]]
\end{align*}
\end{minipage}

\newpage

\item 2-factor type $[4, 3, 2, 2, 2, 2]$: three starters \vspace*{-3mm} \\
\hspace*{-4mm}
\begin{minipage}{8cm}
\begin{align*}
C_0 =& \; [10, [5, -1], 5, [4, -1], 1, [6, -1], 9, [1, 1]] \\
C_1 =& \; [12, [5, 0], 3, [14, 1], 14, [14, 2]] \\
C_2 =& \; [2, [2, -1], 4, [2, 1]] \\
C_3 =& \; [6, [2, 0], 8, [2, 2]] \\
C_4 =& \; [11, [3, 0], 0, [3, 2]] \\
C_5 =& \; [7, [6, 0], 13, [6, 2]]
\end{align*}
\end{minipage}
\begin{minipage}{8cm}
\begin{align*}
C_0' =& \; [1, [4, 1], 5, [5, 1], 10, [1, -1], 9, [6, 1]] \\
C_1' =& \; [14, [14, -1], 11, [5, 2], 2, [14, 0]] \\
C_2' =& \; [4, [4, 0], 8, [4, 2]] \\
C_3' =& \; [6, [3, -1], 3, [3, 1]] \\
C_4' =& \; [12, [1, 0], 13, [1, 2]] \\
C_5' =& \; [0, [7, 0], 7, [7, 2]]
\end{align*}
\end{minipage}

\item 2-factor type $[3, 3, 3, 2, 2, 2]$: three starters \vspace*{-3mm} \\
\hspace*{-11mm}
\begin{minipage}{8cm}
\begin{align*}
C_0 =& \; [3, [6, -1], 9, [5, 1], 4, [1, 1]] \\
C_1 =& \; [2, [4, 0], 12, [2, 2], 0, [2, -1]] \\
C_2 =& \; [1, [14, 1], 14, [14, -1], 11, [4, 1]] \\
C_3 =& \; [6, [1, 0], 7, [1, 2]] \\
C_4 =& \; [10, [3, 0], 13, [3, 2]] \\
C_5 =& \; [8, [3, -1], 5, [3, 1]]
\end{align*}
\end{minipage}
\begin{minipage}{8cm}
\begin{align*}
C_0' =& \; [1, [4, 2], 11, [2, 0], 13, [2, 1]] \\
C_1' =& \; [2, [4, -1], 6, [14, 0], 14, [14, 2]] \\
C_2' =& \; [8, [5, -1], 3, [6, 1], 9, [1, -1]] \\
C_3' =& \; [0, [7, 0], 7, [7, 2]] \\
C_4' =& \; [12, [6, 0], 4, [6, 2]] \\
C_5' =& \; [5, [5, 0], 10, [5, 2]]
\end{align*}
\end{minipage}

\item 2-factor type $[7, 2, 2, 2, 2]$: three starters \vspace*{-3mm} \\
\begin{minipage}{8cm}
\begin{align*}
C_0 =& \; [8, [3, -1], 11, [2, 0], 13, [6, 1], 5, [3, 2], 2, [4, 0], 6, [6, -1], 0, [6, 2]] \\
C_1 =& \; [3, [4, -1], 7, [4, 1]] \\
C_2 =& \; [1, [14, -1], 14, [14, 1]] \\
C_3 =& \; [12, [2, -1], 10, [2, 1]] \\
C_4 =& \; [4, [5, 0], 9, [5, 2]]
\end{align*}
\end{minipage}

\vspace*{-4mm}
\begin{minipage}{8cm}
\begin{align*}
C_0' =& \; [2, [4, 2], 12, [1, 1], 11, [6, 0], 5, [3, 1], 8, [2, 2], 6, [3, 0], 3, [1, -1]] \\
C_1' =& \; [0, [7, 0], 7, [7, 2]] \\
C_2' =& \; [10, [1, 0], 9, [1, 2]] \\
C_3' =& \; [1, [14, 0], 14, [14, 2]] \\
C_4' =& \; [4, [5, -1], 13, [5, 1]]
\end{align*}
\end{minipage}

\item 2-factor type $[6, 3, 2, 2, 2]$: three starters \vspace*{-3mm} \\
\begin{minipage}{8cm}
\begin{align*}
C_0 =& \; [1, [4, -1], 5, [1, 0], 6, [6, 2], 12, [5, 1], 7, [14, -1], 14, [14, 1]] \\
C_1 =& \; [11, [5, -1], 2, [3, 1], 13, [2, 1]] \\
C_2 =& \; [8, [4, 0], 4, [4, 2]] \\
C_3 =& \; [9, [1, -1], 10, [1, 1]] \\
C_4 =& \; [0, [3, 0], 3, [3, 2]]
\end{align*}
\end{minipage}

\vspace*{-4mm}
\begin{minipage}{8cm}
\begin{align*}
C_0' =& \; [14, [14, 0], 10, [5, 2], 1, [2, -1], 3, [6, 0], 9, [3, -1], 6, [14, 2]] \\
C_1' =& \; [8, [4, 1], 12, [1, 2], 13, [5, 0]] \\
C_2' =& \; [0, [7, 0], 7, [7, 2]] \\
C_3' =& \; [2, [2, 0], 4, [2, 2]] \\
C_4' =& \; [11, [6, -1], 5, [6, 1]]
\end{align*}
\end{minipage}

\newpage

\item 2-factor type $[5, 4, 2, 2, 2]$: three starters \vspace*{-3mm} \\
\hspace*{-4mm}
\begin{minipage}{8cm}
\begin{align*}
C_0 =& \; [2, [4, 0], 6, [3, 2], 3, [2, 1], 1, [3, 0], 4, [2, 2]] \\
C_1 =& \; [0, [1, 2], 13, [4, 1], 9, [1, 0], 8, [6, 1]] \\
C_2 =& \; [12, [14, 0], 14, [14, 2]] \\
C_3 =& \; [11, [6, 0], 5, [6, 2]] \\
C_4 =& \; [10, [3, -1], 7, [3, 1]]
\end{align*}
\end{minipage}

\vspace*{-4mm}
\begin{minipage}{8cm}
\begin{align*}
C_0' =& \; [5, [2, 0], 3, [2, -1], 1, [1, 1], 2, [6, -1], 10, [5, 2]] \\
C_1' =& \; [4, [4, 2], 8, [4, -1], 12, [1, -1], 13, [5, 0]] \\
C_2' =& \; [0, [7, 0], 7, [7, 2]] \\
C_3' =& \; [14, [14, -1], 9, [14, 1]] \\
C_4' =& \; [6, [5, -1], 11, [5, 1]]
\end{align*}
\end{minipage}

\item 2-factor type $[5, 3, 3, 2, 2]$: three starters \vspace*{-3mm} \\
\begin{minipage}{8cm}
\begin{align*}
C_0 =& \; [12, [5, 0], 7, [6, 1], 13, [2, 2], 11, [3, 0], 8, [4, 2]] \\
C_1 =& \; [6, [3, 1], 3, [2, 0], 5, [1, 2]] \\
C_2 =& \; [14, [14, -1], 0, [4, -1], 10, [14, 1]] \\
C_3 =& \; [1, [6, 0], 9, [6, 2]] \\
C_4 =& \; [4, [2, -1], 2, [2, 1]]
\end{align*}
\end{minipage}

\vspace*{-4mm}
\begin{minipage}{8cm}
\begin{align*}
C_0' =& \; [2, [1, 1], 1, [3, -1], 4, [6, -1], 10, [14, 0], 14, [14, 2]] \\
C_1' =& \; [13, [4, 1], 3, [5, 2], 12, [1, 0]] \\
C_2' =& \; [8, [1, -1], 9, [4, 0], 5, [3, 2]] \\
C_3' =& \; [0, [7, 0], 7, [7, 2]] \\
C_4' =& \; [11, [5, -1], 6, [5, 1]]
\end{align*}
\end{minipage}

\item 2-factor type $[4, 4, 3, 2, 2]$: three starters \vspace*{-3mm} \\
\hspace*{-5mm}
\begin{minipage}{8cm}
\begin{align*}
C_0 =& \; [7, [14, 2], 14, [14, 0], 2, [1, -1], 1, [6, 1]] \\
C_1 =& \; [0, [4, 0], 4, [4, 1], 8, [5, 1], 13, [1, 2]] \\
C_2 =& \; [5, [6, 0], 11, [2, -1], 9, [4, 2]] \\
C_3 =& \; [10, [2, 0], 12, [2, 2]] \\
C_4 =& \; [6, [3, 0], 3, [3, 2]]
\end{align*}
\end{minipage}
\begin{minipage}{8cm}
\begin{align*}
C_0' =& \; [9, [6, -1], 1, [1, 0], 2, [5, 2], 11, [2, 1]] \\
C_1' =& \; [4, [6, 2], 12, [5, -1], 3, [5, 0], 8, [4, -1]] \\
C_2' =& \; [6, [1, 1], 5, [14, -1], 14, [14, 1]] \\
C_3' =& \; [0, [7, 0], 7, [7, 2]] \\
C_4' =& \; [13, [3, -1], 10, [3, 1]]
\end{align*}
\end{minipage}

\item 2-factor type $[4, 3, 3, 3, 2]$: three starters \vspace*{-3mm} \\
\hspace*{-6mm}
\begin{minipage}{8cm}
\begin{align*}
C_0 =& \; [12, [1, 1], 11, [3, -1], 0, [5, 0], 9, [3, 2]] \\
C_1 =& \; [2, [3, 0], 5, [4, 2], 1, [1, -1]] \\
C_2 =& \; [6, [3, 1], 3, [5, -1], 8, [2, 1]] \\
C_3 =& \; [10, [6, -1], 4, [14, 1], 14, [14, -1]] \\
C_4 =& \; [7, [6, 0], 13, [6, 2]]
\end{align*}
\end{minipage}
\begin{minipage}{8cm}
\begin{align*}
C_0' =& \; [6, [5, 2], 11, [1, 0], 12, [14, 2], 14, [14, 0]] \\
C_1' =& \; [1, [2, 0], 13, [4, 1], 3, [2, 2]] \\
C_2' =& \; [10, [1, 2], 9, [4, 0], 5, [5, 1]] \\
C_3' =& \; [2, [6, 1], 8, [4, -1], 4, [2, -1]] \\
C_4' =& \; [0, [7, 0], 7, [7, 2]]
\end{align*}
\end{minipage}

\newpage

\item 2-factor type $[9, 2, 2, 2]$: three starters \vspace*{-3mm} \\
\begin{minipage}{8cm}
\begin{align*}
C_0 =& \; [14, [14, 2], 4, [6, 0], 12, [2, -1], 10, [3, 1], 13, [4, -1], 9, [6, 2], 3, [4, 0], 7, [4, 2], 11, [14, 0]] \\
C_1 =& \; [5, [3, 0], 8, [3, 2]] \\
C_2 =& \; [6, [5, -1], 1, [5, 1]] \\
C_3 =& \; [0, [2, 0], 2, [2, 2]]
\end{align*}
\end{minipage}

\vspace*{-4mm}
\begin{minipage}{8cm}
\begin{align*}
C_0' =& \; [9, [3, -1], 12, [1, 0], 11, [5, 2], 6, [14, -1], 14, [14, 1], 5, [4, 1], 1, [2, 1], 13, [5, 0], 8, [1, 2]] \\
C_1' =& \; [0, [7, 0], 7, [7, 2]] \\
C_2' =& \; [4, [1, -1], 3, [1, 1]] \\
C_3' =& \; [10, [6, -1], 2, [6, 1]]
\end{align*}
\end{minipage}

\item 2-factor type $[8, 3, 2, 2]$: three starters \vspace*{-3mm} \\
\begin{minipage}{8cm}
\begin{align*}
C_0 =& \; [10, [5, 2], 5, [2, -1], 7, [1, 0], 6, [2, 2], 4, [3, 0], 1, [1, 2], 2, [4, 1], 12, [2, 0]] \\
C_1 =& \; [11, [6, -1], 3, [5, 0], 8, [3, 2]] \\
C_2 =& \; [9, [4, 0], 13, [4, 2]] \\
C_3 =& \; [0, [14, 0], 14, [14, 2]]
\end{align*}
\end{minipage}

\vspace*{-4mm}
\begin{minipage}{8cm}
\begin{align*}
C_0' =& \; [2, [1, -1], 3, [6, 0], 11, [2, 1], 13, [5, -1], 8, [3, -1], 5, [1, 1], 6, [6, 2], 12, [4, -1]] \\
C_1' =& \; [10, [6, 1], 4, [3, 1], 1, [5, 1]] \\
C_2' =& \; [0, [7, 0], 7, [7, 2]] \\
C_3' =& \; [14, [14, -1], 9, [14, 1]]
\end{align*}
\end{minipage}

\item 2-factor type $[7, 4, 2, 2]$: three starters \vspace*{-3mm} \\
\begin{minipage}{8cm}
\begin{align*}
C_0 =& \; [5, [4, 1], 1, [2, 1], 13, [5, 1], 8, [1, 0], 7, [1, -1], 6, [3, 1], 9, [4, 2]] \\
C_1 =& \; [10, [4, 0], 0, [2, 2], 2, [5, 0], 11, [1, 2]] \\
C_2 =& \; [14, [14, -1], 3, [14, 1]] \\
C_3 =& \; [4, [6, -1], 12, [6, 1]]
\end{align*}
\end{minipage}

\vspace*{-4mm}
\begin{minipage}{8cm}
\begin{align*}
C_0' =& \; [3, [1, 1], 2, [6, 0], 10, [2, -1], 8, [4, -1], 4, [5, 2], 13, [2, 0], 11, [6, 2]] \\
C_1' =& \; [6, [3, 0], 9, [3, 2], 12, [3, -1], 1, [5, -1]] \\
C_2' =& \; [0, [7, 0], 7, [7, 2]] \\
C_3' =& \; [14, [14, 0], 5, [14, 2]]
\end{align*}
\end{minipage}

\item 2-factor type $[6, 5, 2, 2]$: three starters \vspace*{-3mm} \\
\begin{minipage}{8cm}
\begin{align*}
C_0 =& \; [14, [14, 2], 2, [3, 0], 5, [6, 1], 11, [5, -1], 6, [2, 2], 8, [14, 0]] \\
C_1 =& \; [3, [2, 1], 1, [5, 0], 10, [4, 2], 0, [4, 0], 4, [1, 2]] \\
C_2 =& \; [13, [6, 0], 7, [6, 2]] \\
C_3 =& \; [12, [3, -1], 9, [3, 1]]
\end{align*}
\end{minipage}

\vspace*{-4mm}
\begin{minipage}{8cm}
\begin{align*}
C_0' =& \; [3, [3, 2], 6, [2, 0], 8, [4, -1], 4, [14, 1], 14, [14, -1], 12, [5, 1]] \\
C_1' =& \; [10, [5, 2], 5, [6, -1], 11, [2, -1], 13, [4, 1], 9, [1, 0]] \\
C_2' =& \; [0, [7, 0], 7, [7, 2]] \\
C_3' =& \; [1, [1, -1], 2, [1, 1]]
\end{align*}
\end{minipage}

\newpage

\item 2-factor type $[7, 3, 3, 2]$: three starters \vspace*{-3mm} \\
\begin{minipage}{8cm}
\begin{align*}
C_0 =& \; [0, [4, 0], 10, [6, 2], 4, [6, 1], 12, [5, 0], 3, [14, 2], 14, [14, 1], 13, [1, -1]] \\
C_1 =& \; [5, [6, 0], 11, [2, 2], 9, [4, 1]] \\
C_2 =& \; [2, [6, -1], 8, [1, 0], 7, [5, 2]] \\
C_3 =& \; [1, [5, -1], 6, [5, 1]]
\end{align*}
\end{minipage}

\vspace*{-4mm}
\begin{minipage}{8cm}
\begin{align*}
C_0' =& \; [6, [2, -1], 8, [3, -1], 11, [14, 0], 14, [14, -1], 9, [4, 2], 5, [3, 1], 2, [4, -1]] \\
C_1' =& \; [3, [1, 2], 4, [3, 0], 1, [2, 1]] \\
C_2' =& \; [10, [3, 2], 13, [1, 1], 12, [2, 0]] \\
C_3' =& \; [0, [7, 0], 7, [7, 2]]
\end{align*}
\end{minipage}

\item 2-factor type $[6, 4, 3, 2]$: three starters \vspace*{-3mm} \\
\begin{minipage}{8cm}
\begin{align*}
C_0 =& \; [9, [14, 1], 14, [14, -1], 10, [3, 1], 13, [2, 2], 11, [4, -1], 1, [6, 0]] \\
C_1 =& \; [2, [2, -1], 4, [1, -1], 5, [2, 0], 3, [1, 2]] \\
C_2 =& \; [7, [1, 0], 8, [4, 2], 12, [5, 1]] \\
C_3 =& \; [6, [6, -1], 0, [6, 1]]
\end{align*}
\end{minipage}

\vspace*{-4mm}
\begin{minipage}{8cm}
\begin{align*}
C_0' =& \; [1, [6, 2], 9, [5, 0], 4, [5, -1], 13, [3, 2], 2, [3, -1], 5, [4, 0]] \\
C_1' =& \; [3, [3, 0], 6, [2, 1], 8, [4, 1], 12, [5, 2]] \\
C_2' =& \; [10, [14, 0], 14, [14, 2], 11, [1, 1]] \\
C_3' =& \; [0, [7, 0], 7, [7, 2]]
\end{align*}
\end{minipage}

\item 2-factor type $[5, 5, 3, 2]$: three starters \vspace*{-3mm} \\
\begin{minipage}{8cm}
\begin{align*}
C_0 =& \; [1, [6, -1], 7, [5, 0], 2, [2, -1], 0, [3, 1], 3, [2, 2]] \\
C_1 =& \; [5, [1, -1], 6, [4, -1], 10, [2, 1], 8, [1, 0], 9, [4, 2]] \\
C_2 =& \; [4, [5, -1], 13, [1, 2], 12, [6, 0]] \\
C_3 =& \; [14, [14, -1], 11, [14, 1]]
\end{align*}
\end{minipage}

\vspace*{-4mm}
\begin{minipage}{8cm}
\begin{align*}
C_0' =& \; [4, [5, 2], 13, [14, 0], 14, [14, 2], 9, [1, 1], 8, [4, 0]] \\
C_1' =& \; [11, [6, 1], 5, [3, 0], 2, [4, 1], 6, [3, -1], 3, [6, 2]] \\
C_2' =& \; [10, [2, 0], 12, [3, 2], 1, [5, 1]] \\
C_3' =& \; [0, [7, 0], 7, [7, 2]]
\end{align*}
\end{minipage}

\item 2-factor type $[5, 4, 4, 2]$: three starters \vspace*{-3mm} \\
\begin{minipage}{8cm}
\begin{align*}
C_0 =& \; [1, [1, 2], 2, [4, 0], 12, [5, 2], 7, [14, -1], 14, [14, 0]] \\
C_1 =& \; [11, [3, 2], 8, [2, 0], 10, [5, -1], 5, [6, 1]] \\
C_2 =& \; [13, [4, 2], 3, [6, -1], 9, [5, 0], 0, [1, -1]] \\
C_3 =& \; [6, [2, -1], 4, [2, 1]]
\end{align*}
\end{minipage}

\vspace*{-4mm}
\begin{minipage}{8cm}
\begin{align*}
C_0' =& \; [10, [4, 1], 6, [5, 1], 1, [3, -1], 4, [6, 0], 12, [2, 2]] \\
C_1' =& \; [13, [6, 2], 5, [3, 1], 2, [1, 0], 3, [4, -1]] \\
C_2' =& \; [8, [3, 0], 11, [14, 2], 14, [14, 1], 9, [1, 1]] \\
C_3' =& \; [0, [7, 0], 7, [7, 2]]
\end{align*}
\end{minipage}

\newpage

\item 2-factor type $[11, 2, 2]$: three starters \vspace*{-3mm} \\
\begin{minipage}{8cm}
\begin{align*}
C_0 =& \; [3, [3, 0], 6, [1, 1], 7, [2, 2], 9, [4, 0], 13, [2, 1], 1, [4, 2], 5, [5, 0], 10, [6, -1], 4, [2, -1], 2, [4, -1], 12, [5, 2]] \\
C_1 =& \; [8, [14, -1], 14, [14, 1]] \\
C_2 =& \; [11, [3, -1], 0, [3, 1]]
\end{align*}
\end{minipage}

\vspace*{-4mm}
\begin{minipage}{8cm}
\begin{align*}
C_0' =& \; [6, [5, 1], 1, [1, -1], 2, [6, 0], 8, [3, 2], 5, [1, 0], 4, [6, 2], 12, [5, -1], 3, [6, 1], 11, [2, 0], 9, [1, 2], 10, [4, 1]] \\
C_1' =& \; [0, [7, 0], 7, [7, 2]] \\
C_2' =& \; [14, [14, 0], 13, [14, 2]]
\end{align*}
\end{minipage}

\item 2-factor type $[10, 3, 2]$: three starters \vspace*{-3mm} \\
\begin{minipage}{8cm}
\begin{align*}
C_0 =& \; [3, [3, 1], 0, [1, 1], 13, [14, 0], 14, [14, 2], 8, [4, 1], 4, [5, 0], 9, [3, -1], 6, [5, 2], 1, [4, -1], 5, [2, 1]] \\
C_1 =& \; [11, [1, -1], 12, [2, 0], 10, [1, 2]] \\
C_2 =& \; [2, [5, -1], 7, [5, 1]]
\end{align*}
\end{minipage}

\vspace*{-4mm}
\begin{minipage}{8cm}
\begin{align*}
C_0' =& \; [14, [14, 1], 12, [3, 0], 1, [2, 2], 3, [6, 1], 11, [6, 0], 5, [3, 2], 2, [6, -1], 8, [1, 0], 9, [4, 2], 13, [14, -1]] \\
C_1' =& \; [10, [6, 2], 4, [2, -1], 6, [4, 0]] \\
C_2' =& \; [0, [7, 0], 7, [7, 2]]
\end{align*}
\end{minipage}

\item 2-factor type $[9, 4, 2]$: three starters \vspace*{-3mm} \\
\begin{minipage}{8cm}
\begin{align*}
C_0 =& \; [14, [14, 1], 0, [1, 0], 1, [2, 1], 3, [2, 2], 5, [2, 0], 7, [6, 1], 13, [1, 2], 12, [1, 1], 11, [14, -1]] \\
C_1 =& \; [8, [6, 2], 2, [6, 0], 10, [6, -1], 4, [4, 1]] \\
C_2 =& \; [9, [3, 0], 6, [3, 2]]
\end{align*}
\end{minipage}

\vspace*{-4mm}
\begin{minipage}{8cm}
\begin{align*}
C_0' =& \; [8, [5, -1], 3, [3, 1], 6, [5, 1], 11, [2, -1], 9, [4, 2], 13, [14, 0], 14, [14, 2], 2, [4, 0], 12, [4, -1]] \\
C_1' =& \; [4, [3, -1], 1, [5, 2], 10, [5, 0], 5, [1, -1]] \\
C_2' =& \; [0, [7, 0], 7, [7, 2]]
\end{align*}
\end{minipage}

\item 2-factor type $[8, 5, 2]$: three starters \vspace*{-3mm} \\
\begin{minipage}{8cm}
\begin{align*}
C_0 =& \; [14, [14, 1], 7, [5, 1], 2, [3, 0], 5, [6, 1], 11, [3, -1], 8, [2, 2], 10, [2, -1], 12, [14, -1]] \\
C_1 =& \; [4, [5, -1], 13, [4, 1], 3, [6, 2], 9, [3, 1], 6, [2, 0]] \\
C_2 =& \; [1, [1, -1], 0, [1, 1]]
\end{align*}
\end{minipage}

\vspace*{-4mm}
\begin{minipage}{8cm}
\begin{align*}
C_0' =& \; [12, [3, 2], 1, [6, 0], 9, [1, 2], 10, [5, 0], 5, [6, -1], 13, [4, 2], 3, [1, 0], 2, [4, -1]] \\
C_1' =& \; [6, [5, 2], 11, [14, 0], 14, [14, 2], 8, [4, 0], 4, [2, 1]] \\
C_2' =& \; [0, [7, 0], 7, [7, 2]
\end{align*}
\end{minipage}

\item 2-factor type $[7, 6, 2]$: three starters \vspace*{-3mm} \\
\begin{minipage}{8cm}
\begin{align*}
C_0 =& \; [8, [1, 1], 9, [5, 2], 4, [5, 1], 13, [1, 0], 0, [6, 2], 6, [14, -1], 14, [14, 0]] \\
C_1 =& \; [5, [4, 2], 1, [4, 1], 11, [5, -1], 2, [1, -1], 3, [4, -1], 7, [2, 0]] \\
C_2 =& \; [10, [2, -1], 12, [2, 1]]
\end{align*}
\end{minipage}

\vspace*{-4mm}
\begin{minipage}{8cm}
\begin{align*}
C_0' =& \; [1, [6, 0], 9, [14, 1], 14, [14, 2], 3, [5, 0], 8, [3, 1], 11, [1, 2], 12, [3, -1]] \\
C_1' =& \; [4, [2, 2], 6, [4, 0], 2, [3, 2], 5, [6, 1], 13, [3, 0], 10, [6, -1]] \\
C_2' =& \; [0, [7, 0], 7, [7, 2]]
\end{align*}
\end{minipage}

\end{itemizenew}

\section{Computational results for $n=16$}\label{app:16}

\begin{itemizenew}
\item 2-factor type $[4, 2, 2, 2, 2, 2, 2]$: one starter \vspace*{-3mm} \\
\hspace*{-8mm}
\begin{minipage}{8cm}
\begin{align*}
C_0 =& \; [5, [6, 0], 11, [7, 1], 4, [6, 1], 13, [7, 0]] \\
C_1 =& \; [15, [15, 0], 3, [15, 1]] \\
C_2 =& \; [8, [1, 0], 9, [1, 1]] \\
C_3 =& \; [6, [4, 0], 10, [4, 1]] \\
C_4 =& \; [0, [3, 0], 12, [3, 1]] \\
C_5 =& \; [2, [5, 0], 7, [5, 1]] \\
C_6 =& \; [1, [2, 0], 14, [2, 1]]
\end{align*}
\end{minipage}

\item 2-factor type $[3, 3, 2, 2, 2, 2, 2]$: two starters \vspace*{-3mm} \\
\hspace*{-12mm}
\begin{minipage}{8cm}
\begin{align*}
C_0 =& \; [6, [2, 0], 8, [2, 2], 10, [4, 1]] \\
C_1 =& \; [13, [4, -1], 2, [6, 1], 11, [2, -1]] \\
C_2 =& \; [14, [1, -1], 0, [1, 1]] \\
C_3 =& \; [7, [3, -1], 4, [3, 1]] \\
C_4 =& \; [3, [15, -1], 15, [15, 1]] \\
C_5 =& \; [12, [7, 0], 5, [7, 2]] \\
C_6 =& \; [1, [7, -1], 9, [7, 1]]
\end{align*}
\end{minipage}
\begin{minipage}{8cm}
\begin{align*}
C_0' =& \; [0, [15, 2], 15, [15, 0], 13, [2, 1]] \\
C_1' =& \; [11, [3, 0], 8, [3, 2], 5, [6, -1]] \\
C_2' =& \; [6, [1, 0], 7, [1, 2]] \\
C_3' =& \; [1, [6, 0], 10, [6, 2]] \\
C_4' =& \; [3, [4, 0], 14, [4, 2]] \\
C_5' =& \; [9, [5, -1], 4, [5, 1]] \\
C_6' =& \; [2, [5, 0], 12, [5, 2]]
\end{align*}
\end{minipage}

\item 2-factor type $[6, 2, 2, 2, 2, 2]$: two starters \vspace*{-3mm} \\
\begin{minipage}{8cm}
\begin{align*}
C_0 =& \; [1, [7, 1], 9, [1, 0], 8, [4, 1], 12, [2, 2], 14, [4, -1], 3, [2, 1]] \\
C_1 =& \; [0, [7, 0], 7, [7, 2]] \\
C_2 =& \; [6, [4, 0], 10, [4, 2]] \\
C_3 =& \; [2, [15, -1], 15, [15, 1]] \\
C_4 =& \; [11, [6, 0], 5, [6, 2]] \\
C_5 =& \; [4, [6, -1], 13, [6, 1]]
\end{align*}
\end{minipage}

\vspace*{-4mm}
\begin{minipage}{8cm}
\begin{align*}
C_0' =& \; [14, [7, -1], 6, [2, -1], 8, [2, 0], 10, [1, 2], 9, [5, 0], 4, [5, 2]] \\
C_1' =& \; [7, [15, 0], 15, [15, 2]] \\
C_2' =& \; [12, [1, -1], 13, [1, 1]] \\
C_3' =& \; [3, [3, -1], 0, [3, 1]] \\
C_4' =& \; [11, [5, -1], 1, [5, 1]] \\
C_5' =& \; [2, [3, 0], 5, [3, 2]]
\end{align*}
\end{minipage}

\item 2-factor type $[5, 3, 2, 2, 2, 2]$: one starter \vspace*{-3mm} \\
\begin{minipage}{8cm}
\begin{align*}
C_0 =& \; [4, [2, 1], 2, [4, 1], 13, [2, 0], 11, [4, 0], 7, [3, 1]] \\
C_1 =& \; [3, [3, 0], 0, [15, 1], 15, [15, 0]] \\
C_2 =& \; [1, [5, 0], 6, [5, 1]] \\
C_3 =& \; [8, [6, 0], 14, [6, 1]] \\
C_4 =& \; [5, [7, 0], 12, [7, 1]] \\
C_5 =& \; [10, [1, 0], 9, [1, 1]]
\end{align*}
\end{minipage}

\newpage

\item 2-factor type $[4, 4, 2, 2, 2, 2]$: one starter \vspace*{-3mm} \\
\hspace*{-9mm}
\begin{minipage}{8cm}
\begin{align*}
C_0 =& \; [2, [5, 0], 12, [7, 0], 5, [1, 1], 6, [4, 1]] \\
C_1 =& \; [3, [5, 1], 8, [7, 1], 0, [1, 0], 14, [4, 0]] \\
C_2 =& \; [9, [2, 0], 7, [2, 1]] \\
C_3 =& \; [15, [15, 0], 11, [15, 1]] \\
C_4 =& \; [13, [3, 0], 1, [3, 1]] \\
C_5 =& \; [10, [6, 0], 4, [6, 1]]
\end{align*}
\end{minipage}

\item 2-factor type $[4, 3, 3, 2, 2, 2]$: two starters \vspace*{-3mm} \\
\hspace*{-9mm}
\begin{minipage}{8cm}
\begin{align*}
C_0 =& \; [12, [5, 1], 2, [7, 1], 9, [6, 2], 0, [3, 0]] \\
C_1 =& \; [15, [15, -1], 6, [1, 2], 5, [15, 0]] \\
C_2 =& \; [1, [3, 1], 13, [6, 0], 4, [3, 2]] \\
C_3 =& \; [8, [5, 0], 3, [5, 2]] \\
C_4 =& \; [14, [7, 0], 7, [7, 2]] \\
C_5 =& \; [10, [1, -1], 11, [1, 1]]
\end{align*}
\end{minipage}
\begin{minipage}{8cm}
\begin{align*}
C_0' =& \; [0, [5, -1], 10, [4, 2], 6, [7, -1], 13, [2, 0]] \\
C_1' =& \; [9, [15, 2], 15, [15, 1], 5, [4, 0]] \\
C_2' =& \; [1, [2, 2], 3, [1, 0], 4, [3, -1]] \\
C_3' =& \; [7, [4, -1], 11, [4, 1]] \\
C_4' =& \; [12, [2, -1], 14, [2, 1]] \\
C_5' =& \; [8, [6, -1], 2, [6, 1]]
\end{align*}
\end{minipage}

\item 2-factor type $[3, 3, 3, 3, 2, 2]$: one starter \vspace*{-3mm} \\
\hspace*{-14mm}
\begin{minipage}{8cm}
\begin{align*}
C_0 =& \; [0, [6, 0], 9, [15, 0], 15, [15, 1]] \\
C_1 =& \; [12, [6, 1], 3, [2, 0], 5, [7, 0]] \\
C_2 =& \; [8, [3, 0], 11, [2, 1], 13, [5, 0]] \\
C_3 =& \; [14, [5, 1], 4, [3, 1], 7, [7, 1]] \\
C_4 =& \; [1, [1, 0], 2, [1, 1]] \\
C_5 =& \; [6, [4, 0], 10, [4, 1]]
\end{align*}
\end{minipage}

\item 2-factor type $[8, 2, 2, 2, 2]$: one starter \vspace*{-3mm} \\
\hspace*{-1mm}
\begin{minipage}{8cm}
\begin{align*}
C_0 =& \; [13, [4, 1], 2, [7, 0], 9, [3, 1], 12, [3, 0], 0, [4, 0], 11, [15, 1], 15, [15, 0], 6, [7, 1]] \\
C_1 =& \; [5, [2, 0], 3, [2, 1]] \\
C_2 =& \; [4, [5, 0], 14, [5, 1]] \\
C_3 =& \; [8, [1, 0], 7, [1, 1]] \\
C_4 =& \; [1, [6, 0], 10, [6, 1]]
\end{align*}
\end{minipage}

\item 2-factor type $[7, 3, 2, 2, 2]$: two starters \vspace*{-3mm} \\
\begin{minipage}{8cm}
\begin{align*}
C_0 =& \; [15, [15, 2], 6, [4, -1], 10, [5, 1], 5, [4, 0], 1, [6, 2], 7, [2, -1], 9, [15, 0]] \\
C_1 =& \; [2, [6, 0], 8, [4, 2], 12, [5, -1]] \\
C_2 =& \; [4, [6, -1], 13, [6, 1]] \\
C_3 =& \; [14, [3, 0], 11, [3, 2]] \\
C_4 =& \; [0, [3, -1], 3, [3, 1]]
\end{align*}
\end{minipage}

\vspace*{-4mm}
\begin{minipage}{8cm}
\begin{align*}
C_0' =& \; [3, [1, 1], 4, [1, 2], 5, [4, 1], 1, [1, 0], 0, [2, 1], 2, [7, 2], 10, [7, 0]] \\
C_1' =& \; [15, [15, 1], 8, [1, -1], 9, [15, -1]] \\
C_2' =& \; [14, [7, -1], 6, [7, 1]] \\
C_3' =& \; [11, [2, 0], 13, [2, 2]] \\
C_4' =& \; [12, [5, 0], 7, [5, 2]]
\end{align*}
\end{minipage}

\item 2-factor type $[6, 4, 2, 2, 2]$: two starters \vspace*{-3mm} \\
\begin{minipage}{8cm}
\begin{align*}
C_0 =& \; [0, [2, 0], 13, [1, 1], 14, [3, -1], 11, [1, -1], 10, [2, -1], 8, [7, 2]] \\
C_1 =& \; [4, [2, 1], 6, [5, 2], 1, [4, -1], 5, [1, 0]] \\
C_2 =& \; [9, [7, -1], 2, [7, 1]] \\
C_3 =& \; [7, [5, -1], 12, [5, 1]] \\
C_4 =& \; [15, [15, 0], 3, [15, 2]]
\end{align*}
\end{minipage}

\vspace*{-4mm}
\begin{minipage}{8cm}
\begin{align*}
C_0' =& \; [6, [5, 0], 11, [1, 2], 12, [4, 1], 8, [7, 0], 0, [3, 1], 3, [3, 2]] \\
C_1' =& \; [4, [15, 1], 15, [15, -1], 9, [2, 2], 7, [3, 0]] \\
C_2' =& \; [5, [6, 0], 14, [6, 2]] \\
C_3' =& \; [10, [6, -1], 1, [6, 1]] \\
C_4' =& \; [2, [4, 0], 13, [4, 2]]
\end{align*}
\end{minipage}

\item 2-factor type $[5, 5, 2, 2, 2]$: two starters \vspace*{-3mm} \\
\begin{minipage}{8cm}
\begin{align*}
C_0 =& \; [10, [2, 1], 8, [2, 0], 6, [7, -1], 14, [2, 2], 1, [6, 1]] \\
C_1 =& \; [4, [1, 1], 5, [5, 2], 0, [4, 0], 11, [15, 2], 15, [15, 0]] \\
C_2 =& \; [9, [6, 0], 3, [6, 2]] \\
C_3 =& \; [7, [5, -1], 2, [5, 1]] \\
C_4 =& \; [13, [1, 0], 12, [1, 2]]
\end{align*}
\end{minipage}

\vspace*{-4mm}
\begin{minipage}{8cm}
\begin{align*}
C_0' =& \; [3, [2, -1], 5, [15, -1], 15, [15, 1], 9, [6, -1], 0, [3, -1]] \\
C_1' =& \; [2, [4, 2], 13, [5, 0], 8, [3, 1], 11, [1, -1], 10, [7, 1]] \\
C_2' =& \; [6, [7, 0], 14, [7, 2]] \\
C_3' =& \; [12, [4, -1], 1, [4, 1]] \\
C_4' =& \; [4, [3, 0], 7, [3, 2]]
\end{align*}
\end{minipage}

\item 2-factor type $[6, 3, 3, 2, 2]$: one starter \vspace*{-3mm} \\
\hspace*{-2mm}
\begin{minipage}{8cm}
\begin{align*}
C_0 =& \; [1, [6, 1], 10, [15, 0], 15, [15, 1], 5, [6, 0], 11, [7, 1], 3, [2, 1]] \\
C_1 =& \; [4, [4, 1], 0, [7, 0], 8, [4, 0]] \\
C_2 =& \; [12, [1, 0], 13, [1, 1], 14, [2, 0]] \\
C_3 =& \; [7, [5, 0], 2, [5, 1]] \\
C_4 =& \; [6, [3, 0], 9, [3, 1]]
\end{align*}
\end{minipage}

\item 2-factor type $[5, 4, 3, 2, 2]$: one starter \vspace*{-3mm} \\
\hspace*{-2mm}
\begin{minipage}{8cm}
\begin{align*}
C_0 =& \; [9, [6, 1], 3, [4, 0], 7, [15, 0], 15, [15, 1], 8, [1, 1]] \\
C_1 =& \; [5, [1, 0], 6, [5, 0], 11, [4, 1], 0, [5, 1]] \\
C_2 =& \; [1, [3, 0], 13, [6, 0], 4, [3, 1]] \\
C_3 =& \; [12, [2, 0], 14, [2, 1]] \\
C_4 =& \; [2, [7, 0], 10, [7, 1]]
\end{align*}
\end{minipage}

\newpage

\item 2-factor type $[4, 4, 4, 2, 2]$: one starter \vspace*{-3mm} \\
\hspace*{-8mm}
\begin{minipage}{8cm}
\begin{align*}
C_0 =& \; [3, [7, 1], 11, [2, 0], 9, [5, 1], 14, [4, 0]] \\
C_1 =& \; [6, [2, 1], 4, [4, 1], 0, [3, 0], 12, [6, 0]] \\
C_2 =& \; [7, [6, 1], 13, [7, 0], 5, [5, 0], 10, [3, 1]] \\
C_3 =& \; [2, [1, 0], 1, [1, 1]] \\
C_4 =& \; [8, [15, 0], 15, [15, 1]]
\end{align*}
\end{minipage}

\item 2-factor type $[5, 3, 3, 3, 2]$: two starters \vspace*{-3mm} \\
\begin{minipage}{8cm}
\begin{align*}
C_0 =& \; [15, [15, 1], 3, [2, 1], 1, [3, -1], 4, [1, -1], 5, [15, -1]] \\
C_1 =& \; [2, [5, 2], 7, [1, 0], 8, [6, -1]] \\
C_2 =& \; [14, [1, 1], 0, [6, 2], 9, [5, 0]] \\
C_3 =& \; [10, [2, -1], 12, [6, 1], 6, [4, -1]] \\
C_4 =& \; [11, [2, 0], 13, [2, 2]]
\end{align*}
\end{minipage}

\vspace*{-4mm}
\begin{minipage}{8cm}
\begin{align*}
C_0' =& \; [1, [4, 0], 12, [5, -1], 7, [4, 1], 11, [3, 2], 8, [7, 1]] \\
C_1' =& \; [0, [6, 0], 6, [7, -1], 14, [1, 2]] \\
C_2' =& \; [3, [7, 2], 10, [3, 0], 13, [5, 1]] \\
C_3' =& \; [5, [4, 2], 9, [7, 0], 2, [3, 1]] \\
C_4' =& \; [4, [15, 0], 15, [15, 2]]
\end{align*}
\end{minipage}

\item 2-factor type $[4, 4, 3, 3, 2]$: two starters \vspace*{-3mm} \\
\hspace*{-4mm}
\begin{minipage}{8cm}
\begin{align*}
C_0 =& \; [10, [3, 0], 13, [6, 1], 4, [1, 1], 5, [5, 2]] \\
C_1 =& \; [3, [5, -1], 8, [15, 0], 15, [15, -1], 0, [3, 2]] \\
C_2 =& \; [9, [7, 1], 2, [6, 0], 11, [2, 2]] \\
C_3 =& \; [7, [1, -1], 6, [7, 2], 14, [7, 0]] \\
C_4 =& \; [12, [4, 0], 1, [4, 2]]
\end{align*}
\end{minipage}
\begin{minipage}{8cm}
\begin{align*}
C_0' =& \; [15, [15, 1], 7, [6, -1], 13, [5, 0], 8, [15, 2]] \\
C_1' =& \; [10, [4, 1], 14, [5, 1], 4, [7, -1], 12, [2, -1]] \\
C_2' =& \; [0, [2, 1], 2, [1, 2], 1, [1, 0]] \\
C_3' =& \; [11, [6, 2], 5, [4, -1], 9, [2, 0]] \\
C_4' =& \; [3, [3, -1], 6, [3, 1]]
\end{align*}
\end{minipage}

\item 2-factor type $[4, 3, 3, 3, 3]$: one starter \vspace*{-3mm} \\
\hspace*{-6mm}
\begin{minipage}{8cm}
\begin{align*}
C_0 =& \; [9, [15, 0], 15, [15, 1], 12, [5, 1], 7, [2, 0]] \\
C_1 =& \; [13, [3, 0], 1, [1, 1], 2, [4, 0]] \\
C_2 =& \; [0, [7, 0], 8, [2, 1], 6, [6, 0]] \\
C_3 =& \; [5, [1, 0], 4, [6, 1], 10, [5, 0]] \\
C_4 =& \; [11, [3, 1], 14, [4, 1], 3, [7, 1]]
\end{align*}
\end{minipage}

\item 2-factor type $[10, 2, 2, 2]$: two starters \vspace*{-3mm} \\
\begin{minipage}{8cm}
\begin{align*}
C_0 =& \; [8, [5, 1], 13, [2, 2], 0, [3, -1], 3, [7, 1], 11, [6, -1], 2, [2, 0], 4, [5, -1], 14, [15, 2], 15, [15, 0], 12, [4, -1]] \\
C_1 =& \; [6, [1, -1], 5, [1, 1]] \\
C_2 =& \; [9, [7, 0], 1, [7, 2]] \\
C_3 =& \; [10, [3, 0], 7, [3, 2]]
\end{align*}
\end{minipage}

\vspace*{-4mm}
\begin{minipage}{8cm}
\begin{align*}
C_0' =& \; [3, [3, 1], 6, [4, 2], 2, [5, 0], 12, [4, 1], 1, [15, 1], 15, [15, -1], 0, [5, 2], 5, [6, 1], 14, [4, 0], 10, [7, -1]] \\
C_1' =& \; [11, [2, -1], 9, [2, 1]] \\
C_2' =& \; [4, [6, 0], 13, [6, 2]] \\
C_3' =& \; [7, [1, 0], 8, [1, 2]]
\end{align*}
\end{minipage}

\item 2-factor type $[9, 3, 2, 2]$: one starter \vspace*{-3mm} \\
\hspace*{-2mm}
\begin{minipage}{8cm}
\begin{align*}
C_0 =& \; [7, [7, 1], 14, [1, 0], 13, [3, 0], 1, [7, 0], 8, [4, 0], 4, [15, 0], 15, [15, 1], 11, [2, 0], 9, [2, 1]] \\
C_1 =& \; [5, [1, 1], 6, [4, 1], 2, [3, 1]] \\
C_2 =& \; [3, [6, 0], 12, [6, 1]] \\
C_3 =& \; [0, [5, 0], 10, [5, 1]]
\end{align*}
\end{minipage}

\item 2-factor type $[8, 4, 2, 2]$: one starter \vspace*{-3mm} \\
\hspace*{-2mm}
\begin{minipage}{8cm}
\begin{align*}
C_0 =& \; [12, [6, 0], 6, [4, 1], 2, [7, 1], 10, [6, 1], 1, [1, 1], 0, [15, 1], 15, [15, 0], 9, [3, 1]] \\
C_1 =& \; [4, [1, 0], 3, [4, 0], 14, [3, 0], 11, [7, 0]] \\
C_2 =& \; [7, [2, 0], 5, [2, 1]] \\
C_3 =& \; [13, [5, 0], 8, [5, 1]]
\end{align*}
\end{minipage}

\item 2-factor type $[7, 5, 2, 2]$: one starter \vspace*{-3mm} \\
\hspace*{-2mm}
\begin{minipage}{8cm}
\begin{align*}
C_0 =& \; [14, [7, 1], 7, [5, 0], 12, [15, 0], 15, [15, 1], 13, [2, 0], 0, [6, 1], 6, [7, 0]] \\
C_1 =& \; [10, [6, 0], 4, [3, 0], 1, [5, 1], 11, [3, 1], 8, [2, 1]] \\
C_2 =& \; [3, [1, 0], 2, [1, 1]] \\
C_3 =& \; [9, [4, 0], 5, [4, 1]]
\end{align*}
\end{minipage}

\item 2-factor type $[6, 6, 2, 2]$: one starter \vspace*{-3mm} \\
\hspace*{-2mm}
\begin{minipage}{8cm}
\begin{align*}
C_0 =& \; [0, [3, 1], 3, [2, 1], 5, [3, 0], 8, [2, 0], 10, [6, 1], 1, [1, 1]] \\
C_1 =& \; [7, [15, 1], 15, [15, 0], 11, [7, 1], 4, [7, 0], 12, [1, 0], 13, [6, 0]] \\
C_2 =& \; [6, [4, 0], 2, [4, 1]] \\
C_3 =& \; [14, [5, 0], 9, [5, 1]]
\end{align*}
\end{minipage}

\item 2-factor type $[8, 3, 3, 2]$: two starters \vspace*{-3mm} \\
\begin{minipage}{8cm}
\begin{align*}
C_0 =& \; [0, [4, 2], 4, [7, -1], 11, [3, 0], 8, [4, 1], 12, [15, 2], 15, [15, 1], 14, [2, 0], 1, [1, -1]] \\
C_1 =& \; [10, [7, 0], 2, [1, 1], 3, [7, 2]] \\
C_2 =& \; [13, [7, 1], 5, [2, 1], 7, [6, 1]] \\
C_3 =& \; [6, [3, -1], 9, [3, 1]]
\end{align*}
\end{minipage}

\vspace*{-4mm}
\begin{minipage}{8cm}
\begin{align*}
C_0' =& \; [7, [4, 0], 11, [3, 2], 14, [5, 0], 4, [1, 2], 3, [6, -1], 9, [6, 0], 0, [2, -1], 13, [6, 2]] \\
C_1' =& \; [2, [1, 0], 1, [4, -1], 12, [5, 2]] \\
C_2' =& \; [6, [15, -1], 15, [15, 0], 8, [2, 2]] \\
C_3' =& \; [5, [5, -1], 10, [5, 1]]
\end{align*}
\end{minipage}

\item 2-factor type $[7, 4, 3, 2]$: two starters \vspace*{-3mm} \\
\begin{minipage}{8cm}
\begin{align*}
C_0 =& \; [15, [15, -1], 14, [7, 2], 7, [1, 1], 6, [6, 1], 0, [7, 0], 8, [6, 2], 2, [15, 0]] \\
C_1 =& \; [4, [7, -1], 11, [2, 0], 13, [4, -1], 9, [5, 2]] \\
C_2 =& \; [1, [6, -1], 10, [2, 2], 12, [4, 0]] \\
C_3 =& \; [5, [2, -1], 3, [2, 1]]
\end{align*}
\end{minipage}

\vspace*{-4mm}
\begin{minipage}{8cm}
\begin{align*}
C_0' =& \; [0, [4, 2], 11, [1, 0], 10, [1, 2], 9, [5, -1], 14, [3, 0], 2, [3, 2], 5, [5, 0]] \\
C_1' =& \; [13, [5, 1], 8, [4, 1], 4, [7, 1], 12, [1, -1]] \\
C_2' =& \; [7, [6, 0], 1, [15, 1], 15, [15, 2]] \\
C_3' =& \; [6, [3, -1], 3, [3, 1]]
\end{align*}
\end{minipage}

\item 2-factor type $[6, 5, 3, 2]$: two starters \vspace*{-3mm} \\
\begin{minipage}{8cm}
\begin{align*}
C_0 =& \; [10, [5, -1], 5, [3, -1], 2, [2, -1], 0, [3, 1], 3, [6, 2], 9, [1, 0]] \\
C_1 =& \; [11, [3, 0], 8, [7, -1], 1, [5, 1], 6, [7, 2], 13, [2, 1]] \\
C_2 =& \; [7, [3, 2], 4, [7, 1], 12, [5, 0]] \\
C_3 =& \; [15, [15, -1], 14, [15, 1]]
\end{align*}
\end{minipage}

\vspace*{-4mm}
\begin{minipage}{8cm}
\begin{align*}
C_0' =& \; [4, [6, -1], 13, [1, 2], 12, [4, -1], 1, [4, 0], 5, [5, 2], 10, [6, 0]] \\
C_1' =& \; [0, [6, 1], 9, [7, 0], 2, [1, 1], 3, [4, 2], 14, [1, -1]] \\
C_2' =& \; [7, [4, 1], 11, [15, 2], 15, [15, 0]] \\
C_3' =& \; [8, [2, 0], 6, [2, 2]]
\end{align*}
\end{minipage}

\item 2-factor type $[6, 4, 4, 2]$: two starters \vspace*{-3mm} \\
\begin{minipage}{8cm}
\begin{align*}
C_0 =& \; [9, [2, 0], 7, [7, 1], 14, [1, 2], 0, [15, -1], 15, [15, 1], 13, [4, 1]] \\
C_1 =& \; [12, [1, 1], 11, [7, 0], 4, [1, -1], 3, [6, 2]] \\
C_2 =& \; [2, [6, 0], 8, [3, -1], 5, [5, 2], 10, [7, -1]] \\
C_3 =& \; [1, [5, -1], 6, [5, 1]]
\end{align*}
\end{minipage}

\vspace*{-4mm}
\begin{minipage}{8cm}
\begin{align*}
C_0' =& \; [4, [6, -1], 13, [2, 1], 0, [6, 1], 6, [4, 2], 2, [1, 0], 1, [3, 1]] \\
C_1' =& \; [10, [2, 2], 12, [5, 0], 7, [7, 2], 14, [4, 0]] \\
C_2' =& \; [9, [4, -1], 5, [3, 2], 8, [3, 0], 11, [2, -1]] \\
C_3' =& \; [15, [15, 0], 3, [15, 2]]
\end{align*}
\end{minipage}

\item 2-factor type $[5, 5, 4, 2]$: two starters \vspace*{-3mm} \\
\begin{minipage}{8cm}
\begin{align*}
C_0 =& \; [8, [5, 0], 13, [7, -1], 6, [6, 2], 12, [1, 1], 11, [3, 1]] \\
C_1 =& \; [15, [15, 1], 7, [6, 0], 1, [1, 2], 2, [3, -1], 5, [15, -1]] \\
C_2 =& \; [4, [4, -1], 0, [5, 2], 10, [7, 1], 3, [1, 0]] \\
C_3 =& \; [14, [5, -1], 9, [5, 1]]
\end{align*}
\end{minipage}

\vspace*{-4mm}
\begin{minipage}{8cm}
\begin{align*}
C_0' =& \; [12, [1, -1], 11, [15, 2], 15, [15, 0], 5, [3, 2], 8, [4, 0]] \\
C_1' =& \; [1, [3, 0], 4, [6, 1], 10, [7, 2], 3, [4, 1], 14, [2, -1]] \\
C_2' =& \; [6, [7, 0], 13, [4, 2], 2, [2, 1], 0, [6, -1]] \\
C_3' =& \; [7, [2, 0], 9, [2, 2]]
\end{align*}
\end{minipage}

\item 2-factor type $[7, 3, 3, 3]$: one starter \vspace*{-3mm} \\
\hspace*{-2mm}
\begin{minipage}{8cm}
\begin{align*}
C_0 =& \; [12, [7, 0], 5, [2, 1], 7, [15, 1], 15, [15, 0], 14, [1, 1], 0, [6, 0], 9, [3, 0]] \\
C_1 =& \; [13, [4, 1], 2, [1, 0], 3, [5, 0]] \\
C_2 =& \; [8, [3, 1], 11, [5, 1], 1, [7, 1]] \\
C_3 =& \; [4, [6, 1], 10, [4, 0], 6, [2, 0]]
\end{align*}
\end{minipage}

\item 2-factor type $[6, 4, 3, 3]$: one starter \vspace*{-3mm} \\
\hspace*{-2mm}
\begin{minipage}{8cm}
\begin{align*}
C_0 =& \; [11, [4, 0], 0, [5, 1], 10, [7, 1], 2, [15, 1], 15, [15, 0], 14, [3, 1]] \\
C_1 =& \; [7, [2, 0], 5, [7, 0], 12, [4, 1], 1, [6, 1]] \\
C_2 =& \; [9, [1, 0], 8, [2, 1], 6, [3, 0]] \\
C_3 =& \; [4, [1, 1], 3, [5, 0], 13, [6, 0]]
\end{align*}
\end{minipage}

\item 2-factor type $[5, 5, 3, 3]$: one starter \vspace*{-3mm} \\
\hspace*{-2mm}
\begin{minipage}{8cm}
\begin{align*}
C_0 =& \; [1, [15, 0], 15, [15, 1], 12, [7, 0], 5, [6, 1], 14, [2, 1]] \\
C_1 =& \; [10, [4, 1], 6, [6, 0], 0, [3, 1], 3, [7, 1], 11, [1, 0]] \\
C_2 =& \; [2, [2, 0], 4, [3, 0], 7, [5, 1]] \\
C_3 =& \; [8, [1, 1], 9, [4, 0], 13, [5, 0]]
\end{align*}
\end{minipage}

\item 2-factor type $[5, 4, 4, 3]$: one starter \vspace*{-3mm} \\
\hspace*{-2mm}
\begin{minipage}{8cm}
\begin{align*}
C_0 =& \; [12, [6, 1], 3, [6, 0], 9, [4, 0], 5, [7, 1], 13, [1, 1]] \\
C_1 =& \; [0, [15, 1], 15, [15, 0], 2, [5, 1], 7, [7, 0]] \\
C_2 =& \; [8, [2, 1], 6, [4, 1], 10, [1, 0], 11, [3, 0]] \\
C_3 =& \; [14, [2, 0], 1, [3, 1], 4, [5, 0]]
\end{align*}
\end{minipage}

\item 2-factor type $[12, 2, 2]$: one starter \vspace*{-3mm} \\
\hspace*{-2mm}
\begin{minipage}{8cm}
\begin{align*}
C_0 =& \; [15, [15, 1], 14, [6, 1], 5, [5, 0], 0, [2, 1], 13, [7, 1], 6, [5, 1], 11, [7, 0], 4, [2, 0], 2, [1, 1], 3, [6, 0], 9, [1, 0], \\ & \; 8, [15, 0]] \\
C_1 =& \; [7, [3, 0], 10, [3, 1]] \\
C_2 =& \; [1, [4, 0], 12, [4, 1]]
\end{align*}
\end{minipage}

\item 2-factor type $[11, 3, 2]$: two starters \vspace*{-3mm} \\
\begin{minipage}{8cm}
\begin{align*}
C_0 =& \; [6, [1, 2], 5, [5, 0], 0, [3, 2], 3, [4, 0], 14, [2, -1], 12, [5, 2], 7, [6, 0], 13, [2, 2], 11, [15, 0], 15, [15, 2], 8, [2, 0]] \\
C_1 =& \; [2, [7, 1], 10, [1, 0], 9, [7, 2]] \\
C_2 =& \; [4, [3, -1], 1, [3, 1]]
\end{align*}
\end{minipage}

\vspace*{-4mm}
\begin{minipage}{8cm}
\begin{align*}
C_0' =& \; [5, [6, 1], 11, [4, 1], 0, [7, -1], 8, [1, -1], 7, [2, 1], 9, [4, 2], 13, [6, -1], 4, [7, 0], 12, [6, 2], 3, [1, 1], 2, [3, 0]] \\
C_1' =& \; [14, [15, -1], 15, [15, 1], 10, [4, -1]] \\
C_2' =& \; [6, [5, -1], 1, [5, 1]]
\end{align*}
\end{minipage}

\item 2-factor type $[10, 4, 2]$: two starters \vspace*{-3mm} \\
\begin{minipage}{8cm}
\begin{align*}
C_0 =& \; [14, [7, 2], 6, [4, 0], 2, [4, -1], 13, [2, 2], 0, [5, -1], 5, [2, 0], 7, [5, 1], 12, [7, 1], 4, [3, -1], 1, [2, -1]] \\
C_1 =& \; [3, [5, 2], 8, [15, -1], 15, [15, 1], 11, [7, 0]] \\
C_2 =& \; [10, [1, 0], 9, [1, 2]]
\end{align*}
\end{minipage}

\vspace*{-4mm}
\begin{minipage}{8cm}
\begin{align*}
C_0' =& \; [0, [4, 1], 4, [6, 2], 13, [5, 0], 3, [6, -1], 12, [4, 2], 1, [1, -1], 2, [7, -1], 9, [2, 1], 7, [1, 1], 6, [6, 0]] \\
C_1' =& \; [14, [6, 1], 5, [3, 1], 8, [3, 2], 11, [3, 0]] \\
C_2' =& \; [10, [15, 0], 15, [15, 2]]
\end{align*}
\end{minipage}

\item 2-factor type $[9, 5, 2]$: two starters \vspace*{-3mm} \\
\begin{minipage}{8cm}
\begin{align*}
C_0 =& \; [6, [3, 0], 3, [4, 2], 14, [5, 0], 4, [7, 2], 12, [4, -1], 1, [1, 0], 2, [15, 2], 15, [15, 1], 7, [1, 1]] \\
C_1 =& \; [0, [7, 1], 8, [3, 1], 5, [7, 0], 13, [3, -1], 10, [5, 2]] \\
C_2 =& \; [11, [2, -1], 9, [2, 1]]
\end{align*}
\end{minipage}

\vspace*{-4mm}
\begin{minipage}{8cm}
\begin{align*}
C_0' =& \; [1, [2, 0], 14, [3, 2], 2, [6, 0], 11, [4, 1], 0, [6, 2], 6, [15, -1], 15, [15, 0], 9, [6, -1], 3, [2, 2]] \\
C_1' =& \; [13, [6, 1], 7, [1, -1], 8, [4, 0], 4, [7, -1], 12, [1, 2]] \\
C_2' =& \; [5, [5, -1], 10, [5, 1]]
\end{align*}
\end{minipage}

\item 2-factor type $[8, 6, 2]$: two starters \vspace*{-3mm} \\
\begin{minipage}{8cm}
\begin{align*}
C_0 =& \; [6, [2, 2], 8, [1, 0], 9, [7, 2], 1, [3, 1], 13, [5, 0], 3, [7, -1], 11, [1, 2], 10, [4, 0]] \\
C_1 =& \; [15, [15, 2], 12, [3, 0], 0, [2, 1], 2, [5, 2], 7, [2, 0], 5, [15, 1]] \\
C_2 =& \; [4, [5, -1], 14, [5, 1]]
\end{align*}
\end{minipage}

\vspace*{-4mm}
\begin{minipage}{8cm}
\begin{align*}
C_0' =& \; [0, [7, 0], 7, [4, 2], 11, [1, 1], 10, [7, 1], 3, [3, -1], 6, [2, -1], 8, [1, -1], 9, [6, -1]] \\
C_1' =& \; [2, [15, 0], 15, [15, -1], 12, [4, 1], 1, [3, 2], 4, [6, 1], 13, [4, -1]] \\
C_2' =& \; [5, [6, 0], 14, [6, 2]]
\end{align*}
\end{minipage}

\item 2-factor type $[7, 7, 2]$: two starters \vspace*{-3mm} \\
\begin{minipage}{8cm}
\begin{align*}
C_0 =& \; [10, [1, -1], 9, [4, -1], 5, [2, 2], 3, [4, 0], 14, [15, 2], 15, [15, 0], 1, [6, -1]] \\
C_1 =& \; [0, [6, 0], 6, [4, 2], 2, [4, 1], 13, [1, 1], 12, [5, 0], 7, [1, 2], 8, [7, -1]] \\
C_2 =& \; [4, [7, 0], 11, [7, 2]]
\end{align*}
\end{minipage}

\vspace*{-4mm}
\begin{minipage}{8cm}
\begin{align*}
C_0' =& \; [2, [5, 2], 12, [5, 1], 7, [2, 0], 9, [6, 2], 3, [3, 0], 0, [5, -1], 10, [7, 1]] \\
C_1' =& \; [13, [1, 0], 14, [6, 1], 5, [3, 2], 8, [3, -1], 11, [15, -1], 15, [15, 1], 1, [3, 1]] \\
C_2' =& \; [4, [2, -1], 6, [2, 1]]
\end{align*}
\end{minipage}

\item 2-factor type $[10, 3, 3]$: one starter \vspace*{-3mm} \\
\hspace*{-2mm}
\begin{minipage}{8cm}
\begin{align*}
C_0 =& \; [1, [4, 0], 5, [1, 1], 6, [4, 1], 10, [7, 1], 2, [6, 1], 11, [7, 0], 4, [15, 1], 15, [15, 0], 0, [3, 0], 3, [2, 1]] \\
C_1 =& \; [9, [3, 1], 12, [5, 0], 7, [2, 0]] \\
C_2 =& \; [13, [5, 1], 8, [6, 0], 14, [1, 0]]
\end{align*}
\end{minipage}

\item 2-factor type $[9, 4, 3]$: one starter \vspace*{-3mm} \\
\hspace*{-2mm}
\begin{minipage}{8cm}
\begin{align*}
C_0 =& \; [9, [3, 1], 6, [6, 1], 12, [2, 0], 14, [1, 0], 0, [4, 1], 11, [7, 1], 4, [4, 0], 8, [15, 1], 15, [15, 0]] \\
C_1 =& \; [1, [1, 1], 2, [5, 1], 7, [3, 0], 10, [6, 0]] \\
C_2 =& \; [3, [5, 0], 13, [7, 0], 5, [2, 1]]
\end{align*}
\end{minipage}

\item 2-factor type $[8, 5, 3]$: one starter \vspace*{-3mm} \\
\hspace*{-2mm}
\begin{minipage}{8cm}
\begin{align*}
C_0 =& \; [1, [7, 1], 8, [5, 0], 3, [15, 0], 15, [15, 1], 10, [1, 1], 9, [6, 1], 0, [2, 1], 13, [3, 1]] \\
C_1 =& \; [5, [3, 0], 2, [4, 0], 6, [1, 0], 7, [4, 1], 11, [6, 0]] \\
C_2 =& \; [14, [5, 1], 4, [7, 0], 12, [2, 0]]
\end{align*}
\end{minipage}

\item 2-factor type $[7, 6, 3]$: one starter \vspace*{-3mm} \\
\hspace*{-2mm}
\begin{minipage}{8cm}
\begin{align*}
C_0 =& \; [14, [2, 0], 12, [15, 0], 15, [15, 1], 8, [1, 0], 7, [3, 1], 4, [6, 0], 13, [1, 1]] \\
C_1 =& \; [10, [5, 1], 0, [2, 1], 2, [7, 0], 9, [4, 0], 5, [4, 1], 1, [6, 1]] \\
C_2 =& \; [6, [5, 0], 11, [7, 1], 3, [3, 0]]
\end{align*}
\end{minipage}

\item 2-factor type $[7, 5, 4]$: one starter \vspace*{-3mm} \\
\hspace*{-2mm}
\begin{minipage}{8cm}
\begin{align*}
C_0 =& \; [10, [3, 1], 13, [15, 0], 15, [15, 1], 5, [7, 0], 12, [6, 1], 3, [3, 0], 6, [4, 0]] \\
C_1 =& \; [0, [1, 0], 14, [2, 0], 1, [5, 1], 11, [4, 1], 7, [7, 1]] \\
C_2 =& \; [4, [2, 1], 2, [6, 0], 8, [1, 1], 9, [5, 0]]
\end{align*}
\end{minipage}

\item 2-factor type $[6, 6, 4]$: one starter \vspace*{-3mm} \\
\hspace*{-2mm}
\begin{minipage}{8cm}
\begin{align*}
C_0 =& \; [4, [7, 1], 12, [3, 1], 0, [6, 1], 6, [2, 0], 8, [1, 0], 9, [5, 1]] \\
C_1 =& \; [15, [15, 1], 10, [5, 0], 5, [6, 0], 14, [7, 0], 7, [4, 1], 11, [15, 0]] \\
C_2 =& \; [13, [4, 0], 2, [1, 1], 3, [2, 1], 1, [3, 0]]
\end{align*}
\end{minipage}

\item 2-factor type $[6, 5, 5]$: one starter \vspace*{-3mm} \\
\hspace*{-2mm}
\begin{minipage}{8cm}
\begin{align*}
C_0 =& \; [2, [15, 1], 15, [15, 0], 3, [2, 0], 1, [5, 0], 6, [3, 0], 9, [7, 1]] \\
C_1 =& \; [8, [1, 0], 7, [7, 0], 14, [4, 1], 10, [6, 0], 4, [4, 0]] \\
C_2 =& \; [5, [6, 1], 11, [2, 1], 13, [1, 1], 12, [3, 1], 0, [5, 1]]
\end{align*}
\end{minipage}

\item 2-factor type $[14, 2]$: two starters \vspace*{-3mm} \\
\begin{minipage}{8cm}
\begin{align*}
C_0 =& \; [10, [4, 0], 14, [2, -1], 12, [3, 2], 9, [2, 1], 7, [2, 0], 5, [7, -1], 13, [15, 1], 15, [15, 2], 11, [6, 0], 2, [6, 1], 8, \\ & \; [7, 2], 1, [1, 1], 0, [4, -1], 4, [6, -1]] \\
C_1 =& \; [3, [3, -1], 6, [3, 1]]
\end{align*}
\end{minipage}

\vspace*{-4mm}
\begin{minipage}{8cm}
\begin{align*}
C_0' =& \; [0, [2, 2], 2, [4, 1], 13, [15, -1], 15, [15, 0], 9, [6, 2], 3, [1, -1], 4, [7, 0], 12, [5, -1], 7, [1, 2], 8, [3, 0], 5, \\ & \; [5, 1], 10, [4, 2], 6, [7, 1], 14, [1, 0]] \\
C_1' =& \; [1, [5, 0], 11, [5, 2]]
\end{align*}
\end{minipage}

\item 2-factor type $[13, 3]$: one starter \vspace*{-3mm} \\
\hspace*{-2mm}
\begin{minipage}{8cm}
\begin{align*}
C_0 =& \; [8, [7, 0], 1, [2, 0], 14, [4, 1], 10, [1, 0], 11, [5, 0], 6, [3, 1], 9, [4, 0], 5, [3, 0], 2, [2, 1], 0, [15, 1], 15, [15, 0], 7, \\ & \; [6, 0], 13, [5, 1]] \\
C_1 =& \; [3, [1, 1], 4, [7, 1], 12, [6, 1]]
\end{align*}
\end{minipage}

\item 2-factor type $[11, 5]$: one starter \vspace*{-3mm} \\
\hspace*{-2mm}
\begin{minipage}{8cm}
\begin{align*}
C_0 =& \; [12, [1, 0], 13, [4, 0], 9, [6, 1], 3, [7, 0], 10, [6, 0], 1, [5, 1], 6, [1, 1], 5, [3, 1], 8, [15, 1], 15, [15, 0], 14, [2, 0]] \\
C_1 =& \; [2, [5, 0], 7, [3, 0], 4, [7, 1], 11, [4, 1], 0, [2, 1]]
\end{align*}
\end{minipage}

\item 2-factor type $[10, 6]$: one starter \vspace*{-3mm} \\
\hspace*{-2mm}
\begin{minipage}{8cm}
\begin{align*}
C_0 =& \; [9, [7, 1], 1, [1, 1], 2, [15, 0], 15, [15, 1], 6, [2, 1], 8, [3, 0], 11, [4, 1], 0, [5, 0], 5, [2, 0], 3, [6, 1]] \\
C_1 =& \; [13, [6, 0], 4, [7, 0], 12, [5, 1], 7, [3, 1], 10, [4, 0], 14, [1, 0]]
\end{align*}
\end{minipage}

\item 2-factor type $[9, 7]$: one starter \vspace*{-3mm} \\
\hspace*{-2mm}
\begin{minipage}{8cm}
\begin{align*}
C_0 =& \; [9, [6, 1], 0, [3, 1], 3, [15, 0], 15, [15, 1], 5, [1, 0], 4, [7, 1], 11, [4, 0], 7, [7, 0], 14, [5, 1]] \\
C_1 =& \; [6, [5, 0], 1, [3, 0], 13, [1, 1], 12, [2, 1], 10, [2, 0], 8, [6, 0], 2, [4, 1]]
\end{align*}
\end{minipage}
\end{itemizenew}

\section{Computational results for $n=17$}\label{app:17}

\begin{itemizenew}
\item 2-factor type $[3, 2, 2, 2, 2, 2, 2, 2]$: three starters \vspace*{-3mm} \\
\hspace*{-14mm}
\begin{minipage}{8cm}
\begin{align*}
C_0 =& \;  [11, [6, 0], 1, [7, 2], 10, [1, -1]] \\
C_1 =& \;  [13, [7, -1], 6, [7, 1]] \\
C_2 =& \;  [14, [5, 0], 9, [5, 2]] \\
C_3 =& \;  [15, [6, -1], 5, [6, 1]] \\
C_4 =& \;  [12, [4, 0], 0, [4, 2]] \\
C_5 =& \;  [7, [4, -1], 3, [4, 1]] \\
C_6 =& \;  [2, [2, 0], 4, [2, 2]] \\
C_7 =& \;  [8, [16, 0], 16, [16, 2]]
\end{align*}
\end{minipage}
\begin{minipage}{8cm}
\begin{align*}
C_0' =& \;  [10, [7, 0], 3, [1, 1], 4, [6, 2]] \\
C_1' =& \;  [13, [2, -1], 15, [2, 1]] \\
C_2' =& \;  [12, [1, 0], 11, [1, 2]] \\
C_3' =& \;  [14, [3, -1], 1, [3, 1]] \\
C_4' =& \;  [16, [16, -1], 5, [16, 1]] \\
C_5' =& \;  [7, [5, -1], 2, [5, 1]] \\
C_6' =& \;  [6, [3, 0], 9, [3, 2]] \\
C_7' =& \;  [0, [8, 0], 8, [8, 2]]
\end{align*}
\end{minipage}

\item 2-factor type $[5, 2, 2, 2, 2, 2, 2]$: three starters \vspace*{-3mm} \\
\begin{minipage}{8cm}
\begin{align*}
C_0 =& \;  [7, [6, -1], 13, [3, 0], 10, [5, 2], 15, [1, 0], 14, [7, 2]] \\
C_1 =& \;  [3, [1, -1], 4, [1, 1]] \\
C_2 =& \;  [12, [6, 0], 6, [6, 2]] \\
C_3 =& \;  [2, [7, -1], 9, [7, 1]] \\
C_4 =& \;  [1, [4, -1], 5, [4, 1]] \\
C_5 =& \;  [0, [5, -1], 11, [5, 1]] \\
C_6 =& \;  [16, [16, -1], 8, [16, 1]]
\end{align*}
\end{minipage}

\vspace*{-4mm}
\begin{minipage}{8cm}
\begin{align*}
C_0' =& \;  [12, [3, 2], 15, [7, 0], 6, [1, 2], 7, [6, 1], 1, [5, 0]] \\
C_1' =& \;  [0, [8, 0], 8, [8, 2]] \\
C_2' =& \;  [3, [2, -1], 5, [2, 1]] \\
C_3' =& \;  [13, [3, -1], 10, [3, 1]] \\
C_4' =& \;  [16, [16, 0], 4, [16, 2]] \\
C_5' =& \;  [14, [4, 0], 2, [4, 2]] \\
C_6' =& \;  [9, [2, 0], 11, [2, 2]]
\end{align*}
\end{minipage}

\item 2-factor type $[4, 3, 2, 2, 2, 2, 2]$: three starters \vspace*{-3mm} \\
\hspace*{-8mm}
\begin{minipage}{8cm}
\begin{align*}
C_0 =& \;  [9, [5, 0], 14, [1, 2], 15, [3, 0], 2, [7, 2]] \\
C_1 =& \;  [12, [16, 0], 16, [16, 2], 6, [6, -1]] \\
C_2 =& \;  [0, [5, -1], 11, [5, 1]] \\
C_3 =& \;  [5, [3, -1], 8, [3, 1]] \\
C_4 =& \;  [4, [6, 0], 10, [6, 2]] \\
C_5 =& \;  [13, [4, -1], 1, [4, 1]] \\
C_6 =& \;  [7, [4, 0], 3, [4, 2]]
\end{align*}
\end{minipage}
\begin{minipage}{8cm}
\begin{align*}
C_0' =& \;  [2, [3, 2], 5, [7, 0], 14, [16, 1], 16, [16, -1]] \\
C_1' =& \;  [7, [1, 0], 6, [6, 1], 12, [5, 2]] \\
C_2' =& \;  [0, [8, 0], 8, [8, 2]] \\
C_3' =& \;  [9, [1, -1], 10, [1, 1]] \\
C_4' =& \;  [3, [2, 0], 1, [2, 2]] \\
C_5' =& \;  [4, [7, -1], 11, [7, 1]] \\
C_6' =& \;  [15, [2, -1], 13, [2, 1]]
\end{align*}
\end{minipage}
\item 2-factor type $[3, 3, 3, 2, 2, 2, 2]$: three starters \vspace*{-3mm} \\
\hspace*{-12mm}
\begin{minipage}{8cm}
\begin{align*}
C_0 =& \;  [16, [16, 2], 4, [5, -1], 9, [16, 0]] \\
C_1 =& \;  [1, [4, 2], 13, [1, 1], 12, [5, 0]] \\
C_2 =& \;  [14, [7, 1], 7, [7, 0], 0, [2, 2]] \\
C_3 =& \;  [6, [4, -1], 10, [4, 1]] \\
C_4 =& \;  [8, [3, 0], 11, [3, 2]] \\
C_5 =& \;  [3, [2, -1], 5, [2, 1]] \\
C_6 =& \;  [15, [3, -1], 2, [3, 1]]
\end{align*}
\end{minipage}
\begin{minipage}{8cm}
\begin{align*}
C_0' =& \;  [10, [16, -1], 16, [16, 1], 9, [1, -1]] \\
C_1' =& \;  [11, [2, 0], 13, [7, -1], 6, [5, 2]] \\
C_2' =& \;  [12, [5, 1], 7, [4, 0], 3, [7, 2]] \\
C_3' =& \;  [0, [8, 0], 8, [8, 2]] \\
C_4' =& \;  [1, [1, 0], 2, [1, 2]] \\
C_5' =& \;  [5, [6, 0], 15, [6, 2]] \\
C_6' =& \;  [4, [6, -1], 14, [6, 1]]
\end{align*}
\end{minipage}

\item 2-factor type $[7, 2, 2, 2, 2, 2]$: three starters \vspace*{-3mm} \\
\begin{minipage}{8cm}
\begin{align*}
C_0 =& \;  [9, [1, 0], 8, [1, -1], 7, [5, -1], 2, [1, 2], 1, [3, -1], 4, [7, 0], 11, [2, 2]] \\
C_1 =& \;  [15, [7, -1], 6, [7, 1]] \\
C_2 =& \;  [16, [16, -1], 14, [16, 1]] \\
C_3 =& \;  [5, [2, -1], 3, [2, 1]] \\
C_4 =& \;  [10, [3, 0], 13, [3, 2]] \\
C_5 =& \;  [12, [4, -1], 0, [4, 1]]
\end{align*}
\end{minipage}

\newpage

\begin{minipage}{8cm}
\begin{align*}
C_0 =& \;  [6, [2, 0], 4, [5, 1], 9, [5, 2], 14, [3, 1], 11, [1, 1], 10, [5, 0], 15, [7, 2]] \\
C_1 =& \;  [0, [8, 0], 8, [8, 2]] \\
C_2 =& \;  [2, [6, 0], 12, [6, 2]] \\
C_3 =& \;  [16, [16, 0], 7, [16, 2]] \\
C_4 =& \;  [5, [4, 0], 1, [4, 2]] \\
C_5 =& \;  [13, [6, -1], 3, [6, 1]]
\end{align*}
\end{minipage}

\item 2-factor type $[6, 3, 2, 2, 2, 2]$: three starters \vspace*{-3mm} \\
\begin{minipage}{8cm}
\begin{align*}
C_0 =& \;  [1, [3, 0], 4, [3, 2], 7, [7, 0], 14, [6, -1], 8, [5, 2], 13, [4, -1]] \\
C_1 =& \;  [2, [2, 0], 0, [6, 2], 6, [4, 1]] \\
C_2 =& \;  [10, [5, -1], 5, [5, 1]] \\
C_3 =& \;  [11, [16, 0], 16, [16, 2]] \\
C_4 =& \;  [15, [4, 0], 3, [4, 2]] \\
C_5 =& \;  [9, [3, -1], 12, [3, 1]]
\end{align*}
\end{minipage}

\vspace*{-4mm}
\begin{minipage}{8cm}
\begin{align*}
C_0' =& \;  [7, [2, 1], 5, [5, 0], 10, [2, 2], 12, [6, 0], 6, [7, 2], 13, [6, 1]] \\
C_1' =& \;  [2, [1, 2], 1, [2, -1], 3, [1, 0]] \\
C_2' =& \;  [0, [8, 0], 8, [8, 2]] \\
C_3' =& \;  [9, [16, -1], 16, [16, 1]] \\
C_4' =& \;  [4, [7, -1], 11, [7, 1]] \\
C_5' =& \;  [15, [1, -1], 14, [1, 1]]
\end{align*}
\end{minipage}

\item 2-factor type $[5, 4, 2, 2, 2, 2]$: three starters \vspace*{-3mm} \\
\begin{minipage}{8cm}
\begin{align*}
C_0 =& \;  [14, [1, 2], 13, [4, 0], 1, [5, 2], 12, [3, 1], 9, [5, 0]] \\
C_1 =& \;  [10, [6, 2], 4, [16, -1], 16, [16, 1], 3, [7, 0]] \\
C_2 =& \;  [5, [6, -1], 11, [6, 1]] \\
C_3 =& \;  [8, [2, -1], 6, [2, 1]] \\
C_4 =& \;  [15, [3, 0], 2, [3, 2]] \\
C_5 =& \;  [0, [7, -1], 7, [7, 1]]
\end{align*}
\end{minipage}

\vspace*{-4mm}
\begin{minipage}{8cm}
\begin{align*}
C_0' =& \;  [7, [5, -1], 12, [16, 0], 16, [16, 2], 9, [6, 0], 3, [4, 2]] \\
C_1' =& \;  [11, [3, -1], 14, [1, 0], 13, [5, 1], 2, [7, 2]] \\
C_2' =& \;  [0, [8, 0], 8, [8, 2]] \\
C_3' =& \;  [10, [4, -1], 6, [4, 1]] \\
C_4' =& \;  [15, [2, 0], 1, [2, 2]] \\
C_5' =& \;  [5, [1, -1], 4, [1, 1]]
\end{align*}
\end{minipage}

\item 2-factor type $[5, 3, 3, 2, 2, 2]$: three starters \vspace*{-3mm} \\
\begin{minipage}{8cm}
\begin{align*}
C_0 =& \;  [7, [2, 0], 5, [6, 1], 11, [2, 2], 13, [1, -1], 14, [7, 1]] \\
C_1 =& \;  [1, [1, 2], 2, [1, 0], 3, [2, -1]] \\
C_2 =& \;  [10, [4, 2], 6, [2, 1], 4, [6, 0]] \\
C_3 =& \;  [0, [4, -1], 12, [4, 1]] \\
C_4 =& \;  [15, [7, 0], 8, [7, 2]] \\
C_5 =& \;  [16, [16, 0], 9, [16, 2]]
\end{align*}
\end{minipage}

\newpage

\begin{minipage}{8cm}
\begin{align*}
C_0' =& \;  [16, [16, -1], 11, [3, 2], 14, [6, -1], 4, [5, 0], 9, [16, 1]] \\
C_1' =& \;  [5, [4, 0], 1, [5, 2], 6, [1, 1]] \\
C_2' =& \;  [13, [6, 2], 3, [7, -1], 10, [3, 0]] \\
C_3' =& \;  [0, [8, 0], 8, [8, 2]] \\
C_4' =& \;  [2, [5, -1], 7, [5, 1]] \\
C_5' =& \;  [12, [3, -1], 15, [3, 1]]
\end{align*}
\end{minipage}

\item 2-factor type $[4, 4, 3, 2, 2, 2]$: three starters \vspace*{-3mm} \\
\hspace*{-6mm}
\begin{minipage}{8cm}
\begin{align*}
C_0 =& \;  [1, [3, 0], 14, [1, 2], 15, [6, 0], 5, [4, 2]] \\
C_1 =& \;  [6, [6, -1], 12, [1, 0], 13, [6, 1], 3, [3, 2]] \\
C_2 =& \;  [7, [5, -1], 2, [7, 2], 9, [2, 0]] \\
C_3 =& \;  [11, [1, -1], 10, [1, 1]] \\
C_4 =& \;  [4, [4, -1], 0, [4, 1]] \\
C_5 =& \;  [8, [16, -1], 16, [16, 1]]
\end{align*}
\end{minipage}
\begin{minipage}{8cm}
\begin{align*}
C_0' =& \;  [7, [2, -1], 9, [3, 1], 6, [2, 1], 4, [3, -1]] \\
C_1' =& \;  [2, [7, -1], 11, [6, 2], 1, [4, 0], 13, [5, 1]] \\
C_2' =& \;  [12, [2, 2], 14, [7, 0], 5, [7, 1]] \\
C_3' =& \;  [0, [8, 0], 8, [8, 2]] \\
C_4' =& \;  [16, [16, 0], 3, [16, 2]] \\
C_5' =& \;  [10, [5, 0], 15, [5, 2]]
\end{align*}
\end{minipage}

\item 2-factor type $[4, 3, 3, 3, 2, 2]$: three starters \vspace*{-3mm} \\
\hspace*{-7mm}
\begin{minipage}{8cm}
\begin{align*}
C_0 =& \;  [8, [1, -1], 9, [5, 1], 4, [6, 0], 14, [6, 2]] \\
C_1 =& \;  [5, [5, -1], 0, [2, 1], 2, [3, 1]] \\
C_2 =& \;  [7, [4, 0], 3, [7, 1], 10, [3, 2]] \\
C_3 =& \;  [6, [16, 0], 16, [16, -1], 13, [7, 2]] \\
C_4 =& \;  [1, [2, 0], 15, [2, 2]] \\
C_5 =& \;  [12, [1, 0], 11, [1, 2]]
\end{align*}
\end{minipage}
\begin{minipage}{8cm}
\begin{align*}
C_0' =& \;  [4, [1, 1], 3, [4, 1], 15, [2, -1], 1, [3, -1]] \\
C_1' =& \;  [10, [4, 2], 6, [7, -1], 13, [3, 0]] \\
C_2' =& \;  [14, [5, 2], 9, [7, 0], 2, [4, -1]] \\
C_3' =& \;  [16, [16, 2], 7, [5, 0], 12, [16, 1]] \\
C_4' =& \;  [0, [8, 0], 8, [8, 2]] \\
C_5' =& \;  [11, [6, -1], 5, [6, 1]]
\end{align*}
\end{minipage}

\item 2-factor type $[3, 3, 3, 3, 3, 2]$: three starters \vspace*{-3mm} \\
\hspace*{-12mm}
\begin{minipage}{8cm}
\begin{align*}
C_0 =& \;  [4, [4, 2], 8, [5, 1], 3, [1, 0]] \\
C_1 =& \;  [9, [1, 1], 10, [3, 2], 13, [4, 0]] \\
C_2 =& \;  [5, [4, -1], 1, [2, -1], 15, [6, 1]] \\
C_3 =& \;  [6, [6, 2], 12, [5, 0], 7, [1, -1]] \\
C_4 =& \;  [2, [16, 2], 16, [16, 1], 11, [7, 0]] \\
C_5 =& \;  [14, [2, 0], 0, [2, 2]]
\end{align*}
\end{minipage}
\begin{minipage}{8cm}
\begin{align*}
C_0' =& \;  [10, [1, 2], 11, [7, 1], 4, [6, 0]] \\
C_1' =& \;  [9, [3, 0], 6, [7, 2], 13, [4, 1]] \\
C_2' =& \;  [7, [5, 2], 2, [16, -1], 16, [16, 0]] \\
C_3' =& \;  [5, [7, -1], 12, [3, -1], 15, [6, -1]] \\
C_4' =& \;  [3, [2, 1], 1, [3, 1], 14, [5, -1]] \\
C_5' =& \;  [0, [8, 0], 8, [8, 2]]
\end{align*}
\end{minipage}

\item 2-factor type $[9, 2, 2, 2, 2]$: three starters \vspace*{-3mm} \\
\begin{minipage}{8cm}
\begin{align*}
C_0 =& \;  [3, [5, 2], 8, [7, 1], 1, [3, 0], 14, [4, -1], 10, [3, 2], 13, [3, -1], 0, [16, -1], 16, [16, 1], 2, [1, 0]] \\
C_1 =& \;  [12, [6, 0], 6, [6, 2]] \\
C_2 =& \;  [7, [2, 0], 9, [2, 2]] \\
C_3 =& \;  [15, [4, 0], 11, [4, 2]] \\
C_4 =& \;  [5, [1, -1], 4, [1, 1]]
\end{align*}
\end{minipage}

\vspace*{-4mm}
\begin{minipage}{8cm}
\begin{align*}
C_0' =& \;  [10, [7, -1], 1, [16, 0], 16, [16, 2], 13, [4, 1], 9, [3, 1], 6, [7, 0], 15, [1, 2], 14, [5, 0], 3, [7, 2]] \\
C_1' =& \;  [0, [8, 0], 8, [8, 2]] \\
C_2' =& \;  [7, [5, -1], 12, [5, 1]] \\
C_3' =& \;  [4, [2, -1], 2, [2, 1]] \\
C_4' =& \;  [5, [6, -1], 11, [6, 1]]
\end{align*}
\end{minipage}

\item 2-factor type $[8, 3, 2, 2, 2]$: three starters \vspace*{-3mm} \\
\begin{minipage}{8cm}
\begin{align*}
C_0 =& \;  [14, [4, 0], 10, [7, 2], 3, [6, 0], 13, [5, 1], 2, [2, 1], 4, [4, 2], 0, [7, 1], 9, [5, -1]] \\
C_1 =& \;  [5, [3, -1], 8, [3, 0], 11, [6, 2]] \\
C_2 =& \;  [16, [16, -1], 12, [16, 1]] \\
C_3 =& \;  [15, [2, 0], 1, [2, 2]] \\
C_4 =& \;  [6, [1, -1], 7, [1, 1]]
\end{align*}
\end{minipage}

\vspace*{-4mm}
\begin{minipage}{8cm}
\begin{align*}
C_0' =& \;  [12, [2, -1], 10, [3, 1], 13, [5, 2], 2, [7, 0], 11, [3, 2], 14, [1, 0], 15, [4, 1], 3, [7, -1]] \\
C_1' =& \;  [4, [1, 2], 5, [4, -1], 9, [5, 0]] \\
C_2' =& \;  [0, [8, 0], 8, [8, 2]] \\
C_3' =& \;  [1, [6, -1], 7, [6, 1]] \\
C_4' =& \;  [16, [16, 0], 6, [16, 2]]
\end{align*}
\end{minipage}

\item 2-factor type $[7, 4, 2, 2, 2]$: three starters \vspace*{-3mm} \\
\begin{minipage}{8cm}
\begin{align*}
C_0 =& \;  [2, [1, 0], 3, [6, 2], 13, [1, -1], 14, [16, -1], 16, [16, 1], 12, [4, 0], 0, [2, 2]] \\
C_1 =& \;  [7, [1, 1], 6, [2, -1], 8, [7, 1], 1, [6, -1]] \\
C_2 =& \;  [10, [5, -1], 15, [5, 1]] \\
C_3 =& \;  [11, [7, 0], 4, [7, 2]] \\
C_4 =& \;  [9, [4, -1], 5, [4, 1]]
\end{align*}
\end{minipage}

\vspace*{-4mm}
\begin{minipage}{8cm}
\begin{align*}
C_0' =& \;  [13, [6, 0], 3, [2, 1], 5, [1, 2], 4, [7, -1], 11, [5, 0], 6, [4, 2], 10, [3, -1]] \\
C_1' =& \;  [2, [5, 2], 7, [6, 1], 1, [2, 0], 15, [3, 1]] \\
C_2' =& \;  [0, [8, 0], 8, [8, 2]] \\
C_3' =& \;  [16, [16, 0], 14, [16, 2]] \\
C_4' =& \;  [12, [3, 0], 9, [3, 2]]
\end{align*}
\end{minipage}

\item 2-factor type $[6, 5, 2, 2, 2]$: three starters \vspace*{-3mm} \\
\begin{minipage}{8cm}
\begin{align*}
C_0 =& \;  [3, [4, 2], 15, [7, 0], 8, [3, 2], 5, [6, 0], 11, [1, -1], 10, [7, -1]] \\
C_1 =& \;  [2, [1, 2], 1, [3, -1], 4, [2, -1], 6, [16, 0], 16, [16, -1]] \\
C_2 =& \;  [12, [5, -1], 7, [5, 1]] \\
C_3 =& \;  [13, [4, -1], 9, [4, 1]] \\
C_4 =& \;  [14, [2, 0], 0, [2, 2]]
\end{align*}
\end{minipage}

\vspace*{-4mm}
\begin{minipage}{8cm}
\begin{align*}
C_0' =& \;  [9, [16, 1], 16, [16, 2], 11, [4, 0], 7, [6, 2], 13, [7, 1], 6, [3, 0]] \\
C_1' =& \;  [3, [1, 1], 2, [3, 1], 15, [1, 0], 14, [7, 2], 5, [2, 1]] \\
C_2' =& \;  [0, [8, 0], 8, [8, 2]] \\
C_3' =& \;  [10, [6, -1], 4, [6, 1]] \\
C_4' =& \;  [1, [5, 0], 12, [5, 2]]
\end{align*}
\end{minipage}

\item 2-factor type $[7, 3, 3, 2, 2]$: three starters \vspace*{-3mm} \\
\begin{minipage}{8cm}
\begin{align*}
C_0 =& \;  [15, [7, 1], 8, [3, -1], 11, [1, 0], 10, [2, 1], 12, [4, 1], 0, [4, 2], 4, [5, 1]] \\
C_1 =& \;  [1, [4, -1], 13, [1, 2], 14, [3, 0]] \\
C_2 =& \;  [2, [4, 0], 6, [3, 1], 9, [7, 2]] \\
C_3 =& \;  [5, [2, 0], 3, [2, 2]] \\
C_4 =& \;  [16, [16, -1], 7, [16, 1]]
\end{align*}
\end{minipage}

\newpage

\begin{minipage}{8cm}
\begin{align*}
C_0' =& \;  [1, [16, 2], 16, [16, 0], 5, [6, -1], 15, [3, 2], 12, [5, 0], 7, [5, -1], 2, [1, -1]] \\
C_1' =& \;  [6, [5, 2], 11, [2, -1], 13, [7, 0]] \\
C_2' =& \;  [3, [6, 1], 9, [1, 1], 10, [7, -1]] \\
C_3' =& \;  [0, [8, 0], 8, [8, 2]] \\
C_4' =& \;  [4, [6, 0], 14, [6, 2]]
\end{align*}
\end{minipage}

\item 2-factor type $[6, 4, 3, 2, 2]$: three starters \vspace*{-3mm} \\
\begin{minipage}{8cm}
\begin{align*}
C_0 =& \;  [13, [3, 2], 10, [1, 0], 9, [2, 2], 11, [1, -1], 12, [5, -1], 1, [4, 0]] \\
C_1 =& \;  [15, [4, 2], 3, [5, 0], 8, [2, -1], 6, [7, -1]] \\
C_2 =& \;  [16, [16, 0], 2, [3, 1], 5, [16, 2]] \\
C_3 =& \;  [7, [7, 0], 0, [7, 2]] \\
C_4 =& \;  [14, [6, -1], 4, [6, 1]]
\end{align*}
\end{minipage}

\vspace*{-4mm}
\begin{minipage}{8cm}
\begin{align*}
C_0' =& \;  [9, [4, 1], 13, [1, 1], 14, [4, -1], 10, [5, 2], 15, [3, 0], 2, [7, 1]] \\
C_1' =& \;  [5, [2, 1], 3, [3, -1], 6, [2, 0], 4, [1, 2]] \\
C_2' =& \;  [16, [16, -1], 7, [5, 1], 12, [16, 1]] \\
C_3' =& \;  [0, [8, 0], 8, [8, 2]] \\
C_4' =& \;  [11, [6, 0], 1, [6, 2]]
\end{align*}
\end{minipage}

\item 2-factor type $[5, 5, 3, 2, 2]$: three starters \vspace*{-3mm} \\
\begin{minipage}{8cm}
\begin{align*}
C_0 =& \;  [3, [4, -1], 15, [4, 2], 11, [16, -1], 16, [16, 1], 2, [1, 0]] \\
C_1 =& \;  [6, [2, -1], 8, [7, 1], 1, [7, 0], 10, [3, 1], 13, [7, 2]] \\
C_2 =& \;  [12, [5, 1], 7, [2, 0], 9, [3, 2]] \\
C_3 =& \;  [5, [5, 0], 0, [5, 2]] \\
C_4 =& \;  [4, [6, -1], 14, [6, 1]]
\end{align*}
\end{minipage}

\vspace*{-4mm}
\begin{minipage}{8cm}
\begin{align*}
C_0' =& \;  [13, [2, 1], 11, [6, 0], 5, [4, 1], 9, [1, 2], 10, [3, -1]] \\
C_1' =& \;  [14, [4, 0], 2, [6, 2], 12, [7, -1], 3, [16, 0], 16, [16, 2]] \\
C_2' =& \;  [15, [5, -1], 4, [3, 0], 1, [2, 2]] \\
C_3' =& \;  [0, [8, 0], 8, [8, 2]] \\
C_4' =& \;  [6, [1, -1], 7, [1, 1]]
\end{align*}
\end{minipage}

\item 2-factor type $[5, 4, 4, 2, 2]$: three starters \vspace*{-3mm} \\
\hspace*{-3mm}
\begin{minipage}{8cm}
\begin{align*}
C_0 =& \;  [2, [1, 0], 1, [7, 2], 10, [1, 1], 9, [3, 1], 6, [4, 1]] \\
C_1 =& \;  [5, [6, 1], 15, [2, 1], 13, [6, 0], 7, [2, 2]] \\
C_2 =& \;  [14, [6, 2], 8, [5, 0], 3, [16, 2], 16, [16, 0]] \\
C_3 =& \;  [12, [4, 0], 0, [4, 2]] \\
C_4 =& \;  [11, [7, -1], 4, [7, 1]]
\end{align*}
\end{minipage}

\vspace*{-4mm}
\begin{minipage}{8cm}
\begin{align*}
C_0' =& \;  [15, [6, -1], 5, [2, 0], 7, [1, 2], 6, [5, -1], 11, [4, -1]] \\
C_1' =& \;  [9, [3, -1], 12, [7, 0], 3, [5, 2], 14, [5, 1]] \\
C_2' =& \;  [2, [1, -1], 1, [16, 1], 16, [16, -1], 4, [2, -1]] \\
C_3' =& \;  [0, [8, 0], 8, [8, 2]] \\
C_4' =& \;  [10, [3, 0], 13, [3, 2]]
\end{align*}
\end{minipage}

\item 2-factor type $[6, 3, 3, 3, 2]$: three starters \vspace*{-3mm} \\
\begin{minipage}{8cm}
\begin{align*}
C_0 =& \;  [11, [6, 2], 5, [7, -1], 12, [3, 0], 15, [1, 2], 14, [5, 1], 9, [2, 0]] \\
C_1 =& \;  [7, [3, -1], 10, [4, 1], 6, [1, -1]] \\
C_2 =& \;  [4, [16, 2], 16, [16, 0], 13, [7, 1]] \\
C_3 =& \;  [3, [3, 2], 0, [2, -1], 2, [1, 0]] \\
C_4 =& \;  [8, [7, 0], 1, [7, 2]]
\end{align*}
\end{minipage}

\vspace*{-4mm}
\begin{minipage}{8cm}
\begin{align*}
C_0' =& \;  [2, [3, 1], 15, [6, -1], 5, [1, 1], 4, [6, 0], 10, [4, -1], 6, [4, 2]] \\
C_1' =& \;  [14, [16, 1], 16, [16, -1], 3, [5, -1]] \\
C_2' =& \;  [9, [2, 1], 11, [2, 2], 13, [4, 0]] \\
C_3' =& \;  [12, [5, 0], 1, [6, 1], 7, [5, 2]] \\
C_4' =& \;  [0, [8, 0], 8, [8, 2]]
\end{align*}
\end{minipage}

\item 2-factor type $[5, 4, 3, 3, 2]$: three starters \vspace*{-3mm} \\
\begin{minipage}{8cm}
\begin{align*}
C_0 =& \;  [9, [1, 0], 8, [5, 1], 13, [2, 1], 15, [5, -1], 10, [1, 2]] \\
C_1 =& \;  [0, [6, -1], 6, [5, 0], 11, [3, 1], 14, [2, 2]] \\
C_2 =& \;  [2, [1, -1], 3, [7, 0], 12, [6, 2]] \\
C_3 =& \;  [16, [16, 0], 4, [3, 2], 7, [16, -1]] \\
C_4 =& \;  [5, [4, 0], 1, [4, 2]]
\end{align*}
\end{minipage}

\vspace*{-4mm}
\begin{minipage}{8cm}
\begin{align*}
C_0' =& \;  [13, [4, 1], 1, [5, 2], 12, [2, 0], 10, [7, 2], 3, [6, 0]] \\
C_1' =& \;  [16, [16, 1], 7, [2, -1], 9, [3, 0], 6, [16, 2]] \\
C_2' =& \;  [11, [4, -1], 15, [3, -1], 2, [7, 1]] \\
C_3' =& \;  [14, [7, -1], 5, [1, 1], 4, [6, 1]] \\
C_4' =& \;  [0, [8, 0], 8, [8, 2]]
\end{align*}
\end{minipage}

\item 2-factor type $[4, 4, 4, 3, 2]$: three starters \vspace*{-3mm} \\
\hspace*{-4mm}
\begin{minipage}{8cm}
\begin{align*}
C_0 =& \;  [5, [5, -1], 0, [2, 2], 2, [7, 1], 11, [6, 0]] \\
C_1 =& \;  [9, [5, 2], 14, [1, 0], 13, [3, -1], 10, [1, -1]] \\
C_2 =& \;  [15, [7, 2], 8, [7, 0], 1, [3, 2], 4, [5, 0]] \\
C_3 =& \;  [3, [7, -1], 12, [6, 2], 6, [3, 0]] \\
C_4 =& \;  [16, [16, 0], 7, [16, 2]]
\end{align*}
\end{minipage}
\begin{minipage}{8cm}
\begin{align*}
C_0' =& \;  [13, [1, 1], 12, [5, 1], 7, [4, 1], 3, [6, 1]] \\
C_1' =& \;  [11, [1, 2], 10, [4, 0], 6, [3, 1], 9, [2, 1]] \\
C_2' =& \;  [5, [4, -1], 1, [2, -1], 15, [16, 1], 16, [16, -1]] \\
C_3' =& \;  [4, [2, 0], 2, [4, 2], 14, [6, -1]] \\
C_4' =& \;  [0, [8, 0], 8, [8, 2]]
\end{align*}
\end{minipage}

\item 2-factor type $[11, 2, 2, 2]$: three starters \vspace*{-3mm} \\
\begin{minipage}{8cm}
\begin{align*}
C_0 =& \;  [8, [1, 0], 7, [6, 2], 13, [6, 0], 3, [4, -1], 15, [5, 2], 10, [2, 0], 12, [5, 1], 1, [1, 2], 0, [2, 1], 14, [3, 0], 11, [3, 2]] \\
C_1 =& \;  [6, [16, 0], 16, [16, 2]] \\
C_2 =& \;  [4, [1, -1], 5, [1, 1]] \\
C_3 =& \;  [2, [7, 0], 9, [7, 2]]
\end{align*}
\end{minipage}

\vspace*{-4mm}
\begin{minipage}{8cm}
\begin{align*}
C_0' =& \;  [10, [5, 0], 15, [3, -1], 12, [7, 1], 3, [3, 1], 6, [2, 2], 4, [7, -1], 11, [2, -1], 13, [5, -1], 2, [4, 1], 14, [16, -1], \\ & \; 16, [16, 1]] \\
C_1' =& \;  [0, [8, 0], 8, [8, 2]] \\
C_2' =& \;  [9, [4, 0], 5, [4, 2]] \\
C_3' =& \;  [7, [6, -1], 1, [6, 1]]
\end{align*}
\end{minipage}

\item 2-factor type $[10, 3, 2, 2]$: three starters \vspace*{-3mm} \\
\begin{minipage}{8cm}
\begin{align*}
C_0 =& \;  [10, [1, 1], 11, [5, 2], 0, [4, 0], 4, [5, 1], 9, [3, 2], 12, [6, 0], 2, [3, 1], 5, [1, 2], 6, [1, -1], 7, [3, 0]] \\
C_1 =& \;  [8, [5, 0], 3, [5, -1], 14, [6, 2]] \\
C_2 =& \;  [15, [2, -1], 1, [2, 1]] \\
C_3 =& \;  [13, [16, -1], 16, [16, 1]]
\end{align*}
\end{minipage}

\vspace*{-4mm}
\begin{minipage}{8cm}
\begin{align*}
C_0' =& \;  [13, [7, -1], 4, [1, 0], 3, [7, 2], 12, [6, 1], 6, [4, 1], 2, [3, -1], 5, [7, 1], 14, [7, 0], 7, [2, 2], 9, [4, -1]] \\
C_1' =& \;  [11, [4, 2], 15, [2, 0], 1, [6, -1]] \\
C_2' =& \;  [0, [8, 0], 8, [8, 2]] \\
C_3' =& \;  [10, [16, 0], 16, [16, 2]]
\end{align*}
\end{minipage}

\item 2-factor type $[9, 4, 2, 2]$: three starters \vspace*{-3mm} \\
\begin{minipage}{8cm}
\begin{align*}
C_0 =& \;  [10, [4, 1], 14, [1, 2], 13, [7, 1], 6, [7, 0], 15, [6, -1], 9, [6, 2], 3, [4, 0], 7, [7, 2], 0, [6, 0]] \\
C_1 =& \;  [4, [1, 0], 5, [3, 2], 2, [16, -1], 16, [16, 1]] \\
C_2 =& \;  [8, [3, -1], 11, [3, 1]] \\
C_3 =& \;  [1, [5, -1], 12, [5, 1]]
\end{align*}
\end{minipage}

\vspace*{-4mm}
\begin{minipage}{8cm}
\begin{align*}
C_0' =& \;  [15, [6, 1], 9, [4, -1], 13, [1, 1], 12, [2, 0], 14, [7, -1], 7, [4, 2], 11, [5, 0], 6, [2, -1], 4, [5, 2]] \\
C_1' =& \;  [5, [3, 0], 2, [1, -1], 1, [2, 1], 3, [2, 2]] \\
C_2' =& \;  [0, [8, 0], 8, [8, 2]] \\
C_3' =& \;  [10, [16, 0], 16, [16, 2]]
\end{align*}
\end{minipage}

\item 2-factor type $[8, 5, 2, 2]$: three starters \vspace*{-3mm} \\
\begin{minipage}{8cm}
\begin{align*}
C_0 =& \;  [4, [1, 1], 3, [3, 0], 0, [5, -1], 11, [16, 2], 16, [16, 1], 8, [2, 0], 6, [7, 1], 13, [7, 2]] \\
C_1 =& \;  [10, [5, 1], 15, [1, 2], 14, [7, -1], 5, [4, 1], 1, [7, 0]] \\
C_2 =& \;  [2, [5, 0], 7, [5, 2]] \\
C_3 =& \;  [12, [3, -1], 9, [3, 1]]
\end{align*}
\end{minipage}

\vspace*{-4mm}
\begin{minipage}{8cm}
\begin{align*}
C_0' =& \;  [11, [6, -1], 1, [2, 1], 15, [4, -1], 3, [16, -1], 16, [16, 0], 2, [4, 2], 14, [1, 0], 13, [2, 2]] \\
C_1' =& \;  [7, [3, 2], 10, [6, 1], 4, [1, -1], 5, [4, 0], 9, [2, -1]] \\
C_2' =& \;  [0, [8, 0], 8, [8, 2]] \\
C_3' =& \;  [12, [6, 0], 6, [6, 2]]
\end{align*}
\end{minipage}

\item 2-factor type $[7, 6, 2, 2]$: three starters \vspace*{-3mm} \\
\begin{minipage}{8cm}
\begin{align*}
C_0 =& \;  [6, [5, 1], 1, [16, 0], 16, [16, 2], 11, [4, 0], 15, [6, 1], 5, [2, 2], 3, [3, -1]] \\
C_1 =& \;  [12, [2, 1], 10, [4, -1], 14, [5, 0], 9, [7, 2], 0, [2, 0], 2, [6, 2]] \\
C_2 =& \;  [4, [7, -1], 13, [7, 1]] \\
C_3 =& \;  [8, [1, 0], 7, [1, 2]]
\end{align*}
\end{minipage}

\vspace*{-4mm}
\begin{minipage}{8cm}
\begin{align*}
C_0' =& \;  [10, [3, 0], 7, [3, 2], 4, [5, -1], 9, [3, 1], 6, [7, 0], 15, [1, -1], 14, [4, 2]] \\
C_1' =& \;  [13, [4, 1], 1, [5, 2], 12, [1, 1], 11, [6, 0], 5, [2, -1], 3, [6, -1]] \\
C_2' =& \;  [0, [8, 0], 8, [8, 2]] \\
C_3' =& \;  [2, [16, -1], 16, [16, 1]]
\end{align*}
\end{minipage}

\newpage

\item 2-factor type $[9, 3, 3, 2]$: three starters \vspace*{-3mm} \\
\begin{minipage}{8cm}
\begin{align*}
C_0 =& \;  [13, [5, 2], 8, [2, -1], 10, [4, -1], 14, [7, 1], 7, [2, 1], 5, [6, 0], 15, [1, 1], 0, [16, 1], 16, [16, -1]] \\
C_1 =& \;  [2, [6, -1], 12, [1, 2], 11, [7, 0]] \\
C_2 =& \;  [3, [3, -1], 6, [5, 0], 1, [2, 2]] \\
C_3 =& \;  [9, [5, -1], 4, [5, 1]]
\end{align*}
\end{minipage}

\vspace*{-4mm}
\begin{minipage}{8cm}
\begin{align*}
C_0' =& \;  [12, [6, 1], 6, [7, -1], 13, [4, 0], 9, [7, 2], 2, [3, 0], 5, [1, -1], 4, [16, 2], 16, [16, 0], 15, [3, 2]] \\
C_1' =& \;  [14, [3, 1], 11, [1, 0], 10, [4, 2]] \\
C_2' =& \;  [7, [4, 1], 3, [2, 0], 1, [6, 2]] \\
C_3' =& \;  [0, [8, 0], 8, [8, 2]]
\end{align*}
\end{minipage}

\item 2-factor type $[8, 4, 3, 2]$: three starters \vspace*{-3mm} \\
\begin{minipage}{8cm}
\begin{align*}
C_0 =& \;  [1, [7, 0], 10, [5, 2], 5, [1, 0], 4, [6, -1], 14, [3, 2], 11, [2, -1], 13, [5, -1], 2, [1, 1]] \\
C_1 =& \;  [6, [7, -1], 15, [1, 2], 0, [7, 1], 9, [3, 0]] \\
C_2 =& \;  [8, [16, 2], 16, [16, 1], 12, [4, 0]] \\
C_3 =& \;  [3, [4, -1], 7, [4, 1]]
\end{align*}
\end{minipage}

\vspace*{-4mm}
\begin{minipage}{8cm}
\begin{align*}
C_0' =& \;  [5, [7, 2], 12, [6, 0], 6, [4, 2], 2, [2, 0], 4, [1, -1], 3, [6, 2], 13, [16, 0], 16, [16, -1]] \\
C_1' =& \;  [1, [2, 1], 15, [5, 1], 10, [3, 1], 7, [6, 1]] \\
C_2' =& \;  [14, [3, -1], 11, [2, 2], 9, [5, 0]] \\
C_3' =& \;  [0, [8, 0], 8, [8, 2]]
\end{align*}
\end{minipage}

\item 2-factor type $[7, 5, 3, 2]$: three starters \vspace*{-3mm} \\
\begin{minipage}{8cm}
\begin{align*}
C_0 =& \;  [10, [3, 1], 7, [1, 1], 6, [7, -1], 13, [4, -1], 1, [1, 0], 0, [1, 2], 15, [5, 1]] \\
C_1 =& \;  [3, [6, 1], 9, [5, 2], 4, [4, 0], 8, [6, 2], 14, [5, 0]] \\
C_2 =& \;  [11, [1, -1], 12, [6, -1], 2, [7, 1]] \\
C_3 =& \;  [16, [16, 0], 5, [16, 2]]
\end{align*}
\end{minipage}

\vspace*{-4mm}
\begin{minipage}{8cm}
\begin{align*}
C_0' =& \;  [14, [3, 2], 1, [16, -1], 16, [16, 1], 3, [3, 0], 6, [4, 1], 10, [2, 2], 12, [2, 0]] \\
C_1' =& \;  [9, [5, -1], 4, [7, 2], 13, [2, -1], 15, [3, -1], 2, [7, 0]] \\
C_2' =& \;  [11, [4, 2], 7, [2, 1], 5, [6, 0]] \\
C_3' =& \;  [0, [8, 0], 8, [8, 2]]
\end{align*}
\end{minipage}

\item 2-factor type $[6, 6, 3, 2]$: three starters \vspace*{-3mm} \\
\begin{minipage}{8cm}
\begin{align*}
C_0 =& \;  [2, [7, -1], 11, [6, 1], 1, [16, 2], 16, [16, 0], 0, [1, -1], 15, [3, 1]] \\
C_1 =& \;  [5, [3, -1], 8, [6, -1], 14, [5, 0], 9, [6, 2], 3, [4, 0], 7, [2, 2]] \\
C_2 =& \;  [10, [2, 0], 12, [1, 1], 13, [3, 2]] \\
C_3 =& \;  [4, [2, -1], 6, [2, 1]]
\end{align*}
\end{minipage}

\vspace*{-4mm}
\begin{minipage}{8cm}
\begin{align*}
C_0' =& \;  [11, [5, 1], 6, [1, 0], 7, [4, 2], 3, [6, 0], 13, [5, 2], 2, [7, 1]] \\
C_1' =& \;  [10, [7, 0], 1, [16, 1], 16, [16, -1], 4, [5, -1], 15, [1, 2], 14, [4, 1]] \\
C_2' =& \;  [12, [7, 2], 5, [4, -1], 9, [3, 0]] \\
C_3' =& \;  [0, [8, 0], 8, [8, 2]]
\end{align*}
\end{minipage}

\newpage

\item 2-factor type $[7, 4, 4, 2]$: three starters \vspace*{-3mm} \\
\begin{minipage}{8cm}
\begin{align*}
C_0 =& \;  [6, [1, 0], 7, [4, 2], 3, [1, -1], 4, [6, 0], 14, [1, 2], 13, [4, 1], 9, [3, 1]] \\
C_1 =& \;  [11, [1, 1], 12, [3, 2], 15, [2, 0], 1, [6, -1]] \\
C_2 =& \;  [8, [3, -1], 5, [5, 1], 10, [16, 1], 16, [16, -1]] \\
C_3 =& \;  [2, [2, -1], 0, [2, 1]]
\end{align*}
\end{minipage}

\vspace*{-4mm}
\begin{minipage}{8cm}
\begin{align*}
C_0' =& \;  [5, [7, -1], 12, [5, 0], 7, [2, 2], 9, [5, -1], 14, [4, -1], 2, [7, 0], 11, [6, 2]] \\
C_1' =& \;  [16, [16, 2], 1, [3, 0], 4, [5, 2], 15, [16, 0]] \\
C_2' =& \;  [3, [6, 1], 13, [7, 1], 6, [4, 0], 10, [7, 2]] \\
C_3' =& \;  [0, [8, 0], 8, [8, 2]]
\end{align*}
\end{minipage}

\item 2-factor type $[6, 5, 4, 2]$: three starters \vspace*{-3mm} \\
\begin{minipage}{8cm}
\begin{align*}
C_0 =& \;  [7, [4, -1], 11, [5, -1], 0, [1, 0], 15, [6, 1], 5, [1, 1], 6, [1, 2]] \\
C_1 =& \;  [14, [1, -1], 13, [16, 2], 16, [16, 1], 10, [7, 1], 3, [5, 0]] \\
C_2 =& \;  [1, [5, 2], 12, [4, 1], 8, [4, 0], 4, [3, -1]] \\
C_3 =& \;  [9, [7, 0], 2, [7, 2]]
\end{align*}
\end{minipage}

\vspace*{-4mm}
\begin{minipage}{8cm}
\begin{align*}
C_0' =& \;  [5, [4, 2], 1, [16, 0], 16, [16, -1], 12, [6, -1], 6, [3, 2], 3, [2, 0]] \\
C_1' =& \;  [2, [7, -1], 9, [2, -1], 11, [3, 0], 14, [6, 2], 4, [2, 1]] \\
C_2' =& \;  [13, [2, 2], 15, [5, 1], 10, [3, 1], 7, [6, 0]] \\
C_3' =& \;  [0, [8, 0], 8, [8, 2]]
\end{align*}
\end{minipage}

\item 2-factor type $[5, 5, 5, 2]$: three starters \vspace*{-3mm} \\
\begin{minipage}{8cm}
\begin{align*}
C_0 =& \;  [13, [7, 1], 4, [5, -1], 15, [3, 1], 2, [4, -1], 14, [1, -1]] \\
C_1 =& \;  [3, [16, 0], 16, [16, 2], 7, [2, 1], 5, [4, 0], 1, [2, 2]] \\
C_2 =& \;  [8, [4, 2], 12, [1, 0], 11, [1, 2], 10, [4, 1], 6, [2, 0]] \\
C_3 =& \;  [0, [7, 0], 9, [7, 2]]
\end{align*}
\end{minipage}

\vspace*{-4mm}
\begin{minipage}{8cm}
\begin{align*}
C_0' =& \;  [13, [5, 1], 2, [3, 2], 15, [6, 1], 9, [6, 0], 3, [6, -1]] \\
C_1' =& \;  [6, [16, 1], 16, [16, -1], 4, [1, 1], 5, [6, 2], 11, [5, 0]] \\
C_2' =& \;  [7, [5, 2], 12, [2, -1], 14, [3, 0], 1, [7, -1], 10, [3, -1]] \\
C_3' =& \;  [0, [8, 0], 8, [8, 2]]
\end{align*}
\end{minipage}

\item 2-factor type $[13, 2, 2]$: three starters \vspace*{-3mm} \\
\begin{minipage}{8cm}
\begin{align*}
C_0 =& \;  [7, [2, -1], 9, [4, 1], 5, [1, -1], 6, [6, -1], 12, [6, 0], 2, [3, 2], 15, [7, 0], 8, [7, -1], 1, [16, 1], 16, [16, -1], \\ & \; 14, [4, -1], 10, [3, 1], 13, [6, 2]] \\
C_1 =& \;  [3, [1, 0], 4, [1, 2]] \\
C_2 =& \;  [0, [5, -1], 11, [5, 1]]
\end{align*}
\end{minipage}

\vspace*{-4mm}
\begin{minipage}{8cm}
\begin{align*}
C_0' =& \;  [1, [3, -1], 4, [16, 0], 16, [16, 2], 5, [2, 1], 3, [7, 1], 12, [3, 0], 15, [5, 2], 10, [4, 0], 6, [4, 2], 2, [5, 0], 13, \\ & \; [1, 1], 14, [7, 2], 7, [6, 1]] \\
C_1' =& \;  [0, [8, 0], 8, [8, 2]] \\
C_2' =& \;  [9, [2, 0], 11, [2, 2]]
\end{align*}
\end{minipage}

\newpage

\item 2-factor type $[12, 3, 2]$: three starters \vspace*{-3mm} \\
\begin{minipage}{8cm}
\begin{align*}
C_0 =& \;  [4, [7, -1], 13, [1, 1], 14, [2, 1], 0, [4, 2], 12, [4, 1], 8, [3, 0], 11, [6, -1], 5, [4, -1], 1, [2, 2], 3, [3, -1], 6, \\ & \; [7, 0], 15, [5, 1]] \\
C_1 =& \;  [7, [16, -1], 16, [16, 1], 2, [5, -1]] \\
C_2 =& \;  [10, [1, 0], 9, [1, 2]]
\end{align*}
\end{minipage}

\vspace*{-4mm}
\begin{minipage}{8cm}
\begin{align*}
C_0' =& \;  [4, [2, 0], 2, [5, 2], 13, [7, 1], 6, [5, 0], 1, [7, 2], 10, [16, 0], 16, [16, 2], 5, [6, 1], 15, [6, 0], 9, [6, 2], 3, [4, 0], \\ & \; 7, [3, 2]] \\
C_1' =& \;  [12, [1, -1], 11, [3, 1], 14, [2, -1]] \\
C_2' =& \;  [0, [8, 0], 8, [8, 2]]
\end{align*}
\end{minipage}

\item 2-factor type $[11, 4, 2]$: three starters \vspace*{-3mm} \\
\begin{minipage}{8cm}
\begin{align*}
C_0 =& \;  [4, [3, 1], 1, [7, -1], 8, [7, 0], 15, [1, 2], 14, [2, 0], 12, [7, 2], 5, [5, -1], 10, [1, 0], 11, [2, 2], 9, [7, 1], 2, \\ & \; [2, -1]] \\
C_1 =& \;  [6, [6, -1], 0, [3, -1], 13, [16, 1], 16, [16, -1]] \\
C_2 =& \;  [7, [4, 0], 3, [4, 2]]
\end{align*}
\end{minipage}

\vspace*{-4mm}
\begin{minipage}{8cm}
\begin{align*}
C_0' =& \;  [14, [1, -1], 13, [3, 2], 10, [5, 1], 5, [1, 1], 4, [2, 1], 2, [4, -1], 6, [3, 0], 3, [6, 2], 9, [6, 0], 15, [16, 2], 16, \\ & \; [16, 0]] \\
C_1' =& \;  [11, [4, 1], 7, [5, 0], 12, [5, 2], 1, [6, 1]] \\
C_2' =& \;  [0, [8, 0], 8, [8, 2]]
\end{align*}
\end{minipage}

\item 2-factor type $[10, 5, 2]$: three starters \vspace*{-3mm} \\
\begin{minipage}{8cm}
\begin{align*}
C_0 =& \;  [14, [3, 1], 11, [4, 1], 7, [16, 0], 16, [16, -1], 3, [5, 2], 8, [2, -1], 10, [6, 1], 4, [5, 0], 15, [6, 2], 5, [7, 1]] \\
C_1 =& \;  [6, [7, 2], 13, [5, -1], 2, [7, -1], 9, [3, 0], 12, [6, -1]] \\
C_2 =& \;  [0, [1, 0], 1, [1, 2]]
\end{align*}
\end{minipage}

\vspace*{-4mm}
\begin{minipage}{8cm}
\begin{align*}
C_0' =& \;  [12, [7, 0], 3, [4, 2], 7, [6, 0], 1, [3, 2], 4, [2, 0], 6, [4, -1], 2, [3, -1], 15, [2, 2], 13, [1, -1], 14, [2, 1]] \\
C_1' =& \;  [9, [16, 1], 16, [16, 2], 11, [1, 1], 10, [5, 1], 5, [4, 0]] \\
C_2' =& \;  [0, [8, 0], 8, [8, 2]]
\end{align*}
\end{minipage}

\item 2-factor type $[9, 6, 2]$: three starters \vspace*{-3mm} \\
\begin{minipage}{8cm}
\begin{align*}
C_0 =& \;  [12, [1, 2], 13, [6, 0], 7, [5, 2], 2, [2, 1], 0, [4, -1], 4, [2, -1], 6, [16, 0], 16, [16, -1], 3, [7, -1]] \\
C_1 =& \;  [10, [5, 0], 15, [1, -1], 14, [6, 2], 8, [3, 0], 11, [6, 1], 1, [7, 2]] \\
C_2 =& \;  [5, [4, 0], 9, [4, 2]]
\end{align*}
\end{minipage}

\vspace*{-4mm}
\begin{minipage}{8cm}
\begin{align*}
C_0' =& \;  [15, [2, 0], 13, [4, 1], 1, [1, 1], 2, [16, 1], 16, [16, 2], 11, [5, 1], 6, [3, 1], 3, [7, 0], 12, [3, 2]] \\
C_1' =& \;  [4, [5, -1], 9, [1, 0], 10, [3, -1], 7, [2, 2], 5, [7, 1], 14, [6, -1]] \\
C_2' =& \;  [0, [8, 0], 8, [8, 2]]
\end{align*}
\end{minipage}

\item 2-factor type $[8, 7, 2]$: three starters \vspace*{-3mm} \\
\begin{minipage}{8cm}
\begin{align*}
C_0 =& \;  [3, [4, -1], 7, [3, 1], 4, [7, 1], 13, [4, 0], 1, [2, -1], 15, [3, 2], 2, [2, 1], 0, [3, -1]] \\
C_1 =& \;  [8, [3, 0], 5, [5, 1], 10, [4, 2], 14, [5, 0], 9, [16, 1], 16, [16, -1], 6, [2, 2]] \\
C_2 =& \;  [12, [1, -1], 11, [1, 1]]
\end{align*}
\end{minipage}

\newpage

\begin{minipage}{8cm}
\begin{align*}
C_0' =& \;  [5, [7, -1], 14, [5, 2], 9, [2, 0], 11, [6, 2], 1, [16, 0], 16, [16, 2], 10, [6, 1], 4, [1, 0]] \\
C_1' =& \;  [2, [5, -1], 13, [6, -1], 7, [1, 2], 6, [7, 0], 15, [4, 1], 3, [7, 2], 12, [6, 0]] \\
C_2' =& \;  [0, [8, 0], 8, [8, 2]]
\end{align*}
\end{minipage}

\end{itemizenew}

\section{Computational results for $n=18$}\label{app:18}

\begin{itemizenew}
\item 2-factor type $[4, 2, 2, 2, 2, 2, 2, 2]$: one starter \vspace*{-3mm} \\
\hspace*{-10mm}
\begin{minipage}{8cm}
\begin{align*}
C_0 =& \; [2, [4, 0], 15, [5, 1], 3, [4, 1], 7, [5, 0]] \\
C_1 =& \; [1, [8, 0], 9, [8, 1]] \\
C_2 =& \; [10, [7, 0], 0, [7, 1]] \\
C_3 =& \; [11, [3, 0], 8, [3, 1]] \\
C_4 =& \; [12, [6, 0], 6, [6, 1]] \\
C_5 =& \; [5, [1, 0], 4, [1, 1]] \\
C_6 =& \; [17, [17, 0], 13, [17, 1]] \\
C_7 =& \; [16, [2, 0], 14, [2, 1]]
\end{align*}
\end{minipage}

\item 2-factor type $[3, 3, 2, 2, 2, 2, 2, 2]$: two starters \vspace*{-3mm} \\
\hspace*{-12mm}
\begin{minipage}{8cm}
\begin{align*}
C_0 =& \; [17, [17, 1], 11, [4, 1], 7, [17, -1]] \\
C_1 =& \; [13, [5, 1], 8, [5, 0], 3, [7, 2]] \\
C_2 =& \; [5, [1, 0], 4, [1, 2]] \\
C_3 =& \; [15, [1, -1], 16, [1, 1]] \\
C_4 =& \; [2, [7, -1], 9, [7, 1]] \\
C_5 =& \; [0, [6, 0], 6, [6, 2]] \\
C_6 =& \; [1, [6, -1], 12, [6, 1]] \\
C_7 =& \; [10, [4, 0], 14, [4, 2]]
\end{align*}
\end{minipage}
\begin{minipage}{8cm}
\begin{align*}
C_0' =& \; [8, [7, 0], 15, [5, 2], 3, [5, -1]] \\
C_1' =& \; [9, [2, 2], 11, [2, 0], 13, [4, -1]] \\
C_2' =& \; [2, [3, 0], 5, [3, 2]] \\
C_3' =& \; [1, [3, -1], 4, [3, 1]] \\
C_4' =& \; [7, [8, -1], 16, [8, 1]] \\
C_5' =& \; [10, [2, -1], 12, [2, 1]] \\
C_6' =& \; [0, [17, 0], 17, [17, 2]] \\
C_7' =& \; [14, [8, 0], 6, [8, 2]]
\end{align*}
\end{minipage}

\item 2-factor type $[6, 2, 2, 2, 2, 2, 2]$: two starters \vspace*{-3mm} \\
\begin{minipage}{8cm}
\begin{align*}
C_0 =& \; [10, [17, -1], 17, [17, 1], 1, [6, 0], 7, [1, -1], 6, [8, 2], 14, [4, 1]] \\
C_1 =& \; [2, [1, 0], 3, [1, 2]] \\
C_2 =& \; [5, [3, 0], 8, [3, 2]] \\
C_3 =& \; [15, [2, -1], 13, [2, 1]] \\
C_4 =& \; [9, [8, -1], 0, [8, 1]] \\
C_5 =& \; [12, [4, 0], 16, [4, 2]] \\
C_6 =& \; [11, [7, 0], 4, [7, 2]]
\end{align*}
\end{minipage}

\vspace*{-4mm}
\begin{minipage}{8cm}
\begin{align*}
C_0' =& \; [0, [1, 1], 16, [6, -1], 5, [4, -1], 9, [6, 1], 3, [8, 0], 11, [6, 2]] \\
C_1' =& \; [10, [3, -1], 13, [3, 1]] \\
C_2' =& \; [4, [17, 0], 17, [17, 2]] \\
C_3' =& \; [15, [7, -1], 8, [7, 1]] \\
C_4' =& \; [1, [5, 0], 6, [5, 2]] \\
C_5' =& \; [14, [2, 0], 12, [2, 2]] \\
C_6' =& \; [7, [5, -1], 2, [5, 1]]
\end{align*}
\end{minipage}

\item 2-factor type $[5, 3, 2, 2, 2, 2, 2]$: one starter \vspace*{-3mm} \\
\hspace*{-1mm}
\begin{minipage}{8cm}
\begin{align*}
C_0 =& \; [15, [6, 1], 4, [17, 1], 17, [17, 0], 9, [7, 0], 2, [4, 1]] \\
C_1 =& \; [14, [7, 1], 7, [4, 0], 3, [6, 0]] \\
C_2 =& \; [6, [2, 0], 8, [2, 1]] \\
C_3 =& \; [13, [3, 0], 16, [3, 1]] \\
C_4 =& \; [5, [5, 0], 0, [5, 1]] \\
C_5 =& \; [12, [1, 0], 11, [1, 1]] \\
C_6 =& \; [1, [8, 0], 10, [8, 1]]
\end{align*}
\end{minipage}

\item 2-factor type $[4, 4, 2, 2, 2, 2, 2]$: one starter \vspace*{-3mm} \\
\hspace*{-9mm}
\begin{minipage}{8cm}
\begin{align*}
C_0 =& \; [14, [8, 1], 5, [5, 1], 0, [6, 1], 11, [3, 1]] \\
C_1 =& \; [16, [8, 0], 7, [6, 0], 1, [3, 0], 4, [5, 0]] \\
C_2 =& \; [6, [17, 0], 17, [17, 1]] \\
C_3 =& \; [13, [4, 0], 9, [4, 1]] \\
C_4 =& \; [8, [7, 0], 15, [7, 1]] \\
C_5 =& \; [2, [1, 0], 3, [1, 1]] \\
C_6 =& \; [12, [2, 0], 10, [2, 1]]
\end{align*}
\end{minipage}

\item 2-factor type $[4, 3, 3, 2, 2, 2, 2]$: two starters \vspace*{-3mm} \\
\hspace*{-5mm}
\begin{minipage}{8cm}
\begin{align*}
C_0 =& \; [1, [17, 2], 17, [17, 0], 4, [8, -1], 13, [5, 1]] \\
C_1 =& \; [12, [2, 0], 10, [7, 2], 3, [8, 1]] \\
C_2 =& \; [6, [1, 0], 5, [4, 2], 9, [3, 1]] \\
C_3 =& \; [2, [6, 0], 8, [6, 2]] \\
C_4 =& \; [0, [3, 0], 14, [3, 2]] \\
C_5 =& \; [16, [1, -1], 15, [1, 1]] \\
C_6 =& \; [11, [4, -1], 7, [4, 1]]
\end{align*}
\end{minipage}
\begin{minipage}{8cm}
\begin{align*}
C_0' =& \; [8, [3, -1], 5, [1, 2], 4, [8, 0], 13, [5, -1]] \\
C_1' =& \; [11, [8, 2], 3, [2, 1], 1, [7, 0]] \\
C_2' =& \; [14, [4, 0], 10, [2, 2], 12, [2, -1]] \\
C_3' =& \; [9, [7, -1], 16, [7, 1]] \\
C_4' =& \; [2, [5, 0], 7, [5, 2]] \\
C_5' =& \; [15, [17, -1], 17, [17, 1]] \\
C_6' =& \; [0, [6, -1], 6, [6, 1]]
\end{align*}
\end{minipage}

\item 2-factor type $[3, 3, 3, 3, 2, 2, 2]$: one starter \vspace*{-3mm} \\
\hspace*{-15mm}
\begin{minipage}{8cm}
\begin{align*}
C_0 =& \; [13, [7, 1], 3, [1, 0], 2, [6, 0]] \\
C_1 =& \; [0, [1, 1], 16, [2, 1], 14, [3, 1]] \\
C_2 =& \; [4, [6, 1], 10, [8, 0], 1, [3, 0]] \\
C_3 =& \; [6, [8, 1], 15, [7, 0], 8, [2, 0]] \\
C_4 =& \; [9, [4, 0], 5, [4, 1]] \\
C_5 =& \; [12, [5, 0], 7, [5, 1]] \\
C_6 =& \; [11, [17, 0], 17, [17, 1]]
\end{align*}
\end{minipage}

\item 2-factor type $[8, 2, 2, 2, 2, 2]$: one starter \vspace*{-3mm} \\
\begin{minipage}{8cm}
\begin{align*}
C_0 =& \; [7, [4, 1], 11, [4, 0], 15, [3, 0], 12, [2, 1], 14, [2, 0], 16, [3, 1], 2, [6, 1], 13, [6, 0]] \\
C_1 =& \; [17, [17, 0], 1, [17, 1]] \\
C_2 =& \; [6, [1, 0], 5, [1, 1]] \\
C_3 =& \; [4, [5, 0], 9, [5, 1]] \\
C_4 =& \; [3, [7, 0], 10, [7, 1]] \\
C_5 =& \; [0, [8, 0], 8, [8, 1]]
\end{align*}
\end{minipage}

\item 2-factor type $[7, 3, 2, 2, 2, 2]$: two starters \vspace*{-3mm} \\
\begin{minipage}{8cm}
\begin{align*}
C_0 =& \; [10, [8, 1], 1, [4, -1], 14, [17, 1], 17, [17, 2], 15, [8, 0], 6, [5, 1], 11, [1, -1]] \\
C_1 =& \; [12, [1, 1], 13, [7, 1], 3, [8, -1]] \\
C_2 =& \; [7, [3, -1], 4, [3, 1]] \\
C_3 =& \; [5, [5, 0], 0, [5, 2]] \\
C_4 =& \; [16, [7, 0], 9, [7, 2]] \\
C_5 =& \; [2, [6, 0], 8, [6, 2]]
\end{align*}
\end{minipage}

\vspace*{-4mm}
\begin{minipage}{8cm}
\begin{align*}
C_0' =& \; [11, [1, 0], 12, [3, 2], 9, [5, -1], 14, [17, 0], 17, [17, -1], 13, [8, 2], 4, [7, -1]] \\
C_1' =& \; [1, [1, 2], 2, [4, 1], 15, [3, 0]] \\
C_2' =& \; [6, [6, -1], 0, [6, 1]] \\
C_3' =& \; [8, [2, 0], 10, [2, 2]] \\
C_4' =& \; [7, [2, -1], 5, [2, 1]] \\
C_5' =& \; [3, [4, 0], 16, [4, 2]]
\end{align*}
\end{minipage}

\item 2-factor type $[6, 4, 2, 2, 2, 2]$: two starters \vspace*{-3mm} \\
\begin{minipage}{8cm}
\begin{align*}
C_0 =& \; [0, [5, 0], 12, [3, 2], 15, [6, 1], 9, [5, 1], 4, [2, 0], 6, [6, 2]] \\
C_1 =& \; [5, [5, 2], 10, [6, 0], 16, [5, -1], 11, [6, -1]] \\
C_2 =& \; [8, [17, 0], 17, [17, 2]] \\
C_3 =& \; [14, [1, -1], 13, [1, 1]] \\
C_4 =& \; [3, [4, -1], 7, [4, 1]] \\
C_5 =& \; [1, [1, 0], 2, [1, 2]]
\end{align*}
\end{minipage}

\vspace*{-4mm}
\begin{minipage}{8cm}
\begin{align*}
C_0' =& \; [9, [8, -1], 1, [2, 2], 16, [8, 0], 7, [17, 1], 17, [17, -1], 12, [3, -1]] \\
C_1' =& \; [5, [8, 1], 14, [3, 0], 0, [8, 2], 8, [3, 1]] \\
C_2' =& \; [4, [2, -1], 2, [2, 1]] \\
C_3' =& \; [13, [7, 0], 6, [7, 2]] \\
C_4' =& \; [11, [4, 0], 15, [4, 2]] \\
C_5' =& \; [10, [7, -1], 3, [7, 1]]
\end{align*}
\end{minipage}

\item 2-factor type $[5, 5, 2, 2, 2, 2]$: two starters \vspace*{-3mm} \\
\begin{minipage}{8cm}
\begin{align*}
C_0 =& \; [6, [4, -1], 10, [3, 1], 7, [2, 0], 9, [4, 2], 5, [1, -1]] \\
C_1 =& \; [11, [7, 0], 4, [5, 2], 16, [4, 0], 3, [1, 2], 2, [8, 1]] \\
C_2 =& \; [13, [5, -1], 1, [5, 1]] \\
C_3 =& \; [12, [2, -1], 14, [2, 1]] \\
C_4 =& \; [15, [7, -1], 8, [7, 1]] \\
C_5 =& \; [17, [17, -1], 0, [17, 1]]
\end{align*}
\end{minipage}

\vspace*{-4mm}
\begin{minipage}{8cm}
\begin{align*}
C_0' =& \; [8, [7, 2], 1, [3, -1], 4, [8, -1], 12, [5, 0], 7, [1, 1]] \\
C_1' =& \; [14, [8, 2], 6, [8, 0], 15, [2, 2], 0, [4, 1], 13, [1, 0]] \\
C_2' =& \; [2, [3, 0], 16, [3, 2]] \\
C_3' =& \; [11, [6, 0], 5, [6, 2]] \\
C_4' =& \; [10, [17, 0], 17, [17, 2]] \\
C_5' =& \; [3, [6, -1], 9, [6, 1]]
\end{align*}
\end{minipage}

\item 2-factor type $[6, 3, 3, 2, 2, 2]$: one starter \vspace*{-3mm} \\
\begin{minipage}{8cm}
\begin{align*}
C_0 =& \; [6, [4, 0], 10, [7, 1], 3, [2, 1], 5, [5, 1], 0, [2, 0], 15, [8, 1]] \\
C_1 =& \; [12, [4, 1], 16, [5, 0], 4, [8, 0]] \\
C_2 =& \; [9, [7, 0], 2, [17, 1], 17, [17, 0]] \\
C_3 =& \; [13, [1, 0], 14, [1, 1]] \\
C_4 =& \; [7, [6, 0], 1, [6, 1]] \\
C_5 =& \; [11, [3, 0], 8, [3, 1]]
\end{align*}
\end{minipage}

\item 2-factor type $[5, 4, 3, 2, 2, 2]$: one starter \vspace*{-3mm} \\
\begin{minipage}{8cm}
\begin{align*}
C_0 =& \; [11, [1, 0], 12, [7, 1], 5, [5, 0], 10, [3, 1], 7, [4, 1]] \\
C_1 =& \; [13, [3, 0], 16, [2, 1], 1, [1, 1], 0, [4, 0]] \\
C_2 =& \; [14, [5, 1], 2, [2, 0], 4, [7, 0]] \\
C_3 =& \; [3, [6, 0], 9, [6, 1]] \\
C_4 =& \; [15, [8, 0], 6, [8, 1]] \\
C_5 =& \; [8, [17, 0], 17, [17, 1]]
\end{align*}
\end{minipage}

\item 2-factor type $[4, 4, 4, 2, 2, 2]$: one starter \vspace*{-3mm} \\
\hspace*{-6mm}
\begin{minipage}{8cm}
\begin{align*}
C_0 =& \; [6, [8, 0], 14, [8, 1], 5, [5, 1], 0, [6, 0]] \\
C_1 =& \; [4, [5, 0], 16, [2, 0], 1, [2, 1], 3, [1, 1]] \\
C_2 =& \; [15, [6, 1], 9, [1, 0], 8, [17, 1], 17, [17, 0]] \\
C_3 =& \; [2, [7, 0], 12, [7, 1]] \\
C_4 =& \; [7, [4, 0], 11, [4, 1]] \\
C_5 =& \; [10, [3, 0], 13, [3, 1]]
\end{align*}
\end{minipage}

\item 2-factor type $[5, 3, 3, 3, 2, 2]$: two starters \vspace*{-3mm} \\
\begin{minipage}{8cm}
\begin{align*}
C_0 =& \; [10, [1, 0], 11, [17, 2], 17, [17, 0], 2, [1, 2], 3, [7, -1]] \\
C_1 =& \; [9, [3, 2], 6, [6, 0], 12, [3, -1]] \\
C_2 =& \; [4, [8, 1], 13, [4, 0], 0, [4, 2]] \\
C_3 =& \; [5, [4, 1], 1, [3, 0], 15, [7, 2]] \\
C_4 =& \; [14, [6, -1], 8, [6, 1]] \\
C_5 =& \; [7, [8, 0], 16, [8, 2]]
\end{align*}
\end{minipage}

\vspace*{-4mm}
\begin{minipage}{8cm}
\begin{align*}
C_0' =& \; [11, [2, 1], 13, [7, 1], 3, [4, -1], 16, [8, -1], 8, [3, 1]] \\
C_1' =& \; [1, [6, 2], 12, [5, 0], 0, [1, 1]] \\
C_2' =& \; [6, [1, -1], 5, [17, 1], 17, [17, -1]] \\
C_3' =& \; [2, [5, 2], 14, [7, 0], 4, [2, -1]] \\
C_4' =& \; [10, [5, -1], 15, [5, 1]] \\
C_5' =& \; [9, [2, 0], 7, [2, 2]]
\end{align*}
\end{minipage}

\newpage

\item 2-factor type $[4, 4, 3, 3, 2, 2]$: two starters \vspace*{-3mm} \\
\hspace*{-5mm}
\begin{minipage}{8cm}
\begin{align*}
C_0 =& \; [10, [8, -1], 1, [2, 0], 3, [4, 2], 7, [3, -1]] \\
C_1 =& \; [6, [8, 0], 14, [6, 2], 8, [17, -1], 17, [17, 1]] \\
C_2 =& \; [4, [4, -1], 0, [8, 2], 9, [5, 0]] \\
C_3 =& \; [5, [7, -1], 15, [4, 1], 2, [3, 1]] \\
C_4 =& \; [16, [5, -1], 11, [5, 1]] \\
C_5 =& \; [13, [1, -1], 12, [1, 1]]
\end{align*}
\end{minipage}
\begin{minipage}{8cm}
\begin{align*}
C_0' =& \; [11, [3, 2], 8, [6, -1], 14, [7, 1], 7, [4, 0]] \\
C_1' =& \; [12, [6, 0], 6, [2, -1], 4, [6, 1], 10, [2, 2]] \\
C_2' =& \; [15, [3, 0], 1, [2, 1], 3, [5, 2]] \\
C_3' =& \; [17, [17, 2], 13, [8, 1], 5, [17, 0]] \\
C_4' =& \; [2, [7, 0], 9, [7, 2]] \\
C_5' =& \; [16, [1, 0], 0, [1, 2]]
\end{align*}
\end{minipage}

\item 2-factor type $[4, 3, 3, 3, 3, 2]$: one starter \vspace*{-3mm} \\
\hspace*{-9mm}
\begin{minipage}{8cm}
\begin{align*}
C_0 =& \; [8, [4, 0], 4, [8, 0], 12, [7, 1], 5, [3, 1]] \\
C_1 =& \; [11, [2, 0], 9, [5, 1], 14, [3, 0]] \\
C_2 =& \; [10, [6, 1], 16, [2, 1], 1, [8, 1]] \\
C_3 =& \; [15, [4, 1], 2, [1, 0], 3, [5, 0]] \\
C_4 =& \; [6, [1, 1], 7, [7, 0], 0, [6, 0]] \\
C_5 =& \; [13, [17, 0], 17, [17, 1]]
\end{align*}
\end{minipage}

\item 2-factor type $[10, 2, 2, 2, 2]$: two starters \vspace*{-3mm} \\
\begin{minipage}{8cm}
\begin{align*}
C_0 =& \; [7, [2, -1], 5, [3, 1], 8, [6, 1], 14, [4, 2], 10, [4, 0], 6, [6, 2], 0, [17, 0], 17, [17, 2], 2, [1, 0], 3, [4, 1]] \\
C_1 =& \; [4, [8, -1], 12, [8, 1]] \\
C_2 =& \; [13, [5, 0], 1, [5, 2]] \\
C_3 =& \; [15, [1, -1], 16, [1, 1]] \\
C_4 =& \; [9, [2, 0], 11, [2, 2]]
\end{align*}
\end{minipage}

\vspace*{-4mm}
\begin{minipage}{8cm}
\begin{align*}
C_0' =& \; [1, [6, 0], 12, [1, 2], 13, [4, -1], 0, [17, -1], 17, [17, 1], 8, [8, 0], 16, [6, -1], 10, [3, -1], 7, [2, 1], \\ & \; 9, [8, 2]] \\
C_1' =& \; [4, [7, -1], 11, [7, 1]] \\
C_2' =& \; [14, [5, -1], 2, [5, 1]] \\
C_3' =& \; [3, [3, 0], 6, [3, 2]] \\
C_4' =& \; [5, [7, 0], 15, [7, 2]]
\end{align*}
\end{minipage}

\item 2-factor type $[9, 3, 2, 2, 2]$: one starter \vspace*{-3mm} \\
\begin{minipage}{8cm}
\begin{align*}
C_0 =& \; [13, [2, 1], 15, [5, 1], 10, [3, 0], 7, [4, 0], 3, [2, 0], 5, [17, 1], 17, [17, 0], 0, [8, 1], 9, [4, 1]] \\
C_1 =& \; [11, [5, 0], 16, [8, 0], 8, [3, 1]] \\
C_2 =& \; [4, [7, 0], 14, [7, 1]] \\
C_3 =& \; [12, [6, 0], 6, [6, 1]] \\
C_4 =& \; [2, [1, 0], 1, [1, 1]]
\end{align*}
\end{minipage}

\item 2-factor type $[8, 4, 2, 2, 2]$: one starter \vspace*{-3mm} \\
\begin{minipage}{8cm}
\begin{align*}
C_0 =& \; [11, [2, 1], 13, [5, 0], 8, [17, 0], 17, [17, 1], 4, [7, 1], 14, [2, 0], 16, [3, 0], 2, [8, 1]] \\
C_1 =& \; [3, [5, 1], 15, [8, 0], 7, [3, 1], 10, [7, 0]] \\
C_2 =& \; [12, [6, 0], 6, [6, 1]] \\
C_3 =& \; [5, [4, 0], 9, [4, 1]] \\
C_4 =& \; [0, [1, 0], 1, [1, 1]]
\end{align*}
\end{minipage}

\item 2-factor type $[7, 5, 2, 2, 2]$: one starter \vspace*{-3mm} \\
\begin{minipage}{8cm}
\begin{align*}
C_0 =& \; [12, [6, 0], 1, [6, 1], 7, [7, 0], 0, [2, 1], 2, [8, 0], 10, [5, 0], 5, [7, 1]] \\
C_1 =& \; [14, [1, 1], 15, [5, 1], 3, [1, 0], 4, [2, 0], 6, [8, 1]] \\
C_2 =& \; [8, [3, 0], 11, [3, 1]] \\
C_3 =& \; [16, [17, 0], 17, [17, 1]] \\
C_4 =& \; [13, [4, 0], 9, [4, 1]]
\end{align*}
\end{minipage}

\item 2-factor type $[6, 6, 2, 2, 2]$: one starter \vspace*{-3mm} \\
\begin{minipage}{8cm}
\begin{align*}
C_0 =& \; [0, [3, 1], 14, [2, 1], 16, [8, 0], 8, [2, 0], 10, [6, 1], 4, [4, 1]] \\
C_1 =& \; [2, [8, 1], 11, [4, 0], 15, [6, 0], 9, [3, 0], 12, [5, 1], 7, [5, 0]] \\
C_2 =& \; [1, [17, 0], 17, [17, 1]] \\
C_3 =& \; [5, [1, 0], 6, [1, 1]] \\
C_4 =& \; [3, [7, 0], 13, [7, 1]]
\end{align*}
\end{minipage}

\item 2-factor type $[8, 3, 3, 2, 2]$: two starters \vspace*{-3mm} \\
\begin{minipage}{8cm}
\begin{align*}
C_0 =& \; [8, [8, 1], 0, [6, 0], 11, [2, -1], 9, [7, 2], 2, [3, -1], 5, [4, 1], 1, [8, 0], 10, [2, 2]] \\
C_1 =& \; [3, [4, 2], 7, [7, 0], 14, [6, 1]] \\
C_2 =& \; [13, [17, -1], 17, [17, 1], 6, [7, -1]] \\
C_3 =& \; [12, [3, 0], 15, [3, 2]] \\
C_4 =& \; [4, [5, -1], 16, [5, 1]]
\end{align*}
\end{minipage}

\vspace*{-4mm}
\begin{minipage}{8cm}
\begin{align*}
C_0' =& \; [17, [17, 0], 6, [4, -1], 2, [2, 1], 4, [8, -1], 13, [8, 2], 5, [4, 0], 9, [6, -1], 3, [17, 2]] \\
C_1' =& \; [0, [3, 1], 14, [2, 0], 16, [1, 2]] \\
C_2' =& \; [11, [7, 1], 1, [6, 2], 12, [1, 0]] \\
C_3' =& \; [15, [5, 0], 10, [5, 2]] \\
C_4' =& \; [7, [1, -1], 8, [1, 1]]
\end{align*}
\end{minipage}

\item 2-factor type $[7, 4, 3, 2, 2]$: two starters \vspace*{-3mm} \\
\begin{minipage}{8cm}
\begin{align*}
C_0 =& \; [8, [4, 0], 12, [5, -1], 7, [3, -1], 4, [5, 2], 9, [1, -1], 10, [6, -1], 16, [8, 1]] \\
C_1 =& \; [13, [7, -1], 3, [2, -1], 5, [1, 0], 6, [7, 2]] \\
C_2 =& \; [2, [8, -1], 11, [7, 1], 1, [1, 1]] \\
C_3 =& \; [15, [2, 0], 0, [2, 2]] \\
C_4 =& \; [17, [17, -1], 14, [17, 1]]
\end{align*}
\end{minipage}

\vspace*{-4mm}
\begin{minipage}{8cm}
\begin{align*}
C_0' =& \; [8, [1, 2], 9, [3, 1], 6, [17, 0], 17, [17, 2], 5, [5, 1], 0, [5, 0], 12, [4, -1]] \\
C_1' =& \; [15, [2, 1], 13, [3, 0], 16, [3, 2], 2, [4, 1]] \\
C_2' =& \; [4, [7, 0], 14, [4, 2], 10, [6, 1]] \\
C_3' =& \; [1, [6, 0], 7, [6, 2]] \\
C_4' =& \; [3, [8, 0], 11, [8, 2]]
\end{align*}
\end{minipage}

\newpage

\item 2-factor type $[6, 5, 3, 2, 2]$: two starters \vspace*{-3mm} \\
\begin{minipage}{8cm}
\begin{align*}
C_0 =& \; [16, [7, 1], 9, [1, 0], 10, [5, 2], 5, [5, 0], 0, [4, 2], 13, [3, -1]] \\
C_1 =& \; [2, [2, 2], 4, [6, 0], 15, [17, 2], 17, [17, 0], 6, [4, -1]] \\
C_2 =& \; [7, [7, -1], 14, [3, 1], 11, [4, 1]] \\
C_3 =& \; [1, [7, 0], 8, [7, 2]] \\
C_4 =& \; [12, [8, -1], 3, [8, 1]]
\end{align*}
\end{minipage}

\vspace*{-4mm}
\begin{minipage}{8cm}
\begin{align*}
C_0' =& \; [4, [2, -1], 6, [3, 0], 3, [5, 1], 8, [3, 2], 11, [4, 0], 15, [6, 2]] \\
C_1' =& \; [7, [5, -1], 12, [1, -1], 13, [17, -1], 17, [17, 1], 9, [2, 1]] \\
C_2' =& \; [0, [2, 0], 2, [1, 2], 1, [1, 1]] \\
C_3' =& \; [14, [8, 0], 5, [8, 2]] \\
C_4' =& \; [10, [6, -1], 16, [6, 1]]
\end{align*}
\end{minipage}

\item 2-factor type $[6, 4, 4, 2, 2]$: two starters \vspace*{-3mm} \\
\begin{minipage}{8cm}
\begin{align*}
C_0 =& \; [0, [7, 0], 7, [17, 2], 17, [17, 0], 11, [1, 1], 12, [4, 1], 16, [1, 2]] \\
C_1 =& \; [4, [5, 2], 9, [7, 1], 2, [6, 0], 8, [4, -1]] \\
C_2 =& \; [14, [1, 0], 15, [7, 2], 5, [8, -1], 13, [1, -1]] \\
C_3 =& \; [1, [2, -1], 3, [2, 1]] \\
C_4 =& \; [10, [4, 0], 6, [4, 2]]
\end{align*}
\end{minipage}

\vspace*{-4mm}
\begin{minipage}{8cm}
\begin{align*}
C_0' =& \; [4, [6, 2], 15, [5, -1], 3, [6, 1], 14, [2, 0], 12, [3, 2], 9, [5, 0]] \\
C_1' =& \; [6, [17, -1], 17, [17, 1], 11, [3, 0], 8, [2, 2]] \\
C_2' =& \; [7, [7, -1], 0, [5, 1], 5, [8, 1], 13, [6, -1]] \\
C_3' =& \; [1, [8, 0], 10, [8, 2]] \\
C_4' =& \; [2, [3, -1], 16, [3, 1]]
\end{align*}
\end{minipage}

\item 2-factor type $[5, 5, 4, 2, 2]$: two starters \vspace*{-3mm} \\
\begin{minipage}{8cm}
\begin{align*}
C_0 =& \; [5, [1, 0], 4, [3, -1], 1, [3, 2], 15, [7, 1], 8, [3, 1]] \\
C_1 =& \; [17, [17, -1], 14, [1, -1], 13, [6, 2], 2, [2, 0], 0, [17, 1]] \\
C_2 =& \; [7, [8, 2], 16, [5, 0], 11, [5, -1], 6, [1, 1]] \\
C_3 =& \; [9, [6, -1], 3, [6, 1]] \\
C_4 =& \; [10, [2, -1], 12, [2, 1]]
\end{align*}
\end{minipage}

\vspace*{-4mm}
\begin{minipage}{8cm}
\begin{align*}
C_0' =& \; [14, [5, 1], 9, [3, 0], 12, [7, 2], 2, [6, 0], 13, [1, 2]] \\
C_1' =& \; [16, [5, 2], 11, [8, 0], 3, [2, 2], 1, [7, 0], 8, [8, 1]] \\
C_2' =& \; [5, [7, -1], 15, [8, -1], 7, [17, 2], 17, [17, 0]] \\
C_3' =& \; [0, [4, 0], 4, [4, 2]] \\
C_4' =& \; [10, [4, -1], 6, [4, 1]]
\end{align*}
\end{minipage}

\item 2-factor type $[7, 3, 3, 3, 2]$: one starter \vspace*{-3mm} \\
\begin{minipage}{8cm}
\begin{align*}
C_0 =& \; [13, [2, 1], 11, [2, 0], 9, [17, 1], 17, [17, 0], 8, [6, 0], 2, [4, 1], 6, [7, 0]] \\
C_1 =& \; [0, [3, 1], 14, [4, 0], 1, [1, 0]] \\
C_2 =& \; [12, [8, 1], 4, [1, 1], 5, [7, 1]] \\
C_3 =& \; [10, [3, 0], 7, [8, 0], 16, [6, 1]] \\
C_4 =& \; [15, [5, 0], 3, [5, 1]]
\end{align*}
\end{minipage}

\item 2-factor type $[6, 4, 3, 3, 2]$: one starter \vspace*{-3mm} \\
\begin{minipage}{8cm}
\begin{align*}
C_0 =& \; [14, [1, 1], 13, [4, 1], 9, [3, 0], 6, [17, 1], 17, [17, 0], 4, [7, 1]] \\
C_1 =& \; [3, [5, 0], 8, [2, 0], 10, [6, 0], 16, [4, 0]] \\
C_2 =& \; [5, [3, 1], 2, [2, 1], 0, [5, 1]] \\
C_3 =& \; [11, [7, 0], 1, [6, 1], 12, [1, 0]] \\
C_4 =& \; [7, [8, 0], 15, [8, 1]]
\end{align*}
\end{minipage}

\item 2-factor type $[5, 5, 3, 3, 2]$: one starter \vspace*{-3mm} \\
\begin{minipage}{8cm}
\begin{align*}
C_0 =& \; [17, [17, 1], 4, [7, 1], 11, [3, 1], 8, [7, 0], 15, [17, 0]] \\
C_1 =& \; [9, [2, 0], 7, [2, 1], 5, [1, 0], 6, [8, 1], 14, [5, 1]] \\
C_2 =& \; [1, [1, 1], 0, [5, 0], 12, [6, 0]] \\
C_3 =& \; [2, [6, 1], 13, [3, 0], 10, [8, 0]] \\
C_4 =& \; [3, [4, 0], 16, [4, 1]]
\end{align*}
\end{minipage}

\item 2-factor type $[5, 4, 4, 3, 2]$: one starter \vspace*{-3mm} \\
\begin{minipage}{8cm}
\begin{align*}
C_0 =& \; [11, [4, 1], 15, [6, 0], 4, [17, 1], 17, [17, 0], 1, [7, 1]] \\
C_1 =& \; [5, [6, 1], 16, [8, 1], 8, [5, 0], 13, [8, 0]] \\
C_2 =& \; [2, [4, 0], 6, [3, 0], 3, [3, 1], 0, [2, 1]] \\
C_3 =& \; [14, [2, 0], 12, [5, 1], 7, [7, 0]] \\
C_4 =& \; [9, [1, 0], 10, [1, 1]]
\end{align*}
\end{minipage}

\item 2-factor type $[4, 4, 4, 4, 2]$: one starter \vspace*{-3mm} \\
\hspace*{-6mm}
\begin{minipage}{8cm}
\begin{align*}
C_0 =& \; [16, [3, 1], 13, [6, 0], 7, [17, 0], 17, [17, 1]] \\
C_1 =& \; [9, [5, 0], 14, [8, 0], 6, [5, 1], 11, [2, 1]] \\
C_2 =& \; [10, [7, 1], 0, [8, 1], 8, [4, 0], 12, [2, 0]] \\
C_3 =& \; [5, [4, 1], 1, [3, 0], 4, [6, 1], 15, [7, 0]] \\
C_4 =& \; [3, [1, 0], 2, [1, 1]]
\end{align*}
\end{minipage}

\item 2-factor type $[6, 3, 3, 3, 3]$: two starters \vspace*{-3mm} \\
\begin{minipage}{8cm}
\begin{align*}
C_0 =& \; [12, [7, -1], 2, [8, -1], 10, [4, -1], 14, [3, 0], 0, [4, 1], 4, [8, 2]] \\
C_1 =& \; [7, [2, 2], 5, [2, 1], 3, [4, 0]] \\
C_2 =& \; [16, [2, 0], 1, [5, -1], 13, [3, 2]] \\
C_3 =& \; [8, [3, -1], 11, [17, -1], 17, [17, 1]] \\
C_4 =& \; [15, [8, 0], 6, [3, 1], 9, [6, 2]]
\end{align*}
\end{minipage}

\vspace*{-4mm}
\begin{minipage}{8cm}
\begin{align*}
C_0' =& \; [1, [17, 0], 17, [17, 2], 10, [6, 1], 4, [2, -1], 6, [7, 0], 13, [5, 2]] \\
C_1' =& \; [9, [5, 1], 14, [6, -1], 8, [1, 1]] \\
C_2' =& \; [16, [8, 1], 7, [7, 2], 0, [1, 0]] \\
C_3' =& \; [15, [5, 0], 3, [1, -1], 2, [4, 2]] \\
C_4' =& \; [5, [6, 0], 11, [1, 2], 12, [7, 1]]
\end{align*}
\end{minipage}

\newpage

\item 2-factor type $[5, 4, 3, 3, 3]$: two starters \vspace*{-3mm} \\
\begin{minipage}{8cm}
\begin{align*}
C_0 =& \; [5, [8, 1], 14, [6, -1], 3, [5, -1], 8, [5, 0], 13, [8, 2]] \\
C_1 =& \; [12, [1, -1], 11, [5, 1], 16, [6, 2], 10, [2, 0]] \\
C_2 =& \; [0, [1, 1], 1, [3, -1], 15, [2, 1]] \\
C_3 =& \; [17, [17, 2], 4, [3, 0], 7, [17, 1]] \\
C_4 =& \; [9, [7, -1], 2, [4, 1], 6, [3, 1]]
\end{align*}
\end{minipage}

\vspace*{-4mm}
\begin{minipage}{8cm}
\begin{align*}
C_0' =& \; [5, [4, 2], 9, [17, -1], 17, [17, 0], 13, [2, 2], 11, [6, 0]] \\
C_1' =& \; [12, [5, 2], 0, [2, -1], 2, [1, 0], 3, [8, -1]] \\
C_2' =& \; [7, [8, 0], 15, [1, 2], 14, [7, 1]] \\
C_3' =& \; [8, [7, 0], 1, [3, 2], 4, [4, -1]] \\
C_4' =& \; [16, [7, 2], 6, [4, 0], 10, [6, 1]]
\end{align*}
\end{minipage}

\item 2-factor type $[4, 4, 4, 3, 3]$: two starters \vspace*{-3mm} \\
\hspace*{-6mm}
\begin{minipage}{8cm}
\begin{align*}
C_0 =& \; [9, [7, 1], 16, [8, 1], 7, [3, -1], 4, [5, 1]] \\
C_1 =& \; [12, [8, 0], 3, [8, 2], 11, [2, -1], 13, [1, 1]] \\
C_2 =& \; [15, [7, -1], 5, [3, 0], 8, [2, 1], 10, [5, 2]] \\
C_3 =& \; [1, [4, 0], 14, [17, 2], 17, [17, 1]] \\
C_4 =& \; [0, [2, 2], 2, [4, -1], 6, [6, 0]]
\end{align*}
\end{minipage}
\begin{minipage}{8cm}
\begin{align*}
C_0' =& \; [4, [5, -1], 9, [17, -1], 17, [17, 0], 11, [7, 2]] \\
C_1' =& \; [10, [3, 2], 13, [1, 0], 12, [6, -1], 6, [4, 1]] \\
C_2' =& \; [1, [1, 2], 0, [3, 1], 14, [6, 1], 8, [7, 0]] \\
C_3' =& \; [5, [6, 2], 16, [8, -1], 7, [2, 0]] \\
C_4' =& \; [15, [5, 0], 3, [1, -1], 2, [4, 2]]
\end{align*}
\end{minipage}

\item 2-factor type $[12, 2, 2, 2]$: one starter \vspace*{-3mm} \\
\begin{minipage}{8cm}
\begin{align*}
C_0 =& \; [15, [5, 0], 3, [3, 1], 6, [5, 1], 11, [1, 0], 10, [17, 0], 17, [17, 1], 14, [3, 0], 0, [1, 1], 16, [6, 0], 5, [2, 0], 7, \\ & \; [6, 1], 13, [2, 1]] \\
C_1 =& \; [12, [7, 0], 2, [7, 1]] \\
C_2 =& \; [1, [8, 0], 9, [8, 1]] \\
C_3 =& \; [8, [4, 0], 4, [4, 1]]
\end{align*}
\end{minipage}

\item 2-factor type $[11, 3, 2, 2]$: two starters \vspace*{-3mm} \\
\begin{minipage}{8cm}
\begin{align*}
C_0 =& \; [12, [6, 1], 1, [5, -1], 13, [7, 2], 3, [4, -1], 7, [2, 1], 5, [4, 0], 9, [5, 2], 4, [7, 0], 14, [6, 2], 8, [6, 0], 2, [7, -1]] \\
C_1 =& \; [16, [17, 2], 17, [17, 0], 0, [1, -1]] \\
C_2 =& \; [6, [8, -1], 15, [8, 1]] \\
C_3 =& \; [10, [1, 0], 11, [1, 2]]
\end{align*}
\end{minipage}

\vspace*{-4mm}
\begin{minipage}{8cm}
\begin{align*}
C_0' =& \; [8, [2, -1], 10, [5, 0], 15, [3, 2], 1, [1, 1], 0, [8, 0], 9, [4, 1], 13, [17, 1], 17, [17, -1], 16, [7, 1], 6, [4, 2], 2, \\ & \; [6, -1]] \\
C_1' =& \; [4, [8, 2], 12, [5, 1], 7, [3, 0]] \\
C_2' =& \; [3, [2, 0], 5, [2, 2]] \\
C_3' =& \; [11, [3, -1], 14, [3, 1]]
\end{align*}
\end{minipage}

\item 2-factor type $[10, 4, 2, 2]$: two starters \vspace*{-3mm} \\
\begin{minipage}{8cm}
\begin{align*}
C_0 =& \; [6, [3, 1], 9, [2, -1], 7, [17, 1], 17, [17, 2], 0, [6, 0], 11, [8, 1], 2, [6, 2], 8, [7, -1], 15, [5, 1], 10, [4, 0]] \\
C_1 =& \; [5, [1, 2], 4, [5, 0], 16, [2, 1], 1, [4, 1]] \\
C_2 =& \; [14, [1, -1], 13, [1, 1]] \\
C_3 =& \; [3, [8, 0], 12, [8, 2]]
\end{align*}
\end{minipage}

\newpage

\begin{minipage}{8cm}
\begin{align*}
C_0' =& \; [7, [6, -1], 13, [5, -1], 1, [7, 1], 11, [3, -1], 14, [6, 1], 8, [4, -1], 12, [17, -1], 17, [17, 0], 6, [7, 2], 16, \\ & \; [8, -1]] \\
C_1' =& \; [5, [4, 2], 9, [1, 0], 10, [5, 2], 15, [7, 0]] \\
C_2' =& \; [2, [2, 0], 4, [2, 2]] \\
C_3' =& \; [0, [3, 0], 3, [3, 2]]
\end{align*}
\end{minipage}

\item 2-factor type $[9, 5, 2, 2]$: two starters \vspace*{-3mm} \\
\begin{minipage}{8cm}
\begin{align*}
C_0 =& \; [13, [17, 2], 17, [17, 0], 12, [3, -1], 9, [7, 1], 16, [8, 2], 8, [6, 1], 2, [5, -1], 7, [1, 1], 6, [7, 0]] \\
C_1 =& \; [1, [3, 0], 4, [6, -1], 15, [4, -1], 11, [1, 2], 10, [8, -1]] \\
C_2 =& \; [14, [6, 0], 3, [6, 2]] \\
C_3 =& \; [5, [5, 0], 0, [5, 2]]
\end{align*}
\end{minipage}

\vspace*{-4mm}
\begin{minipage}{8cm}
\begin{align*}
C_0' =& \; [11, [4, 0], 7, [7, 2], 14, [1, 0], 13, [4, 2], 9, [8, 0], 0, [8, 1], 8, [4, 1], 12, [2, 2], 10, [1, -1]] \\
C_1' =& \; [15, [3, 1], 1, [5, 1], 6, [3, 2], 3, [2, 0], 5, [7, -1]] \\
C_2' =& \; [17, [17, -1], 16, [17, 1]] \\
C_3' =& \; [4, [2, -1], 2, [2, 1]]
\end{align*}
\end{minipage}

\item 2-factor type $[8, 6, 2, 2]$: two starters \vspace*{-3mm} \\
\begin{minipage}{8cm}
\begin{align*}
C_0 =& \; [12, [7, 1], 2, [4, 2], 6, [7, 0], 16, [2, -1], 14, [1, -1], 13, [4, -1], 9, [5, -1], 4, [8, 1]] \\
C_1 =& \; [15, [2, 2], 0, [5, 0], 5, [2, 1], 7, [6, -1], 1, [8, -1], 10, [5, 1]] \\
C_2 =& \; [8, [3, 0], 11, [3, 2]] \\
C_3 =& \; [17, [17, 0], 3, [17, 2]]
\end{align*}
\end{minipage}

\vspace*{-4mm}
\begin{minipage}{8cm}
\begin{align*}
C_0' =& \; [5, [7, -1], 15, [8, 2], 7, [2, 0], 9, [7, 2], 2, [1, 0], 3, [5, 2], 8, [4, 0], 4, [1, 1]] \\
C_1' =& \; [0, [1, 2], 1, [8, 0], 10, [6, 2], 16, [4, 1], 12, [6, 1], 6, [6, 0]] \\
C_2' =& \; [11, [3, -1], 14, [3, 1]] \\
C_3' =& \; [13, [17, -1], 17, [17, 1]]
\end{align*}
\end{minipage}

\item 2-factor type $[7, 7, 2, 2]$: two starters \vspace*{-3mm} \\
\begin{minipage}{8cm}
\begin{align*}
C_0 =& \; [11, [5, 2], 16, [4, -1], 3, [5, 0], 15, [8, -1], 7, [6, 2], 1, [8, 0], 10, [1, 1]] \\
C_1 =& \; [17, [17, 2], 4, [1, -1], 5, [3, 1], 2, [7, 1], 12, [4, 0], 8, [1, 2], 9, [17, 0]] \\
C_2 =& \; [6, [7, 0], 13, [7, 2]] \\
C_3 =& \; [0, [3, 0], 14, [3, 2]]
\end{align*}
\end{minipage}

\vspace*{-4mm}
\begin{minipage}{8cm}
\begin{align*}
C_0' =& \; [9, [6, 1], 3, [4, 1], 16, [1, 0], 15, [8, 2], 7, [2, 0], 5, [8, 1], 13, [4, 2]] \\
C_1' =& \; [1, [3, -1], 4, [17, -1], 17, [17, 1], 11, [6, 0], 0, [7, -1], 10, [2, 2], 12, [6, -1]] \\
C_2' =& \; [14, [5, -1], 2, [5, 1]] \\
C_3' =& \; [6, [2, -1], 8, [2, 1]]
\end{align*}
\end{minipage}

\item 2-factor type $[10, 3, 3, 2]$: one starter \vspace*{-3mm} \\
\begin{minipage}{8cm}
\begin{align*}
C_0 =& \; [8, [7, 1], 1, [8, 1], 10, [4, 1], 6, [3, 0], 3, [1, 0], 2, [17, 1], 17, [17, 0], 13, [2, 0], 15, [6, 1], 9, [1, 1]] \\
C_1 =& \; [7, [4, 0], 11, [7, 0], 4, [3, 1]] \\
C_2 =& \; [16, [6, 0], 5, [8, 0], 14, [2, 1]] \\
C_3 =& \; [12, [5, 0], 0, [5, 1]]
\end{align*}
\end{minipage}

\item 2-factor type $[9, 4, 3, 2]$: one starter \vspace*{-3mm} \\
\begin{minipage}{8cm}
\begin{align*}
C_0 =& \; [5, [8, 1], 14, [2, 1], 16, [7, 0], 9, [17, 0], 17, [17, 1], 2, [8, 0], 10, [4, 0], 6, [5, 1], 11, [6, 1]] \\
C_1 =& \; [7, [6, 0], 1, [7, 1], 8, [4, 1], 4, [3, 0]] \\
C_2 =& \; [3, [3, 1], 0, [2, 0], 15, [5, 0]] \\
C_3 =& \; [13, [1, 0], 12, [1, 1]]
\end{align*}
\end{minipage}

\item 2-factor type $[8, 5, 3, 2]$: one starter \vspace*{-3mm} \\
\begin{minipage}{8cm}
\begin{align*}
C_0 =& \; [1, [7, 0], 8, [4, 0], 12, [4, 1], 16, [5, 0], 4, [1, 1], 3, [6, 1], 9, [17, 1], 17, [17, 0]] \\
C_1 =& \; [13, [2, 0], 15, [1, 0], 14, [3, 0], 0, [6, 0], 6, [7, 1]] \\
C_2 =& \; [10, [5, 1], 5, [2, 1], 7, [3, 1]] \\
C_3 =& \; [11, [8, 0], 2, [8, 1]]
\end{align*}
\end{minipage}

\item 2-factor type $[7, 6, 3, 2]$: one starter \vspace*{-3mm} \\
\begin{minipage}{8cm}
\begin{align*}
C_0 =& \; [1, [5, 0], 6, [5, 1], 11, [7, 1], 4, [4, 0], 8, [1, 1], 9, [8, 0], 0, [1, 0]] \\
C_1 =& \; [13, [3, 1], 10, [17, 1], 17, [17, 0], 5, [2, 1], 7, [8, 1], 15, [2, 0]] \\
C_2 =& \; [16, [4, 1], 12, [7, 0], 2, [3, 0]] \\
C_3 =& \; [3, [6, 0], 14, [6, 1]]
\end{align*}
\end{minipage}

\item 2-factor type $[8, 4, 4, 2]$: one starter \vspace*{-3mm} \\
\begin{minipage}{8cm}
\begin{align*}
C_0 =& \; [4, [8, 0], 12, [6, 1], 1, [7, 0], 8, [17, 0], 17, [17, 1], 16, [5, 1], 11, [4, 1], 15, [6, 0]] \\
C_1 =& \; [3, [3, 1], 6, [1, 0], 7, [3, 0], 10, [7, 1]] \\
C_2 =& \; [14, [5, 0], 9, [4, 0], 5, [8, 1], 13, [1, 1]] \\
C_3 =& \; [2, [2, 0], 0, [2, 1]]
\end{align*}
\end{minipage}

\item 2-factor type $[7, 5, 4, 2]$: one starter \vspace*{-3mm} \\
\begin{minipage}{8cm}
\begin{align*}
C_0 =& \; [2, [4, 0], 6, [7, 1], 16, [6, 1], 10, [2, 1], 12, [8, 0], 4, [3, 1], 1, [1, 1]] \\
C_1 =& \; [5, [7, 0], 15, [17, 0], 17, [17, 1], 13, [2, 0], 11, [6, 0]] \\
C_2 =& \; [8, [1, 0], 7, [4, 1], 3, [3, 0], 0, [8, 1]] \\
C_3 =& \; [14, [5, 0], 9, [5, 1]]
\end{align*}
\end{minipage}

\item 2-factor type $[6, 6, 4, 2]$: one starter \vspace*{-3mm} \\
\begin{minipage}{8cm}
\begin{align*}
C_0 =& \; [12, [2, 1], 14, [3, 1], 11, [7, 1], 4, [1, 1], 3, [7, 0], 13, [1, 0]] \\
C_1 =& \; [0, [5, 1], 5, [4, 1], 1, [6, 1], 7, [17, 0], 17, [17, 1], 2, [2, 0]] \\
C_2 =& \; [6, [3, 0], 9, [6, 0], 15, [5, 0], 10, [4, 0]] \\
C_3 =& \; [16, [8, 0], 8, [8, 1]]
\end{align*}
\end{minipage}

\item 2-factor type $[6, 5, 5, 2]$: one starter \vspace*{-3mm} \\
\begin{minipage}{8cm}
\begin{align*}
C_0 =& \; [10, [17, 0], 17, [17, 1], 7, [4, 1], 3, [8, 1], 11, [2, 0], 9, [1, 1]] \\
C_1 =& \; [8, [4, 0], 4, [8, 0], 12, [6, 0], 1, [5, 0], 6, [2, 1]] \\
C_2 =& \; [14, [5, 1], 2, [6, 1], 13, [3, 0], 16, [1, 0], 0, [3, 1]] \\
C_3 =& \; [5, [7, 0], 15, [7, 1]]
\end{align*}
\end{minipage}

\newpage

\item 2-factor type $[9, 3, 3, 3]$: two starters \vspace*{-3mm} \\
\begin{minipage}{8cm}
\begin{align*}
C_0 =& \; [10, [3, 2], 7, [6, 0], 13, [8, -1], 5, [1, 1], 6, [8, 2], 15, [3, 0], 1, [1, 2], 2, [1, -1], 3, [7, 0]] \\
C_1 =& \; [16, [17, 0], 17, [17, 2], 11, [5, -1]] \\
C_2 =& \; [9, [5, 0], 14, [7, 1], 4, [5, 2]] \\
C_3 =& \; [0, [5, 1], 12, [4, 1], 8, [8, 1]]
\end{align*}
\end{minipage}

\vspace*{-4mm}
\begin{minipage}{8cm}
\begin{align*}
C_0' =& \; [7, [7, -1], 0, [2, 1], 2, [3, -1], 16, [4, 2], 3, [2, 0], 1, [7, 2], 11, [1, 0], 12, [17, 1], 17, [17, -1]] \\
C_1' =& \; [4, [2, 2], 6, [4, 0], 10, [6, -1]] \\
C_2' =& \; [9, [6, 1], 15, [2, -1], 13, [4, -1]] \\
C_3' =& \; [5, [8, 0], 14, [6, 2], 8, [3, 1]]
\end{align*}
\end{minipage}

\item 2-factor type $[8, 4, 3, 3]$: two starters \vspace*{-3mm} \\
\begin{minipage}{8cm}
\begin{align*}
C_0 =& \; [12, [2, 1], 14, [5, 2], 2, [6, 0], 13, [8, 2], 5, [1, -1], 6, [1, 0], 7, [1, 2], 8, [4, 0]] \\
C_1 =& \; [9, [17, 2], 17, [17, 0], 16, [4, 1], 3, [6, 1]] \\
C_2 =& \; [4, [6, -1], 10, [7, -1], 0, [4, -1]] \\
C_3 =& \; [11, [4, 2], 15, [3, -1], 1, [7, 0]]
\end{align*}
\end{minipage}

\vspace*{-4mm}
\begin{minipage}{8cm}
\begin{align*}
C_0' =& \; [11, [2, 2], 9, [5, 0], 14, [2, -1], 12, [7, 1], 2, [8, 1], 10, [7, 2], 0, [8, 0], 8, [3, 1]] \\
C_1' =& \; [7, [6, 2], 13, [5, -1], 1, [3, 0], 15, [8, -1]] \\
C_2' =& \; [3, [3, 2], 6, [1, 1], 5, [2, 0]] \\
C_3' =& \; [16, [17, -1], 17, [17, 1], 4, [5, 1]]
\end{align*}
\end{minipage}

\item 2-factor type $[7, 5, 3, 3]$: two starters \vspace*{-3mm} \\
\begin{minipage}{8cm}
\begin{align*}
C_0 =& \; [15, [6, -1], 4, [4, 0], 8, [6, 2], 2, [6, 1], 13, [7, -1], 3, [8, 1], 12, [3, -1]] \\
C_1 =& \; [11, [1, 2], 10, [4, 1], 6, [7, 0], 16, [17, 2], 17, [17, 0]] \\
C_2 =& \; [5, [2, 2], 7, [2, 0], 9, [4, -1]] \\
C_3 =& \; [14, [3, 0], 0, [1, 1], 1, [4, 2]]
\end{align*}
\end{minipage}

\vspace*{-4mm}
\begin{minipage}{8cm}
\begin{align*}
C_0' =& \; [16, [6, 0], 5, [5, 2], 0, [8, 0], 9, [7, 2], 2, [7, 1], 12, [5, 0], 7, [8, 2]] \\
C_1' =& \; [17, [17, 1], 11, [5, 1], 6, [8, -1], 14, [1, -1], 15, [17, -1]] \\
C_2' =& \; [8, [5, -1], 13, [3, 1], 10, [2, 1]] \\
C_3' =& \; [4, [1, 0], 3, [2, -1], 1, [3, 2]]
\end{align*}
\end{minipage}

\item 2-factor type $[6, 6, 3, 3]$: two starters \vspace*{-3mm} \\
\begin{minipage}{8cm}
\begin{align*}
C_0 =& \; [10, [8, -1], 1, [2, 0], 16, [7, -1], 9, [17, 2], 17, [17, 1], 12, [2, 1]] \\
C_1 =& \; [14, [7, 2], 7, [7, 0], 0, [8, 2], 8, [5, -1], 13, [6, -1], 2, [5, 0]] \\
C_2 =& \; [15, [4, 0], 11, [6, 2], 5, [7, 1]] \\
C_3 =& \; [4, [1, -1], 3, [3, 1], 6, [2, -1]]
\end{align*}
\end{minipage}

\vspace*{-4mm}
\begin{minipage}{8cm}
\begin{align*}
C_0' =& \; [9, [3, 0], 12, [4, -1], 8, [2, 2], 6, [5, 1], 1, [1, 1], 0, [8, 1]] \\
C_1' =& \; [11, [4, 1], 15, [1, 2], 14, [17, 0], 17, [17, -1], 7, [4, 2], 3, [8, 0]] \\
C_2' =& \; [13, [3, 2], 16, [3, -1], 2, [6, 0]] \\
C_3' =& \; [4, [6, 1], 10, [5, 2], 5, [1, 0]]
\end{align*}
\end{minipage}

\newpage

\item 2-factor type $[7, 4, 4, 3]$: two starters \vspace*{-3mm} \\
\begin{minipage}{8cm}
\begin{align*}
C_0 =& \; [0, [5, 0], 12, [1, -1], 11, [5, 2], 6, [7, 0], 13, [3, 2], 10, [1, 0], 9, [8, 2]] \\
C_1 =& \; [4, [3, 1], 1, [17, 0], 17, [17, -1], 3, [1, 2]] \\
C_2 =& \; [2, [4, 1], 15, [1, 1], 14, [7, 1], 7, [5, 1]] \\
C_3 =& \; [16, [6, 2], 5, [3, 0], 8, [8, 1]]
\end{align*}
\end{minipage}

\vspace*{-4mm}
\begin{minipage}{8cm}
\begin{align*}
C_0' =& \; [11, [4, 0], 7, [2, 1], 9, [17, 2], 17, [17, 1], 10, [4, -1], 14, [6, -1], 3, [8, -1]] \\
C_1' =& \; [8, [7, 2], 1, [3, -1], 4, [2, 0], 2, [6, 1]] \\
C_2' =& \; [15, [2, -1], 0, [5, -1], 5, [8, 0], 13, [2, 2]] \\
C_3' =& \; [12, [4, 2], 16, [7, -1], 6, [6, 0]]
\end{align*}
\end{minipage}

\item 2-factor type $[6, 5, 4, 3]$: two starters \vspace*{-3mm} \\
\begin{minipage}{8cm}
\begin{align*}
C_0 =& \; [7, [7, 2], 14, [2, 0], 16, [4, 2], 12, [8, 1], 4, [17, -1], 17, [17, 0]] \\
C_1 =& \; [15, [5, 0], 10, [1, 2], 11, [7, 0], 1, [8, 2], 9, [6, -1]] \\
C_2 =& \; [6, [6, 1], 0, [5, -1], 5, [8, -1], 13, [7, 1]] \\
C_3 =& \; [2, [1, 0], 3, [5, 1], 8, [6, 2]]
\end{align*}
\end{minipage}

\vspace*{-4mm}
\begin{minipage}{8cm}
\begin{align*}
C_0' =& \; [7, [2, -1], 9, [3, -1], 12, [8, 0], 4, [1, 1], 5, [2, 2], 3, [4, -1]] \\
C_1' =& \; [13, [3, 0], 16, [1, -1], 15, [2, 1], 0, [7, -1], 10, [3, 2]] \\
C_2' =& \; [8, [6, 0], 2, [4, 1], 6, [5, 2], 11, [3, 1]] \\
C_3' =& \; [17, [17, 2], 14, [4, 0], 1, [17, 1]]
\end{align*}
\end{minipage}

\item 2-factor type $[5, 5, 5, 3]$: two starters \vspace*{-3mm} \\
\begin{minipage}{8cm}
\begin{align*}
C_0 =& \; [3, [7, 0], 10, [1, 1], 11, [2, 2], 9, [5, 0], 4, [1, 2]] \\
C_1 =& \; [15, [17, 0], 17, [17, 2], 7, [6, 1], 1, [4, 1], 14, [1, -1]] \\
C_2 =& \; [16, [3, 2], 13, [7, 1], 6, [6, -1], 12, [7, -1], 2, [3, 0]] \\
C_3 =& \; [0, [5, 1], 5, [3, 1], 8, [8, -1]]
\end{align*}
\end{minipage}

\vspace*{-4mm}
\begin{minipage}{8cm}
\begin{align*}
C_0' =& \; [9, [2, 1], 7, [1, 0], 6, [6, 2], 0, [17, -1], 17, [17, 1]] \\
C_1' =& \; [2, [3, -1], 16, [5, 2], 11, [6, 0], 5, [7, 2], 15, [4, 0]] \\
C_2' =& \; [12, [8, 1], 3, [2, -1], 1, [8, 2], 10, [4, -1], 14, [2, 0]] \\
C_3' =& \; [4, [4, 2], 8, [5, -1], 13, [8, 0]]
\end{align*}
\end{minipage}

\item 2-factor type $[6, 4, 4, 4]$: two starters \vspace*{-3mm} \\
\begin{minipage}{8cm}
\begin{align*}
C_0 =& \; [13, [2, -1], 15, [7, 0], 8, [3, 2], 11, [5, 0], 16, [5, 2], 4, [8, 1]] \\
C_1 =& \; [6, [6, 1], 12, [7, 2], 5, [3, 0], 2, [4, 1]] \\
C_2 =& \; [7, [17, -1], 17, [17, 1], 3, [2, 1], 1, [6, -1]] \\
C_3 =& \; [14, [3, 1], 0, [7, -1], 10, [1, -1], 9, [5, 1]]
\end{align*}
\end{minipage}

\vspace*{-4mm}
\begin{minipage}{8cm}
\begin{align*}
C_0' =& \; [6, [7, 1], 13, [3, -1], 10, [1, 2], 11, [1, 0], 12, [17, 2], 17, [17, 0]] \\
C_1' =& \; [16, [2, 0], 14, [6, 2], 3, [4, 0], 7, [8, 2]] \\
C_2' =& \; [9, [8, 0], 1, [4, 2], 5, [1, 1], 4, [5, -1]] \\
C_3' =& \; [8, [8, -1], 0, [2, 2], 15, [4, -1], 2, [6, 0]]
\end{align*}
\end{minipage}

\newpage

\item 2-factor type $[5, 5, 4, 4]$: two starters \vspace*{-3mm} \\
\begin{minipage}{8cm}
\begin{align*}
C_0 =& \; [7, [3, 1], 4, [8, -1], 12, [8, 0], 3, [5, 2], 15, [8, 1]] \\
C_1 =& \; [17, [17, 0], 14, [6, 2], 8, [6, 0], 2, [3, 2], 5, [17, -1]] \\
C_2 =& \; [16, [5, 0], 11, [2, 2], 9, [1, -1], 10, [6, -1]] \\
C_3 =& \; [13, [4, 0], 0, [1, 1], 1, [5, 1], 6, [7, 2]]
\end{align*}
\end{minipage}

\vspace*{-4mm}
\begin{minipage}{8cm}
\begin{align*}
C_0' =& \; [10, [7, 0], 0, [5, -1], 12, [8, 2], 4, [2, -1], 6, [4, -1]] \\
C_1' =& \; [15, [6, 1], 9, [2, 0], 7, [4, 1], 11, [17, 1], 17, [17, 2]] \\
C_2' =& \; [3, [4, 2], 16, [2, 1], 14, [1, 0], 13, [7, 1]] \\
C_3' =& \; [2, [1, 2], 1, [7, -1], 8, [3, 0], 5, [3, -1]]
\end{align*}
\end{minipage}

\item 2-factor type $[14, 2, 2]$: two starters \vspace*{-3mm} \\
\begin{minipage}{8cm}
\begin{align*}
C_0 =& \; [6, [5, 0], 1, [7, 2], 8, [17, -1], 17, [17, 0], 16, [1, 2], 0, [2, 1], 15, [6, 1], 9, [4, -1], 13, [1, 1], 12, [8, -1], 3,  \\ & \; [2, -1], 5, [5, -1], 10, [1, -1], 11, [5, 1]] \\
C_1 =& \; [2, [2, 0], 4, [2, 2]] \\
C_2 =& \; [14, [7, -1], 7, [7, 1]]
\end{align*}
\end{minipage}

\vspace*{-4mm}
\begin{minipage}{8cm}
\begin{align*}
C_0' =& \; [15, [17, 1], 17, [17, 2], 16, [4, 0], 3, [8, 1], 11, [5, 2], 6, [6, -1], 12, [7, 0], 5, [3, 2], 2, [6, 0], 8, [4, 2], 4, \\ & \; [3, 0], 7, [6, 2], 1, [4, 1], 14, [1, 0]] \\
C_1' =& \; [9, [8, 0], 0, [8, 2]] \\
C_2' =& \; [10, [3, -1], 13, [3, 1]]
\end{align*}
\end{minipage}

\item 2-factor type $[13, 3, 2]$: one starter \vspace*{-3mm} \\
\begin{minipage}{8cm}
\begin{align*}
C_0 =& \; [1, [17, 0], 17, [17, 1], 7, [1, 1], 8, [3, 0], 11, [8, 0], 2, [2, 1], 4, [5, 0], 16, [7, 1], 6, [4, 1], 10, [7, 0], 0, [5, 1], \\ & \; 5, [4, 0], 9, [8, 1]] \\
C_1 =& \; [12, [1, 0], 13, [2, 0], 15, [3, 1]] \\
C_2 =& \; [14, [6, 0], 3, [6, 1]]
\end{align*}
\end{minipage}

\item 2-factor type $[12, 4, 2]$: one starter \vspace*{-3mm} \\
\begin{minipage}{8cm}
\begin{align*}
C_0 =& \; [4, [8, 1], 13, [4, 0], 9, [7, 0], 16, [8, 0], 7, [7, 1], 14, [3, 0], 11, [4, 1], 15, [17, 1], 17, [17, 0], 3, [1, 0], 2, \\ & \; [3, 1], 5, [1, 1]] \\
C_1 =& \; [0, [5, 0], 12, [6, 1], 1, [5, 1], 6, [6, 0]] \\
C_2 =& \; [8, [2, 0], 10, [2, 1]]
\end{align*}
\end{minipage}

\item 2-factor type $[11, 5, 2]$: one starter \vspace*{-3mm} \\
\begin{minipage}{8cm}
\begin{align*}
C_0 =& \; [17, [17, 1], 11, [3, 0], 8, [7, 1], 15, [6, 0], 4, [5, 1], 16, [1, 1], 0, [1, 0], 1, [4, 0], 5, [5, 0], 10, [3, 1], 13, [17, 0]] \\
C_1 =& \; [12, [7, 0], 2, [4, 1], 6, [8, 1], 14, [6, 1], 3, [8, 0]] \\
C_2 =& \; [7, [2, 0], 9, [2, 1]]
\end{align*}
\end{minipage}

\item 2-factor type $[10, 6, 2]$: one starter \vspace*{-3mm} \\
\begin{minipage}{8cm}
\begin{align*}
C_0 =& \; [10, [7, 0], 0, [8, 0], 8, [5, 1], 13, [1, 1], 12, [6, 0], 6, [1, 0], 7, [8, 1], 16, [2, 0], 1, [4, 1], 14, [4, 0]] \\
C_1 =& \; [11, [7, 1], 4, [17, 1], 17, [17, 0], 3, [5, 0], 15, [6, 1], 9, [2, 1]] \\
C_2 =& \; [2, [3, 0], 5, [3, 1]]
\end{align*}
\end{minipage}

\item 2-factor type $[9, 7, 2]$: one starter \vspace*{-3mm} \\
\begin{minipage}{8cm}
\begin{align*}
C_0 =& \; [9, [3, 1], 12, [4, 0], 8, [3, 0], 5, [6, 0], 11, [5, 1], 16, [8, 0], 7, [8, 1], 15, [5, 0], 10, [1, 0]] \\
C_1 =& \; [1, [1, 1], 2, [2, 0], 0, [6, 1], 6, [2, 1], 4, [17, 1], 17, [17, 0], 14, [4, 1]] \\
C_2 =& \; [3, [7, 0], 13, [7, 1]]
\end{align*}
\end{minipage}

\item 2-factor type $[8, 8, 2]$: one starter \vspace*{-3mm} \\
\begin{minipage}{8cm}
\begin{align*}
C_0 =& \; [4, [7, 0], 14, [3, 1], 0, [4, 0], 13, [1, 0], 12, [7, 1], 5, [3, 0], 2, [4, 1], 6, [2, 1]] \\
C_1 =& \; [7, [6, 0], 1, [8, 0], 10, [17, 0], 17, [17, 1], 11, [2, 0], 9, [6, 1], 15, [1, 1], 16, [8, 1]] \\
C_2 =& \; [3, [5, 0], 8, [5, 1]]
\end{align*}
\end{minipage}

\item 2-factor type $[12, 3, 3]$: two starters \vspace*{-3mm} \\
\begin{minipage}{8cm}
\begin{align*}
C_0 =& \; [12, [1, -1], 11, [6, 1], 0, [7, 1], 7, [17, 2], 17, [17, 1], 14, [6, -1], 3, [6, 0], 9, [1, 2], 8, [3, 1], 5, [7, 0], 15, \\ & \; [6, 2], 4, [8, 0]] \\
C_1 =& \; [13, [5, 0], 1, [2, -1], 16, [3, 2]] \\
C_2 =& \; [6, [4, -1], 2, [8, 2], 10, [4, 0]]
\end{align*}
\end{minipage}

\vspace*{-4mm}
\begin{minipage}{8cm}
\begin{align*}
C_0' =& \; [0, [8, -1], 8, [1, 0], 7, [7, 2], 14, [2, 0], 12, [4, 1], 16, [3, -1], 13, [8, 1], 4, [5, 2], 9, [3, 0], 6, [5, -1], 1, \\ & \; [1, 1], 2, [2, 2]] \\
C_1' =& \; [17, [17, 0], 11, [4, 2], 15, [17, -1]] \\
C_2' =& \; [10, [5, 1], 5, [2, 1], 3, [7, -1]]
\end{align*}
\end{minipage}

\item 2-factor type $[11, 4, 3]$: two starters \vspace*{-3mm} \\
\begin{minipage}{8cm}
\begin{align*}
C_0 =& \; [16, [2, 1], 14, [17, -1], 17, [17, 1], 13, [3, 0], 10, [8, 1], 1, [6, 2], 12, [7, 0], 2, [3, 2], 5, [6, 0], 11, [2, -1], \\ & \; 9, [7, 2]] \\
C_1 =& \; [0, [7, -1], 7, [4, 0], 3, [5, 1], 8, [8, 2]] \\
C_2 =& \; [6, [2, 2], 4, [6, 1], 15, [8, 0]]
\end{align*}
\end{minipage}

\vspace*{-4mm}
\begin{minipage}{8cm}
\begin{align*}
C_0' =& \; [4, [7, 1], 11, [5, -1], 6, [4, 1], 10, [1, -1], 9, [3, 1], 12, [4, 2], 16, [3, -1], 2, [17, 0], 17, [17, 2], 14, [6, -1], \\ & \; 3, [1, 0]] \\
C_1' =& \; [13, [2, 0], 15, [8, -1], 7, [1, 1], 8, [5, 2]] \\
C_2' =& \; [0, [5, 0], 5, [4, -1], 1, [1, 2]]
\end{align*}
\end{minipage}

\item 2-factor type $[10, 5, 3]$: two starters \vspace*{-3mm} \\
\begin{minipage}{8cm}
\begin{align*}
C_0 =& \; [5, [17, -1], 17, [17, 1], 8, [5, 0], 3, [8, -1], 12, [5, -1], 7, [2, 1], 9, [4, 1], 13, [7, 2], 6, [6, 0], 0, [5, 2]] \\
C_1 =& \; [10, [6, -1], 4, [7, 1], 11, [8, 2], 2, [1, 1], 1, [8, 0]] \\
C_2 =& \; [14, [2, 0], 16, [1, -1], 15, [1, 2]]
\end{align*}
\end{minipage}

\vspace*{-4mm}
\begin{minipage}{8cm}
\begin{align*}
C_0' =& \; [10, [3, 2], 13, [17, 0], 17, [17, 2], 4, [5, 1], 16, [2, -1], 1, [3, 0], 15, [3, -1], 12, [7, -1], 5, [4, 2], 9, [1, 0]] \\
C_1' =& \; [0, [3, 1], 14, [7, 0], 7, [4, -1], 3, [8, 1], 11, [6, 2]] \\
C_2' =& \; [2, [4, 0], 6, [2, 2], 8, [6, 1]]
\end{align*}
\end{minipage}

\newpage

\item 2-factor type $[9, 6, 3]$: two starters \vspace*{-3mm} \\
\begin{minipage}{8cm}
\begin{align*}
C_0 =& \; [16, [1, 2], 0, [6, 0], 6, [1, 1], 7, [2, 1], 9, [4, 1], 13, [7, 1], 3, [8, 1], 11, [7, 2], 4, [5, 0]] \\
C_1 =& \; [1, [3, -1], 15, [4, 2], 2, [7, 0], 12, [2, 2], 14, [8, -1], 5, [4, 0]] \\
C_2 =& \; [10, [2, -1], 8, [17, 2], 17, [17, 0]]
\end{align*}
\end{minipage}

\vspace*{-4mm}
\begin{minipage}{8cm}
\begin{align*}
C_0' =& \; [2, [6, 2], 8, [5, 1], 3, [3, 1], 0, [7, -1], 7, [4, -1], 11, [1, 0], 12, [6, -1], 6, [17, 1], 17, [17, -1]] \\
C_1' =& \; [14, [2, 0], 16, [1, -1], 15, [6, 1], 4, [3, 2], 1, [8, 0], 9, [5, 2]] \\
C_2' =& \; [10, [3, 0], 13, [8, 2], 5, [5, -1]]
\end{align*}
\end{minipage}

\item 2-factor type $[8, 7, 3]$: two starters \vspace*{-3mm} \\
\begin{minipage}{8cm}
\begin{align*}
C_0 =& \; [6, [7, 2], 16, [8, -1], 7, [6, -1], 13, [17, 0], 17, [17, -1], 5, [5, 2], 0, [3, 0], 3, [3, 1]] \\
C_1 =& \; [4, [4, 2], 8, [7, -1], 15, [1, 1], 14, [5, 0], 2, [8, 1], 10, [1, 2], 11, [7, 0]] \\
C_2 =& \; [1, [6, 0], 12, [3, -1], 9, [8, 2]]
\end{align*}
\end{minipage}

\vspace*{-4mm}
\begin{minipage}{8cm}
\begin{align*}
C_0' =& \; [16, [1, -1], 0, [1, 0], 1, [6, 1], 7, [7, 1], 14, [3, 2], 11, [8, 0], 2, [2, 2], 4, [5, 1]] \\
C_1' =& \; [10, [4, 0], 6, [2, 1], 8, [4, 1], 12, [17, 2], 17, [17, 1], 3, [2, -1], 5, [5, -1]] \\
C_2' =& \; [9, [4, -1], 13, [2, 0], 15, [6, 2]]
\end{align*}
\end{minipage}

\item 2-factor type $[10, 4, 4]$: two starters \vspace*{-3mm} \\
\begin{minipage}{8cm}
\begin{align*}
C_0 =& \; [3, [7, 2], 13, [3, 1], 10, [4, 0], 14, [8, 1], 5, [1, -1], 4, [2, 1], 6, [5, 2], 1, [1, 0], 2, [8, 2], 11, [8, 0]] \\
C_1 =& \; [12, [4, -1], 8, [1, 2], 7, [7, 0], 0, [5, -1]] \\
C_2 =& \; [15, [1, 1], 16, [17, 1], 17, [17, -1], 9, [6, 1]]
\end{align*}
\end{minipage}

\vspace*{-4mm}
\begin{minipage}{8cm}
\begin{align*}
C_0' =& \; [8, [4, 2], 12, [17, 0], 17, [17, 2], 16, [5, 0], 4, [6, 2], 15, [4, 1], 11, [7, -1], 1, [5, 1], 13, [3, 0], 10, \\ & \; [2, -1]] \\
C_1' =& \; [3, [6, 0], 14, [8, -1], 6, [6, -1], 0, [3, 2]] \\
C_2' =& \; [2, [3, -1], 5, [2, 0], 7, [2, 2], 9, [7, 1]]
\end{align*}
\end{minipage}

\item 2-factor type $[9, 5, 4]$: two starters \vspace*{-3mm} \\
\begin{minipage}{8cm}
\begin{align*}
C_0 =& \; [1, [1, 0], 2, [4, 2], 15, [1, -1], 16, [7, -1], 6, [2, 0], 4, [1, 2], 5, [5, -1], 10, [2, -1], 12, [6, -1]] \\
C_1 =& \; [17, [17, 1], 9, [4, -1], 13, [7, 0], 3, [5, 2], 8, [17, -1]] \\
C_2 =& \; [11, [4, 0], 7, [7, 1], 14, [3, 2], 0, [6, 1]]
\end{align*}
\end{minipage}

\vspace*{-4mm}
\begin{minipage}{8cm}
\begin{align*}
C_0' =& \; [8, [8, 1], 0, [8, 0], 9, [17, 2], 17, [17, 0], 16, [6, 2], 5, [3, 0], 2, [7, 2], 12, [1, 1], 11, [3, 1]] \\
C_1' =& \; [6, [2, 1], 4, [3, -1], 7, [4, 1], 3, [6, 0], 14, [8, 2]] \\
C_2' =& \; [1, [8, -1], 10, [5, 1], 15, [2, 2], 13, [5, 0]]
\end{align*}
\end{minipage}

\item 2-factor type $[8, 6, 4]$: two starters \vspace*{-3mm} \\
\begin{minipage}{8cm}
\begin{align*}
C_0 =& \; [11, [7, -1], 1, [5, 1], 13, [2, 0], 15, [4, 2], 2, [6, 0], 8, [1, 1], 9, [2, 2], 7, [4, -1]] \\
C_1 =& \; [3, [4, 0], 16, [7, 1], 6, [1, 2], 5, [5, 0], 0, [3, 2], 14, [6, -1]] \\
C_2 =& \; [12, [8, 1], 4, [17, 0], 17, [17, 2], 10, [2, -1]]
\end{align*}
\end{minipage}

\vspace*{-4mm}
\begin{minipage}{8cm}
\begin{align*}
C_0' =& \; [1, [3, 0], 15, [4, 1], 2, [5, 2], 7, [1, 0], 8, [3, 1], 11, [2, 1], 13, [3, -1], 10, [8, 2]] \\
C_1' =& \; [16, [7, 0], 9, [5, -1], 4, [1, -1], 3, [8, -1], 12, [7, 2], 5, [6, 1]] \\
C_2' =& \; [0, [6, 2], 6, [8, 0], 14, [17, 1], 17, [17, -1]]
\end{align*}
\end{minipage}

\item 2-factor type $[7, 7, 4]$: two starters \vspace*{-3mm} \\
\begin{minipage}{8cm}
\begin{align*}
C_0 =& \; [3, [7, -1], 10, [4, 0], 14, [8, 1], 5, [6, 2], 16, [4, 1], 12, [3, -1], 15, [5, -1]] \\
C_1 =& \; [13, [7, 0], 6, [17, 2], 17, [17, 0], 11, [7, 2], 4, [2, 1], 2, [1, 1], 1, [5, 1]] \\
C_2 =& \; [7, [2, 0], 9, [8, 2], 0, [8, 0], 8, [1, 2]]
\end{align*}
\end{minipage}

\vspace*{-4mm}
\begin{minipage}{8cm}
\begin{align*}
C_0' =& \; [1, [1, 0], 0, [3, 2], 14, [6, 0], 8, [5, 2], 13, [3, 0], 10, [2, 2], 12, [6, -1]] \\
C_1' =& \; [4, [17, 1], 17, [17, -1], 16, [3, 1], 2, [8, -1], 11, [4, -1], 7, [1, -1], 6, [2, -1]] \\
C_2' =& \; [9, [4, 2], 5, [7, 1], 15, [5, 0], 3, [6, 1]]
\end{align*}
\end{minipage}

\item 2-factor type $[8, 5, 5]$: two starters \vspace*{-3mm} \\
\begin{minipage}{8cm}
\begin{align*}
C_0 =& \; [16, [6, -1], 10, [5, -1], 5, [8, 1], 13, [8, 2], 4, [2, 0], 6, [5, 2], 1, [6, 0], 12, [4, 1]] \\
C_1 =& \; [3, [8, -1], 11, [3, -1], 14, [3, 0], 0, [7, 1], 7, [4, 2]] \\
C_2 =& \; [15, [4, 0], 2, [7, 2], 9, [17, -1], 17, [17, 1], 8, [7, -1]]
\end{align*}
\end{minipage}

\vspace*{-4mm}
\begin{minipage}{8cm}
\begin{align*}
C_0' =& \; [16, [1, 0], 15, [5, 1], 3, [3, 2], 6, [6, 1], 0, [1, -1], 1, [4, -1], 5, [8, 0], 14, [2, 2]] \\
C_1' =& \; [10, [6, 2], 4, [2, 1], 2, [7, 0], 9, [2, -1], 7, [3, 1]] \\
C_2' =& \; [17, [17, 2], 8, [5, 0], 13, [1, 2], 12, [1, 1], 11, [17, 0]]
\end{align*}
\end{minipage}

\item 2-factor type $[7, 6, 5]$: two starters \vspace*{-3mm} \\
\begin{minipage}{8cm}
\begin{align*}
C_0 =& \; [14, [1, 1], 15, [4, -1], 11, [3, -1], 8, [8, -1], 0, [1, -1], 16, [2, 2], 1, [4, 0]] \\
C_1 =& \; [13, [7, 1], 6, [2, 0], 4, [1, 2], 5, [5, 0], 10, [3, 2], 7, [6, -1]] \\
C_2 =& \; [9, [7, -1], 2, [17, 1], 17, [17, 2], 3, [8, 1], 12, [3, 0]]
\end{align*}
\end{minipage}

\vspace*{-4mm}
\begin{minipage}{8cm}
\begin{align*}
C_0' =& \; [15, [17, 0], 17, [17, -1], 16, [7, 2], 6, [8, 0], 14, [2, -1], 12, [6, 2], 1, [3, 1]] \\
C_1' =& \; [9, [2, 1], 11, [8, 2], 3, [1, 0], 4, [6, 1], 10, [5, -1], 5, [4, 1]] \\
C_2' =& \; [13, [5, 1], 8, [6, 0], 2, [5, 2], 7, [7, 0], 0, [4, 2]]
\end{align*}
\end{minipage}

\item 2-factor type $[16, 2]$: one starter \vspace*{-3mm} \\
\begin{minipage}{8cm}
\begin{align*}
 C_0 =& \; [11, [5, 0], 16, [4, 1], 12, [2, 0], 10, [5, 1], 5, [1, 1], 4, [3, 1], 7, [4, 0], 3, [6, 1], 14, [17, 0], 17, [17, 1], 2, \\ & \; [1, 0], 1, [8, 0], 9, [6, 0], 15, [2, 1], 0, [8, 1], 8, [3, 0]] \\
C_1 =& \; [6, [7, 0], 13, [7, 1]]
\end{align*}
\end{minipage}

\item 2-factor type $[15, 3]$: two starters \vspace*{-3mm} \\
\begin{minipage}{8cm}
\begin{align*}
C_0 =& \; [15, [8, 1], 6, [17, 1], 17, [17, -1], 16, [8, 2], 8, [1, 0], 7, [5, 2], 2, [2, -1], 4, [3, 1], 1, [5, 0], 13, [4, -1], 9, \\ & \; [3, 2], 12, [7, 0], 5, [8, -1], 14, [6, 1], 3, [5, -1]] \\
C_1 =& \; [11, [1, 2], 10, [7, -1], 0, [6, 0]]
\end{align*}
\end{minipage}

\vspace*{-4mm}
\begin{minipage}{8cm}
\begin{align*}
C_0' =& \; [0, [5, 1], 5, [4, 2], 9, [3, 0], 6, [2, 1], 8, [17, 2], 17, [17, 0], 13, [6, -1], 7, [6, 2], 1, [4, 0], 14, [3, -1], 11, \\ & \; [1, -1], 10, [7, 2], 3, [4, 1], 16, [1, 1], 15, [2, 0]] \\
C_1' =& \; [4, [8, 0], 12, [7, 1], 2, [2, 2]]
\end{align*}
\end{minipage}

\item 2-factor type $[14, 4]$: two starters \vspace*{-3mm} \\
\begin{minipage}{8cm}
\begin{align*}
C_0 =& \; [7, [6, 1], 1, [3, 0], 15, [2, 1], 0, [8, 2], 9, [1, 0], 10, [4, 2], 6, [2, 0], 8, [3, 1], 11, [6, 2], 5, [17, 0], 17, [17, 2], \\ & \; 16, [3, -1], 2, [1, -1], 3, [4, -1]] \\
C_1 =& \; [13, [1, 2], 14, [7, -1], 4, [8, 0], 12, [1, 1]]
\end{align*}
\end{minipage}

\newpage

\begin{minipage}{8cm}
\begin{align*}
C_0' =& \; [9, [3, 2], 12, [8, -1], 3, [5, 1], 15, [17, -1], 17, [17, 1], 2, [6, -1], 8, [8, 1], 0, [5, 0], 5, [2, 2], 7, [6, 0], 13, \\ & \; [7, 2], 6, [4, 0], 10, [4, 1], 14, [5, -1]] \\
C_1' =& \; [4, [7, 0], 11, [7, 1], 1, [2, -1], 16, [5, 2]]
\end{align*}
\end{minipage}

\item 2-factor type $[13, 5]$: two starters \vspace*{-3mm} \\
\begin{minipage}{8cm}
\begin{align*}
C_0 =& \; [8, [2, 1], 6, [4, 1], 2, [2, 0], 4, [6, -1], 15, [4, 2], 11, [5, 0], 16, [6, 1], 5, [2, 2], 3, [7, -1], 10, [7, 0], 0, [4, -1], \\ & \; 13, [6, 2], 7, [1, -1]] \\
C_1 =& \; [9, [3, 1], 12, [17, 2], 17, [17, 1], 1, [4, 0], 14, [5, -1]]
\end{align*}
\end{minipage}

\vspace*{-4mm}
\begin{minipage}{8cm}
\begin{align*}
C_0' =& \; [11, [7, 1], 1, [8, 2], 9, [17, 0], 17, [17, -1], 13, [5, 2], 8, [8, -1], 16, [3, 0], 2, [2, -1], 0, [3, 2], 14, [8, 0], 5, \\ & \; [7, 2], 15, [5, 1], 10, [1, 0]] \\
C_1' =& \; [4, [3, -1], 7, [1, 1], 6, [6, 0], 12, [8, 1], 3, [1, 2]]
\end{align*}
\end{minipage}

\item 2-factor type $[12, 6]$: two starters \vspace*{-3mm} \\
\begin{minipage}{8cm}
\begin{align*}
C_0 =& \; [5, [8, -1], 13, [1, 0], 14, [4, 2], 10, [17, 0], 17, [17, 2], 0, [3, 0], 3, [8, 2], 11, [7, 1], 4, [3, 1], 1, [2, 0], 16, \\ & \; [1, -1], 15, [7, 2]] \\
C_1 =& \; [9, [7, -1], 2, [5, 2], 7, [5, 0], 12, [6, -1], 6, [2, 1], 8, [1, 1]]
\end{align*}
\end{minipage}

\vspace*{-4mm}
\begin{minipage}{8cm}
\begin{align*}
C_0' =& \; [1, [17, 1], 17, [17, -1], 0, [3, 2], 14, [2, -1], 16, [7, 0], 9, [1, 2], 10, [5, 1], 5, [8, 0], 13, [8, 1], 4, [2, 2], 6, \\ & \; [6, 0], 12, [6, 1]] \\
C_1' =& \; [3, [5, -1], 15, [4, 1], 2, [6, 2], 8, [3, -1], 11, [4, 0], 7, [4, -1]]
\end{align*}
\end{minipage}

\item 2-factor type $[11, 7]$: two starters \vspace*{-3mm} \\
\begin{minipage}{8cm}
\begin{align*}
C_0 =& \; [4, [17, 1], 17, [17, 2], 15, [3, -1], 1, [4, 0], 14, [7, 2], 7, [4, -1], 11, [8, 1], 3, [6, 0], 9, [7, -1], 2, [3, 2], 5, \\ & \; [1, 0]] \\
C_1 =& \; [16, [7, 0], 6, [6, 2], 0, [5, 1], 12, [2, 0], 10, [3, 1], 13, [5, 2], 8, [8, -1]]
\end{align*}
\end{minipage}

\vspace*{-4mm}
\begin{minipage}{8cm}
\begin{align*}
C_0' =& \; [5, [1, 2], 6, [2, 1], 4, [5, 0], 9, [1, -1], 8, [8, 2], 0, [2, -1], 2, [17, -1], 17, [17, 0], 10, [6, 1], 16, [5, -1], 11, \\ & \; [6, -1]] \\
C_1' =& \; [13, [2, 2], 15, [3, 0], 1, [4, 2], 14, [7, 1], 7, [4, 1], 3, [8, 0], 12, [1, 1]]
\end{align*}
\end{minipage}

\item 2-factor type $[10, 8]$: two starters \vspace*{-3mm} \\
\begin{minipage}{8cm}
\begin{align*}
C_0 =& \; [3, [5, 1], 8, [7, 2], 15, [6, 1], 9, [7, -1], 16, [2, 0], 1, [1, -1], 0, [4, 2], 4, [3, -1], 7, [1, 1], 6, [3, 0]] \\
C_1 =& \; [11, [6, -1], 5, [8, 2], 13, [17, 0], 17, [17, 2], 2, [5, 0], 14, [2, -1], 12, [2, 2], 10, [1, 0]]
\end{align*}
\end{minipage}

\vspace*{-4mm}
\begin{minipage}{8cm}
\begin{align*}
C_0' =& \; [7, [2, 1], 9, [17, 1], 17, [17, -1], 6, [8, -1], 15, [4, -1], 11, [5, 2], 16, [6, 0], 5, [5, -1], 0, [3, 2], \\ & \; 3, [4, 0]] \\
C_1' =& \; [14, [7, 0], 4, [8, 1], 12, [1, 2], 13, [3, 1], 10, [8, 0], 2, [6, 2], 8, [7, 1], 1, [4, 1]]
\end{align*}
\end{minipage}

\end{itemizenew}


\section{Computational results for $n=19$}\label{app:19}

\begin{itemizenew}
\item 2-factor type $[3, 2, 2, 2, 2, 2, 2, 2, 2]$: three starters \vspace*{-3mm} \\
\hspace*{-14mm}
\begin{minipage}{8cm}
\begin{align*}
C_0 =& \;  [12, [1, 0], 11, [4, 2], 15, [3, 1]] \\
C_1 =& \;  [5, [4, -1], 1, [4, 1]] \\
C_2 =& \;  [3, [8, 0], 13, [8, 2]] \\
C_3 =& \;  [7, [7, 0], 0, [7, 2]] \\
C_4 =& \;  [4, [5, -1], 17, [5, 1]] \\
C_5 =& \;  [8, [2, 0], 6, [2, 2]] \\
C_6 =& \;  [2, [18, -1], 18, [18, 1]] \\
C_7 =& \;  [10, [6, 0], 16, [6, 2]] \\
C_8 =& \;  [14, [5, 0], 9, [5, 2]]
\end{align*}
\end{minipage}
\begin{minipage}{8cm}
\begin{align*}
C_0' =& \;  [6, [3, -1], 3, [1, 2], 2, [4, 0]] \\
C_1' =& \;  [12, [8, -1], 4, [8, 1]] \\
C_2' =& \;  [16, [1, -1], 15, [1, 1]] \\
C_3' =& \;  [17, [7, -1], 10, [7, 1]] \\
C_4' =& \;  [14, [3, 0], 11, [3, 2]] \\
C_5' =& \;  [5, [2, -1], 7, [2, 1]] \\
C_6' =& \;  [13, [6, -1], 1, [6, 1]] \\
C_7' =& \;  [18, [18, 0], 8, [18, 2]] \\
C_8' =& \;  [0, [9, 0], 9, [9, 2]]
\end{align*}
\end{minipage}

\item 2-factor type $[5, 2, 2, 2, 2, 2, 2, 2]$: three starters \vspace*{-3mm} \\
\begin{minipage}{8cm}
\begin{align*}
C_0 =& \;  [7, [5, 2], 2, [8, -1], 10, [18, -1], 18, [18, 1], 6, [1, 0]] \\
C_1 =& \;  [4, [6, -1], 16, [6, 1]] \\
C_2 =& \;  [15, [8, 0], 5, [8, 2]] \\
C_3 =& \;  [0, [5, -1], 13, [5, 1]] \\
C_4 =& \;  [11, [3, 0], 8, [3, 2]] \\
C_5 =& \;  [17, [2, 0], 1, [2, 2]] \\
C_6 =& \;  [12, [3, -1], 9, [3, 1]] \\
C_7 =& \;  [14, [7, -1], 3, [7, 1]]
\end{align*}
\end{minipage}

\vspace*{-4mm}
\begin{minipage}{8cm}
\begin{align*}
C_0' =& \;  [7, [7, 2], 14, [7, 0], 3, [1, 2], 2, [8, 1], 12, [5, 0]] \\
C_1' =& \;  [0, [9, 0], 9, [9, 2]] \\
C_2' =& \;  [4, [6, 0], 10, [6, 2]] \\
C_3' =& \;  [17, [1, -1], 16, [1, 1]] \\
C_4' =& \;  [13, [18, 0], 18, [18, 2]] \\
C_5' =& \;  [8, [2, -1], 6, [2, 1]] \\
C_6' =& \;  [11, [4, 0], 15, [4, 2]] \\
C_7' =& \;  [5, [4, -1], 1, [4, 1]]
\end{align*}
\end{minipage}

\item 2-factor type $[4, 3, 2, 2, 2, 2, 2, 2]$: three starters \vspace*{-3mm} \\
\hspace*{-6mm}
\begin{minipage}{8cm}
\begin{align*}
C_0 =& \;  [10, [2, -1], 12, [6, 1], 6, [2, 0], 4, [6, 2]] \\
C_1 =& \;  [0, [1, 2], 17, [1, 0], 16, [2, 1]] \\
C_2 =& \;  [14, [3, -1], 11, [3, 1]] \\
C_3 =& \;  [13, [5, 0], 8, [5, 2]] \\
C_4 =& \;  [9, [18, -1], 18, [18, 1]] \\
C_5 =& \;  [2, [3, 0], 5, [3, 2]] \\
C_6 =& \;  [3, [4, -1], 7, [4, 1]] \\
C_7 =& \;  [1, [4, 0], 15, [4, 2]]
\end{align*}
\end{minipage}
\begin{minipage}{8cm}
\begin{align*}
C_0' =& \;  [3, [2, 2], 1, [7, -1], 8, [6, 0], 2, [1, 1]] \\
C_1' =& \;  [11, [6, -1], 17, [7, 1], 10, [1, -1]] \\
C_2' =& \;  [0, [9, 0], 9, [9, 2]] \\
C_3' =& \;  [13, [8, -1], 5, [8, 1]] \\
C_4' =& \;  [18, [18, 0], 14, [18, 2]] \\
C_5' =& \;  [6, [8, 0], 16, [8, 2]] \\
C_6' =& \;  [12, [5, -1], 7, [5, 1]] \\
C_7' =& \;  [15, [7, 0], 4, [7, 2]]
\end{align*}
\end{minipage}

\newpage

\item 2-factor type $[3, 3, 3, 2, 2, 2, 2, 2]$: three starters \vspace*{-3mm} \\
\hspace*{-14mm}
\begin{minipage}{8cm}
\begin{align*}
C_0 =& \;  [17, [3, -1], 2, [1, 1], 1, [2, 1]] \\
C_1 =& \;  [4, [4, -1], 0, [3, 2], 3, [1, 0]] \\
C_2 =& \;  [13, [7, -1], 6, [3, 1], 9, [4, 1]] \\
C_3 =& \;  [10, [2, 0], 12, [2, 2]] \\
C_4 =& \;  [18, [18, 0], 5, [18, 2]] \\
C_5 =& \;  [11, [4, 0], 15, [4, 2]] \\
C_6 =& \;  [16, [8, 0], 8, [8, 2]] \\
C_7 =& \;  [14, [7, 0], 7, [7, 2]]
\end{align*}
\end{minipage}
\begin{minipage}{8cm}
\begin{align*}
C_0' =& \;  [1, [7, 1], 12, [1, -1], 13, [6, -1]] \\
C_1' =& \;  [15, [1, 2], 16, [6, 1], 10, [5, 0]] \\
C_2' =& \;  [2, [2, -1], 4, [3, 0], 7, [5, 2]] \\
C_3' =& \;  [0, [9, 0], 9, [9, 2]] \\
C_4' =& \;  [17, [6, 0], 11, [6, 2]] \\
C_5' =& \;  [5, [18, -1], 18, [18, 1]] \\
C_6' =& \;  [3, [5, -1], 8, [5, 1]] \\
C_7' =& \;  [6, [8, -1], 14, [8, 1]]
\end{align*}
\end{minipage}

\item 2-factor type $[7, 2, 2, 2, 2, 2, 2]$: three starters \vspace*{-3mm} \\
\begin{minipage}{8cm}
\begin{align*}
C_0 =& \;  [9, [1, 2], 8, [5, -1], 13, [1, 0], 14, [5, 1], 1, [3, 2], 4, [5, 0], 17, [8, -1]] \\
C_1 =& \;  [5, [2, -1], 3, [2, 1]] \\
C_2 =& \;  [6, [1, -1], 7, [1, 1]] \\
C_3 =& \;  [18, [18, 0], 11, [18, 2]] \\
C_4 =& \;  [16, [2, 0], 0, [2, 2]] \\
C_5 =& \;  [15, [3, -1], 12, [3, 1]] \\
C_6 =& \;  [10, [8, 0], 2, [8, 2]]
\end{align*}
\end{minipage}

\vspace*{-4mm}
\begin{minipage}{8cm}
\begin{align*}
C_0' =& \;  [16, [3, 0], 13, [8, 1], 3, [18, 1], 18, [18, -1], 11, [4, -1], 7, [5, 2], 2, [4, 1]] \\
C_1' =& \;  [0, [9, 0], 9, [9, 2]] \\
C_2' =& \;  [17, [7, 0], 10, [7, 2]] \\
C_3' =& \;  [6, [6, -1], 12, [6, 1]] \\
C_4' =& \;  [5, [4, 0], 1, [4, 2]] \\
C_5' =& \;  [8, [6, 0], 14, [6, 2]] \\
C_6' =& \;  [15, [7, -1], 4, [7, 1]]
\end{align*}
\end{minipage}

\item 2-factor type $[6, 3, 2, 2, 2, 2, 2]$: three starters \vspace*{-3mm} \\
\begin{minipage}{8cm}
\begin{align*}
C_0 =& \;  [13, [2, -1], 11, [3, -1], 8, [18, 2], 18, [18, 0], 0, [1, 1], 1, [6, -1]] \\
C_1 =& \;  [7, [7, 1], 14, [3, 1], 17, [8, 1]] \\
C_2 =& \;  [2, [1, 0], 3, [1, 2]] \\
C_3 =& \;  [16, [4, -1], 12, [4, 1]] \\
C_4 =& \;  [15, [6, 0], 9, [6, 2]] \\
C_5 =& \;  [4, [2, 0], 6, [2, 2]] \\
C_6 =& \;  [5, [5, -1], 10, [5, 1]]
\end{align*}
\end{minipage}

\vspace*{-4mm}
\begin{minipage}{8cm}
\begin{align*}
C_0' =& \;  [10, [7, -1], 17, [4, 0], 3, [1, -1], 2, [4, 2], 6, [18, -1], 18, [18, 1]] \\
C_1' =& \;  [7, [6, 1], 13, [2, 1], 15, [8, -1]] \\
C_2' =& \;  [0, [9, 0], 9, [9, 2]] \\
C_3' =& \;  [5, [7, 0], 16, [7, 2]] \\
C_4' =& \;  [12, [8, 0], 4, [8, 2]] \\
C_5' =& \;  [8, [3, 0], 11, [3, 2]] \\
C_6' =& \;  [1, [5, 0], 14, [5, 2]]
\end{align*}
\end{minipage}

\item 2-factor type $[5, 4, 2, 2, 2, 2, 2]$: three starters \vspace*{-3mm} \\
\begin{minipage}{8cm}
\begin{align*}
C_0 =& \;  [6, [6, 2], 12, [6, 0], 0, [1, 2], 1, [6, -1], 7, [1, 0]] \\
C_1 =& \;  [15, [7, 1], 8, [18, 0], 18, [18, -1], 4, [7, 2]] \\
C_2 =& \;  [11, [1, -1], 10, [1, 1]] \\
C_3 =& \;  [16, [5, 0], 3, [5, 2]] \\
C_4 =& \;  [14, [3, -1], 17, [3, 1]] \\
C_5 =& \;  [9, [4, -1], 13, [4, 1]] \\
C_6 =& \;  [5, [3, 0], 2, [3, 2]]
\end{align*}
\end{minipage}

\vspace*{-4mm}
\begin{minipage}{8cm}
\begin{align*}
C_0' =& \;  [17, [7, -1], 10, [8, 2], 2, [7, 0], 13, [2, -1], 11, [6, 1]] \\
C_1' =& \;  [14, [18, 1], 18, [18, 2], 4, [8, 0], 12, [2, 1]] \\
C_2' =& \;  [0, [9, 0], 9, [9, 2]] \\
C_3' =& \;  [8, [2, 0], 6, [2, 2]] \\
C_4' =& \;  [1, [4, 0], 5, [4, 2]] \\
C_5' =& \;  [7, [8, -1], 15, [8, 1]] \\
C_6' =& \;  [3, [5, -1], 16, [5, 1]]
\end{align*}
\end{minipage}

\item 2-factor type $[5, 3, 3, 2, 2, 2, 2]$: three starters \vspace*{-3mm} \\
\begin{minipage}{8cm}
\begin{align*}
C_0 =& \;  [12, [5, 2], 17, [6, -1], 5, [1, 0], 6, [3, 2], 9, [3, 0]] \\
C_1 =& \;  [15, [2, 0], 13, [18, 2], 18, [18, 1]] \\
C_2 =& \;  [0, [4, -1], 14, [6, 1], 2, [2, -1]] \\
C_3 =& \;  [16, [5, -1], 3, [5, 1]] \\
C_4 =& \;  [1, [7, 0], 8, [7, 2]] \\
C_5 =& \;  [4, [7, -1], 11, [7, 1]] \\
C_6 =& \;  [10, [3, -1], 7, [3, 1]]
\end{align*}
\end{minipage}

\vspace*{-4mm}
\begin{minipage}{8cm}
\begin{align*}
C_0' =& \;  [11, [5, 0], 6, [6, 2], 12, [4, 0], 8, [2, 1], 10, [1, 2]] \\
C_1' =& \;  [2, [4, 1], 16, [6, 0], 4, [2, 2]] \\
C_2' =& \;  [1, [4, 2], 5, [18, 0], 18, [18, -1]] \\
C_3' =& \;  [0, [9, 0], 9, [9, 2]] \\
C_4' =& \;  [17, [8, -1], 7, [8, 1]] \\
C_5' =& \;  [14, [1, -1], 15, [1, 1]] \\
C_6' =& \;  [3, [8, 0], 13, [8, 2]]
\end{align*}
\end{minipage}

\item 2-factor type $[4, 4, 3, 2, 2, 2, 2]$: three starters \vspace*{-3mm} \\
\hspace*{-4mm}
\begin{minipage}{8cm}
\begin{align*}
C_0 =& \;  [6, [1, 0], 5, [1, -1], 4, [8, -1], 14, [8, 2]] \\
C_1 =& \;  [12, [1, 2], 13, [3, -1], 16, [2, -1], 0, [6, 0]] \\
C_2 =& \;  [1, [8, 0], 11, [6, 1], 17, [2, 2]] \\
C_3 =& \;  [3, [7, -1], 10, [7, 1]] \\
C_4 =& \;  [9, [18, 0], 18, [18, 2]] \\
C_5 =& \;  [2, [5, 0], 7, [5, 2]] \\
C_6 =& \;  [8, [7, 0], 15, [7, 2]]
\end{align*}
\end{minipage}
\begin{minipage}{8cm}
\begin{align*}
C_0' =& \;  [1, [2, 1], 17, [4, 1], 13, [18, -1], 18, [18, 1]] \\
C_1' =& \;  [16, [2, 0], 14, [8, 1], 4, [6, 2], 10, [6, -1]] \\
C_2' =& \;  [2, [1, 1], 3, [3, 1], 6, [4, -1]] \\
C_3' =& \;  [0, [9, 0], 9, [9, 2]] \\
C_4' =& \;  [5, [3, 0], 8, [3, 2]] \\
C_5' =& \;  [12, [5, -1], 7, [5, 1]] \\
C_6' =& \;  [15, [4, 0], 11, [4, 2]
\end{align*}
\end{minipage}

\newpage

\item 2-factor type $[4, 3, 3, 3, 2, 2, 2]$: three starters \vspace*{-3mm} \\
\hspace*{-7mm}
\begin{minipage}{8cm}
\begin{align*}
C_0 =& \;  [9, [7, 0], 16, [4, 2], 12, [5, 1], 7, [2, -1]] \\
C_1 =& \;  [10, [8, 2], 0, [3, 1], 15, [5, 0]] \\
C_2 =& \;  [1, [2, 2], 3, [1, 1], 2, [1, 0]] \\
C_3 =& \;  [18, [18, 1], 13, [4, 0], 17, [18, 2]] \\
C_4 =& \;  [6, [8, -1], 14, [8, 1]] \\
C_5 =& \;  [11, [7, -1], 4, [7, 1]] \\
C_6 =& \;  [5, [3, 0], 8, [3, 2]]
\end{align*}
\end{minipage}
\begin{minipage}{8cm}
\begin{align*}
C_0' =& \;  [14, [1, 2], 15, [2, 0], 17, [2, 1], 1, [5, -1]] \\
C_1' =& \;  [10, [18, 0], 18, [18, -1], 4, [6, 2]] \\
C_2' =& \;  [5, [7, 2], 12, [6, 0], 6, [1, -1]] \\
C_3' =& \;  [11, [8, 0], 3, [5, 2], 8, [3, -1]] \\
C_4' =& \;  [0, [9, 0], 9, [9, 2]] \\
C_5' =& \;  [16, [4, -1], 2, [4, 1]] \\
C_6' =& \;  [7, [6, -1], 13, [6, 1]]
\end{align*}
\end{minipage}

\item 2-factor type $[3, 3, 3, 3, 3, 2, 2]$: three starters \vspace*{-3mm} \\
\hspace*{-10mm}
\begin{minipage}{8cm}
\begin{align*}
C_0 =& \;  [10, [5, 1], 5, [3, 1], 2, [8, -1]] \\
C_1 =& \;  [11, [7, 1], 4, [4, -1], 8, [3, -1]] \\
C_2 =& \;  [0, [7, 0], 7, [6, 2], 1, [1, 1]] \\
C_3 =& \;  [15, [1, 2], 14, [18, 0], 18, [18, -1]] \\
C_4 =& \;  [13, [4, 1], 9, [8, 0], 17, [4, 2]] \\
C_5 =& \;  [3, [5, 0], 16, [5, 2]] \\
C_6 =& \;  [6, [6, -1], 12, [6, 1]]
\end{align*}
\end{minipage}
\begin{minipage}{8cm}
\begin{align*}
C_0' =& \;  [8, [4, 0], 12, [7, 2], 1, [7, -1]] \\
C_1' =& \;  [10, [1, -1], 11, [2, 0], 13, [3, 2]] \\
C_2' =& \;  [15, [1, 0], 14, [18, 2], 18, [18, 1]] \\
C_3' =& \;  [3, [5, -1], 16, [8, 2], 6, [3, 0]] \\
C_4' =& \;  [5, [2, 2], 7, [8, 1], 17, [6, 0]] \\
C_5' =& \;  [0, [9, 0], 9, [9, 2]] \\
C_6' =& \;  [4, [2, -1], 2, [2, 1]]
\end{align*}
\end{minipage}

\item 2-factor type $[9, 2, 2, 2, 2, 2]$: three starters \vspace*{-3mm} \\
\begin{minipage}{8cm}
\begin{align*}
C_0 =& \;  [1, [6, -1], 7, [7, 0], 0, [2, 2], 16, [5, 0], 11, [5, -1], 6, [3, 2], 9, [7, 1], 2, [6, 0], 8, [7, 2]] \\
C_1 =& \;  [10, [4, 0], 14, [4, 2]] \\
C_2 =& \;  [4, [1, -1], 5, [1, 1]] \\
C_3 =& \;  [3, [8, 0], 13, [8, 2]] \\
C_4 =& \;  [12, [3, -1], 15, [3, 1]] \\
C_5 =& \;  [18, [18, -1], 17, [18, 1]]
\end{align*}
\end{minipage}

\vspace*{-4mm}
\begin{minipage}{8cm}
\begin{align*}
C_0' =& \;  [11, [3, 0], 14, [18, 2], 18, [18, 0], 6, [5, 2], 1, [6, 1], 13, [7, -1], 2, [2, 0], 4, [6, 2], 16, [5, 1]] \\
C_1' =& \;  [0, [9, 0], 9, [9, 2]] \\
C_2' =& \;  [8, [1, 0], 7, [1, 2]] \\
C_3' =& \;  [5, [8, -1], 15, [8, 1]] \\
C_4' =& \;  [10, [2, -1], 12, [2, 1]] \\
C_5' =& \;  [3, [4, -1], 17, [4, 1]]
\end{align*}
\end{minipage}

\item 2-factor type $[8, 3, 2, 2, 2, 2]$: three starters \vspace*{-3mm} \\
\begin{minipage}{8cm}
\begin{align*}
C_0 =& \;  [10, [1, 1], 11, [3, 2], 14, [4, 0], 0, [1, 2], 1, [3, -1], 4, [1, 0], 3, [4, 2], 17, [7, 0]] \\
C_1 =& \;  [9, [3, 1], 12, [4, 1], 16, [7, -1]] \\
C_2 =& \;  [8, [2, 0], 6, [2, 2]] \\
C_3 =& \;  [5, [8, 0], 13, [8, 2]] \\
C_4 =& \;  [7, [8, -1], 15, [8, 1]] \\
C_5 =& \;  [18, [18, -1], 2, [18, 1]]
\end{align*}
\end{minipage}

\newpage

\begin{minipage}{8cm}
\begin{align*}
C_0' =& \;  [8, [2, 1], 6, [3, 0], 3, [4, -1], 17, [1, -1], 16, [5, 2], 11, [7, 1], 4, [18, 0], 18, [18, 2]] \\
C_1' =& \;  [5, [2, -1], 7, [5, 0], 12, [7, 2]] \\
C_2' =& \;  [0, [9, 0], 9, [9, 2]] \\
C_3' =& \;  [10, [5, -1], 15, [5, 1]] \\
C_4' =& \;  [2, [6, 0], 14, [6, 2]] \\
C_5' =& \;  [1, [6, -1], 13, [6, 1]]
\end{align*}
\end{minipage}

\item 2-factor type $[7, 4, 2, 2, 2, 2]$: three starters \vspace*{-3mm} \\
\begin{minipage}{8cm}
\begin{align*}
C_0 =& \;  [9, [2, -1], 7, [7, 2], 14, [18, 0], 18, [18, -1], 17, [1, 1], 0, [2, 1], 2, [7, 1]] \\
C_1 =& \;  [12, [7, 0], 1, [8, 2], 11, [6, 1], 5, [7, -1]] \\
C_2 =& \;  [15, [2, 0], 13, [2, 2]] \\
C_3 =& \;  [6, [8, -1], 16, [8, 1]] \\
C_4 =& \;  [4, [6, 0], 10, [6, 2]] \\
C_5 =& \;  [3, [5, -1], 8, [5, 1]]
\end{align*}
\end{minipage}

\vspace*{-4mm}
\begin{minipage}{8cm}
\begin{align*}
C_0' =& \;  [2, [1, 0], 3, [1, 2], 4, [1, -1], 5, [8, 0], 15, [4, -1], 11, [18, 1], 18, [18, 2]] \\
C_1' =& \;  [12, [4, 1], 8, [5, 0], 13, [6, -1], 7, [5, 2]] \\
C_2' =& \;  [0, [9, 0], 9, [9, 2]] \\
C_3' =& \;  [10, [4, 0], 6, [4, 2]] \\
C_4' =& \;  [14, [3, -1], 17, [3, 1]] \\
C_5' =& \;  [16, [3, 0], 1, [3, 2]]
\end{align*}
\end{minipage}

\item 2-factor type $[6, 5, 2, 2, 2, 2]$: three starters \vspace*{-3mm} \\
\begin{minipage}{8cm}
\begin{align*}
C_0 =& \;  [11, [5, 0], 16, [2, 2], 0, [4, 1], 14, [3, -1], 17, [7, 0], 6, [5, 2]] \\
C_1 =& \;  [1, [8, 1], 9, [2, -1], 7, [3, 2], 10, [2, 0], 8, [7, -1]] \\
C_2 =& \;  [18, [18, -1], 2, [18, 1]] \\
C_3 =& \;  [4, [1, 0], 3, [1, 2]] \\
C_4 =& \;  [15, [8, 0], 5, [8, 2]] \\
C_5 =& \;  [12, [1, -1], 13, [1, 1]]
\end{align*}
\end{minipage}

\vspace*{-4mm}
\begin{minipage}{8cm}
\begin{align*}
C_0' =& \;  [7, [5, 1], 12, [18, 2], 18, [18, 0], 10, [3, 1], 13, [2, 1], 15, [8, -1]] \\
C_1' =& \;  [3, [4, -1], 17, [7, 2], 6, [5, -1], 11, [3, 0], 14, [7, 1]] \\
C_2' =& \;  [0, [9, 0], 9, [9, 2]] \\
C_3' =& \;  [16, [6, -1], 4, [6, 1]] \\
C_4' =& \;  [1, [4, 0], 5, [4, 2]] \\
C_5' =& \;  [8, [6, 0], 2, [6, 2]]
\end{align*}
\end{minipage}

\item 2-factor type $[7, 3, 3, 2, 2, 2]$: three starters \vspace*{-3mm} \\
\begin{minipage}{8cm}
\begin{align*}
C_0 =& \;  [3, [8, 2], 11, [3, 0], 8, [18, 1], 18, [18, 2], 0, [8, 0], 10, [8, -1], 2, [1, 1]] \\
C_1 =& \;  [1, [3, -1], 16, [1, 2], 15, [4, 0]] \\
C_2 =& \;  [7, [6, 1], 13, [4, 1], 17, [8, 1]] \\
C_3 =& \;  [5, [7, -1], 12, [7, 1]] \\
C_4 =& \;  [14, [5, 0], 9, [5, 2]] \\
C_5 =& \;  [6, [2, -1], 4, [2, 1]]
\end{align*}
\end{minipage}

\newpage

\begin{minipage}{8cm}
\begin{align*}
C_0' =& \;  [15, [3, 2], 12, [1, 0], 11, [1, -1], 10, [6, -1], 4, [3, 1], 7, [6, 2], 13, [2, 0]] \\
C_1' =& \;  [14, [6, 0], 2, [4, -1], 16, [2, 2]] \\
C_2' =& \;  [18, [18, -1], 5, [4, 2], 1, [18, 0]] \\
C_3' =& \;  [0, [9, 0], 9, [9, 2]] \\
C_4' =& \;  [8, [5, -1], 3, [5, 1]] \\
C_5' =& \;  [17, [7, 0], 6, [7, 2]]
\end{align*}
\end{minipage}

\item 2-factor type $[6, 4, 3, 2, 2, 2]$: three starters \vspace*{-3mm} \\
\begin{minipage}{8cm}
\begin{align*}
C_0 =& \;  [13, [7, 0], 6, [1, 1], 7, [5, -1], 2, [1, 2], 3, [4, 0], 17, [4, 2]] \\
C_1 =& \;  [4, [8, -1], 14, [2, 2], 12, [2, 0], 10, [6, -1]] \\
C_2 =& \;  [11, [5, 2], 16, [7, -1], 5, [6, 0]] \\
C_3 =& \;  [18, [18, 0], 8, [18, 2]] \\
C_4 =& \;  [0, [3, 0], 15, [3, 2]] \\
C_5 =& \;  [1, [8, 0], 9, [8, 2]]
\end{align*}
\end{minipage}

\vspace*{-4mm}
\begin{minipage}{8cm}
\begin{align*}
C_0' =& \;  [11, [3, 1], 8, [1, 0], 7, [8, 1], 15, [18, 1], 18, [18, -1], 5, [6, 2]] \\
C_1' =& \;  [17, [3, -1], 14, [7, 2], 3, [1, -1], 4, [5, 0]] \\
C_2' =& \;  [6, [7, 1], 13, [6, 1], 1, [5, 1]] \\
C_3' =& \;  [0, [9, 0], 9, [9, 2]] \\
C_4' =& \;  [2, [4, -1], 16, [4, 1]] \\
C_5' =& \;  [12, [2, -1], 10, [2, 1]]
\end{align*}
\end{minipage}

\item 2-factor type $[5, 5, 3, 2, 2, 2]$: three starters \vspace*{-3mm} \\
\begin{minipage}{8cm}
\begin{align*}
C_0 =& \;  [12, [8, -1], 2, [2, 1], 0, [7, -1], 7, [3, -1], 10, [2, -1]] \\
C_1 =& \;  [13, [7, 2], 6, [2, 0], 8, [8, 1], 16, [1, 2], 17, [4, 0]] \\
C_2 =& \;  [3, [7, 1], 14, [5, 0], 1, [2, 2]] \\
C_3 =& \;  [5, [6, 0], 11, [6, 2]] \\
C_4 =& \;  [4, [5, -1], 9, [5, 1]] \\
C_5 =& \;  [18, [18, 0], 15, [18, 2]]
\end{align*}
\end{minipage}

\vspace*{-4mm}
\begin{minipage}{8cm}
\begin{align*}
C_0' =& \;  [12, [7, 0], 1, [4, 2], 15, [1, 0], 14, [3, 2], 11, [1, -1]] \\
C_1' =& \;  [2, [1, 1], 3, [18, 1], 18, [18, -1], 7, [8, 2], 17, [3, 0]] \\
C_2' =& \;  [8, [5, 2], 13, [8, 0], 5, [3, 1]] \\
C_3' =& \;  [0, [9, 0], 9, [9, 2]] \\
C_4' =& \;  [6, [4, -1], 10, [4, 1]] \\
C_5' =& \;  [16, [6, -1], 4, [6, 1]]
\end{align*}
\end{minipage}

\item 2-factor type $[5, 4, 4, 2, 2, 2]$: three starters \vspace*{-3mm} \\
\begin{minipage}{8cm}
\begin{align*}
C_0 =& \;  [15, [18, 1], 18, [18, -1], 8, [1, 2], 7, [5, 0], 2, [5, -1]] \\
C_1 =& \;  [14, [5, 2], 9, [4, 1], 5, [1, 0], 6, [8, 1]] \\
C_2 =& \;  [1, [2, -1], 3, [8, 0], 13, [3, 1], 16, [3, 2]] \\
C_3 =& \;  [4, [4, 0], 0, [4, 2]] \\
C_4 =& \;  [10, [7, 0], 17, [7, 2]] \\
C_5 =& \;  [11, [1, -1], 12, [1, 1]]
\end{align*}
\end{minipage}

\newpage

\begin{minipage}{8cm}
\begin{align*}
C_0' =& \;  [7, [4, -1], 3, [7, 1], 10, [8, -1], 2, [18, 2], 18, [18, 0]] \\
C_1' =& \;  [15, [3, 0], 12, [7, -1], 5, [6, 2], 17, [2, 1]] \\
C_2' =& \;  [4, [3, -1], 1, [5, 1], 6, [8, 2], 16, [6, 0]] \\
C_3' =& \;  [0, [9, 0], 9, [9, 2]] \\
C_4' =& \;  [14, [6, -1], 8, [6, 1]] \\
C_5' =& \;  [13, [2, 0], 11, [2, 2]]
\end{align*}
\end{minipage}

\item 2-factor type $[6, 3, 3, 3, 2, 2]$: three starters \vspace*{-3mm} \\
\begin{minipage}{8cm}
\begin{align*}
C_0 =& \;  [1, [7, 1], 12, [3, 0], 9, [7, 2], 16, [1, 0], 17, [2, -1], 15, [4, 2]] \\
C_1 =& \;  [13, [1, -1], 14, [4, 1], 10, [3, -1]] \\
C_2 =& \;  [7, [3, 2], 4, [18, 0], 18, [18, -1]] \\
C_3 =& \;  [0, [6, 2], 6, [1, 1], 5, [5, 0]] \\
C_4 =& \;  [2, [6, -1], 8, [6, 1]] \\
C_5 =& \;  [11, [8, -1], 3, [8, 1]]
\end{align*}
\end{minipage}

\vspace*{-4mm}
\begin{minipage}{8cm}
\begin{align*}
C_0' =& \;  [10, [5, 2], 15, [4, 0], 11, [5, 1], 16, [4, -1], 12, [8, 2], 4, [6, 0]] \\
C_1' =& \;  [7, [7, -1], 14, [8, 0], 6, [1, 2]] \\
C_2' =& \;  [8, [5, -1], 3, [2, 1], 5, [3, 1]] \\
C_3' =& \;  [13, [18, 1], 18, [18, 2], 2, [7, 0]] \\
C_4' =& \;  [0, [9, 0], 9, [9, 2]] \\
C_5' =& \;  [1, [2, 0], 17, [2, 2]]
\end{align*}
\end{minipage}

\item 2-factor type $[5, 4, 3, 3, 2, 2]$: three starters \vspace*{-3mm} \\
\begin{minipage}{8cm}
\begin{align*}
C_0 =& \;  [17, [8, 2], 7, [7, 0], 0, [4, 2], 14, [4, 0], 10, [7, 1]] \\
C_1 =& \;  [9, [4, -1], 5, [18, 1], 18, [18, 2], 12, [3, 0]] \\
C_2 =& \;  [15, [4, 1], 11, [8, 1], 3, [6, 1]] \\
C_3 =& \;  [1, [7, -1], 8, [2, 1], 6, [5, 1]] \\
C_4 =& \;  [16, [3, -1], 13, [3, 1]] \\
C_5 =& \;  [2, [2, 0], 4, [2, 2]]
\end{align*}
\end{minipage}

\vspace*{-4mm}
\begin{minipage}{8cm}
\begin{align*}
C_0' =& \;  [14, [18, -1], 18, [18, 0], 17, [1, 2], 16, [5, 0], 11, [3, 2]] \\
C_1' =& \;  [2, [1, 0], 3, [6, -1], 15, [8, -1], 7, [5, 2]] \\
C_2' =& \;  [1, [7, 2], 8, [5, -1], 13, [6, 0]] \\
C_3' =& \;  [10, [2, -1], 12, [8, 0], 4, [6, 2]] \\
C_4' =& \;  [0, [9, 0], 9, [9, 2]] \\
C_5' =& \;  [6, [1, -1], 5, [1, 1]]
\end{align*}
\end{minipage}

\item 2-factor type $[4, 4, 4, 3, 2, 2]$: three starters \vspace*{-3mm} \\
\hspace*{-3mm}
\begin{minipage}{8cm}
\begin{align*}
C_0 =& \;  [8, [2, 0], 10, [5, 1], 15, [8, 2], 5, [3, -1]] \\
C_1 =& \;  [6, [18, 0], 18, [18, -1], 1, [3, 2], 16, [8, -1]] \\
C_2 =& \;  [3, [1, -1], 2, [7, 1], 9, [2, 2], 11, [8, 0]] \\
C_3 =& \;  [4, [8, 1], 14, [3, 0], 17, [5, 2]] \\
C_4 =& \;  [7, [7, 0], 0, [7, 2]] \\
C_5 =& \;  [12, [1, 0], 13, [1, 2]]
\end{align*}
\end{minipage}
\begin{minipage}{8cm}
\begin{align*}
C_0' =& \;  [3, [5, -1], 16, [2, -1], 14, [1, 1], 15, [6, 1]] \\
C_1' =& \;  [2, [5, 0], 7, [6, -1], 1, [18, 2], 18, [18, 1]] \\
C_2' =& \;  [11, [7, -1], 4, [4, 2], 8, [3, 1], 5, [6, 0]] \\
C_3' =& \;  [6, [6, 2], 12, [2, 1], 10, [4, 0]] \\
C_4' =& \;  [0, [9, 0], 9, [9, 2]] \\
C_5' =& \;  [13, [4, -1], 17, [4, 1]]
\end{align*}
\end{minipage}

\item 2-factor type $[5, 3, 3, 3, 3, 2]$: three starters \vspace*{-3mm} \\
\begin{minipage}{8cm}
\begin{align*}
C_0 =& \;  [9, [5, 2], 14, [4, 0], 0, [3, 2], 15, [5, -1], 2, [7, 0]] \\
C_1 =& \;  [4, [1, 1], 5, [2, 2], 7, [3, 0]] \\
C_2 =& \;  [6, [3, 1], 3, [4, 1], 17, [7, -1]] \\
C_3 =& \;  [12, [4, -1], 8, [5, 1], 13, [1, -1]] \\
C_4 =& \;  [18, [18, 0], 10, [6, 1], 16, [18, 2]] \\
C_5 =& \;  [1, [8, 0], 11, [8, 2]]
\end{align*}
\end{minipage}

\vspace*{-4mm}
\begin{minipage}{8cm}
\begin{align*}
C_0' =& \;  [8, [3, -1], 11, [1, 0], 12, [1, 2], 13, [2, -1], 15, [7, 1]] \\
C_1' =& \;  [6, [4, 2], 2, [6, 0], 14, [8, -1]] \\
C_2' =& \;  [10, [6, -1], 4, [18, 1], 18, [18, -1]] \\
C_3' =& \;  [1, [6, 2], 7, [8, 1], 17, [2, 0]] \\
C_4' =& \;  [3, [2, 1], 5, [7, 2], 16, [5, 0]] \\
C_5' =& \;  [0, [9, 0], 9, [9, 2]]
\end{align*}
\end{minipage}

\item 2-factor type $[4, 4, 3, 3, 3, 2]$: three starters \vspace*{-3mm} \\
\hspace*{-4mm}
\begin{minipage}{8cm}
\begin{align*}
C_0 =& \;  [11, [2, 0], 13, [1, -1], 12, [2, 2], 14, [3, 1]] \\
C_1 =& \;  [10, [7, 1], 3, [4, 0], 7, [2, 1], 9, [1, 2]] \\
C_2 =& \;  [4, [2, -1], 2, [18, 2], 18, [18, 0]] \\
C_3 =& \;  [15, [7, -1], 8, [8, 2], 0, [3, 0]] \\
C_4 =& \;  [16, [1, 1], 17, [6, 2], 5, [7, 0]] \\
C_5 =& \;  [6, [5, 0], 1, [5, 2]]
\end{align*}
\end{minipage}
\begin{minipage}{8cm}
\begin{align*}
C_0' =& \;  [10, [7, 2], 3, [6, 0], 15, [5, 1], 2, [8, 1]] \\
C_1' =& \;  [14, [3, -1], 11, [18, 1], 18, [18, -1], 1, [5, -1]] \\
C_2' =& \;  [6, [6, 1], 12, [4, 2], 16, [8, 0]] \\
C_3' =& \;  [17, [8, -1], 7, [6, -1], 13, [4, -1]] \\
C_4' =& \;  [8, [4, 1], 4, [1, 0], 5, [3, 2]] \\
C_5' =& \;  [0, [9, 0], 9, [9, 2]]
\end{align*}
\end{minipage}

\item 2-factor type $[11, 2, 2, 2, 2]$: three starters \vspace*{-3mm} \\
\begin{minipage}{8cm}
\begin{align*}
C_0 =& \;  [8, [4, 0], 12, [2, 1], 14, [1, 2], 13, [7, 0], 6, [7, 2], 17, [1, -1], 0, [3, 0], 15, [5, 2], 2, [1, 0], 3, [5, -1], 16, \\ & \; [8, 2]] \\
C_1 =& \;  [1, [8, -1], 9, [8, 1]] \\
C_2 =& \;  [18, [18, -1], 11, [18, 1]] \\
C_3 =& \;  [4, [6, 0], 10, [6, 2]] \\
C_4 =& \;  [7, [2, 0], 5, [2, 2]]
\end{align*}
\end{minipage}

\vspace*{-4mm}
\begin{minipage}{8cm}
\begin{align*}
C_0' =& \;  [17, [3, 2], 14, [1, 1], 13, [8, 0], 5, [6, 1], 11, [5, 1], 16, [4, 2], 2, [5, 0], 7, [6, -1], 1, [7, -1], 12, [2, -1], 10, \\ & \; [7, 1]] \\
C_1' =& \;  [0, [9, 0], 9, [9, 2]] \\
C_2' =& \;  [15, [18, 0], 18, [18, 2]] \\
C_3' =& \;  [8, [4, -1], 4, [4, 1]] \\
C_4' =& \;  [6, [3, -1], 3, [3, 1]]
\end{align*}
\end{minipage}

\item 2-factor type $[10, 3, 2, 2, 2]$: three starters \vspace*{-3mm} \\
\begin{minipage}{8cm}
\begin{align*}
C_0 =& \;  [2, [4, 2], 16, [3, -1], 1, [18, 0], 18, [18, 2], 13, [6, 1], 7, [4, 0], 3, [3, 2], 6, [8, 0], 14, [6, 2], 8, [6, 0]] \\
C_1 =& \;  [11, [7, 1], 0, [6, -1], 12, [1, -1]] \\
C_2 =& \;  [15, [8, -1], 5, [8, 1]] \\
C_3 =& \;  [4, [5, 0], 9, [5, 2]] \\
C_4 =& \;  [10, [7, 0], 17, [7, 2]]
\end{align*}
\end{minipage}

\newpage

\begin{minipage}{8cm}
\begin{align*}
C_0' =& \;  [14, [2, 0], 16, [1, 1], 17, [4, 1], 3, [3, 1], 6, [5, -1], 1, [7, -1], 12, [8, 2], 4, [4, -1], 8, [18, -1], 18, [18, 1]] \\
C_1' =& \;  [2, [3, 0], 5, [2, 2], 7, [5, 1]] \\
C_2' =& \;  [0, [9, 0], 9, [9, 2]] \\
C_3' =& \;  [13, [2, -1], 15, [2, 1]] \\
C_4' =& \;  [11, [1, 0], 10, [1, 2]]
\end{align*}
\end{minipage}

\item 2-factor type $[9, 4, 2, 2, 2]$: three starters \vspace*{-3mm} \\
\begin{minipage}{8cm}
\begin{align*}
C_0 =& \;  [13, [3, 2], 16, [5, 0], 11, [4, -1], 7, [3, 1], 10, [5, 2], 15, [4, 0], 1, [5, -1], 14, [6, 2], 2, [7, 0]] \\
C_1 =& \;  [17, [5, 1], 12, [3, 0], 9, [18, 2], 18, [18, 1]] \\
C_2 =& \;  [0, [8, -1], 8, [8, 1]] \\
C_3 =& \;  [6, [1, -1], 5, [1, 1]] \\
C_4 =& \;  [4, [1, 0], 3, [1, 2]]
\end{align*}
\end{minipage}

\vspace*{-4mm}
\begin{minipage}{8cm}
\begin{align*}
C_0' =& \;  [11, [6, 0], 5, [2, 2], 3, [2, 0], 1, [4, 2], 15, [18, 0], 18, [18, -1], 7, [7, 1], 14, [8, 2], 4, [7, -1]] \\
C_1' =& \;  [17, [7, 2], 6, [8, 0], 16, [3, -1], 13, [4, 1]] \\
C_2' =& \;  [0, [9, 0], 9, [9, 2]] \\
C_3' =& \;  [8, [6, -1], 2, [6, 1]] \\
C_4' =& \;  [10, [2, -1], 12, [2, 1]]
\end{align*}
\end{minipage}

\item 2-factor type $[8, 5, 2, 2, 2]$: three starters \vspace*{-3mm} \\
\begin{minipage}{8cm}
\begin{align*}
C_0 =& \;  [7, [3, 0], 4, [4, -1], 0, [2, -1], 16, [3, 2], 1, [5, 0], 14, [6, 2], 8, [1, 0], 9, [2, 2]] \\
C_1 =& \;  [6, [18, -1], 18, [18, 0], 15, [2, 1], 17, [3, 1], 2, [4, 2]] \\
C_2 =& \;  [12, [1, -1], 13, [1, 1]] \\
C_3 =& \;  [5, [5, -1], 10, [5, 1]] \\
C_4 =& \;  [3, [8, -1], 11, [8, 1]]
\end{align*}
\end{minipage}

\vspace*{-4mm}
\begin{minipage}{8cm}
\begin{align*}
C_0' =& \;  [11, [6, 0], 17, [1, 2], 16, [2, 0], 14, [18, 2], 18, [18, 1], 7, [6, -1], 13, [8, 0], 3, [8, 2]] \\
C_1' =& \;  [2, [6, 1], 8, [3, -1], 5, [5, 2], 10, [4, 1], 6, [4, 0]] \\
C_2' =& \;  [0, [9, 0], 9, [9, 2]] \\
C_3' =& \;  [1, [7, 0], 12, [7, 2]] \\
C_4' =& \;  [4, [7, -1], 15, [7, 1]]
\end{align*}
\end{minipage}

\item 2-factor type $[7, 6, 2, 2, 2]$: three starters \vspace*{-3mm} \\
\begin{minipage}{8cm}
\begin{align*}
C_0 =& \;  [6, [7, -1], 13, [5, 0], 0, [1, 2], 1, [2, 0], 17, [3, 2], 2, [3, -1], 5, [1, -1]] \\
C_1 =& \;  [12, [3, 0], 9, [5, 2], 14, [2, -1], 16, [6, 0], 10, [2, 2], 8, [4, -1]] \\
C_2 =& \;  [4, [18, 0], 18, [18, 2]] \\
C_3 =& \;  [3, [8, -1], 11, [8, 1]] \\
C_4 =& \;  [7, [8, 0], 15, [8, 2]]
\end{align*}
\end{minipage}

\vspace*{-4mm}
\begin{minipage}{8cm}
\begin{align*}
C_0' =& \;  [11, [7, 1], 4, [6, 1], 16, [1, 0], 15, [7, 2], 8, [1, 1], 7, [6, -1], 13, [2, 1]] \\
C_1' =& \;  [2, [18, -1], 18, [18, 1], 10, [4, 1], 6, [3, 1], 3, [7, 0], 14, [6, 2]] \\
C_2' =& \;  [0, [9, 0], 9, [9, 2]] \\
C_3' =& \;  [5, [4, 0], 1, [4, 2]] \\
C_4' =& \;  [12, [5, -1], 17, [5, 1]]
\end{align*}
\end{minipage}

\item 2-factor type $[9, 3, 3, 2, 2]$: three starters \vspace*{-3mm} \\
\begin{minipage}{8cm}
\begin{align*}
C_0 =& \;  [2, [2, -1], 4, [3, -1], 7, [6, 0], 1, [8, 2], 9, [3, 0], 6, [6, -1], 0, [2, 2], 16, [5, 0], 3, [1, 2]] \\
C_1 =& \;  [10, [2, 0], 12, [2, 1], 14, [4, 2]] \\
C_2 =& \;  [8, [3, 1], 11, [18, 2], 18, [18, 0]] \\
C_3 =& \;  [15, [8, -1], 5, [8, 1]] \\
C_4 =& \;  [17, [4, -1], 13, [4, 1]]
\end{align*}
\end{minipage}

\vspace*{-4mm}
\begin{minipage}{8cm}
\begin{align*}
C_0' =& \;  [14, [18, 1], 18, [18, -1], 1, [3, 2], 4, [7, 1], 15, [5, 1], 10, [1, 0], 11, [6, 2], 17, [7, -1], 6, [8, 0]] \\
C_1' =& \;  [16, [5, 2], 3, [1, 1], 2, [4, 0]] \\
C_2' =& \;  [13, [5, -1], 8, [1, -1], 7, [6, 1]] \\
C_3' =& \;  [0, [9, 0], 9, [9, 2]] \\
C_4' =& \;  [12, [7, 0], 5, [7, 2]]
\end{align*}
\end{minipage}

\item 2-factor type $[8, 4, 3, 2, 2]$: three starters \vspace*{-3mm} \\
\begin{minipage}{8cm}
\begin{align*}
C_0 =& \;  [17, [1, 1], 0, [7, -1], 11, [3, 1], 14, [6, 1], 2, [4, 1], 6, [4, 2], 10, [2, 0], 12, [5, 1]] \\
C_1 =& \;  [4, [4, 0], 8, [18, 1], 18, [18, -1], 3, [1, 2]] \\
C_2 =& \;  [13, [6, 0], 1, [6, 2], 7, [6, -1]] \\
C_3 =& \;  [9, [7, 0], 16, [7, 2]] \\
C_4 =& \;  [15, [8, -1], 5, [8, 1]]
\end{align*}
\end{minipage}

\vspace*{-4mm}
\begin{minipage}{8cm}
\begin{align*}
C_0' =& \;  [16, [1, -1], 17, [4, -1], 3, [18, 0], 18, [18, 2], 4, [1, 0], 5, [7, 1], 12, [2, 1], 14, [2, 2]] \\
C_1' =& \;  [6, [5, 2], 1, [8, 0], 11, [3, -1], 8, [2, -1]] \\
C_2' =& \;  [7, [8, 2], 15, [5, -1], 2, [5, 0]] \\
C_3' =& \;  [0, [9, 0], 9, [9, 2]] \\
C_4' =& \;  [13, [3, 0], 10, [3, 2]]
\end{align*}
\end{minipage}

\item 2-factor type $[7, 5, 3, 2, 2]$: three starters \vspace*{-3mm} \\
\begin{minipage}{8cm}
\begin{align*}
C_0 =& \;  [17, [4, 2], 3, [6, 1], 15, [3, 0], 0, [6, 2], 6, [1, 0], 7, [2, 2], 5, [6, 0]] \\
C_1 =& \;  [4, [18, 1], 18, [18, -1], 10, [1, 2], 9, [7, 1], 2, [2, 0]] \\
C_2 =& \;  [11, [3, -1], 14, [5, 0], 1, [8, 2]] \\
C_3 =& \;  [12, [4, -1], 16, [4, 1]] \\
C_4 =& \;  [13, [5, -1], 8, [5, 1]]
\end{align*}
\end{minipage}

\vspace*{-4mm}
\begin{minipage}{8cm}
\begin{align*}
C_0' =& \;  [4, [6, -1], 16, [1, -1], 15, [7, -1], 8, [5, 2], 13, [1, 1], 12, [8, 0], 2, [2, 1]] \\
C_1' =& \;  [7, [4, 0], 11, [8, 1], 1, [2, -1], 17, [7, 2], 10, [3, 1]] \\
C_2' =& \;  [3, [3, 2], 6, [8, -1], 14, [7, 0]] \\
C_3' =& \;  [0, [9, 0], 9, [9, 2]] \\
C_4' =& \;  [5, [18, 0], 18, [18, 2]]
\end{align*}
\end{minipage}

\item 2-factor type $[6, 6, 3, 2, 2]$: three starters \vspace*{-3mm} \\
\begin{minipage}{8cm}
\begin{align*}
C_0 =& \;  [9, [4, 0], 5, [6, 1], 11, [7, 2], 4, [2, 1], 2, [5, 0], 7, [2, 2]] \\
C_1 =& \;  [16, [1, 0], 17, [2, -1], 15, [7, -1], 8, [8, -1], 0, [6, 2], 12, [4, -1]] \\
C_2 =& \;  [18, [18, -1], 10, [3, 2], 13, [18, 0]] \\
C_3 =& \;  [1, [5, -1], 14, [5, 1]] \\
C_4 =& \;  [6, [3, -1], 3, [3, 1]]
\end{align*}
\end{minipage}

\newpage

\begin{minipage}{8cm}
\begin{align*}
C_0' =& \;  [8, [18, 1], 18, [18, 2], 1, [8, 1], 11, [7, 0], 4, [1, 2], 5, [3, 0]] \\
C_1' =& \;  [3, [7, 1], 14, [6, 0], 2, [5, 2], 15, [2, 0], 17, [8, 2], 7, [4, 1]] \\
C_2' =& \;  [16, [8, 0], 6, [4, 2], 10, [6, -1]] \\
C_3' =& \;  [0, [9, 0], 9, [9, 2]] \\
C_4' =& \;  [13, [1, -1], 12, [1, 1]]
\end{align*}
\end{minipage}

\item 2-factor type $[7, 4, 4, 2, 2]$: three starters \vspace*{-3mm} \\
\begin{minipage}{8cm}
\begin{align*}
C_0 =& \;  [0, [3, 1], 3, [5, 2], 8, [3, -1], 11, [4, 0], 7, [5, 1], 12, [2, 1], 14, [4, 1]] \\
C_1 =& \;  [17, [2, 2], 15, [4, -1], 1, [1, -1], 2, [3, 0]] \\
C_2 =& \;  [5, [7, 2], 16, [6, 0], 4, [18, 2], 18, [18, 0]] \\
C_3 =& \;  [9, [1, 0], 10, [1, 2]] \\
C_4 =& \;  [13, [7, -1], 6, [7, 1]]
\end{align*}
\end{minipage}

\vspace*{-4mm}
\begin{minipage}{8cm}
\begin{align*}
C_0' =& \;  [14, [3, 2], 17, [6, -1], 5, [8, 0], 13, [6, 2], 7, [5, 0], 2, [8, 2], 12, [2, 0]] \\
C_1' =& \;  [11, [4, 2], 15, [7, 0], 4, [6, 1], 10, [1, 1]] \\
C_2' =& \;  [18, [18, -1], 8, [5, -1], 3, [2, -1], 1, [18, 1]] \\
C_3' =& \;  [0, [9, 0], 9, [9, 2]] \\
C_4' =& \;  [6, [8, -1], 16, [8, 1]]
\end{align*}
\end{minipage}

\item 2-factor type $[6, 5, 4, 2, 2]$: three starters \vspace*{-3mm} \\
\begin{minipage}{8cm}
\begin{align*}
C_0 =& \;  [8, [1, 1], 7, [3, -1], 10, [8, 0], 0, [3, 1], 3, [2, -1], 1, [7, 2]] \\
C_1 =& \;  [9, [7, 0], 2, [7, -1], 13, [3, 2], 16, [5, 0], 11, [2, 2]] \\
C_2 =& \;  [6, [18, 0], 18, [18, 2], 15, [3, 0], 12, [6, 2]] \\
C_3 =& \;  [5, [6, -1], 17, [6, 1]] \\
C_4 =& \;  [14, [8, -1], 4, [8, 1]]
\end{align*}
\end{minipage}

\vspace*{-4mm}
\begin{minipage}{8cm}
\begin{align*}
C_0' =& \;  [15, [7, 1], 4, [5, 2], 17, [1, 0], 16, [4, 2], 12, [1, -1], 13, [2, 0]] \\
C_1' =& \;  [3, [5, -1], 8, [6, 0], 2, [18, 1], 18, [18, -1], 11, [8, 2]] \\
C_2' =& \;  [7, [2, 1], 5, [4, 0], 1, [5, 1], 6, [1, 2]] \\
C_3' =& \;  [0, [9, 0], 9, [9, 2]] \\
C_4' =& \;  [14, [4, -1], 10, [4, 1]]
\end{align*}
\end{minipage}

\item 2-factor type $[5, 5, 5, 2, 2]$: three starters \vspace*{-3mm} \\
\begin{minipage}{8cm}
\begin{align*}
C_0 =& \;  [3, [4, -1], 17, [8, 1], 7, [5, 1], 12, [8, 2], 2, [1, 0]] \\
C_1 =& \;  [6, [5, -1], 1, [3, 2], 4, [8, 0], 14, [18, 1], 18, [18, -1]] \\
C_2 =& \;  [10, [2, 1], 8, [8, -1], 16, [7, 0], 5, [6, 2], 11, [1, 1]] \\
C_3 =& \;  [9, [6, -1], 15, [6, 1]] \\
C_4 =& \;  [0, [5, 0], 13, [5, 2]]
\end{align*}
\end{minipage}

\vspace*{-4mm}
\begin{minipage}{8cm}
\begin{align*}
C_0' =& \;  [15, [18, 0], 18, [18, 2], 17, [4, 0], 3, [2, -1], 1, [4, 2]] \\
C_1' =& \;  [16, [3, 1], 13, [7, -1], 2, [6, 0], 8, [3, -1], 5, [7, 2]] \\
C_2' =& \;  [11, [1, 2], 10, [4, 1], 6, [1, -1], 7, [3, 0], 4, [7, 1]] \\
C_3' =& \;  [0, [9, 0], 9, [9, 2]] \\
C_4' =& \;  [12, [2, 0], 14, [2, 2]]
\end{align*}
\end{minipage}

\item 2-factor type $[8, 3, 3, 3, 2]$: three starters \vspace*{-3mm} \\
\begin{minipage}{8cm}
\begin{align*}
C_0 =& \;  [14, [6, 1], 8, [5, 0], 13, [7, 2], 6, [3, 0], 9, [2, -1], 7, [3, 2], 10, [1, -1], 11, [3, -1]] \\
C_1 =& \;  [0, [6, -1], 12, [5, 2], 17, [1, 0]] \\
C_2 =& \;  [3, [1, 1], 4, [18, 1], 18, [18, -1]] \\
C_3 =& \;  [1, [1, 2], 2, [4, 0], 16, [3, 1]] \\
C_4 =& \;  [15, [8, 0], 5, [8, 2]]
\end{align*}
\end{minipage}

\vspace*{-4mm}
\begin{minipage}{8cm}
\begin{align*}
C_0' =& \;  [5, [7, 0], 12, [5, -1], 7, [18, 2], 18, [18, 0], 6, [8, -1], 16, [4, 2], 2, [7, 1], 13, [8, 1]] \\
C_1' =& \;  [10, [6, 2], 4, [4, -1], 8, [2, 0]] \\
C_2' =& \;  [15, [4, 1], 11, [6, 0], 17, [2, 2]] \\
C_3' =& \;  [1, [2, 1], 3, [7, -1], 14, [5, 1]] \\
C_4' =& \;  [0, [9, 0], 9, [9, 2]]
\end{align*}
\end{minipage}

\item 2-factor type $[7, 4, 3, 3, 2]$: three starters \vspace*{-3mm} \\
\begin{minipage}{8cm}
\begin{align*}
C_0 =& \;  [2, [4, 1], 16, [2, 0], 0, [6, -1], 12, [3, -1], 9, [4, -1], 5, [3, 1], 8, [6, 2]] \\
C_1 =& \;  [3, [4, 0], 17, [4, 2], 13, [2, 1], 11, [8, 1]] \\
C_2 =& \;  [6, [5, 2], 1, [5, 0], 14, [8, -1]] \\
C_3 =& \;  [7, [8, 2], 15, [7, -1], 4, [3, 0]] \\
C_4 =& \;  [18, [18, 0], 10, [18, 2]]
\end{align*}
\end{minipage}

\vspace*{-4mm}
\begin{minipage}{8cm}
\begin{align*}
C_0' =& \;  [13, [6, 0], 1, [7, 1], 8, [6, 1], 14, [7, 2], 3, [1, -1], 4, [7, 0], 11, [2, 2]] \\
C_1' =& \;  [12, [5, 1], 7, [1, 1], 6, [18, -1], 18, [18, 1]] \\
C_2' =& \;  [17, [1, 0], 16, [1, 2], 15, [2, -1]] \\
C_3' =& \;  [2, [8, 0], 10, [5, -1], 5, [3, 2]] \\
C_4' =& \;  [0, [9, 0], 9, [9, 2]]
\end{align*}
\end{minipage}

\item 2-factor type $[6, 5, 3, 3, 2]$: three starters \vspace*{-3mm} \\
\begin{minipage}{8cm}
\begin{align*}
C_0 =& \;  [3, [3, 0], 0, [5, 2], 5, [6, 0], 17, [2, -1], 15, [8, -1], 7, [4, 2]] \\
C_1 =& \;  [13, [7, 0], 6, [5, -1], 1, [3, 2], 16, [2, 1], 14, [1, 1]] \\
C_2 =& \;  [12, [3, -1], 9, [1, 2], 8, [4, 0]] \\
C_3 =& \;  [10, [8, 2], 2, [18, 0], 18, [18, -1]] \\
C_4 =& \;  [4, [7, -1], 11, [7, 1]]
\end{align*}
\end{minipage}

\vspace*{-4mm}
\begin{minipage}{8cm}
\begin{align*}
C_0' =& \;  [15, [8, 1], 5, [6, 2], 11, [6, -1], 17, [3, 1], 14, [4, 1], 10, [5, 0]] \\
C_1' =& \;  [4, [2, 0], 6, [2, 2], 8, [8, 0], 16, [18, 2], 18, [18, 1]] \\
C_2' =& \;  [7, [5, 1], 2, [1, -1], 3, [4, -1]] \\
C_3' =& \;  [13, [1, 0], 12, [7, 2], 1, [6, 1]] \\
C_4' =& \;  [0, [9, 0], 9, [9, 2]]
\end{align*}
\end{minipage}

\item 2-factor type $[6, 4, 4, 3, 2]$: three starters \vspace*{-3mm} \\
\begin{minipage}{8cm}
\begin{align*}
C_0 =& \;  [1, [7, 1], 8, [4, 1], 12, [2, 2], 10, [8, 0], 2, [3, -1], 17, [2, 1]] \\
C_1 =& \;  [14, [7, 2], 7, [6, -1], 13, [18, -1], 18, [18, 0]] \\
C_2 =& \;  [0, [3, 1], 3, [6, 2], 15, [1, -1], 16, [2, 0]] \\
C_3 =& \;  [5, [6, 0], 11, [5, -1], 6, [1, 2]] \\
C_4 =& \;  [4, [5, 0], 9, [5, 2]]
\end{align*}
\end{minipage}

\newpage

\begin{minipage}{8cm}
\begin{align*}
C_0' =& \;  [3, [8, -1], 11, [6, 1], 5, [1, 0], 6, [8, 2], 16, [8, 1], 8, [5, 1]] \\
C_1' =& \;  [2, [2, -1], 4, [3, 0], 7, [18, 1], 18, [18, 2]] \\
C_2' =& \;  [13, [1, 1], 14, [4, 2], 10, [7, -1], 17, [4, 0]] \\
C_3' =& \;  [12, [7, 0], 1, [4, -1], 15, [3, 2]] \\
C_4' =& \;  [0, [9, 0], 9, [9, 2]]
\end{align*}
\end{minipage}

\item 2-factor type $[5, 5, 4, 3, 2]$: three starters \vspace*{-3mm} \\
\begin{minipage}{8cm}
\begin{align*}
C_0 =& \;  [12, [2, 1], 10, [8, 0], 0, [4, -1], 14, [5, 2], 1, [7, 1]] \\
C_1 =& \;  [2, [7, 0], 9, [18, 1], 18, [18, -1], 4, [1, 2], 3, [1, 1]] \\
C_2 =& \;  [5, [7, 2], 16, [1, 0], 17, [6, 2], 11, [6, 0]] \\
C_3 =& \;  [7, [8, 1], 15, [2, -1], 13, [6, -1]] \\
C_4 =& \;  [6, [2, 0], 8, [2, 2]]
\end{align*}
\end{minipage}

\vspace*{-4mm}
\begin{minipage}{8cm}
\begin{align*}
C_0' =& \;  [3, [8, -1], 13, [8, 2], 5, [4, 0], 1, [3, -1], 16, [5, 1]] \\
C_1' =& \;  [15, [3, 0], 12, [1, -1], 11, [3, 2], 14, [6, 1], 8, [7, -1]] \\
C_2' =& \;  [7, [5, 0], 2, [4, 1], 6, [4, 2], 10, [3, 1]] \\
C_3' =& \;  [17, [18, 2], 18, [18, 0], 4, [5, -1]] \\
C_4' =& \;  [0, [9, 0], 9, [9, 2]]
\end{align*}
\end{minipage}

\item 2-factor type $[5, 4, 4, 4, 2]$: three starters \vspace*{-3mm} \\
\begin{minipage}{8cm}
\begin{align*}
C_0 =& \;  [6, [6, 2], 12, [3, 0], 15, [1, 1], 16, [18, 1], 18, [18, -1]] \\
C_1 =& \;  [8, [3, 1], 5, [6, 0], 17, [4, 2], 13, [5, 1]] \\
C_2 =& \;  [7, [7, 1], 14, [4, 1], 0, [3, 2], 3, [4, 0]] \\
C_3 =& \;  [10, [1, 2], 9, [8, 1], 1, [8, 0], 11, [1, -1]] \\
C_4 =& \;  [2, [2, 0], 4, [2, 2]]
\end{align*}
\end{minipage}

\vspace*{-4mm}
\begin{minipage}{8cm}
\begin{align*}
C_0' =& \;  [5, [5, -1], 10, [7, 0], 3, [2, -1], 1, [7, -1], 12, [7, 2]] \\
C_1' =& \;  [14, [3, -1], 17, [5, 0], 4, [2, 1], 6, [8, 2]] \\
C_2' =& \;  [15, [5, 2], 2, [6, -1], 8, [8, -1], 16, [1, 0]] \\
C_3' =& \;  [18, [18, 2], 13, [6, 1], 7, [4, -1], 11, [18, 0]] \\
C_4' =& \;  [0, [9, 0], 9, [9, 2]]
\end{align*}
\end{minipage}

\item 2-factor type $[13, 2, 2, 2]$: three starters \vspace*{-3mm} \\
\begin{minipage}{8cm}
\begin{align*}
C_0 =& \;  [6, [8, 2], 14, [3, 0], 11, [4, 2], 15, [7, 0], 8, [1, 2], 9, [8, 1], 1, [1, 0], 0, [3, 2], 3, [2, 0], 5, [7, 2], 12, [2, 1], 10, \\ & \; [8, 0], 2, [4, 1]] \\
C_1 =& \;  [18, [18, 0], 16, [18, 2]] \\
C_2 =& \;  [4, [5, 0], 17, [5, 2]] \\
C_3 =& \;  [7, [6, -1], 13, [6, 1]]
\end{align*}
\end{minipage}

\vspace*{-4mm}
\begin{minipage}{8cm}
\begin{align*}
C_0' =& \;  [16, [4, 0], 12, [1, -1], 11, [5, -1], 6, [2, 2], 4, [7, 1], 15, [4, -1], 1, [7, -1], 8, [5, 1], 3, [2, -1], 5, [8, -1], \\ & \; 13,  [18, -1], 18, [18, 1], 17, [1, 1]] \\
C_1' =& \;  [0, [9, 0], 9, [9, 2]] \\
C_2' =& \;  [2, [6, 0], 14, [6, 2]] \\
C_3' =& \;  [7, [3, -1], 10, [3, 1]]
\end{align*}
\end{minipage}

\item 2-factor type $[12, 3, 2, 2]$: three starters \vspace*{-3mm} \\
\begin{minipage}{8cm}
\begin{align*}
C_0 =& \;  [10, [3, 2], 7, [2, 0], 9, [4, 2], 5, [6, 1], 17, [3, 1], 14, [6, 0], 8, [5, 2], 13, [2, -1], 15, [1, 0], 16, [4, 1], 2, [2, 2], \\ & \; 0, [8, 0]] \\
C_1 =& \;  [12, [6, 2], 6, [5, 1], 1, [7, 0]] \\
C_2 =& \;  [4, [18, 0], 18, [18, 2]] \\
C_3 =& \;  [11, [8, -1], 3, [8, 1]]
\end{align*}
\end{minipage}

\vspace*{-4mm}
\begin{minipage}{8cm}
\begin{align*}
C_0' =& \;  [5, [8, 2], 15, [2, 1], 13, [3, 0], 10, [1, 2], 11, [6, -1], 17, [3, -1], 2, [5, 0], 7, [1, -1], 6, [5, -1], 1, [7, 2], 8, \\ & \; [4, 0], 4, [1, 1]] \\
C_1' =& \;  [12, [18, 1], 18, [18, -1], 16, [4, -1]] \\
C_2' =& \;  [0, [9, 0], 9, [9, 2]] \\
C_3' =& \;  [3, [7, -1], 14, [7, 1]]
\end{align*}
\end{minipage}

\item 2-factor type $[11, 4, 2, 2]$: three starters \vspace*{-3mm} \\
\begin{minipage}{8cm}
\begin{align*}
C_0 =& \;  [1, [1, 0], 0, [8, -1], 10, [4, 2], 14, [1, -1], 15, [4, 0], 11, [5, 2], 6, [2, 0], 4, [6, 2], 16, [4, 1], 12, [5, 1], 7, \\ & \; [6, 1]] \\
C_1 =& \;  [9, [4, -1], 13, [7, -1], 2, [3, 1], 17, [8, 1]] \\
C_2 =& \;  [5, [3, 0], 8, [3, 2]] \\
C_3 =& \;  [18, [18, 0], 3, [18, 2]]
\end{align*}
\end{minipage}

\vspace*{-4mm}
\begin{minipage}{8cm}
\begin{align*}
C_0' =& \;  [1, [1, 2], 2, [2, -1], 4, [5, 0], 17, [2, 2], 15, [7, 1], 8, [6, 0], 14, [3, -1], 11, [6, -1], 5, [7, 2], 16, [8, 0], 6, \\ & \; [5, -1]] \\
C_1' =& \;  [12, [2, 1], 10, [7, 0], 3, [8, 2], 13, [1, 1]] \\
C_2' =& \;  [0, [9, 0], 9, [9, 2]] \\
C_3' =& \;  [18, [18, -1], 7, [18, 1]]
\end{align*}
\end{minipage}

\item 2-factor type $[10, 5, 2, 2]$: three starters \vspace*{-3mm} \\
\begin{minipage}{8cm}
\begin{align*}
C_0 =& \;  [13, [4, 2], 9, [8, -1], 17, [4, -1], 3, [7, 0], 14, [6, 2], 2, [2, 0], 4, [3, 1], 7, [7, -1], 0, [8, 1], 8, [5, 1]] \\
C_1 =& \;  [18, [18, 2], 11, [5, 0], 6, [5, 2], 1, [7, 1], 12, [18, 0]] \\
C_2 =& \;  [10, [6, -1], 16, [6, 1]] \\
C_3 =& \;  [15, [8, 0], 5, [8, 2]]
\end{align*}
\end{minipage}

\vspace*{-4mm}
\begin{minipage}{8cm}
\begin{align*}
C_0' =& \;  [4, [7, 2], 15, [4, 1], 11, [2, -1], 13, [6, 0], 7, [18, 1], 18, [18, -1], 10, [2, 2], 12, [4, 0], 8, [5, -1], 3, [1, 1]] \\
C_1' =& \;  [1, [3, 2], 16, [2, 1], 14, [3, -1], 17, [3, 0], 2, [1, -1]] \\
C_2' =& \;  [0, [9, 0], 9, [9, 2]] \\
C_3' =& \;  [5, [1, 0], 6, [1, 2]]
\end{align*}
\end{minipage}

\item 2-factor type $[9, 6, 2, 2]$: three starters \vspace*{-3mm} \\
\begin{minipage}{8cm}
\begin{align*}
C_0 =& \;  [0, [5, 1], 5, [7, 2], 16, [3, -1], 1, [18, -1], 18, [18, 1], 7, [6, -1], 13, [4, 0], 17, [5, -1], 12, [6, 1]] \\
C_1 =& \;  [2, [5, 0], 15, [1, 2], 14, [3, 1], 11, [2, 0], 9, [1, 1], 10, [8, 2]] \\
C_2 =& \;  [6, [3, 0], 3, [3, 2]] \\
C_3 =& \;  [8, [4, -1], 4, [4, 1]]
\end{align*}
\end{minipage}

\newpage

\begin{minipage}{8cm}
\begin{align*}
C_0' =& \;  [14, [7, -1], 7, [6, 2], 13, [7, 0], 6, [2, -1], 4, [2, 2], 2, [8, 0], 12, [7, 1], 1, [4, 2], 15, [1, 0]] \\
C_1' =& \;  [10, [2, 1], 8, [18, 0], 18, [18, 2], 16, [1, -1], 17, [6, 0], 5, [5, 2]] \\
C_2' =& \;  [0, [9, 0], 9, [9, 2]] \\
C_3' =& \;  [11, [8, -1], 3, [8, 1]]
\end{align*}
\end{minipage}

\item 2-factor type $[8, 7, 2, 2]$: three starters \vspace*{-3mm} \\
\begin{minipage}{8cm}
\begin{align*}
C_0 =& \;  [3, [3, 1], 0, [5, 0], 13, [18, 1], 18, [18, -1], 14, [2, 2], 12, [8, 1], 4, [6, 0], 16, [5, 2]] \\
C_1 =& \;  [6, [5, -1], 11, [1, 0], 10, [7, 1], 17, [8, 2], 9, [2, 0], 7, [1, 2], 8, [2, 1]] \\
C_2 =& \;  [1, [4, 0], 15, [4, 2]] \\
C_3 =& \;  [2, [3, 0], 5, [3, 2]]
\end{align*}
\end{minipage}

\vspace*{-4mm}
\begin{minipage}{8cm}
\begin{align*}
C_0' =& \;  [4, [1, -1], 5, [4, 1], 1, [6, -1], 7, [3, -1], 10, [2, -1], 12, [4, -1], 16, [1, 1], 15, [7, -1]] \\
C_1' =& \;  [3, [5, 1], 8, [6, 1], 14, [8, -1], 6, [18, 2], 18, [18, 0], 17, [6, 2], 11, [8, 0]] \\
C_2' =& \;  [0, [9, 0], 9, [9, 2]] \\
C_3' =& \;  [13, [7, 0], 2, [7, 2]]
\end{align*}
\end{minipage}

\item 2-factor type $[11, 3, 3, 2]$: three starters \vspace*{-3mm} \\
\begin{minipage}{8cm}
\begin{align*}
C_0 =& \;  [15, [8, -1], 7, [6, 2], 13, [3, 0], 10, [8, 2], 2, [18, 0], 18, [18, -1], 9, [3, 1], 12, [5, 1], 17, [5, 2], 4, [7, -1], \\ & \; 11, [4, 0]] \\
C_1 =& \;  [8, [6, -1], 14, [2, -1], 16, [8, 1]] \\
C_2 =& \;  [0, [3, -1], 3, [2, 0], 1, [1, 2]] \\
C_3 =& \;  [6, [1, -1], 5, [1, 1]]
\end{align*}
\end{minipage}

\vspace*{-4mm}
\begin{minipage}{8cm}
\begin{align*}
C_0' =& \;  [12, [2, 2], 10, [1, 0], 11, [5, -1], 6, [7, 2], 13, [5, 0], 8, [18, 1], 18, [18, 2], 3, [4, -1], 7, [8, 0], 17, [6, 1], 5, \\ & \; [7, 1]] \\
C_1' =& \;  [2, [4, 2], 16, [2, 1], 14, [6, 0]] \\
C_2' =& \;  [4, [3, 2], 1, [4, 1], 15, [7, 0]] \\
C_3' =& \;  [0, [9, 0], 9, [9, 2]]
\end{align*}
\end{minipage}

\item 2-factor type $[10, 4, 3, 2]$: three starters \vspace*{-3mm} \\
\begin{minipage}{8cm}
\begin{align*}
C_0 =& \;  [6, [4, -1], 10, [3, 0], 7, [8, 1], 15, [4, 2], 1, [3, 1], 16, [5, 0], 11, [2, 2], 13, [1, 0], 14, [18, 2], 18, [18, 1]] \\
C_1 =& \;  [5, [7, 1], 12, [6, 2], 0, [1, 1], 17, [6, 0]] \\
C_2 =& \;  [4, [4, 0], 8, [5, -1], 3, [1, 2]] \\
C_3 =& \;  [2, [7, 0], 9, [7, 2]]
\end{align*}
\end{minipage}

\vspace*{-4mm}
\begin{minipage}{8cm}
\begin{align*}
C_0' =& \;  [5, [8, -1], 15, [6, 1], 3, [3, 2], 6, [8, 0], 14, [2, 1], 16, [6, -1], 10, [2, -1], 8, [5, 1], 13, [1, -1], 12, [7, -1]] \\
C_1' =& \;  [1, [8, 2], 11, [4, 1], 7, [18, 0], 18, [18, -1]] \\
C_2' =& \;  [17, [3, -1], 2, [2, 0], 4, [5, 2]] \\
C_3' =& \;  [0, [9, 0], 9, [9, 2]]
\end{align*}
\end{minipage}

\item 2-factor type $[9, 5, 3, 2]$: three starters \vspace*{-3mm} \\
\begin{minipage}{8cm}
\begin{align*}
C_0 =& \;  [10, [5, 1], 15, [1, -1], 14, [6, 2], 8, [4, 0], 4, [4, 2], 0, [3, 0], 3, [8, 2], 13, [1, 1], 12, [2, 0]] \\
C_1 =& \;  [9, [2, -1], 7, [1, 2], 6, [8, 0], 16, [3, 1], 1, [8, 1]] \\
C_2 =& \;  [2, [3, 2], 5, [18, -1], 18, [18, 0]] \\
C_3 =& \;  [11, [6, -1], 17, [6, 1]]
\end{align*}
\end{minipage}

\newpage

\begin{minipage}{8cm}
\begin{align*}
C_0' =& \;  [1, [18, 2], 18, [18, 1], 16, [1, 0], 15, [5, -1], 10, [2, 2], 12, [4, 1], 8, [3, -1], 11, [6, 0], 5, [4, -1]] \\
C_1' =& \;  [7, [5, 2], 2, [7, 1], 13, [8, -1], 3, [7, 0], 14, [7, -1]] \\
C_2' =& \;  [4, [2, 1], 6, [7, 2], 17, [5, 0]] \\
C_3' =& \;  [0, [9, 0], 9, [9, 2]]
\end{align*}
\end{minipage}

\item 2-factor type $[8, 6, 3, 2]$: three starters \vspace*{-3mm} \\
\begin{minipage}{8cm}
\begin{align*}
C_0 =& \;  [7, [18, 2], 18, [18, 0], 12, [7, 1], 1, [1, 2], 2, [6, 0], 14, [1, 1], 15, [3, 2], 0, [7, 0]] \\
C_1 =& \;  [10, [7, 2], 3, [5, 1], 16, [6, -1], 4, [5, 0], 17, [8, 2], 9, [1, 0]] \\
C_2 =& \;  [6, [2, 1], 8, [3, 1], 11, [5, -1]] \\
C_3 =& \;  [13, [8, -1], 5, [8, 1]]
\end{align*}
\end{minipage}

\vspace*{-4mm}
\begin{minipage}{8cm}
\begin{align*}
C_0' =& \;  [13, [6, 2], 7, [4, 0], 3, [2, 2], 5, [18, -1], 18, [18, 1], 4, [8, 0], 12, [4, 2], 16, [3, 0]] \\
C_1' =& \;  [11, [1, -1], 10, [4, 1], 14, [6, 1], 2, [3, -1], 17, [2, -1], 15, [4, -1]] \\
C_2' =& \;  [1, [7, -1], 8, [2, 0], 6, [5, 2]] \\
C_3' =& \;  [0, [9, 0], 9, [9, 2]]
\end{align*}
\end{minipage}

\item 2-factor type $[7, 7, 3, 2]$: three starters \vspace*{-3mm} \\
\begin{minipage}{8cm}
\begin{align*}
C_0 =& \;  [12, [4, 2], 16, [4, -1], 2, [18, 0], 18, [18, 2], 9, [3, 1], 6, [1, 0], 7, [5, 1]] \\
C_1 =& \;  [17, [6, 0], 11, [6, -1], 5, [2, -1], 3, [2, 2], 1, [1, 1], 0, [4, 0], 4, [5, 2]] \\
C_2 =& \;  [10, [3, 0], 13, [1, 2], 14, [4, 1]] \\
C_3 =& \;  [15, [7, 0], 8, [7, 2]]
\end{align*}
\end{minipage}

\vspace*{-4mm}
\begin{minipage}{8cm}
\begin{align*}
C_0' =& \;  [15, [2, 0], 13, [7, 1], 2, [1, -1], 1, [3, -1], 16, [8, 2], 8, [18, -1], 18, [18, 1]] \\
C_1' =& \;  [12, [6, 2], 6, [8, 0], 14, [8, 1], 4, [5, -1], 17, [7, -1], 10, [3, 2], 7, [5, 0]] \\
C_2' =& \;  [3, [8, -1], 11, [6, 1], 5, [2, 1]] \\
C_3' =& \;  [0, [9, 0], 9, [9, 2]]
\end{align*}
\end{minipage}

\item 2-factor type $[9, 4, 4, 2]$: three starters \vspace*{-3mm} \\
\begin{minipage}{8cm}
\begin{align*}
C_0 =& \;  [7, [2, 2], 9, [3, -1], 12, [1, -1], 13, [1, 0], 14, [4, 2], 0, [3, 0], 15, [5, 2], 2, [1, 1], 1, [6, 0]] \\
C_1 =& \;  [4, [18, -1], 18, [18, 1], 5, [6, -1], 11, [7, 1]] \\
C_2 =& \;  [8, [2, 1], 10, [7, -1], 3, [5, -1], 16, [8, -1]] \\
C_3 =& \;  [17, [7, 0], 6, [7, 2]]
\end{align*}
\end{minipage}

\vspace*{-4mm}
\begin{minipage}{8cm}
\begin{align*}
C_0' =& \;  [4, [6, 2], 16, [8, 0], 6, [3, 2], 3, [4, 0], 17, [4, -1], 13, [2, -1], 11, [1, 2], 10, [8, 1], 2, [2, 0]] \\
C_1' =& \;  [12, [18, 2], 18, [18, 0], 15, [8, 2], 7, [5, 0]] \\
C_2' =& \;  [5, [3, 1], 8, [6, 1], 14, [5, 1], 1, [4, 1]] \\
C_3' =& \;  [0, [9, 0], 9, [9, 2]]
\end{align*}
\end{minipage}

\item 2-factor type $[8, 5, 4, 2]$: three starters \vspace*{-3mm} \\
\begin{minipage}{8cm}
\begin{align*}
C_0 =& \;  [3, [1, -1], 2, [5, 1], 7, [18, 2], 18, [18, 0], 8, [6, 2], 14, [3, 1], 11, [1, 0], 10, [7, -1]] \\
C_1 =& \;  [4, [6, -1], 16, [4, 2], 12, [6, 0], 0, [1, 2], 1, [3, 0]] \\
C_2 =& \;  [5, [8, 1], 15, [6, 1], 9, [4, 0], 13, [8, 2]] \\
C_3 =& \;  [17, [7, 0], 6, [7, 2]]
\end{align*}
\end{minipage}

\newpage

\begin{minipage}{8cm}
\begin{align*}
C_0' =& \;  [7, [5, -1], 2, [3, -1], 17, [2, 2], 15, [4, -1], 1, [7, 1], 12, [8, 0], 4, [1, 1], 5, [2, 1]] \\
C_1' =& \;  [10, [3, 2], 13, [8, -1], 3, [5, 0], 8, [2, -1], 6, [4, 1]] \\
C_2' =& \;  [16, [2, 0], 14, [18, 1], 18, [18, -1], 11, [5, 2]] \\
C_3' =& \;  [0, [9, 0], 9, [9, 2]]
\end{align*}
\end{minipage}

\item 2-factor type $[7, 6, 4, 2]$: three starters \vspace*{-3mm} \\
\begin{minipage}{8cm}
\begin{align*}
C_0 =& \;  [2, [1, 1], 1, [2, 0], 17, [7, -1], 10, [2, 2], 8, [6, 0], 14, [8, 2], 4, [2, 1]] \\
C_1 =& \;  [13, [5, 1], 0, [18, 1], 18, [18, 2], 15, [3, 0], 12, [7, 2], 5, [8, 0]] \\
C_2 =& \;  [3, [5, 0], 16, [8, 1], 6, [3, 1], 9, [6, 2]] \\
C_3 =& \;  [7, [4, -1], 11, [4, 1]]
\end{align*}
\end{minipage}

\vspace*{-4mm}
\begin{minipage}{8cm}
\begin{align*}
C_0' =& \;  [18, [18, 0], 10, [5, 2], 15, [1, -1], 16, [6, -1], 4, [1, 0], 3, [4, 2], 17, [18, -1]] \\
C_1' =& \;  [13, [2, -1], 11, [1, 2], 12, [4, 0], 8, [3, -1], 5, [3, 2], 2, [7, 0]] \\
C_2' =& \;  [7, [6, 1], 1, [5, -1], 6, [8, -1], 14, [7, 1]] \\
C_3' =& \;  [0, [9, 0], 9, [9, 2]]
\end{align*}
\end{minipage}

\item 2-factor type $[7, 5, 5, 2]$: three starters \vspace*{-3mm} \\
\begin{minipage}{8cm}
\begin{align*}
C_0 =& \;  [0, [1, 1], 17, [6, 1], 11, [5, -1], 16, [6, 0], 10, [4, -1], 6, [7, 1], 13, [5, 2]] \\
C_1 =& \;  [9, [1, 2], 8, [6, -1], 14, [7, 0], 7, [3, 2], 4, [5, 0]] \\
C_2 =& \;  [2, [8, 1], 12, [7, -1], 1, [4, 0], 15, [18, 1], 18, [18, 2]] \\
C_3 =& \;  [3, [2, -1], 5, [2, 1]]
\end{align*}
\end{minipage}

\vspace*{-4mm}
\begin{minipage}{8cm}
\begin{align*}
C_0' =& \;  [18, [18, -1], 15, [7, 2], 4, [3, 1], 1, [2, 0], 3, [4, 2], 17, [5, 1], 12, [18, 0]] \\
C_1' =& \;  [11, [4, 1], 7, [3, 0], 10, [8, -1], 2, [6, 2], 8, [3, -1]] \\
C_2' =& \;  [14, [1, 0], 13, [8, 2], 5, [1, -1], 6, [8, 0], 16, [2, 2]] \\
C_3' =& \;  [0, [9, 0], 9, [9, 2]]
\end{align*}
\end{minipage}

\item 2-factor type $[6, 6, 5, 2]$: three starters \vspace*{-3mm} \\
\begin{minipage}{8cm}
\begin{align*}
C_0 =& \;  [3, [6, 1], 15, [1, 0], 14, [2, -1], 12, [5, -1], 7, [6, 2], 13, [8, -1]] \\
C_1 =& \;  [10, [1, -1], 9, [8, 2], 1, [2, 1], 17, [6, 0], 5, [3, -1], 2, [8, 1]] \\
C_2 =& \;  [8, [4, -1], 4, [2, 2], 6, [18, -1], 18, [18, 1], 16, [8, 0]] \\
C_3 =& \;  [0, [7, 0], 11, [7, 2]]
\end{align*}
\end{minipage}

\vspace*{-4mm}
\begin{minipage}{8cm}
\begin{align*}
C_0' =& \;  [4, [2, 0], 6, [3, 2], 3, [5, 0], 8, [7, -1], 1, [4, 2], 5, [1, 1]] \\
C_1' =& \;  [2, [7, 1], 13, [3, 1], 10, [6, -1], 16, [4, 0], 12, [5, 2], 7, [5, 1]] \\
C_2' =& \;  [17, [3, 0], 14, [1, 2], 15, [4, 1], 11, [18, 0], 18, [18, 2]] \\
C_3' =& \;  [0, [9, 0], 9, [9, 2]]
\end{align*}
\end{minipage}

\item 2-factor type $[15, 2, 2]$: three starters \vspace*{-3mm} \\
\begin{minipage}{8cm}
\begin{align*}
C_0 =& \;  [11, [1, 1], 12, [7, 2], 5, [3, 0], 2, [4, 2], 6, [18, 0], 18, [18, 2], 13, [6, 1], 7, [2, 0], 9, [6, 2], 15, [1, -1], 16, \\ & \; [2, 1], 14, [4, 0], 10, [7, -1], 3, [4, -1], 17, [6, -1]] \\
C_1 =& \;  [0, [8, -1], 8, [8, 1]] \\
C_2 =& \;  [1, [3, -1], 4, [3, 1]]
\end{align*}
\end{minipage}

\newpage

\begin{minipage}{8cm}
\begin{align*}
C_0' =& \;  [1, [6, 0], 13, [8, 2], 5, [2, -1], 7, [5, 0], 12, [2, 2], 14, [4, 1], 10, [7, 0], 17, [1, 2], 16, [8, 0], 8, [5, 2], 3, [1, 0], \\ & \; 2, [18, 1], 18, [18, -1], 15, [7, 1], 4, [3, 2]] \\
C_1' =& \;  [0, [9, 0], 9, [9, 2]] \\
C_2' =& \;  [6, [5, -1], 11, [5, 1]]
\end{align*}
\end{minipage}

\item 2-factor type $[14, 3, 2]$: three starters \vspace*{-3mm} \\
\begin{minipage}{8cm}
\begin{align*}
C_0 =& \;  [17, [3, 2], 14, [6, 0], 2, [6, 1], 8, [3, -1], 5, [4, 1], 9, [4, 2], 13, [1, 1], 12, [4, 0], 16, [5, -1], 11, [7, 2], 4, \\ & \; [1, 0], 3, [2, -1], 1, [6, 2], 7, [8, 0]] \\
C_1 =& \;  [0, [8, -1], 10, [4, -1], 6, [6, -1]] \\
C_2 =& \;  [18, [18, 0], 15, [18, 2]]
\end{align*}
\end{minipage}

\vspace*{-4mm}
\begin{minipage}{8cm}
\begin{align*}
C_0' =& \;  [10, [3, 0], 7, [2, 2], 5, [1, -1], 6, [2, 0], 4, [8, 1], 12, [2, 1], 14, [1, 2], 13, [7, -1], 2, [3, 1], 17, [18, -1], 18, \\ & \; [18, 1], 1, [7, 0], 8, [7, 1], 15, [5, 2]] \\
C_1' =& \;  [3, [8, 2], 11, [5, 0], 16, [5, 1]] \\
C_2' =& \;  [0, [9, 0], 9, [9, 2]]
\end{align*}
\end{minipage}

\item 2-factor type $[13, 4, 2]$: three starters \vspace*{-3mm} \\
\begin{minipage}{8cm}
\begin{align*}
C_0 =& \;  [3, [1, 2], 2, [7, 1], 13, [5, 1], 8, [2, 1], 6, [5, 0], 1, [2, -1], 17, [5, 2], 4, [7, -1], 11, [6, 1], 5, [7, 0], 12, [4, 2], \\ & \; 16, [18, -1], 18, [18, 0]] \\
C_1 =& \;  [0, [8, 0], 10, [4, 1], 14, [1, 1], 15, [3, 2]] \\
C_2 =& \;  [7, [2, 0], 9, [2, 2]]
\end{align*}
\end{minipage}

\vspace*{-4mm}
\begin{minipage}{8cm}
\begin{align*}
C_0' =& \;  [16, [8, 2], 8, [6, -1], 14, [1, -1], 15, [3, 0], 12, [5, -1], 7, [6, 2], 13, [8, 1], 5, [3, 1], 2, [4, 0], 6, [18, 1], 18, \\ & \; [18, 2], 1, [3, -1], 4, [6, 0]] \\
C_1' =& \;  [11, [8, -1], 3, [4, -1], 17, [7, 2], 10, [1, 0]] \\
C_2' =& \;  [0, [9, 0], 9, [9, 2]]
\end{align*}
\end{minipage}

\item 2-factor type $[12, 5, 2]$: three starters \vspace*{-3mm} \\
\begin{minipage}{8cm}
\begin{align*}
C_0 =& \;  [7, [1, 1], 6, [7, 1], 17, [2, 0], 15, [4, -1], 11, [1, 2], 12, [6, 0], 0, [8, -1], 10, [6, 1], 16, [8, 2], 8, [3, 0], 5, \\ & \; [18, 2], 18, [18, 1]] \\
C_1 =& \;  [13, [8, 1], 3, [1, -1], 2, [6, -1], 14, [5, -1], 9, [4, 1]] \\
C_2 =& \;  [1, [3, -1], 4, [3, 1]]
\end{align*}
\end{minipage}

\vspace*{-4mm}
\begin{minipage}{8cm}
\begin{align*}
C_0' =& \;  [7, [5, 0], 2, [5, 2], 15, [4, 0], 11, [2, 2], 13, [8, 0], 3, [4, 2], 17, [1, 0], 16, [6, 2], 10, [5, 1], 5, [7, 0], 12, [2, 1], \\ & \; 14, [7, 2]] \\
C_1' =& \;  [18, [18, 0], 8, [7, -1], 1, [3, 2], 4, [2, -1], 6, [18, -1]] \\
C_2' =& \;  [0, [9, 0], 9, [9, 2]]
\end{align*}
\end{minipage}

\item 2-factor type $[11, 6, 2]$: three starters \vspace*{-3mm} \\
\begin{minipage}{8cm}
\begin{align*}
C_0 =& \;  [7, [6, 0], 13, [7, 1], 2, [2, -1], 0, [3, 2], 15, [18, 0], 18, [18, -1], 4, [5, 2], 9, [1, -1], 10, [4, 0], 6, [3, -1], 3, \\ & \; [4, 2]] \\
C_1 =& \;  [8, [3, 1], 11, [5, 1], 16, [7, 2], 5, [7, 0], 12, [2, 1], 14, [6, -1]] \\
C_2 =& \;  [17, [2, 0], 1, [2, 2]]
\end{align*}
\end{minipage}

\newpage

\begin{minipage}{8cm}
\begin{align*}
C_0' =& \;  [8, [8, 0], 16, [6, 1], 4, [7, -1], 15, [6, 2], 3, [3, 0], 6, [4, -1], 2, [8, -1], 12, [1, 1], 13, [4, 1], 17, [8, 1], 7, \\ & \; [1, 2]] \\
C_1' =& \;  [18, [18, 1], 14, [5, 0], 1, [8, 2], 11, [1, 0], 10, [5, -1], 5, [18, 2]] \\
C_2' =& \;  [0, [9, 0], 9, [9, 2]]
\end{align*}
\end{minipage}

\item 2-factor type $[10, 7, 2]$: three starters \vspace*{-3mm} \\
\begin{minipage}{8cm}
\begin{align*}
C_0 =& \;  [17, [8, 2], 9, [6, 1], 3, [2, 0], 5, [18, 2], 18, [18, 1], 7, [5, 1], 2, [6, 0], 14, [2, 2], 12, [7, 0], 1, [2, -1]] \\
C_1 =& \;  [11, [1, -1], 10, [3, 1], 13, [7, 2], 6, [2, 1], 4, [4, 0], 0, [3, -1], 15, [4, -1]] \\
C_2 =& \;  [8, [8, -1], 16, [8, 1]]
\end{align*}
\end{minipage}

\vspace*{-4mm}
\begin{minipage}{8cm}
\begin{align*}
C_0' =& \;  [18, [18, 0], 14, [5, 2], 1, [7, 1], 12, [4, 1], 8, [5, 0], 13, [7, -1], 6, [1, 2], 5, [3, 0], 2, [4, 2], 16, [18, -1]] \\
C_1' =& \;  [11, [6, 2], 17, [8, 0], 7, [3, 2], 4, [1, 0], 3, [6, -1], 15, [5, -1], 10, [1, 1]] \\
C_2' =& \;  [0, [9, 0], 9, [9, 2]]
\end{align*}
\end{minipage}

\item 2-factor type $[9, 8, 2]$: three starters \vspace*{-3mm} \\
\begin{minipage}{8cm}
\begin{align*}
C_0 =& \;  [18, [18, 2], 15, [7, 0], 8, [4, 1], 12, [2, -1], 10, [3, 1], 13, [6, 2], 1, [5, -1], 6, [6, 1], 0, [18, 0]] \\
C_1 =& \;  [4, [1, 1], 3, [5, 0], 16, [7, -1], 9, [4, -1], 5, [2, 1], 7, [7, 1], 14, [3, -1], 11, [7, 2]] \\
C_2 =& \;  [2, [3, 0], 17, [3, 2]]
\end{align*}
\end{minipage}

\vspace*{-4mm}
\begin{minipage}{8cm}
\begin{align*}
C_0' =& \;  [18, [18, -1], 10, [5, 2], 5, [1, -1], 6, [2, 0], 4, [8, 2], 14, [1, 0], 15, [2, 2], 13, [6, 0], 1, [18, 1]] \\
C_1' =& \;  [8, [6, -1], 2, [1, 2], 3, [8, 0], 11, [4, 2], 7, [8, 1], 17, [5, 1], 12, [4, 0], 16, [8, -1]] \\
C_2' =& \;  [0, [9, 0], 9, [9, 2]]
\end{align*}
\end{minipage}

\end{itemizenew}

\section{Computational results for $n=20$}\label{app:20}

\begin{itemizenew}
\item 2-factor type $[4, 2, 2, 2, 2, 2, 2, 2, 2]$: one starter \vspace*{-3mm} \\
\hspace*{-8mm}
\begin{minipage}{8cm}
\begin{align*}
C_0 =& \;  [4, [8, 0], 12, [1, 1], 11, [8, 1], 3, [1, 0]] \\
C_1 =& \;  [6, [19, 0], 19, [19, 1]] \\
C_2 =& \;  [0, [5, 0], 5, [5, 1]] \\
C_3 =& \;  [15, [2, 0], 13, [2, 1]] \\
C_4 =& \;  [18, [3, 0], 2, [3, 1]] \\
C_5 =& \;  [1, [6, 0], 7, [6, 1]] \\
C_6 =& \;  [8, [9, 0], 17, [9, 1]] \\
C_7 =& \;  [9, [7, 0], 16, [7, 1]] \\
C_8 =& \;  [14, [4, 0], 10, [4, 1]]
\end{align*}
\end{minipage}

\newpage

\item 2-factor type $[3, 3, 2, 2, 2, 2, 2, 2, 2]$: two starters \vspace*{-3mm} \\
\hspace*{-11mm}
\begin{minipage}{8cm}
\begin{align*}
C_0 =& \;  [2, [1, -1], 1, [19, 1], 19, [19, -1]] \\
C_1 =& \;  [9, [7, 0], 16, [7, 1], 4, [5, 2]] \\
C_2 =& \;  [14, [9, -1], 5, [9, 1]] \\
C_3 =& \;  [8, [1, 0], 7, [1, 2]] \\
C_4 =& \;  [18, [4, 0], 3, [4, 2]] \\
C_5 =& \;  [10, [3, 0], 13, [3, 2]] \\
C_6 =& \;  [17, [2, -1], 0, [2, 1]] \\
C_7 =& \;  [12, [3, -1], 15, [3, 1]] \\
C_8 =& \;  [11, [5, -1], 6, [5, 1]]
\end{align*}
\end{minipage}
\begin{minipage}{8cm}
\begin{align*}
C_0' =& \;  [17, [7, 2], 10, [5, 0], 5, [7, -1]] \\
C_1' =& \;  [0, [1, 1], 1, [19, 2], 19, [19, 0]] \\
C_2' =& \;  [15, [6, -1], 2, [6, 1]] \\
C_3' =& \;  [16, [8, 0], 8, [8, 2]] \\
C_4' =& \;  [18, [6, 0], 12, [6, 2]] \\
C_5' =& \;  [4, [9, 0], 13, [9, 2]] \\
C_6' =& \;  [11, [2, 0], 9, [2, 2]] \\
C_7' =& \;  [7, [4, -1], 3, [4, 1]] \\
C_8' =& \;  [14, [8, -1], 6, [8, 1]]
\end{align*}
\end{minipage}

\item 2-factor type $[6, 2, 2, 2, 2, 2, 2, 2]$: two starters \vspace*{-3mm} \\
\begin{minipage}{8cm}
\begin{align*}
C_0 =& \;  [10, [8, 2], 18, [9, 0], 9, [19, 1], 19, [19, -1], 2, [4, 2], 6, [4, 0]] \\
C_1 =& \;  [12, [1, -1], 13, [1, 1]] \\
C_2 =& \;  [7, [9, -1], 17, [9, 1]] \\
C_3 =& \;  [1, [4, -1], 16, [4, 1]] \\
C_4 =& \;  [8, [7, -1], 15, [7, 1]] \\
C_5 =& \;  [14, [3, -1], 11, [3, 1]] \\
C_6 =& \;  [5, [1, 0], 4, [1, 2]] \\
C_7 =& \;  [3, [3, 0], 0, [3, 2]]
\end{align*}
\end{minipage}

\vspace*{-4mm}
\begin{minipage}{8cm}
\begin{align*}
C_0' =& \;  [1, [7, 0], 13, [9, 2], 4, [19, 0], 19, [19, 2], 16, [8, 0], 8, [7, 2]] \\
C_1' =& \;  [3, [5, 0], 17, [5, 2]] \\
C_2' =& \;  [0, [5, -1], 14, [5, 1]] \\
C_3' =& \;  [5, [2, 0], 7, [2, 2]] \\
C_4' =& \;  [6, [6, -1], 12, [6, 1]] \\
C_5' =& \;  [11, [2, -1], 9, [2, 1]] \\
C_6' =& \;  [15, [6, 0], 2, [6, 2]] \\
C_7' =& \;  [10, [8, -1], 18, [8, 1]]
\end{align*}
\end{minipage}

\item 2-factor type $[5, 3, 2, 2, 2, 2, 2, 2]$: one starter \vspace*{-3mm} \\
\begin{minipage}{8cm}
\begin{align*}
C_0 =& \;  [6, [3, 0], 3, [4, 1], 18, [9, 1], 8, [3, 1], 11, [5, 0]] \\
C_1 =& \;  [10, [4, 0], 14, [5, 1], 0, [9, 0]] \\
C_2 =& \;  [16, [1, 0], 17, [1, 1]] \\
C_3 =& \;  [12, [8, 0], 4, [8, 1]] \\
C_4 =& \;  [9, [19, 0], 19, [19, 1]] \\
C_5 =& \;  [15, [6, 0], 2, [6, 1]] \\
C_6 =& \;  [1, [7, 0], 13, [7, 1]] \\
C_7 =& \;  [7, [2, 0], 5, [2, 1]]
\end{align*}
\end{minipage}

\newpage

\item 2-factor type $[4, 4, 2, 2, 2, 2, 2, 2]$: one starter \vspace*{-3mm} \\
\hspace*{-8mm}
\begin{minipage}{8cm}
\begin{align*}
C_0 =& \;  [15, [3, 0], 12, [8, 0], 4, [2, 1], 2, [6, 1]] \\
C_1 =& \;  [3, [6, 0], 16, [2, 0], 14, [8, 1], 6, [3, 1]] \\
C_2 =& \;  [11, [1, 0], 10, [1, 1]] \\
C_3 =& \;  [19, [19, 0], 0, [19, 1]] \\
C_4 =& \;  [8, [7, 0], 1, [7, 1]] \\
C_5 =& \;  [18, [5, 0], 13, [5, 1]] \\
C_6 =& \;  [17, [9, 0], 7, [9, 1]] \\
C_7 =& \;  [5, [4, 0], 9, [4, 1]]
\end{align*}
\end{minipage}

\item 2-factor type $[4, 3, 3, 2, 2, 2, 2, 2]$: two starters \vspace*{-3mm} \\
\hspace*{-5mm}
\begin{minipage}{8cm}
\begin{align*}
C_0 =& \;  [16, [3, 1], 13, [9, 0], 3, [7, -1], 15, [1, 2]] \\
C_1 =& \;  [18, [7, 1], 11, [9, -1], 1, [2, 1]] \\
C_2 =& \;  [4, [1, 0], 5, [7, 2], 12, [8, 1]] \\
C_3 =& \;  [19, [19, 0], 10, [19, 2]] \\
C_4 =& \;  [14, [6, 0], 8, [6, 2]] \\
C_5 =& \;  [2, [5, -1], 7, [5, 1]] \\
C_6 =& \;  [0, [6, -1], 6, [6, 1]] \\
C_7 =& \;  [17, [8, 0], 9, [8, 2]]
\end{align*}
\end{minipage}
\begin{minipage}{8cm}
\begin{align*}
C_0' =& \;  [2, [9, 2], 11, [8, -1], 0, [5, 0], 5, [3, -1]] \\
C_1' =& \;  [7, [9, 1], 17, [7, 0], 10, [3, 2]] \\
C_2' =& \;  [4, [3, 0], 1, [2, -1], 18, [5, 2]] \\
C_3' =& \;  [8, [2, 0], 6, [2, 2]] \\
C_4' =& \;  [13, [4, -1], 9, [4, 1]] \\
C_5' =& \;  [3, [19, -1], 19, [19, 1]] \\
C_6' =& \;  [14, [1, -1], 15, [1, 1]] \\
C_7' =& \;  [16, [4, 0], 12, [4, 2]]
\end{align*}
\end{minipage}

\item 2-factor type $[3, 3, 3, 3, 2, 2, 2, 2]$: one starter \vspace*{-3mm} \\
\hspace*{-14mm}
\begin{minipage}{8cm}
\begin{align*}
C_0 =& \;  [2, [6, 0], 8, [2, 0], 6, [4, 1]] \\
C_1 =& \;  [3, [9, 0], 12, [2, 1], 10, [7, 0]] \\
C_2 =& \;  [17, [7, 1], 5, [4, 0], 1, [3, 0]] \\
C_3 =& \;  [4, [9, 1], 13, [6, 1], 7, [3, 1]] \\
C_4 =& \;  [18, [19, 0], 19, [19, 1]] \\
C_5 =& \;  [16, [1, 0], 15, [1, 1]] \\
C_6 =& \;  [9, [5, 0], 14, [5, 1]] \\
C_7 =& \;  [11, [8, 0], 0, [8, 1]]
\end{align*}
\end{minipage}

\item 2-factor type $[8, 2, 2, 2, 2, 2, 2]$: one starter \vspace*{-3mm} \\
\begin{minipage}{8cm}
\begin{align*}
C_0 =& \;  [18, [4, 1], 14, [9, 1], 4, [4, 0], 8, [2, 1], 6, [9, 0], 16, [1, 0], 17, [2, 0], 0, [1, 1]] \\
C_1 =& \;  [15, [7, 0], 3, [7, 1]] \\
C_2 =& \;  [10, [3, 0], 13, [3, 1]] \\
C_3 =& \;  [5, [6, 0], 11, [6, 1]] \\
C_4 =& \;  [19, [19, 0], 12, [19, 1]] \\
C_5 =& \;  [2, [5, 0], 7, [5, 1]] \\
C_6 =& \;  [9, [8, 0], 1, [8, 1]]
\end{align*}
\end{minipage}

\newpage

\item 2-factor type $[7, 3, 2, 2, 2, 2, 2]$: two starters \vspace*{-3mm} \\
\begin{minipage}{8cm}
\begin{align*}
C_0 =& \;  [10, [2, 1], 12, [9, 2], 3, [5, -1], 8, [7, -1], 15, [1, 0], 14, [7, 1], 2, [8, 1]] \\
C_1 =& \;  [17, [1, -1], 16, [9, 1], 6, [8, -1]] \\
C_2 =& \;  [18, [2, 0], 1, [2, 2]] \\
C_3 =& \;  [5, [19, 0], 19, [19, 2]] \\
C_4 =& \;  [0, [6, 0], 13, [6, 2]] \\
C_5 =& \;  [11, [4, -1], 7, [4, 1]] \\
C_6 =& \;  [9, [5, 0], 4, [5, 2]]
\end{align*}
\end{minipage}

\vspace*{-4mm}
\begin{minipage}{8cm}
\begin{align*}
C_0' =& \;  [19, [19, 1], 5, [9, 0], 15, [5, 1], 1, [2, -1], 18, [8, 2], 7, [1, 1], 6, [19, -1]] \\
C_1' =& \;  [4, [8, 0], 12, [9, -1], 3, [1, 2]] \\
C_2' =& \;  [0, [3, -1], 16, [3, 1]] \\
C_3' =& \;  [9, [4, 0], 13, [4, 2]] \\
C_4' =& \;  [10, [7, 0], 17, [7, 2]] \\
C_5' =& \;  [11, [3, 0], 14, [3, 2]] \\
C_6' =& \;  [2, [6, -1], 8, [6, 1]]
\end{align*}
\end{minipage}

\item 2-factor type $[6, 4, 2, 2, 2, 2, 2]$: two starters \vspace*{-3mm} \\
\begin{minipage}{8cm}
\begin{align*}
C_0 =& \;  [1, [9, 0], 10, [4, -1], 6, [1, 2], 7, [4, 0], 3, [1, 1], 4, [3, 2]] \\
C_1 =& \;  [11, [2, 0], 9, [7, 1], 16, [1, -1], 15, [4, 2]] \\
C_2 =& \;  [5, [5, 0], 0, [5, 2]] \\
C_3 =& \;  [13, [8, -1], 2, [8, 1]] \\
C_4 =& \;  [19, [19, -1], 14, [19, 1]] \\
C_5 =& \;  [8, [9, -1], 17, [9, 1]] \\
C_6 =& \;  [12, [6, -1], 18, [6, 1]]
\end{align*}
\end{minipage}

\vspace*{-4mm}
\begin{minipage}{8cm}
\begin{align*}
C_0' =& \;  [1, [1, 0], 2, [7, 2], 9, [4, 1], 5, [6, 0], 18, [7, -1], 11, [9, 2]] \\
C_1' =& \;  [10, [6, 2], 16, [3, 0], 0, [2, 2], 17, [7, 0]] \\
C_2' =& \;  [7, [8, 0], 15, [8, 2]] \\
C_3' =& \;  [6, [3, -1], 3, [3, 1]] \\
C_4' =& \;  [19, [19, 0], 4, [19, 2]] \\
C_5' =& \;  [12, [2, -1], 14, [2, 1]] \\
C_6' =& \;  [13, [5, -1], 8, [5, 1]]
\end{align*}
\end{minipage}

\item 2-factor type $[5, 5, 2, 2, 2, 2, 2]$: two starters \vspace*{-3mm} \\
\begin{minipage}{8cm}
\begin{align*}
C_0 =& \;  [5, [1, 2], 4, [5, -1], 9, [9, -1], 18, [7, 0], 6, [1, -1]] \\
C_1 =& \;  [12, [1, 1], 13, [2, -1], 11, [8, -1], 3, [7, 2], 10, [2, 0]] \\
C_2 =& \;  [2, [4, 0], 17, [4, 2]] \\
C_3 =& \;  [0, [8, 0], 8, [8, 2]] \\
C_4 =& \;  [14, [7, -1], 7, [7, 1]] \\
C_5 =& \;  [1, [5, 0], 15, [5, 2]] \\
C_6 =& \;  [19, [19, -1], 16, [19, 1]]
\end{align*}
\end{minipage}

\newpage

\begin{minipage}{8cm}
\begin{align*}
C_0' =& \;  [16, [3, 1], 0, [2, 1], 2, [9, 1], 11, [3, -1], 8, [8, 1]] \\
C_1' =& \;  [15, [1, 0], 14, [2, 2], 12, [5, 1], 7, [19, 0], 19, [19, 2]] \\
C_2' =& \;  [9, [4, -1], 13, [4, 1]] \\
C_3' =& \;  [3, [3, 0], 6, [3, 2]] \\
C_4' =& \;  [1, [9, 0], 10, [9, 2]] \\
C_5' =& \;  [18, [6, -1], 5, [6, 1]] \\
C_6' =& \;  [4, [6, 0], 17, [6, 2]]
\end{align*}
\end{minipage}

\item 2-factor type $[6, 3, 3, 2, 2, 2, 2]$: one starter \vspace*{-3mm} \\
\begin{minipage}{8cm}
\begin{align*}
C_0 =& \;  [6, [9, 1], 16, [4, 0], 12, [4, 1], 8, [9, 0], 18, [8, 0], 7, [1, 0]] \\
C_1 =& \;  [10, [8, 1], 2, [7, 1], 9, [1, 1]] \\
C_2 =& \;  [11, [6, 0], 17, [7, 0], 5, [6, 1]] \\
C_3 =& \;  [4, [3, 0], 1, [3, 1]] \\
C_4 =& \;  [3, [19, 0], 19, [19, 1]] \\
C_5 =& \;  [0, [5, 0], 14, [5, 1]] \\
C_6 =& \;  [15, [2, 0], 13, [2, 1]]
\end{align*}
\end{minipage}

\item 2-factor type $[5, 4, 3, 2, 2, 2, 2]$: one starter \vspace*{-3mm} \\
\begin{minipage}{8cm}
\begin{align*}
C_0 =& \;  [9, [19, 0], 19, [19, 1], 1, [5, 0], 15, [9, 1], 5, [4, 1]] \\
C_1 =& \;  [11, [4, 0], 7, [5, 1], 2, [6, 0], 8, [3, 1]] \\
C_2 =& \;  [6, [6, 1], 12, [9, 0], 3, [3, 0]] \\
C_3 =& \;  [17, [2, 0], 0, [2, 1]] \\
C_4 =& \;  [18, [8, 0], 10, [8, 1]] \\
C_5 =& \;  [16, [7, 0], 4, [7, 1]] \\
C_6 =& \;  [13, [1, 0], 14, [1, 1]]
\end{align*}
\end{minipage}

\item 2-factor type $[4, 4, 4, 2, 2, 2, 2]$: one starter \vspace*{-3mm} \\
\hspace*{-4mm}
\begin{minipage}{8cm}
\begin{align*}
C_0 =& \;  [8, [8, 0], 0, [6, 0], 6, [7, 0], 13, [5, 0]] \\
C_1 =& \;  [5, [2, 1], 7, [8, 1], 18, [7, 1], 11, [6, 1]] \\
C_2 =& \;  [10, [19, 1], 19, [19, 0], 17, [2, 0], 15, [5, 1]] \\
C_3 =& \;  [9, [3, 0], 12, [3, 1]] \\
C_4 =& \;  [14, [9, 0], 4, [9, 1]] \\
C_5 =& \;  [2, [1, 0], 3, [1, 1]] \\
C_6 =& \;  [1, [4, 0], 16, [4, 1]]
\end{align*}
\end{minipage}

\item 2-factor type $[5, 3, 3, 3, 2, 2, 2]$: two starters \vspace*{-3mm} \\
\begin{minipage}{8cm}
\begin{align*}
C_0 =& \;  [4, [9, 1], 14, [1, 0], 15, [3, 2], 12, [6, 1], 6, [2, 1]] \\
C_1 =& \;  [8, [8, 2], 0, [9, 0], 10, [2, -1]] \\
C_2 =& \;  [16, [6, -1], 3, [19, 0], 19, [19, 2]] \\
C_3 =& \;  [9, [8, -1], 1, [9, 2], 11, [2, 0]] \\
C_4 =& \;  [17, [4, -1], 13, [4, 1]] \\
C_5 =& \;  [18, [6, 0], 5, [6, 2]] \\
C_6 =& \;  [7, [5, -1], 2, [5, 1]]
\end{align*}
\end{minipage}

\newpage

\begin{minipage}{8cm}
\begin{align*}
C_0' =& \;  [11, [9, -1], 2, [3, -1], 18, [7, 1], 6, [19, 1], 19, [19, -1]] \\
C_1' =& \;  [17, [5, 2], 3, [8, 1], 14, [3, 0]] \\
C_2' =& \;  [4, [1, 2], 5, [7, -1], 12, [8, 0]] \\
C_3' =& \;  [10, [5, 0], 15, [2, 2], 13, [3, 1]] \\
C_4' =& \;  [1, [4, 0], 16, [4, 2]] \\
C_5' =& \;  [0, [7, 0], 7, [7, 2]] \\
C_6' =& \;  [9, [1, -1], 8, [1, 1]]
\end{align*}
\end{minipage}

\item 2-factor type $[4, 4, 3, 3, 2, 2, 2]$: two starters \vspace*{-3mm} \\
\hspace*{-7mm}
\begin{minipage}{8cm}
\begin{align*}
C_0 =& \;  [16, [7, 2], 9, [8, 0], 17, [2, 2], 15, [1, 0]] \\
C_1 =& \;  [7, [4, 1], 3, [1, -1], 4, [9, -1], 13, [6, 1]] \\
C_2 =& \;  [5, [9, 1], 14, [19, 2], 19, [19, 0]] \\
C_3 =& \;  [1, [9, 2], 10, [9, 0], 0, [1, 1]] \\
C_4 =& \;  [11, [5, 0], 6, [5, 2]] \\
C_5 =& \;  [18, [3, -1], 2, [3, 1]] \\
C_6 =& \;  [12, [4, 0], 8, [4, 2]]
\end{align*}
\end{minipage}
\begin{minipage}{8cm}
\begin{align*}
C_0' =& \;  [18, [1, 2], 0, [6, 0], 13, [8, -1], 5, [6, -1]] \\
C_1' =& \;  [3, [19, -1], 19, [19, 1], 8, [3, 0], 11, [8, 2]] \\
C_2' =& \;  [17, [4, -1], 2, [7, 0], 14, [3, 2]] \\
C_3' =& \;  [15, [8, 1], 7, [2, 0], 9, [6, 2]] \\
C_4' =& \;  [4, [7, -1], 16, [7, 1]] \\
C_5' =& \;  [1, [5, -1], 6, [5, 1]] \\
C_6' =& \;  [10, [2, -1], 12, [2, 1]]
\end{align*}
\end{minipage}

\item 2-factor type $[4, 3, 3, 3, 3, 2, 2]$: one starter \vspace*{-3mm} \\
\hspace*{-8mm}
\begin{minipage}{8cm}
\begin{align*}
C_0 =& \;  [7, [19, 0], 19, [19, 1], 6, [2, 1], 8, [1, 0]] \\
C_1 =& \;  [12, [3, 1], 9, [9, 1], 18, [6, 1]] \\
C_2 =& \;  [16, [5, 0], 2, [6, 0], 15, [1, 1]] \\
C_3 =& \;  [4, [3, 0], 1, [7, 1], 13, [9, 0]] \\
C_4 =& \;  [17, [7, 0], 5, [2, 0], 3, [5, 1]] \\
C_5 =& \;  [10, [4, 0], 14, [4, 1]] \\
C_6 =& \;  [0, [8, 0], 11, [8, 1]]
\end{align*}
\end{minipage}

\item 2-factor type $[3, 3, 3, 3, 3, 3, 2]$: two starters \vspace*{-3mm} \\
\hspace*{-13mm}
\begin{minipage}{8cm}
\begin{align*}
C_0 =& \;  [7, [5, 1], 12, [19, 2], 19, [19, 0]] \\
C_1 =& \;  [4, [2, 2], 2, [7, 0], 14, [9, 1]] \\
C_2 =& \;  [0, [3, 2], 16, [3, 0], 13, [6, 1]] \\
C_3 =& \;  [5, [5, -1], 10, [5, 0], 15, [9, 2]] \\
C_4 =& \;  [1, [2, -1], 18, [4, 2], 3, [2, 0]] \\
C_5 =& \;  [8, [2, 1], 6, [3, -1], 9, [1, 1]] \\
C_6 =& \;  [11, [6, 0], 17, [6, 2]]
\end{align*}
\end{minipage}
\begin{minipage}{8cm}
\begin{align*}
C_0' =& \;  [4, [3, 1], 7, [7, 2], 14, [9, 0]] \\
C_1' =& \;  [1, [5, 2], 15, [1, 0], 16, [4, 1]] \\
C_2' =& \;  [18, [4, -1], 3, [9, -1], 12, [6, -1]] \\
C_3' =& \;  [13, [4, 0], 9, [7, -1], 2, [8, 2]] \\
C_4' =& \;  [5, [1, 2], 6, [8, 0], 17, [7, 1]] \\
C_5' =& \;  [10, [1, -1], 11, [19, -1], 19, [19, 1]] \\
C_6' =& \;  [0, [8, -1], 8, [8, 1]]
\end{align*}
\end{minipage}

\item 2-factor type $[10, 2, 2, 2, 2, 2]$: two starters \vspace*{-3mm} \\
\begin{minipage}{8cm}
\begin{align*}
C_0 =& \;  [1, [4, 1], 5, [2, 1], 7, [1, 1], 8, [3, 1], 11, [8, 1], 0, [3, 2], 3, [5, 0], 17, [8, 2], 6, [3, 0], 9, [8, -1]] \\
C_1 =& \;  [13, [1, 0], 14, [1, 2]] \\
C_2 =& \;  [19, [19, -1], 16, [19, 1]] \\
C_3 =& \;  [18, [6, -1], 12, [6, 1]] \\
C_4 =& \;  [15, [5, -1], 10, [5, 1]] \\
C_5 =& \;  [4, [2, 0], 2, [2, 2]]
\end{align*}
\end{minipage}

\newpage

\begin{minipage}{8cm}
\begin{align*}
C_0' =& \;  [8, [8, 0], 16, [6, 2], 10, [3, -1], 13, [6, 0], 0, [19, 2], 19, [19, 0], 6, [5, 2], 1, [1, -1], 2, [2, -1], 4, [4, -1]] \\
C_1' =& \;  [14, [7, -1], 7, [7, 1]] \\
C_2' =& \;  [12, [9, 0], 3, [9, 2]] \\
C_3' =& \;  [17, [7, 0], 5, [7, 2]] \\
C_4' =& \;  [9, [9, -1], 18, [9, 1]] \\
C_5' =& \;  [15, [4, 0], 11, [4, 2]]
\end{align*}
\end{minipage}

\item 2-factor type $[9, 3, 2, 2, 2, 2]$: one starter \vspace*{-3mm} \\
\begin{minipage}{8cm}
\begin{align*}
C_0 =& \;  [1, [2, 1], 18, [1, 0], 17, [19, 1], 19, [19, 0], 15, [8, 0], 4, [1, 1], 5, [9, 1], 14, [2, 0], 12, [8, 1]] \\
C_1 =& \;  [2, [5, 0], 16, [9, 0], 7, [5, 1]] \\
C_2 =& \;  [8, [3, 0], 11, [3, 1]] \\
C_3 =& \;  [9, [4, 0], 13, [4, 1]] \\
C_4 =& \;  [3, [7, 0], 10, [7, 1]] \\
C_5 =& \;  [0, [6, 0], 6, [6, 1]]
\end{align*}
\end{minipage}

\item 2-factor type $[8, 4, 2, 2, 2, 2]$: one starter \vspace*{-3mm} \\
\begin{minipage}{8cm}
\begin{align*}
C_0 =& \;  [7, [19, 1], 19, [19, 0], 14, [6, 1], 8, [2, 1], 6, [4, 1], 10, [7, 0], 17, [4, 0], 13, [6, 0]] \\
C_1 =& \;  [18, [7, 1], 11, [2, 0], 9, [5, 1], 4, [5, 0]] \\
C_2 =& \;  [16, [3, 0], 0, [3, 1]] \\
C_3 =& \;  [1, [8, 0], 12, [8, 1]] \\
C_4 =& \;  [15, [9, 0], 5, [9, 1]] \\
C_5 =& \;  [3, [1, 0], 2, [1, 1]]
\end{align*}
\end{minipage}

\item 2-factor type $[7, 5, 2, 2, 2, 2]$: one starter \vspace*{-3mm} \\
\begin{minipage}{8cm}
\begin{align*}
C_0 =& \;  [13, [1, 0], 12, [2, 0], 10, [5, 0], 15, [1, 1], 16, [4, 0], 1, [2, 1], 3, [9, 1]] \\
C_1 =& \;  [14, [5, 1], 9, [9, 0], 18, [7, 1], 6, [4, 1], 2, [7, 0]] \\
C_2 =& \;  [7, [19, 0], 19, [19, 1]] \\
C_3 =& \;  [4, [6, 0], 17, [6, 1]] \\
C_4 =& \;  [0, [8, 0], 11, [8, 1]] \\
C_5 =& \;  [8, [3, 0], 5, [3, 1]]
\end{align*}
\end{minipage}

\item 2-factor type $[6, 6, 2, 2, 2, 2]$: one starter \vspace*{-3mm} \\
\begin{minipage}{8cm}
\begin{align*}
C_0 =& \;  [17, [2, 1], 0, [1, 1], 1, [4, 1], 5, [4, 0], 9, [5, 1], 14, [3, 0]] \\
C_1 =& \;  [3, [7, 1], 10, [2, 0], 8, [3, 1], 11, [5, 0], 16, [7, 0], 4, [1, 0]] \\
C_2 =& \;  [13, [8, 0], 2, [8, 1]] \\
C_3 =& \;  [15, [9, 0], 6, [9, 1]] \\
C_4 =& \;  [7, [19, 0], 19, [19, 1]] \\
C_5 =& \;  [12, [6, 0], 18, [6, 1]]
\end{align*}
\end{minipage}

\newpage

\item 2-factor type $[8, 3, 3, 2, 2, 2]$: two starters \vspace*{-3mm} \\
\begin{minipage}{8cm}
\begin{align*}
C_0 =& \;  [17, [6, 1], 11, [6, 0], 5, [5, 1], 10, [3, 2], 13, [7, 0], 1, [3, 1], 4, [4, -1], 0, [2, 2]] \\
C_1 =& \;  [9, [1, -1], 8, [4, 1], 12, [3, -1]] \\
C_2 =& \;  [16, [2, 0], 18, [8, 2], 7, [9, -1]] \\
C_3 =& \;  [6, [9, 0], 15, [9, 2]] \\
C_4 =& \;  [14, [7, -1], 2, [7, 1]] \\
C_5 =& \;  [3, [19, 0], 19, [19, 2]]
\end{align*}
\end{minipage}

\vspace*{-4mm}
\begin{minipage}{8cm}
\begin{align*}
C_0' =& \;  [2, [3, 0], 18, [1, 2], 17, [8, 0], 9, [4, 2], 5, [1, 1], 4, [2, -1], 6, [5, -1], 11, [9, 1]] \\
C_1' =& \;  [1, [7, 2], 13, [6, -1], 0, [1, 0]] \\
C_2' =& \;  [10, [6, 2], 16, [2, 1], 14, [4, 0]] \\
C_3' =& \;  [3, [5, 0], 8, [5, 2]] \\
C_4' =& \;  [12, [19, -1], 19, [19, 1]] \\
C_5' =& \;  [15, [8, -1], 7, [8, 1]]
\end{align*}
\end{minipage}

\item 2-factor type $[7, 4, 3, 2, 2, 2]$: two starters \vspace*{-3mm} \\
\begin{minipage}{8cm}
\begin{align*}
C_0 =& \;  [5, [2, 1], 3, [9, 1], 13, [7, 0], 6, [3, 1], 9, [6, 1], 15, [4, -1], 11, [6, 2]] \\
C_1 =& \;  [0, [5, 0], 14, [19, 2], 19, [19, 0], 16, [3, 2]] \\
C_2 =& \;  [4, [8, 0], 12, [4, 2], 8, [4, 1]] \\
C_3 =& \;  [2, [1, 0], 1, [1, 2]] \\
C_4 =& \;  [10, [7, -1], 17, [7, 1]] \\
C_5 =& \;  [18, [8, -1], 7, [8, 1]]
\end{align*}
\end{minipage}

\vspace*{-4mm}
\begin{minipage}{8cm}
\begin{align*}
C_0' =& \;  [5, [6, -1], 11, [9, -1], 1, [2, -1], 3, [4, 0], 18, [8, 2], 7, [5, 1], 2, [3, -1]] \\
C_1' =& \;  [0, [7, 2], 12, [5, -1], 17, [1, 1], 16, [3, 0]] \\
C_2' =& \;  [15, [6, 0], 9, [5, 2], 14, [1, -1]] \\
C_3' =& \;  [6, [2, 0], 8, [2, 2]] \\
C_4' =& \;  [10, [19, -1], 19, [19, 1]] \\
C_5' =& \;  [13, [9, 0], 4, [9, 2]]
\end{align*}
\end{minipage}

\item 2-factor type $[6, 5, 3, 2, 2, 2]$: two starters \vspace*{-3mm} \\
\begin{minipage}{8cm}
\begin{align*}
C_0 =& \;  [13, [9, -1], 4, [2, -1], 2, [7, 2], 9, [1, 1], 8, [2, 0], 10, [3, 1]] \\
C_1 =& \;  [14, [2, 1], 16, [4, 2], 1, [9, 1], 11, [8, 1], 3, [8, 0]] \\
C_2 =& \;  [6, [7, -1], 18, [6, 2], 5, [1, 0]] \\
C_3 =& \;  [19, [19, 0], 0, [19, 2]] \\
C_4 =& \;  [7, [9, 0], 17, [9, 2]] \\
C_5 =& \;  [15, [3, 0], 12, [3, 2]]
\end{align*}
\end{minipage}

\vspace*{-4mm}
\begin{minipage}{8cm}
\begin{align*}
C_0' =& \;  [1, [1, 2], 0, [4, 0], 4, [8, 2], 12, [4, 1], 8, [6, 0], 2, [1, -1]] \\
C_1' =& \;  [6, [7, 0], 13, [8, -1], 5, [2, 2], 7, [19, -1], 19, [19, 1]] \\
C_2' =& \;  [11, [4, -1], 15, [3, -1], 18, [7, 1]] \\
C_3' =& \;  [16, [6, -1], 10, [6, 1]] \\
C_4' =& \;  [14, [5, 0], 9, [5, 2]] \\
C_5' =& \;  [3, [5, -1], 17, [5, 1]]
\end{align*}
\end{minipage}

\item 2-factor type $[6, 4, 4, 2, 2, 2]$: two starters \vspace*{-3mm} \\
\begin{minipage}{8cm}
\begin{align*}
C_0 =& \;  [16, [7, 1], 4, [1, 2], 3, [19, 0], 19, [19, -1], 14, [2, -1], 12, [4, 1]] \\
C_1 =& \;  [18, [7, 0], 6, [2, 1], 8, [7, 2], 15, [3, -1]] \\
C_2 =& \;  [9, [8, 0], 17, [3, 2], 1, [1, -1], 2, [7, -1]] \\
C_3 =& \;  [13, [6, 0], 0, [6, 2]] \\
C_4 =& \;  [11, [4, 0], 7, [4, 2]] \\
C_5 =& \;  [5, [5, 0], 10, [5, 2]]
\end{align*}
\end{minipage}

\vspace*{-4mm}
\begin{minipage}{8cm}
\begin{align*}
C_0' =& \;  [1, [3, 1], 17, [1, 0], 16, [4, -1], 12, [19, 1], 19, [19, 2], 2, [1, 1]] \\
C_1' =& \;  [11, [3, 0], 14, [6, -1], 8, [5, -1], 3, [8, 2]] \\
C_2' =& \;  [0, [9, 2], 9, [5, 1], 4, [9, 0], 13, [6, 1]] \\
C_3' =& \;  [7, [2, 0], 5, [2, 2]] \\
C_4' =& \;  [15, [9, -1], 6, [9, 1]] \\
C_5' =& \;  [10, [8, -1], 18, [8, 1]]
\end{align*}
\end{minipage}

\item 2-factor type $[5, 5, 4, 2, 2, 2]$: two starters \vspace*{-3mm} \\
\begin{minipage}{8cm}
\begin{align*}
C_0 =& \;  [7, [3, 1], 10, [19, 1], 19, [19, -1], 18, [3, 2], 15, [8, 0]] \\
C_1 =& \;  [17, [3, 0], 14, [9, -1], 5, [1, -1], 4, [2, -1], 2, [4, 2]] \\
C_2 =& \;  [0, [3, -1], 3, [9, 0], 12, [4, -1], 8, [8, 2]] \\
C_3 =& \;  [1, [5, 0], 6, [5, 2]] \\
C_4 =& \;  [13, [2, 0], 11, [2, 2]] \\
C_5 =& \;  [9, [7, -1], 16, [7, 1]]
\end{align*}
\end{minipage}

\vspace*{-4mm}
\begin{minipage}{8cm}
\begin{align*}
C_0' =& \;  [3, [1, 1], 2, [6, -1], 8, [6, 0], 14, [9, 1], 4, [1, 2]] \\
C_1' =& \;  [18, [2, 1], 1, [6, 2], 7, [4, 0], 11, [19, 2], 19, [19, 0]] \\
C_2' =& \;  [9, [6, 1], 15, [4, 1], 0, [9, 2], 10, [1, 0]] \\
C_3' =& \;  [12, [5, -1], 17, [5, 1]] \\
C_4' =& \;  [13, [7, 0], 6, [7, 2]] \\
C_5' =& \;  [16, [8, -1], 5, [8, 1]]
\end{align*}
\end{minipage}

\item 2-factor type $[7, 3, 3, 3, 2, 2]$: one starter \vspace*{-3mm} \\
\begin{minipage}{8cm}
\begin{align*}
C_0 =& \;  [13, [6, 0], 7, [1, 0], 6, [7, 1], 18, [2, 0], 16, [6, 1], 10, [1, 1], 9, [4, 0]] \\
C_1 =& \;  [5, [9, 1], 14, [7, 0], 2, [3, 0]] \\
C_2 =& \;  [8, [4, 1], 12, [9, 0], 3, [5, 0]] \\
C_3 =& \;  [17, [3, 1], 1, [5, 1], 15, [2, 1]] \\
C_4 =& \;  [11, [8, 0], 0, [8, 1]] \\
C_5 =& \;  [4, [19, 0], 19, [19, 1]]
\end{align*}
\end{minipage}

\newpage

\item 2-factor type $[6, 4, 3, 3, 2, 2]$: one starter \vspace*{-3mm} \\
\begin{minipage}{8cm}
\begin{align*}
C_0 =& \;  [16, [5, 0], 2, [8, 0], 13, [6, 0], 0, [9, 0], 9, [1, 1], 8, [8, 1]] \\
C_1 =& \;  [15, [9, 1], 6, [2, 1], 4, [5, 1], 18, [3, 1]] \\
C_2 =& \;  [7, [4, 0], 3, [2, 0], 1, [6, 1]] \\
C_3 =& \;  [14, [3, 0], 11, [1, 0], 10, [4, 1]] \\
C_4 =& \;  [12, [7, 0], 5, [7, 1]] \\
C_5 =& \;  [17, [19, 0], 19, [19, 1]]
\end{align*}
\end{minipage}

\item 2-factor type $[5, 5, 3, 3, 2, 2]$: one starter \vspace*{-3mm} \\
\begin{minipage}{8cm}
\begin{align*}
C_0 =& \;  [15, [6, 0], 2, [3, 1], 18, [7, 1], 11, [4, 0], 7, [8, 1]] \\
C_1 =& \;  [13, [3, 0], 16, [4, 1], 1, [19, 0], 19, [19, 1], 4, [9, 1]] \\
C_2 =& \;  [0, [6, 1], 6, [8, 0], 17, [2, 0]] \\
C_3 =& \;  [5, [9, 0], 14, [2, 1], 12, [7, 0]] \\
C_4 =& \;  [10, [1, 0], 9, [1, 1]] \\
C_5 =& \;  [3, [5, 0], 8, [5, 1]]
\end{align*}
\end{minipage}

\item 2-factor type $[5, 4, 4, 3, 2, 2]$: one starter \vspace*{-3mm} \\
\hspace*{-2mm}
\begin{minipage}{8cm}
\begin{align*}
C_0 =& \;  [9, [8, 0], 1, [1, 0], 0, [6, 1], 6, [8, 1], 14, [5, 1]] \\
C_1 =& \;  [18, [6, 0], 12, [3, 1], 15, [1, 1], 16, [2, 0]] \\
C_2 =& \;  [5, [3, 0], 2, [9, 0], 11, [4, 1], 7, [2, 1]] \\
C_3 =& \;  [8, [9, 1], 17, [4, 0], 13, [5, 0]] \\
C_4 =& \;  [4, [19, 0], 19, [19, 1]] \\
C_5 =& \;  [10, [7, 0], 3, [7, 1]]
\end{align*}
\end{minipage}

\item 2-factor type $[4, 4, 4, 4, 2, 2]$: one starter \vspace*{-3mm} \\
\hspace*{-6mm}
\begin{minipage}{8cm}
\begin{align*}
C_0 =& \;  [5, [8, 0], 16, [6, 1], 10, [7, 1], 17, [7, 0]] \\
C_1 =& \;  [12, [6, 0], 18, [8, 1], 7, [5, 1], 2, [9, 0]] \\
C_2 =& \;  [9, [5, 0], 4, [1, 1], 3, [9, 1], 13, [4, 0]] \\
C_3 =& \;  [19, [19, 1], 1, [1, 0], 0, [4, 1], 15, [19, 0]] \\
C_4 =& \;  [11, [3, 0], 14, [3, 1]] \\
C_5 =& \;  [6, [2, 0], 8, [2, 1]]
\end{align*}
\end{minipage}

\item 2-factor type $[6, 3, 3, 3, 3, 2]$: two starters \vspace*{-3mm} \\
\hspace*{-1mm}
\begin{minipage}{8cm}
\begin{align*}
C_0 =& \;  [4, [2, 1], 2, [8, 1], 13, [1, 1], 12, [6, -1], 18, [7, 0], 11, [7, 2]] \\
C_1 =& \;  [16, [8, -1], 5, [9, -1], 14, [2, -1]] \\
C_2 =& \;  [8, [2, 0], 6, [5, 2], 1, [7, -1]] \\
C_3 =& \;  [7, [4, -1], 3, [7, 1], 10, [3, -1]] \\
C_4 =& \;  [15, [4, 0], 0, [9, 1], 9, [6, 2]] \\
C_5 =& \;  [17, [19, -1], 19, [19, 1]]
\end{align*}
\end{minipage}

\newpage

\begin{minipage}{8cm}
\begin{align*}
C_0' =& \;  [16, [8, 0], 8, [5, 1], 13, [9, 2], 4, [6, 0], 10, [4, 1], 14, [2, 2]] \\
C_1' =& \;  [19, [19, 2], 3, [3, 1], 0, [19, 0]] \\
C_2' =& \;  [1, [8, 2], 9, [3, 0], 6, [5, -1]] \\
C_3' =& \;  [15, [6, 1], 2, [3, 2], 5, [9, 0]] \\
C_4' =& \;  [11, [4, 2], 7, [5, 0], 12, [1, -1]] \\
C_5' =& \;  [18, [1, 0], 17, [1, 2]]
\end{align*}
\end{minipage}

\item 2-factor type $[5, 4, 3, 3, 3, 2]$: two starters \vspace*{-3mm} \\
\begin{minipage}{8cm}
\begin{align*}
C_0 =& \;  [11, [4, 2], 7, [1, 0], 8, [6, 2], 14, [7, -1], 2, [9, 0]] \\
C_1 =& \;  [16, [19, 0], 19, [19, 2], 10, [2, -1], 12, [4, -1]] \\
C_2 =& \;  [15, [9, 1], 6, [6, 1], 0, [4, 1]] \\
C_3 =& \;  [17, [6, -1], 4, [5, 1], 18, [1, 1]] \\
C_4 =& \;  [13, [7, 1], 1, [8, 2], 9, [4, 0]] \\
C_5 =& \;  [3, [2, 0], 5, [2, 2]]
\end{align*}
\end{minipage}

\vspace*{-4mm}
\begin{minipage}{8cm}
\begin{align*}
C_0' =& \;  [18, [19, 1], 19, [19, -1], 8, [3, -1], 5, [9, -1], 15, [3, 1]] \\
C_1' =& \;  [16, [7, 0], 9, [5, 2], 4, [3, 0], 7, [9, 2]] \\
C_2' =& \;  [14, [2, 1], 12, [5, 0], 17, [3, 2]] \\
C_3' =& \;  [1, [1, 2], 0, [6, 0], 6, [5, -1]] \\
C_4' =& \;  [3, [7, 2], 10, [1, -1], 11, [8, 0]] \\
C_5' =& \;  [13, [8, -1], 2, [8, 1]]
\end{align*}
\end{minipage}

\item 2-factor type $[4, 4, 4, 3, 3, 2]$: two starters \vspace*{-3mm} \\
\hspace*{-3mm}
\begin{minipage}{8cm}
\begin{align*}
C_0 =& \;  [16, [3, 0], 13, [8, 2], 2, [1, 0], 3, [6, 2]] \\
C_1 =& \;  [4, [9, 2], 14, [3, 1], 11, [2, 1], 9, [5, 0]] \\
C_2 =& \;  [7, [9, 0], 17, [5, -1], 12, [19, 2], 19, [19, 1]] \\
C_3 =& \;  [10, [5, 2], 15, [7, 1], 8, [2, 0]] \\
C_4 =& \;  [0, [6, 1], 6, [7, -1], 18, [1, 1]] \\
C_5 =& \;  [5, [4, -1], 1, [4, 1]]
\end{align*}
\end{minipage}
\begin{minipage}{8cm}
\begin{align*}
C_0' =& \;  [12, [1, -1], 13, [6, 0], 7, [2, -1], 5, [7, 2]] \\
C_1' =& \;  [2, [6, -1], 15, [1, 2], 14, [3, -1], 17, [4, 0]] \\
C_2' =& \;  [4, [19, -1], 19, [19, 0], 3, [5, 1], 8, [4, 2]] \\
C_3' =& \;  [1, [2, 2], 18, [8, 0], 10, [9, -1]] \\
C_4' =& \;  [6, [3, 2], 9, [7, 0], 16, [9, 1]] \\
C_5' =& \;  [0, [8, -1], 11, [8, 1]]
\end{align*}
\end{minipage}

\item 2-factor type $[5, 3, 3, 3, 3, 3]$: one starter \vspace*{-3mm} \\
\begin{minipage}{8cm}
\begin{align*}
C_0 =& \;  [1, [3, 0], 4, [6, 0], 10, [7, 0], 3, [9, 1], 12, [8, 0]] \\
C_1 =& \;  [5, [7, 1], 17, [8, 1], 9, [4, 1]] \\
C_2 =& \;  [14, [19, 0], 19, [19, 1], 16, [2, 0]] \\
C_3 =& \;  [18, [3, 1], 2, [2, 1], 0, [1, 1]] \\
C_4 =& \;  [15, [9, 0], 6, [5, 1], 11, [4, 0]] \\
C_5 =& \;  [8, [5, 0], 13, [6, 1], 7, [1, 0]]
\end{align*}
\end{minipage}

\newpage

\item 2-factor type $[4, 4, 3, 3, 3, 3]$: one starter \vspace*{-3mm} \\
\hspace*{-8mm}
\begin{minipage}{8cm}
\begin{align*}
C_0 =& \;  [7, [2, 1], 5, [9, 1], 15, [5, 0], 10, [3, 0]] \\
C_1 =& \;  [3, [1, 0], 4, [9, 0], 14, [5, 1], 9, [6, 1]] \\
C_2 =& \;  [12, [7, 1], 0, [1, 1], 1, [8, 1]] \\
C_3 =& \;  [2, [3, 1], 18, [7, 0], 6, [4, 0]] \\
C_4 =& \;  [17, [4, 1], 13, [2, 0], 11, [6, 0]] \\
C_5 =& \;  [8, [8, 0], 16, [19, 0], 19, [19, 1]]
\end{align*}
\end{minipage}

\item 2-factor type $[12, 2, 2, 2, 2]$: one starter \vspace*{-3mm} \\
\begin{minipage}{8cm}
\begin{align*}
C_0 =& \;  [15, [8, 0], 7, [19, 1], 19, [19, 0], 4, [5, 0], 9, [4, 0], 13, [8, 1], 2, [3, 0], 18, [4, 1], 3, [3, 1], 6, [5, 1], 1, [7, 1], \\ & \; 8, [7, 0]] \\
C_1 =& \;  [11, [1, 0], 12, [1, 1]] \\
C_2 =& \;  [5, [9, 0], 14, [9, 1]] \\
C_3 =& \;  [0, [2, 0], 17, [2, 1]] \\
C_4 =& \;  [10, [6, 0], 16, [6, 1]]
\end{align*}
\end{minipage}

\item 2-factor type $[11, 3, 2, 2, 2]$: two starters \vspace*{-3mm} \\
\begin{minipage}{8cm}
\begin{align*}
C_0 =& \;  [18, [5, 1], 13, [3, -1], 16, [8, 0], 5, [19, 2], 19, [19, 1], 11, [6, 0], 17, [5, 2], 3, [7, 1], 15, [7, 0], 8, [2, 1], 10, \\ & \; [8, 2]] \\
C_1 =& \;  [14, [5, -1], 9, [9, 2], 0, [5, 0]] \\
C_2 =& \;  [7, [6, -1], 1, [6, 1]] \\
C_3 =& \;  [2, [9, -1], 12, [9, 1]] \\
C_4 =& \;  [6, [2, 0], 4, [2, 2]]
\end{align*}
\end{minipage}

\vspace*{-4mm}
\begin{minipage}{8cm}
\begin{align*}
C_0' =& \;  [0, [3, 1], 16, [9, 0], 6, [4, 2], 10, [8, 1], 2, [4, 0], 17, [7, 2], 5, [19, 0], 19, [19, -1], 11, [2, -1], 9, [8, -1], 1, \\ & \; [1, 2]] \\
C_1' =& \;  [15, [7, -1], 8, [6, 2], 14, [1, 0]] \\
C_2' =& \;  [3, [4, -1], 18, [4, 1]] \\
C_3' =& \;  [7, [3, 0], 4, [3, 2]] \\
C_4' =& \;  [13, [1, -1], 12, [1, 1]]
\end{align*}
\end{minipage}

\item 2-factor type $[10, 4, 2, 2, 2]$: two starters \vspace*{-3mm} \\
\begin{minipage}{8cm}
\begin{align*}
C_0 =& \;  [18, [7, 0], 11, [8, 2], 0, [5, 1], 14, [7, -1], 2, [5, 0], 16, [8, 1], 5, [1, 2], 4, [9, -1], 13, [19, -1], 19, [19, 1]] \\
C_1 =& \;  [10, [9, 2], 1, [8, 0], 12, [4, 2], 8, [2, 0]] \\
C_2 =& \;  [17, [2, -1], 15, [2, 1]] \\
C_3 =& \;  [9, [6, -1], 3, [6, 1]] \\
C_4 =& \;  [6, [1, -1], 7, [1, 1]]
\end{align*}
\end{minipage}

\vspace*{-4mm}
\begin{minipage}{8cm}
\begin{align*}
C_0' =& \;  [6, [6, 0], 12, [2, 2], 10, [9, 0], 0, [6, 2], 13, [9, 1], 4, [19, 0], 19, [19, 2], 14, [7, 1], 7, [4, 0], 11, [5, 2]] \\
C_1' =& \;  [16, [8, -1], 8, [5, -1], 3, [7, 2], 15, [1, 0]] \\
C_2' =& \;  [17, [3, -1], 1, [3, 1]] \\
C_3' =& \;  [2, [3, 0], 18, [3, 2]] \\
C_4' =& \;  [9, [4, -1], 5, [4, 1]]
\end{align*}
\end{minipage}

\item 2-factor type $[9, 5, 2, 2, 2]$: two starters \vspace*{-3mm} \\
\begin{minipage}{8cm}
\begin{align*}
C_0 =& \;  [4, [1, 2], 3, [5, -1], 8, [9, 0], 17, [9, 2], 7, [6, 0], 13, [6, 1], 0, [7, 2], 12, [9, -1], 2, [2, 0]] \\
C_1 =& \;  [10, [5, 1], 5, [8, 1], 16, [1, 0], 15, [6, 2], 9, [1, -1]] \\
C_2 =& \;  [14, [4, -1], 18, [4, 1]] \\
C_3 =& \;  [1, [19, -1], 19, [19, 1]] \\
C_4 =& \;  [6, [5, 0], 11, [5, 2]]
\end{align*}
\end{minipage}

\vspace*{-4mm}
\begin{minipage}{8cm}
\begin{align*}
C_0' =& \;  [3, [4, 0], 18, [19, 2], 19, [19, 0], 5, [9, 1], 14, [3, 2], 17, [8, 0], 6, [6, -1], 0, [1, 1], 1, [2, 2]] \\
C_1' =& \;  [8, [8, -1], 16, [7, 0], 4, [8, 2], 15, [3, 0], 12, [4, 2]] \\
C_2' =& \;  [13, [2, -1], 11, [2, 1]] \\
C_3' =& \;  [2, [7, -1], 9, [7, 1]] \\
C_4' =& \;  [7, [3, -1], 10, [3, 1]]
\end{align*}
\end{minipage}

\item 2-factor type $[8, 6, 2, 2, 2]$: two starters \vspace*{-3mm} \\
\begin{minipage}{8cm}
\begin{align*}
C_0 =& \;  [7, [4, 0], 11, [6, -1], 5, [7, 1], 12, [4, 2], 16, [4, -1], 1, [6, 0], 14, [7, 2], 2, [5, -1]] \\
C_1 =& \;  [4, [19, 0], 19, [19, 2], 6, [6, 1], 0, [8, 0], 8, [9, 1], 17, [6, 2]] \\
C_2 =& \;  [15, [3, -1], 18, [3, 1]] \\
C_3 =& \;  [9, [1, 0], 10, [1, 2]] \\
C_4 =& \;  [13, [9, 0], 3, [9, 2]]
\end{align*}
\end{minipage}

\vspace*{-4mm}
\begin{minipage}{8cm}
\begin{align*}
C_0' =& \;  [7, [19, 1], 19, [19, -1], 13, [7, -1], 6, [2, 1], 8, [3, 2], 11, [7, 0], 4, [8, 2], 12, [5, 0]] \\
C_1' =& \;  [0, [5, 1], 5, [9, -1], 15, [5, 2], 1, [3, 0], 17, [4, 1], 2, [2, -1]] \\
C_2' =& \;  [3, [8, -1], 14, [8, 1]] \\
C_3' =& \;  [18, [2, 0], 16, [2, 2]] \\
C_4' =& \;  [10, [1, -1], 9, [1, 1]]
\end{align*}
\end{minipage}

\item 2-factor type $[7, 7, 2, 2, 2]$: two starters \vspace*{-3mm} \\
\begin{minipage}{8cm}
\begin{align*}
C_0 =& \;  [13, [7, 2], 1, [8, 0], 9, [2, 1], 11, [19, 1], 19, [19, 2], 16, [8, 1], 8, [5, 0]] \\
C_1 =& \;  [6, [4, -1], 10, [3, 0], 7, [4, 2], 3, [1, 0], 4, [8, 2], 15, [2, -1], 17, [8, -1]] \\
C_2 =& \;  [14, [5, -1], 0, [5, 1]] \\
C_3 =& \;  [18, [3, -1], 2, [3, 1]] \\
C_4 =& \;  [5, [7, -1], 12, [7, 1]]
\end{align*}
\end{minipage}

\vspace*{-4mm}
\begin{minipage}{8cm}
\begin{align*}
C_0' =& \;  [13, [1, -1], 14, [4, 0], 10, [5, 2], 15, [7, 0], 3, [4, 1], 7, [9, 1], 16, [3, 2]] \\
C_1' =& \;  [6, [2, 2], 4, [19, -1], 19, [19, 0], 1, [1, 2], 0, [9, -1], 9, [1, 1], 8, [2, 0]] \\
C_2' =& \;  [18, [6, 0], 5, [6, 2]] \\
C_3' =& \;  [17, [6, -1], 11, [6, 1]] \\
C_4' =& \;  [12, [9, 0], 2, [9, 2]]
\end{align*}
\end{minipage}

\item 2-factor type $[10, 3, 3, 2, 2]$: one starter \vspace*{-3mm} \\
\begin{minipage}{8cm}
\begin{align*}
C_0 =& \;  [6, [5, 1], 1, [2, 1], 18, [8, 1], 10, [2, 0], 8, [1, 1], 9, [6, 0], 3, [7, 1], 15, [9, 0], 5, [5, 0], 0, [6, 1]] \\
C_1 =& \;  [13, [8, 0], 2, [9, 1], 12, [1, 0]] \\
C_2 =& \;  [4, [19, 1], 19, [19, 0], 16, [7, 0]] \\
C_3 =& \;  [11, [4, 0], 7, [4, 1]] \\
C_4 =& \;  [17, [3, 0], 14, [3, 1]]
\end{align*}
\end{minipage}

\item 2-factor type $[9, 4, 3, 2, 2]$: one starter \vspace*{-3mm} \\
\begin{minipage}{8cm}
\begin{align*}
C_0 =& \;  [19, [19, 1], 5, [7, 0], 17, [7, 1], 10, [3, 1], 13, [2, 1], 11, [1, 1], 12, [6, 0], 18, [2, 0], 16, [19, 0]] \\
C_1 =& \;  [8, [4, 0], 4, [3, 0], 7, [4, 1], 3, [5, 1]] \\
C_2 =& \;  [15, [5, 0], 1, [1, 0], 2, [6, 1]] \\
C_3 =& \;  [9, [9, 0], 0, [9, 1]] \\
C_4 =& \;  [6, [8, 0], 14, [8, 1]]
\end{align*}
\end{minipage}

\item 2-factor type $[8, 5, 3, 2, 2]$: one starter \vspace*{-3mm} \\
\begin{minipage}{8cm}
\begin{align*}
C_0 =& \;  [17, [1, 0], 18, [6, 0], 5, [4, 1], 9, [7, 1], 16, [3, 1], 13, [6, 1], 7, [5, 0], 2, [4, 0]] \\
C_1 =& \;  [0, [8, 1], 11, [1, 1], 10, [7, 0], 3, [5, 1], 8, [8, 0]] \\
C_2 =& \;  [1, [19, 1], 19, [19, 0], 4, [3, 0]] \\
C_3 =& \;  [6, [9, 0], 15, [9, 1]] \\
C_4 =& \;  [14, [2, 0], 12, [2, 1]]
\end{align*}
\end{minipage}

\item 2-factor type $[7, 6, 3, 2, 2]$: one starter \vspace*{-3mm} \\
\begin{minipage}{8cm}
\begin{align*}
C_0 =& \;  [12, [6, 0], 18, [8, 0], 7, [9, 0], 16, [7, 0], 4, [4, 1], 8, [9, 1], 17, [5, 1]] \\
C_1 =& \;  [15, [1, 0], 14, [4, 0], 10, [5, 0], 5, [6, 1], 11, [2, 0], 13, [2, 1]] \\
C_2 =& \;  [9, [7, 1], 2, [1, 1], 1, [8, 1]] \\
C_3 =& \;  [6, [19, 0], 19, [19, 1]] \\
C_4 =& \;  [0, [3, 0], 3, [3, 1]]
\end{align*}
\end{minipage}

\item 2-factor type $[8, 4, 4, 2, 2]$: one starter \vspace*{-3mm} \\
\begin{minipage}{8cm}
\begin{align*}
C_0 =& \;  [3, [2, 0], 1, [1, 0], 0, [9, 1], 9, [9, 0], 18, [5, 1], 13, [1, 1], 12, [4, 1], 8, [5, 0]] \\
C_1 =& \;  [6, [2, 1], 4, [8, 1], 15, [6, 0], 2, [4, 0]] \\
C_2 =& \;  [10, [7, 1], 17, [7, 0], 5, [8, 0], 16, [6, 1]] \\
C_3 =& \;  [7, [19, 0], 19, [19, 1]] \\
C_4 =& \;  [11, [3, 0], 14, [3, 1]]
\end{align*}
\end{minipage}

\item 2-factor type $[7, 5, 4, 2, 2]$: one starter \vspace*{-3mm} \\
\begin{minipage}{8cm}
\begin{align*}
C_0 =& \;  [8, [3, 1], 11, [19, 0], 19, [19, 1], 4, [5, 1], 9, [7, 0], 2, [4, 1], 17, [9, 1]] \\
C_1 =& \;  [7, [3, 0], 10, [9, 0], 1, [1, 1], 0, [4, 0], 15, [8, 0]] \\
C_2 =& \;  [6, [1, 0], 5, [8, 1], 13, [5, 0], 18, [7, 1]] \\
C_3 =& \;  [3, [6, 0], 16, [6, 1]] \\
C_4 =& \;  [14, [2, 0], 12, [2, 1]]
\end{align*}
\end{minipage}

\item 2-factor type $[6, 6, 4, 2, 2]$: one starter \vspace*{-3mm} \\
\begin{minipage}{8cm}
\begin{align*}
C_0 =& \;  [5, [5, 0], 0, [9, 1], 10, [6, 0], 16, [8, 0], 8, [9, 0], 18, [6, 1]] \\
C_1 =& \;  [15, [19, 1], 19, [19, 0], 9, [5, 1], 4, [3, 1], 7, [4, 1], 11, [4, 0]] \\
C_2 =& \;  [6, [7, 0], 13, [7, 1], 1, [3, 0], 17, [8, 1]] \\
C_3 =& \;  [3, [1, 0], 2, [1, 1]] \\
C_4 =& \;  [12, [2, 0], 14, [2, 1]]
\end{align*}
\end{minipage}

\newpage

\item 2-factor type $[6, 5, 5, 2, 2]$: one starter \vspace*{-3mm} \\
\begin{minipage}{8cm}
\begin{align*}
C_0 =& \;  [10, [8, 0], 2, [2, 1], 0, [6, 1], 6, [8, 1], 14, [1, 1], 15, [5, 0]] \\
C_1 =& \;  [9, [4, 1], 5, [3, 1], 8, [3, 0], 11, [7, 0], 4, [5, 1]] \\
C_2 =& \;  [13, [1, 0], 12, [6, 0], 18, [4, 0], 3, [2, 0], 1, [7, 1]] \\
C_3 =& \;  [7, [9, 0], 17, [9, 1]] \\
C_4 =& \;  [16, [19, 0], 19, [19, 1]]
\end{align*}
\end{minipage}

\item 2-factor type $[9, 3, 3, 3, 2]$: two starters \vspace*{-3mm} \\
\begin{minipage}{8cm}
\begin{align*}
C_0 =& \;  [19, [19, 2], 18, [8, -1], 7, [6, -1], 13, [3, 1], 10, [9, 0], 0, [5, 2], 5, [1, 1], 4, [4, 0], 8, [19, 1]] \\
C_1 =& \;  [1, [5, 0], 6, [9, -1], 16, [4, 2]] \\
C_2 =& \;  [15, [6, 1], 9, [8, 0], 17, [2, 2]] \\
C_3 =& \;  [11, [9, 1], 2, [1, 0], 3, [8, 2]] \\
C_4 =& \;  [12, [2, -1], 14, [2, 1]]
\end{align*}
\end{minipage}

\vspace*{-4mm}
\begin{minipage}{8cm}
\begin{align*}
C_0' =& \;  [4, [19, -1], 19, [19, 0], 16, [7, -1], 9, [3, -1], 6, [4, -1], 2, [1, -1], 1, [9, 2], 10, [7, 0], 17, [6, 2]] \\
C_1' =& \;  [7, [5, -1], 12, [2, 0], 14, [7, 2]] \\
C_2' =& \;  [15, [4, 1], 11, [8, 1], 3, [7, 1]] \\
C_3' =& \;  [0, [6, 0], 13, [5, 1], 18, [1, 2]] \\
C_4' =& \;  [8, [3, 0], 5, [3, 2]]
\end{align*}
\end{minipage}

\item 2-factor type $[8, 4, 3, 3, 2]$: two starters \vspace*{-3mm} \\
\begin{minipage}{8cm}
\begin{align*}
C_0 =& \;  [7, [1, 0], 6, [5, 1], 11, [8, 1], 0, [5, 2], 14, [1, 1], 13, [19, 0], 19, [19, -1], 15, [8, 2]] \\
C_1 =& \;  [8, [9, 1], 18, [5, 0], 4, [2, 2], 2, [6, -1]] \\
C_2 =& \;  [12, [4, 0], 16, [7, -1], 9, [3, 2]] \\
C_3 =& \;  [17, [3, 0], 1, [4, 2], 5, [7, 1]] \\
C_4 =& \;  [10, [7, 0], 3, [7, 2]]
\end{align*}
\end{minipage}

\vspace*{-4mm}
\begin{minipage}{8cm}
\begin{align*}
C_0' =& \;  [10, [1, -1], 11, [5, -1], 16, [9, 0], 7, [1, 2], 6, [4, 1], 2, [8, 0], 13, [9, -1], 4, [6, 2]] \\
C_1' =& \;  [15, [19, 2], 19, [19, 1], 0, [2, 0], 17, [2, -1]] \\
C_2' =& \;  [9, [6, 0], 3, [4, -1], 18, [9, 2]] \\
C_3' =& \;  [12, [2, 1], 14, [6, 1], 1, [8, -1]] \\
C_4' =& \;  [8, [3, -1], 5, [3, 1]]
\end{align*}
\end{minipage}

\item 2-factor type $[7, 5, 3, 3, 2]$: two starters \vspace*{-3mm} \\
\begin{minipage}{8cm}
\begin{align*}
C_0 =& \;  [3, [1, 0], 4, [4, 1], 8, [7, 2], 15, [9, 1], 6, [8, -1], 14, [1, 1], 13, [9, -1]] \\
C_1 =& \;  [19, [19, 0], 2, [8, 1], 10, [1, -1], 9, [2, 2], 11, [19, -1]] \\
C_2 =& \;  [7, [2, 1], 5, [6, 0], 18, [8, 2]] \\
C_3 =& \;  [0, [3, 2], 16, [4, 0], 12, [7, 1]] \\
C_4 =& \;  [1, [3, -1], 17, [3, 1]]
\end{align*}
\end{minipage}

\vspace*{-4mm}
\begin{minipage}{8cm}
\begin{align*}
C_0' =& \;  [17, [2, -1], 15, [1, 2], 14, [7, -1], 2, [3, 0], 18, [6, 1], 5, [5, 2], 10, [7, 0]] \\
C_1' =& \;  [9, [6, -1], 3, [9, 2], 12, [5, 0], 7, [4, 2], 11, [2, 0]] \\
C_2' =& \;  [13, [6, 2], 0, [4, -1], 4, [9, 0]] \\
C_3' =& \;  [16, [8, 0], 8, [19, 2], 19, [19, 1]] \\
C_4' =& \;  [6, [5, -1], 1, [5, 1]]
\end{align*}
\end{minipage}

\item 2-factor type $[6, 6, 3, 3, 2]$: two starters \vspace*{-3mm} \\
\begin{minipage}{8cm}
\begin{align*}
C_0 =& \;  [4, [9, -1], 14, [8, 2], 3, [4, 1], 18, [7, 0], 11, [2, 2], 9, [5, 0]] \\
C_1 =& \;  [13, [7, 1], 1, [7, 2], 8, [4, -1], 12, [2, 0], 10, [5, 2], 5, [8, 0]] \\
C_2 =& \;  [2, [6, 2], 15, [4, 0], 0, [2, 1]] \\
C_3 =& \;  [17, [1, -1], 16, [19, 2], 19, [19, 0]] \\
C_4 =& \;  [6, [1, 0], 7, [1, 2]]
\end{align*}
\end{minipage}

\vspace*{-4mm}
\begin{minipage}{8cm}
\begin{align*}
C_0' =& \;  [10, [4, 2], 6, [5, 1], 1, [3, 0], 17, [6, 1], 4, [5, -1], 18, [8, -1]] \\
C_1' =& \;  [9, [7, -1], 2, [19, 1], 19, [19, -1], 15, [1, 1], 16, [3, 1], 0, [9, 1]] \\
C_2' =& \;  [7, [6, -1], 13, [8, 1], 5, [2, -1]] \\
C_3' =& \;  [14, [3, -1], 11, [3, 2], 8, [6, 0]] \\
C_4' =& \;  [3, [9, 0], 12, [9, 2]]
\end{align*}
\end{minipage}

\item 2-factor type $[7, 4, 4, 3, 2]$: two starters \vspace*{-3mm} \\
\begin{minipage}{8cm}
\begin{align*}
C_0 =& \;  [16, [1, 0], 17, [7, -1], 10, [19, 1], 19, [19, -1], 4, [9, 1], 13, [6, 2], 0, [3, 1]] \\
C_1 =& \;  [14, [7, 1], 7, [2, 1], 5, [6, -1], 11, [3, -1]] \\
C_2 =& \;  [8, [4, 1], 12, [9, 2], 2, [1, -1], 3, [5, 0]] \\
C_3 =& \;  [6, [9, -1], 15, [6, 0], 9, [3, 2]] \\
C_4 =& \;  [1, [2, 0], 18, [2, 2]]
\end{align*}
\end{minipage}

\vspace*{-4mm}
\begin{minipage}{8cm}
\begin{align*}
C_0' =& \;  [17, [1, 2], 18, [4, -1], 3, [6, 1], 16, [8, 0], 8, [8, 2], 0, [9, 0], 9, [8, 1]] \\
C_1' =& \;  [6, [8, -1], 14, [1, 1], 13, [7, 0], 1, [5, 2]] \\
C_2' =& \;  [12, [7, 2], 5, [5, -1], 10, [3, 0], 7, [5, 1]] \\
C_3' =& \;  [19, [19, 0], 4, [2, -1], 2, [19, 2]] \\
C_4' =& \;  [15, [4, 0], 11, [4, 2]]
\end{align*}
\end{minipage}

\item 2-factor type $[6, 5, 4, 3, 2]$: two starters \vspace*{-3mm} \\
\begin{minipage}{8cm}
\begin{align*}
C_0 =& \;  [6, [4, 1], 10, [19, 2], 19, [19, 0], 13, [6, 2], 7, [7, 1], 0, [6, 0]] \\
C_1 =& \;  [5, [4, -1], 1, [3, 1], 4, [5, 2], 18, [3, -1], 2, [3, 0]] \\
C_2 =& \;  [12, [3, 2], 15, [7, 0], 3, [8, 2], 11, [1, 0]] \\
C_3 =& \;  [14, [5, 1], 9, [7, -1], 16, [2, 1]] \\
C_4 =& \;  [8, [9, -1], 17, [9, 1]]
\end{align*}
\end{minipage}

\vspace*{-4mm}
\begin{minipage}{8cm}
\begin{align*}
C_0' =& \;  [1, [1, 2], 0, [4, 0], 4, [1, 1], 5, [2, -1], 3, [7, 2], 15, [5, 0]] \\
C_1' =& \;  [18, [5, -1], 13, [8, 1], 2, [6, 1], 8, [1, -1], 7, [8, -1]] \\
C_2' =& \;  [12, [2, 0], 14, [4, 2], 10, [19, -1], 19, [19, 1]] \\
C_3' =& \;  [17, [8, 0], 9, [2, 2], 11, [6, -1]] \\
C_4' =& \;  [6, [9, 0], 16, [9, 2]]
\end{align*}
\end{minipage}

\item 2-factor type $[5, 5, 5, 3, 2]$: two starters \vspace*{-3mm} \\
\begin{minipage}{8cm}
\begin{align*}
C_0 =& \;  [5, [1, -1], 6, [4, 1], 2, [5, 0], 7, [2, 1], 9, [4, 2]] \\
C_1 =& \;  [10, [7, -1], 17, [1, 1], 16, [2, 0], 14, [5, 1], 0, [9, 2]] \\
C_2 =& \;  [13, [5, -1], 8, [9, -1], 18, [19, 1], 19, [19, 2], 4, [9, 0]] \\
C_3 =& \;  [1, [8, 1], 12, [3, 0], 15, [5, 2]] \\
C_4 =& \;  [11, [8, 0], 3, [8, 2]]
\end{align*}
\end{minipage}

\newpage

\begin{minipage}{8cm}
\begin{align*}
C_0' =& \;  [14, [19, 0], 19, [19, -1], 10, [1, 2], 11, [6, -1], 17, [3, 1]] \\
C_1' =& \;  [3, [4, 0], 18, [3, 2], 2, [2, -1], 4, [1, 0], 5, [2, 2]] \\
C_2' =& \;  [1, [8, -1], 9, [4, -1], 13, [3, -1], 16, [9, 1], 7, [6, 1]] \\
C_3' =& \;  [12, [6, 0], 6, [6, 2], 0, [7, 1]] \\
C_4' =& \;  [15, [7, 0], 8, [7, 2]]
\end{align*}
\end{minipage}

\item 2-factor type $[6, 4, 4, 4, 2]$: two starters \vspace*{-3mm} \\
\begin{minipage}{8cm}
\begin{align*}
C_0 =& \;  [15, [19, 0], 19, [19, 2], 12, [8, 0], 4, [3, 2], 1, [1, 0], 2, [6, 2]] \\
C_1 =& \;  [14, [5, -1], 0, [2, 1], 17, [9, -1], 7, [7, -1]] \\
C_2 =& \;  [6, [5, 2], 11, [2, -1], 13, [3, -1], 16, [9, 0]] \\
C_3 =& \;  [18, [6, 1], 5, [4, 2], 9, [6, 0], 3, [4, -1]] \\
C_4 =& \;  [10, [2, 0], 8, [2, 2]]
\end{align*}
\end{minipage}

\vspace*{-4mm}
\begin{minipage}{8cm}
\begin{align*}
C_0' =& \;  [13, [8, 2], 5, [4, 0], 9, [8, -1], 1, [7, 2], 8, [9, 1], 18, [5, 0]] \\
C_1' =& \;  [17, [6, -1], 4, [7, 1], 16, [19, -1], 19, [19, 1]] \\
C_2' =& \;  [10, [9, 2], 0, [3, 0], 3, [3, 1], 6, [4, 1]] \\
C_3' =& \;  [14, [7, 0], 2, [5, 1], 7, [8, 1], 15, [1, 2]] \\
C_4' =& \;  [11, [1, -1], 12, [1, 1]]
\end{align*}
\end{minipage}

\item 2-factor type $[5, 5, 4, 4, 2]$: two starters \vspace*{-3mm} \\
\begin{minipage}{8cm}
\begin{align*}
C_0 =& \;  [0, [8, 1], 8, [1, 1], 9, [4, -1], 5, [4, 2], 1, [1, 0]] \\
C_1 =& \;  [4, [7, 0], 11, [5, 2], 6, [1, -1], 7, [19, -1], 19, [19, 1]] \\
C_2 =& \;  [10, [4, 0], 14, [1, 2], 13, [8, 0], 2, [8, 2]] \\
C_3 =& \;  [18, [2, 1], 16, [4, 1], 12, [3, 0], 15, [3, 2]] \\
C_4 =& \;  [3, [5, -1], 17, [5, 1]]
\end{align*}
\end{minipage}

\vspace*{-4mm}
\begin{minipage}{8cm}
\begin{align*}
C_0' =& \;  [4, [7, 2], 11, [19, 0], 19, [19, 2], 12, [6, 0], 6, [2, -1]] \\
C_1' =& \;  [16, [5, 0], 2, [2, 2], 0, [8, -1], 8, [2, 0], 10, [6, 2]] \\
C_2' =& \;  [18, [3, -1], 15, [7, 1], 3, [6, 1], 9, [9, 1]] \\
C_3' =& \;  [7, [9, -1], 17, [3, 1], 1, [7, -1], 13, [6, -1]] \\
C_4' =& \;  [5, [9, 0], 14, [9, 2]]
\end{align*}
\end{minipage}

\item 2-factor type $[8, 3, 3, 3, 3]$: one starter \vspace*{-3mm} \\
\begin{minipage}{8cm}
\begin{align*}
C_0 =& \;  [13, [5, 0], 8, [8, 0], 0, [7, 1], 7, [8, 1], 15, [1, 0], 14, [3, 0], 11, [6, 0], 17, [4, 0]] \\
C_1 =& \;  [3, [19, 0], 19, [19, 1], 1, [2, 0]] \\
C_2 =& \;  [9, [1, 1], 10, [5, 1], 5, [4, 1]] \\
C_3 =& \;  [4, [2, 1], 6, [9, 0], 16, [7, 0]] \\
C_4 =& \;  [18, [6, 1], 12, [9, 1], 2, [3, 1]]
\end{align*}
\end{minipage}

\newpage

\item 2-factor type $[7, 4, 3, 3, 3]$: one starter \vspace*{-3mm} \\
\begin{minipage}{8cm}
\begin{align*}
C_0 =& \;  [12, [4, 1], 8, [8, 0], 16, [5, 0], 2, [1, 1], 3, [1, 0], 4, [19, 0], 19, [19, 1]] \\
C_1 =& \;  [18, [9, 0], 9, [9, 1], 0, [6, 1], 6, [7, 0]] \\
C_2 =& \;  [7, [7, 1], 14, [4, 0], 10, [3, 0]] \\
C_3 =& \;  [15, [5, 1], 1, [3, 1], 17, [2, 1]] \\
C_4 =& \;  [5, [6, 0], 11, [2, 0], 13, [8, 1]]
\end{align*}
\end{minipage}

\item 2-factor type $[6, 5, 3, 3, 3]$: one starter \vspace*{-3mm} \\
\begin{minipage}{8cm}
\begin{align*}
C_0 =& \;  [12, [7, 1], 0, [6, 0], 13, [9, 1], 4, [3, 0], 1, [5, 1], 6, [6, 1]] \\
C_1 =& \;  [14, [4, 0], 18, [3, 1], 15, [4, 1], 11, [8, 1], 3, [8, 0]] \\
C_2 =& \;  [7, [2, 1], 5, [7, 0], 17, [9, 0]] \\
C_3 =& \;  [2, [19, 1], 19, [19, 0], 16, [5, 0]] \\
C_4 =& \;  [10, [2, 0], 8, [1, 1], 9, [1, 0]]
\end{align*}
\end{minipage}

\item 2-factor type $[6, 4, 4, 3, 3]$: one starter \vspace*{-3mm} \\
\begin{minipage}{8cm}
\begin{align*}
C_0 =& \;  [16, [7, 0], 9, [8, 0], 1, [2, 1], 18, [5, 0], 13, [2, 0], 15, [1, 1]] \\
C_1 =& \;  [4, [6, 1], 10, [5, 1], 5, [3, 0], 8, [4, 0]] \\
C_2 =& \;  [6, [4, 1], 2, [19, 0], 19, [19, 1], 0, [6, 0]] \\
C_3 =& \;  [17, [3, 1], 14, [7, 1], 7, [9, 1]] \\
C_4 =& \;  [12, [1, 0], 11, [8, 1], 3, [9, 0]]
\end{align*}
\end{minipage}

\item 2-factor type $[5, 5, 4, 3, 3]$: one starter \vspace*{-3mm} \\
\begin{minipage}{8cm}
\begin{align*}
C_0 =& \;  [0, [2, 1], 17, [19, 1], 19, [19, 0], 9, [4, 0], 13, [6, 1]] \\
C_1 =& \;  [16, [6, 0], 3, [4, 1], 7, [3, 0], 10, [8, 1], 2, [5, 1]] \\
C_2 =& \;  [8, [7, 1], 1, [5, 0], 6, [7, 0], 18, [9, 1]] \\
C_3 =& \;  [11, [3, 1], 14, [2, 0], 12, [1, 0]] \\
C_4 =& \;  [15, [9, 0], 5, [1, 1], 4, [8, 0]]
\end{align*}
\end{minipage}

\item 2-factor type $[5, 4, 4, 4, 3]$: one starter \vspace*{-3mm} \\
\begin{minipage}{8cm}
\begin{align*}
C_0 =& \;  [11, [9, 1], 1, [1, 1], 2, [6, 1], 8, [7, 1], 15, [4, 1]] \\
C_1 =& \;  [14, [2, 1], 12, [8, 1], 4, [3, 0], 7, [7, 0]] \\
C_2 =& \;  [16, [2, 0], 18, [4, 0], 3, [5, 0], 17, [1, 0]] \\
C_3 =& \;  [10, [19, 0], 19, [19, 1], 13, [8, 0], 5, [5, 1]] \\
C_4 =& \;  [9, [3, 1], 6, [6, 0], 0, [9, 0]]
\end{align*}
\end{minipage}

\item 2-factor type $[14, 2, 2, 2]$: two starters \vspace*{-3mm} \\
\begin{minipage}{8cm}
\begin{align*}
C_0 =& \;  [7, [2, 2], 5, [19, 0], 19, [19, 2], 12, [9, 0], 3, [4, 2], 18, [1, 0], 17, [3, 2], 14, [5, -1], 0, [6, -1], 6, [2, 0], 4, \\ & \; [4, 1], 8, [8, 1], 16, [1, 2], 15, [8, 0]] \\
C_1 =& \;  [10, [3, -1], 13, [3, 1]] \\
C_2 =& \;  [11, [2, -1], 9, [2, 1]] \\
C_3 =& \;  [1, [1, -1], 2, [1, 1]]
\end{align*}
\end{minipage}

\newpage

\begin{minipage}{8cm}
\begin{align*}
C_0' =& \;  [18, [4, -1], 14, [19, 1], 19, [19, -1], 12, [5, 2], 7, [4, 0], 11, [6, 1], 17, [9, 2], 8, [8, -1], 16, [7, 0], 4, [5, 1], \\ & \; 9, [7, 2], 2, [3, 0], 5, [8, 2], 13, [5, 0]] \\
C_1' =& \;  [6, [6, 0], 0, [6, 2]] \\
C_2' =& \;  [15, [7, -1], 3, [7, 1]] \\
C_3' =& \;  [10, [9, -1], 1, [9, 1]]
\end{align*}
\end{minipage}

\item 2-factor type $[13, 3, 2, 2]$: one starter \vspace*{-3mm} \\
\begin{minipage}{8cm}
\begin{align*}
C_0 =& \;  [11, [3, 0], 14, [1, 1], 13, [5, 1], 18, [4, 1], 3, [3, 1], 0, [7, 0], 7, [8, 1], 15, [5, 0], 1, [4, 0], 5, [7, 1], 12, [9, 0], \\ & \; 2, [19, 1], 19, [19, 0]] \\
C_1 =& \;  [9, [1, 0], 8, [9, 1], 17, [8, 0]] \\
C_2 =& \;  [10, [6, 0], 16, [6, 1]] \\
C_3 =& \;  [4, [2, 0], 6, [2, 1]]
\end{align*}
\end{minipage}

\item 2-factor type $[12, 4, 2, 2]$: one starter \vspace*{-3mm} \\
\begin{minipage}{8cm}
\begin{align*}
C_0 =& \;  [10, [2, 0], 12, [3, 1], 9, [6, 0], 3, [5, 0], 8, [2, 1], 6, [4, 0], 2, [5, 1], 16, [8, 0], 5, [19, 1], 19, [19, 0], 11, [6, 1], \\ & \; 17, [7, 1]] \\
C_1 =& \;  [14, [7, 0], 7, [8, 1], 15, [3, 0], 18, [4, 1]] \\
C_2 =& \;  [4, [9, 0], 13, [9, 1]] \\
C_3 =& \;  [1, [1, 0], 0, [1, 1]]
\end{align*}
\end{minipage}

\item 2-factor type $[11, 5, 2, 2]$: one starter \vspace*{-3mm} \\
\begin{minipage}{8cm}
\begin{align*}
C_0 =& \;  [11, [9, 1], 2, [7, 1], 14, [8, 1], 6, [7, 0], 18, [3, 0], 15, [9, 0], 5, [5, 0], 10, [2, 1], 12, [3, 1], 9, [4, 0], 13, [2, 0]] \\
C_1 =& \;  [17, [6, 0], 4, [4, 1], 8, [8, 0], 16, [6, 1], 3, [5, 1]] \\
C_2 =& \;  [1, [1, 0], 0, [1, 1]] \\
C_3 =& \;  [7, [19, 0], 19, [19, 1]]
\end{align*}
\end{minipage}

\item 2-factor type $[10, 6, 2, 2]$: one starter \vspace*{-3mm} \\
\begin{minipage}{8cm}
\begin{align*}
C_0 =& \;  [4, [8, 0], 12, [9, 0], 3, [8, 1], 14, [7, 0], 2, [4, 0], 17, [1, 0], 16, [4, 1], 1, [2, 1], 18, [19, 0], 19, [19, 1]] \\
C_1 =& \;  [9, [6, 0], 15, [2, 0], 13, [6, 1], 7, [7, 1], 0, [9, 1], 10, [1, 1]] \\
C_2 =& \;  [8, [3, 0], 5, [3, 1]] \\
C_3 =& \;  [6, [5, 0], 11, [5, 1]]
\end{align*}
\end{minipage}

\item 2-factor type $[9, 7, 2, 2]$: one starter \vspace*{-3mm} \\
\begin{minipage}{8cm}
\begin{align*}
C_0 =& \;  [0, [9, 1], 10, [5, 0], 5, [6, 1], 11, [9, 0], 1, [3, 0], 4, [5, 1], 9, [7, 0], 2, [3, 1], 18, [1, 1]] \\
C_1 =& \;  [7, [1, 0], 6, [19, 0], 19, [19, 1], 3, [6, 0], 16, [8, 0], 8, [7, 1], 15, [8, 1]] \\
C_2 =& \;  [13, [4, 0], 17, [4, 1]] \\
C_3 =& \;  [14, [2, 0], 12, [2, 1]]
\end{align*}
\end{minipage}

\item 2-factor type $[8, 8, 2, 2]$: one starter \vspace*{-3mm} \\
\begin{minipage}{8cm}
\begin{align*}
C_0 =& \;  [11, [8, 1], 3, [2, 0], 1, [6, 1], 7, [5, 0], 12, [6, 0], 6, [7, 0], 13, [2, 1], 15, [4, 1]] \\
C_1 =& \;  [18, [19, 0], 19, [19, 1], 0, [3, 0], 16, [7, 1], 9, [5, 1], 14, [4, 0], 10, [8, 0], 2, [3, 1]] \\
C_2 =& \;  [8, [9, 0], 17, [9, 1]] \\
C_3 =& \;  [4, [1, 0], 5, [1, 1]]
\end{align*}
\end{minipage}

\item 2-factor type $[12, 3, 3, 2]$: two starters \vspace*{-3mm} \\
\begin{minipage}{8cm}
\begin{align*}
C_0 =& \;  [7, [5, 2], 12, [2, 0], 10, [4, 2], 14, [1, 0], 15, [8, 1], 4, [7, -1], 16, [3, 2], 13, [7, 0], 6, [4, -1], 2, [6, 2], 8, \\ & \; [9, 0], 17, [9, 1]] \\
C_1 =& \;  [1, [4, 0], 5, [6, -1], 18, [2, 2]] \\
C_2 =& \;  [9, [2, -1], 11, [8, 0], 0, [9, 2]] \\
C_3 =& \;  [19, [19, 0], 3, [19, 2]]
\end{align*}
\end{minipage}

\vspace*{-4mm}
\begin{minipage}{8cm}
\begin{align*}
C_0' =& \;  [16, [4, 1], 1, [5, -1], 15, [7, 2], 8, [9, -1], 17, [5, 0], 3, [1, 2], 4, [1, -1], 5, [6, 0], 11, [8, 2], 0, [1, 1], 18, \\ & \; [5, 1], 13, [3, 0]] \\
C_1' =& \;  [12, [6, 1], 6, [8, -1], 14, [2, 1]] \\
C_2' =& \;  [19, [19, -1], 2, [7, 1], 9, [19, 1]] \\
C_3' =& \;  [7, [3, -1], 10, [3, 1]]
\end{align*}
\end{minipage}

\item 2-factor type $[11, 4, 3, 2]$: two starters \vspace*{-3mm} \\
\begin{minipage}{8cm}
\begin{align*}
C_0 =& \;  [1, [1, 2], 2, [3, 0], 5, [1, -1], 4, [2, 2], 6, [5, -1], 11, [3, 1], 8, [9, 0], 17, [19, 1], 19, [19, 2], 18, [5, 1], 13, \\ & \; [7, 0]] \\
C_1 =& \;  [10, [2, -1], 12, [7, -1], 0, [5, 0], 14, [4, 2]] \\
C_2 =& \;  [9, [7, 2], 16, [6, -1], 3, [6, 0]] \\
C_3 =& \;  [15, [8, -1], 7, [8, 1]]
\end{align*}
\end{minipage}

\vspace*{-4mm}
\begin{minipage}{8cm}
\begin{align*}
C_0' =& \;  [2, [5, 2], 7, [1, 0], 8, [4, -1], 4, [7, 1], 11, [8, 2], 0, [6, 1], 13, [3, -1], 16, [2, 0], 14, [1, 1], 15, [2, 1], 17, \\ & \; [4, 1]] \\
C_1' =& \;  [6, [4, 0], 10, [9, 2], 1, [8, 0], 9, [3, 2]] \\
C_2' =& \;  [18, [6, 2], 5, [19, 0], 19, [19, -1]] \\
C_3' =& \;  [12, [9, -1], 3, [9, 1]]
\end{align*}
\end{minipage}

\item 2-factor type $[10, 5, 3, 2]$: two starters \vspace*{-3mm} \\
\begin{minipage}{8cm}
\begin{align*}
C_0 =& \;  [9, [5, -1], 14, [4, 0], 18, [1, -1], 17, [8, 1], 6, [2, 1], 8, [7, 2], 1, [2, 0], 3, [3, 2], 0, [3, 1], 16, [7, 1]] \\
C_1 =& \;  [4, [3, 0], 7, [5, 2], 12, [3, -1], 15, [5, 0], 10, [6, 2]] \\
C_2 =& \;  [2, [8, -1], 13, [19, 1], 19, [19, -1]] \\
C_3 =& \;  [5, [6, -1], 11, [6, 1]]
\end{align*}
\end{minipage}

\vspace*{-4mm}
\begin{minipage}{8cm}
\begin{align*}
C_0' =& \;  [2, [4, 1], 17, [8, 0], 6, [1, 2], 5, [9, 0], 15, [4, 2], 11, [7, 0], 4, [9, -1], 14, [2, 2], 16, [9, 1], 7, [5, 1]] \\
C_1' =& \;  [3, [4, -1], 18, [9, 2], 8, [2, -1], 10, [1, 1], 9, [6, 0]] \\
C_2' =& \;  [0, [1, 0], 1, [8, 2], 12, [7, -1]] \\
C_3' =& \;  [13, [19, 0], 19, [19, 2]]
\end{align*}
\end{minipage}

\item 2-factor type $[9, 6, 3, 2]$: two starters \vspace*{-3mm} \\
\begin{minipage}{8cm}
\begin{align*}
C_0 =& \;  [13, [9, 0], 4, [1, 1], 5, [5, 2], 0, [7, 1], 12, [8, -1], 1, [2, 1], 18, [4, 0], 3, [7, 2], 10, [3, -1]] \\
C_1 =& \;  [7, [1, 0], 8, [8, 1], 16, [9, 1], 6, [5, 1], 11, [19, 2], 19, [19, 1]] \\
C_2 =& \;  [14, [5, 0], 9, [6, 1], 15, [1, 2]] \\
C_3 =& \;  [17, [4, -1], 2, [4, 1]]
\end{align*}
\end{minipage}

\newpage

\begin{minipage}{8cm}
\begin{align*}
C_0' =& \;  [14, [5, -1], 9, [6, 2], 15, [9, -1], 5, [7, 0], 12, [4, 2], 16, [3, 0], 0, [6, -1], 13, [9, 2], 3, [8, 0]] \\
C_1' =& \;  [8, [6, 0], 2, [8, 2], 10, [7, -1], 17, [19, 0], 19, [19, -1], 11, [3, 2]] \\
C_2' =& \;  [7, [3, 1], 4, [2, -1], 6, [1, -1]] \\
C_3' =& \;  [18, [2, 0], 1, [2, 2]]
\end{align*}
\end{minipage}

\item 2-factor type $[8, 7, 3, 2]$: two starters \vspace*{-3mm} \\
\begin{minipage}{8cm}
\begin{align*}
C_0 =& \;  [19, [19, 0], 17, [4, -1], 13, [8, 2], 5, [9, 0], 15, [5, 1], 1, [2, -1], 18, [1, 1], 0, [19, 2]] \\
C_1 =& \;  [11, [1, 0], 10, [3, 2], 7, [5, -1], 12, [8, 0], 4, [5, 2], 9, [6, 1], 3, [8, -1]] \\
C_2 =& \;  [8, [6, 0], 14, [2, 2], 16, [8, 1]] \\
C_3 =& \;  [2, [4, 0], 6, [4, 2]]
\end{align*}
\end{minipage}

\vspace*{-4mm}
\begin{minipage}{8cm}
\begin{align*}
C_0' =& \;  [10, [6, 2], 16, [6, -1], 3, [7, 0], 15, [7, 2], 8, [9, -1], 17, [1, -1], 18, [5, 0], 13, [3, -1]] \\
C_1' =& \;  [2, [4, 1], 6, [1, 2], 7, [2, 1], 5, [7, -1], 12, [3, 1], 9, [9, 1], 0, [2, 0]] \\
C_2' =& \;  [1, [9, 2], 11, [7, 1], 4, [3, 0]] \\
C_3' =& \;  [14, [19, -1], 19, [19, 1]]
\end{align*}
\end{minipage}

\item 2-factor type $[10, 4, 4, 2]$: two starters \vspace*{-3mm} \\
\begin{minipage}{8cm}
\begin{align*}
C_0 =& \;  [3, [6, 0], 9, [4, -1], 5, [3, 2], 2, [8, 0], 10, [2, -1], 8, [9, 1], 17, [5, -1], 12, [8, 2], 1, [5, 0], 15, [7, 2]] \\
C_1 =& \;  [13, [2, 0], 11, [4, 2], 7, [7, 0], 14, [1, 2]] \\
C_2 =& \;  [6, [19, 2], 19, [19, 0], 18, [5, 1], 4, [2, 1]] \\
C_3 =& \;  [0, [3, -1], 16, [3, 1]]
\end{align*}
\end{minipage}

\vspace*{-4mm}
\begin{minipage}{8cm}
\begin{align*}
C_0' =& \;  [15, [3, 0], 12, [8, -1], 4, [6, 2], 10, [4, 0], 14, [7, -1], 7, [1, -1], 6, [19, 1], 19, [19, -1], 16, [8, 1], 5, [9, 2]] \\
C_1' =& \;  [11, [9, -1], 1, [1, 0], 2, [7, 1], 9, [2, 2]] \\
C_2' =& \;  [17, [1, 1], 18, [4, 1], 3, [5, 2], 8, [9, 0]] \\
C_3' =& \;  [13, [6, -1], 0, [6, 1]]
\end{align*}
\end{minipage}

\item 2-factor type $[9, 5, 4, 2]$: two starters \vspace*{-3mm} \\
\begin{minipage}{8cm}
\begin{align*}
C_0 =& \;  [13, [9, 0], 4, [7, -1], 16, [8, -1], 8, [3, 2], 11, [1, -1], 12, [3, -1], 15, [1, 1], 14, [6, -1], 1, [7, 1]] \\
C_1 =& \;  [3, [3, 0], 6, [8, 2], 17, [4, 0], 2, [3, 1], 5, [2, 2]] \\
C_2 =& \;  [18, [9, -1], 9, [2, -1], 7, [7, 2], 0, [1, 0]] \\
C_3 =& \;  [10, [19, -1], 19, [19, 1]]
\end{align*}
\end{minipage}

\vspace*{-4mm}
\begin{minipage}{8cm}
\begin{align*}
C_0' =& \;  [15, [9, 1], 5, [4, -1], 1, [5, 1], 6, [6, 2], 0, [6, 0], 13, [1, 2], 14, [2, 1], 12, [8, 1], 4, [8, 0]] \\
C_1' =& \;  [16, [5, 2], 11, [7, 0], 18, [9, 2], 8, [2, 0], 10, [6, 1]] \\
C_2' =& \;  [17, [5, 0], 3, [4, 1], 7, [5, -1], 2, [4, 2]] \\
C_3' =& \;  [19, [19, 0], 9, [19, 2]]
\end{align*}
\end{minipage}

\item 2-factor type $[8, 6, 4, 2]$: two starters \vspace*{-3mm} \\
\begin{minipage}{8cm}
\begin{align*}
C_0 =& \;  [0, [1, 2], 1, [3, 1], 17, [7, 0], 5, [8, 2], 16, [19, 0], 19, [19, 2], 13, [2, 0], 15, [4, 1]] \\
C_1 =& \;  [18, [7, -1], 11, [1, -1], 10, [8, -1], 2, [4, 2], 6, [3, -1], 9, [9, 0]] \\
C_2 =& \;  [12, [9, -1], 3, [5, 2], 8, [6, 0], 14, [2, -1]] \\
C_3 =& \;  [7, [3, 0], 4, [3, 2]]
\end{align*}
\end{minipage}

\newpage

\begin{minipage}{8cm}
\begin{align*}
C_0' =& \;  [13, [7, 2], 6, [19, -1], 19, [19, 1], 16, [2, 1], 14, [5, 0], 0, [8, 1], 8, [4, -1], 4, [9, 1]] \\
C_1' =& \;  [2, [9, 2], 11, [1, 0], 10, [5, 1], 15, [6, 2], 9, [8, 0], 1, [1, 1]] \\
C_2' =& \;  [3, [4, 0], 7, [2, 2], 5, [7, 1], 17, [5, -1]] \\
C_3' =& \;  [12, [6, -1], 18, [6, 1]]
\end{align*}
\end{minipage}

\item 2-factor type $[7, 7, 4, 2]$: two starters \vspace*{-3mm} \\
\begin{minipage}{8cm}
\begin{align*}
C_0 =& \;  [17, [5, -1], 3, [3, 0], 6, [9, 1], 15, [6, 2], 2, [19, -1], 19, [19, 0], 8, [9, 2]] \\
C_1 =& \;  [16, [5, 1], 11, [4, 1], 7, [7, 0], 0, [9, -1], 10, [3, 2], 13, [5, 0], 18, [2, 2]] \\
C_2 =& \;  [1, [8, 2], 12, [8, 1], 4, [9, 0], 14, [6, 1]] \\
C_3 =& \;  [5, [4, 0], 9, [4, 2]]
\end{align*}
\end{minipage}

\vspace*{-4mm}
\begin{minipage}{8cm}
\begin{align*}
C_0' =& \;  [11, [8, -1], 3, [7, 2], 15, [6, 0], 9, [2, -1], 7, [6, -1], 1, [3, 1], 4, [7, 1]] \\
C_1' =& \;  [12, [7, -1], 5, [1, 1], 6, [4, -1], 2, [3, -1], 18, [5, 2], 13, [1, -1], 14, [2, 0]] \\
C_2' =& \;  [19, [19, 1], 0, [8, 0], 8, [2, 1], 10, [19, 2]] \\
C_3' =& \;  [17, [1, 0], 16, [1, 2]]
\end{align*}
\end{minipage}

\item 2-factor type $[8, 5, 5, 2]$: two starters \vspace*{-3mm} \\
\begin{minipage}{8cm}
\begin{align*}
C_0 =& \;  [17, [3, -1], 14, [5, 2], 9, [7, 1], 2, [6, -1], 8, [4, 0], 4, [3, 1], 7, [8, 1], 15, [2, 1]] \\
C_1 =& \;  [13, [9, -1], 3, [6, 0], 16, [8, 2], 5, [19, -1], 19, [19, 1]] \\
C_2 =& \;  [10, [9, 1], 0, [1, 1], 1, [8, -1], 12, [1, 2], 11, [1, 0]] \\
C_3 =& \;  [18, [7, 0], 6, [7, 2]]
\end{align*}
\end{minipage}

\vspace*{-4mm}
\begin{minipage}{8cm}
\begin{align*}
C_0' =& \;  [0, [5, -1], 5, [7, -1], 12, [5, 0], 17, [9, 2], 7, [8, 0], 15, [4, 2], 11, [5, 1], 6, [6, 1]] \\
C_1' =& \;  [8, [9, 0], 18, [2, 2], 1, [3, 0], 4, [2, -1], 2, [6, 2]] \\
C_2' =& \;  [16, [2, 0], 14, [4, -1], 10, [1, -1], 9, [4, 1], 13, [3, 2]] \\
C_3' =& \;  [19, [19, 0], 3, [19, 2]]
\end{align*}
\end{minipage}

\item 2-factor type $[7, 6, 5, 2]$: two starters \vspace*{-3mm} \\
\begin{minipage}{8cm}
\begin{align*}
C_0 =& \;  [19, [19, -1], 18, [7, 1], 6, [5, -1], 1, [3, 2], 4, [5, 0], 9, [3, 1], 12, [19, 1]] \\
C_1 =& \;  [10, [6, -1], 16, [6, 0], 3, [3, -1], 0, [2, 2], 17, [4, 0], 2, [8, 2]] \\
C_2 =& \;  [5, [2, 1], 7, [1, 2], 8, [3, 0], 11, [4, 2], 15, [9, 0]] \\
C_3 =& \;  [14, [1, -1], 13, [1, 1]]
\end{align*}
\end{minipage}

\vspace*{-4mm}
\begin{minipage}{8cm}
\begin{align*}
C_0' =& \;  [6, [6, 1], 12, [5, 2], 7, [19, 0], 19, [19, 2], 15, [4, -1], 0, [8, -1], 8, [2, 0]] \\
C_1' =& \;  [2, [9, 1], 11, [2, -1], 9, [5, 1], 14, [8, 1], 3, [9, 2], 13, [8, 0]] \\
C_2' =& \;  [5, [6, 2], 18, [1, 0], 17, [7, -1], 10, [9, -1], 1, [4, 1]] \\
C_3' =& \;  [4, [7, 0], 16, [7, 2]]
\end{align*}
\end{minipage}

\item 2-factor type $[6, 6, 6, 2]$: two starters \vspace*{-3mm} \\
\begin{minipage}{8cm}
\begin{align*}
C_0 =& \;  [8, [4, 2], 4, [7, 0], 16, [4, -1], 12, [1, 2], 13, [9, -1], 3, [5, 0]] \\
C_1 =& \;  [19, [19, -1], 5, [7, 2], 17, [7, 1], 10, [8, 0], 18, [7, -1], 11, [19, 1]] \\
C_2 =& \;  [15, [4, 1], 0, [5, -1], 14, [8, 2], 6, [1, 0], 7, [5, 2], 2, [6, 0]] \\
C_3 =& \;  [1, [8, -1], 9, [8, 1]]
\end{align*}
\end{minipage}

\newpage

\begin{minipage}{8cm}
\begin{align*}
C_0' =& \;  [8, [6, 1], 2, [2, 1], 0, [6, -1], 6, [1, -1], 7, [2, -1], 9, [1, 1]] \\
C_1' =& \;  [1, [9, 2], 11, [3, -1], 14, [9, 0], 5, [5, 1], 10, [6, 2], 4, [3, 0]] \\
C_2' =& \;  [3, [4, 0], 18, [2, 2], 16, [3, 1], 13, [2, 0], 15, [3, 2], 12, [9, 1]] \\
C_3' =& \;  [19, [19, 0], 17, [19, 2]]
\end{align*}
\end{minipage}

\item 2-factor type $[11, 3, 3, 3]$: one starter \vspace*{-3mm} \\
\begin{minipage}{8cm}
\begin{align*}
C_0 =& \;  [19, [19, 1], 7, [9, 1], 16, [8, 1], 8, [7, 0], 15, [6, 1], 9, [3, 0], 12, [1, 0], 13, [7, 1], 1, [2, 1], 18, [1, 1], 0, [19, 0]] \\
C_1 =& \;  [11, [8, 0], 3, [3, 1], 6, [5, 0]] \\
C_2 =& \;  [4, [6, 0], 17, [4, 1], 2, [2, 0]] \\
C_3 =& \;  [5, [9, 0], 14, [4, 0], 10, [5, 1]]
\end{align*}
\end{minipage}

\item 2-factor type $[10, 4, 3, 3]$: one starter \vspace*{-3mm} \\
\begin{minipage}{8cm}
\begin{align*}
C_0 =& \;  [7, [2, 1], 5, [4, 1], 9, [6, 1], 15, [1, 0], 14, [4, 0], 10, [2, 0], 12, [8, 0], 4, [1, 1], 3, [5, 1], 17, [9, 1]] \\
C_1 =& \;  [19, [19, 1], 0, [6, 0], 6, [7, 1], 18, [19, 0]] \\
C_2 =& \;  [8, [7, 0], 1, [9, 0], 11, [3, 1]] \\
C_3 =& \;  [2, [5, 0], 16, [3, 0], 13, [8, 1]]
\end{align*}
\end{minipage}

\item 2-factor type $[9, 5, 3, 3]$: one starter \vspace*{-3mm} \\
\begin{minipage}{8cm}
\begin{align*}
C_0 =& \;  [2, [6, 1], 8, [7, 0], 1, [6, 0], 14, [1, 1], 15, [3, 1], 18, [2, 0], 16, [7, 1], 4, [2, 1], 6, [4, 0]] \\
C_1 =& \;  [5, [5, 1], 10, [9, 0], 0, [3, 0], 3, [9, 1], 13, [8, 1]] \\
C_2 =& \;  [11, [4, 1], 7, [5, 0], 12, [1, 0]] \\
C_3 =& \;  [19, [19, 0], 17, [8, 0], 9, [19, 1]]
\end{align*}
\end{minipage}

\item 2-factor type $[8, 6, 3, 3]$: one starter \vspace*{-3mm} \\
\begin{minipage}{8cm}
\begin{align*}
C_0 =& \;  [2, [9, 0], 11, [3, 0], 8, [8, 0], 0, [2, 1], 17, [7, 1], 5, [9, 1], 14, [6, 0], 1, [1, 1]] \\
C_1 =& \;  [6, [4, 0], 10, [3, 1], 13, [2, 0], 15, [8, 1], 7, [5, 1], 12, [6, 1]] \\
C_2 =& \;  [9, [7, 0], 16, [19, 1], 19, [19, 0]] \\
C_3 =& \;  [4, [5, 0], 18, [4, 1], 3, [1, 0]]
\end{align*}
\end{minipage}

\item 2-factor type $[7, 7, 3, 3]$: one starter \vspace*{-3mm} \\
\begin{minipage}{8cm}
\begin{align*}
C_0 =& \;  [17, [1, 0], 16, [3, 0], 0, [9, 1], 9, [4, 1], 13, [6, 0], 7, [4, 0], 11, [6, 1]] \\
C_1 =& \;  [5, [2, 0], 3, [8, 1], 14, [7, 0], 2, [9, 0], 12, [19, 0], 19, [19, 1], 4, [1, 1]] \\
C_2 =& \;  [15, [3, 1], 18, [8, 0], 10, [5, 0]] \\
C_3 =& \;  [6, [2, 1], 8, [7, 1], 1, [5, 1]]
\end{align*}
\end{minipage}

\item 2-factor type $[9, 4, 4, 3]$: one starter \vspace*{-3mm} \\
\begin{minipage}{8cm}
\begin{align*}
C_0 =& \;  [11, [9, 0], 2, [6, 1], 8, [7, 1], 1, [8, 1], 12, [3, 0], 9, [5, 0], 14, [8, 0], 6, [1, 1], 7, [4, 1]] \\
C_1 =& \;  [0, [9, 1], 10, [19, 1], 19, [19, 0], 15, [4, 0]] \\
C_2 =& \;  [5, [7, 0], 17, [6, 0], 4, [1, 0], 3, [2, 0]] \\
C_3 =& \;  [13, [5, 1], 18, [2, 1], 16, [3, 1]]
\end{align*}
\end{minipage}

\newpage

\item 2-factor type $[8, 5, 4, 3]$: one starter \vspace*{-3mm} \\
\begin{minipage}{8cm}
\begin{align*}
C_0 =& \;  [15, [9, 1], 5, [8, 1], 16, [9, 0], 7, [1, 0], 6, [6, 1], 0, [5, 0], 14, [6, 0], 8, [7, 1]] \\
C_1 =& \;  [17, [4, 0], 13, [19, 0], 19, [19, 1], 2, [8, 0], 10, [7, 0]] \\
C_2 =& \;  [1, [2, 1], 3, [4, 1], 18, [5, 1], 4, [3, 1]] \\
C_3 =& \;  [12, [1, 1], 11, [2, 0], 9, [3, 0]]
\end{align*}
\end{minipage}

\item 2-factor type $[7, 6, 4, 3]$: one starter \vspace*{-3mm} \\
\begin{minipage}{8cm}
\begin{align*}
C_0 =& \;  [0, [5, 0], 14, [3, 1], 11, [5, 1], 6, [6, 1], 12, [3, 0], 9, [4, 0], 13, [6, 0]] \\
C_1 =& \;  [18, [9, 1], 8, [7, 1], 1, [4, 1], 16, [1, 0], 17, [7, 0], 10, [8, 1]] \\
C_2 =& \;  [4, [1, 1], 3, [19, 1], 19, [19, 0], 2, [2, 0]] \\
C_3 =& \;  [15, [9, 0], 5, [2, 1], 7, [8, 0]]
\end{align*}
\end{minipage}

\item 2-factor type $[7, 5, 5, 3]$: one starter \vspace*{-3mm} \\
\begin{minipage}{8cm}
\begin{align*}
C_0 =& \;  [15, [5, 0], 10, [19, 1], 19, [19, 0], 4, [5, 1], 9, [7, 1], 16, [3, 1], 13, [2, 1]] \\
C_1 =& \;  [18, [6, 0], 5, [4, 1], 1, [1, 0], 0, [2, 0], 2, [3, 0]] \\
C_2 =& \;  [12, [4, 0], 8, [9, 0], 17, [6, 1], 11, [8, 1], 3, [9, 1]] \\
C_3 =& \;  [7, [7, 0], 14, [8, 0], 6, [1, 1]]
\end{align*}
\end{minipage}

\item 2-factor type $[6, 6, 5, 3]$: one starter \vspace*{-3mm} \\
\begin{minipage}{8cm}
\begin{align*}
C_0 =& \;  [14, [8, 0], 3, [6, 1], 16, [3, 0], 13, [19, 0], 19, [19, 1], 10, [4, 0]] \\
C_1 =& \;  [4, [7, 0], 11, [2, 0], 9, [2, 1], 7, [9, 1], 17, [5, 1], 12, [8, 1]] \\
C_2 =& \;  [18, [7, 1], 6, [1, 1], 5, [3, 1], 2, [6, 0], 8, [9, 0]] \\
C_3 =& \;  [0, [4, 1], 15, [5, 0], 1, [1, 0]]
\end{align*}
\end{minipage}

\item 2-factor type $[7, 5, 4, 4]$: one starter \vspace*{-3mm} \\
\begin{minipage}{8cm}
\begin{align*}
C_0 =& \;  [19, [19, 1], 6, [5, 0], 1, [4, 0], 5, [3, 1], 8, [4, 1], 12, [2, 0], 10, [19, 0]] \\
C_1 =& \;  [16, [3, 0], 13, [9, 0], 3, [1, 1], 2, [6, 0], 15, [1, 0]] \\
C_2 =& \;  [17, [6, 1], 4, [7, 0], 11, [8, 0], 0, [2, 1]] \\
C_3 =& \;  [14, [7, 1], 7, [8, 1], 18, [9, 1], 9, [5, 1]]
\end{align*}
\end{minipage}

\item 2-factor type $[6, 6, 4, 4]$: one starter \vspace*{-3mm} \\
\begin{minipage}{8cm}
\begin{align*}
C_0 =& \;  [10, [19, 0], 19, [19, 1], 7, [8, 0], 18, [9, 0], 8, [1, 1], 9, [1, 0]] \\
C_1 =& \;  [11, [8, 1], 0, [4, 1], 15, [2, 1], 13, [9, 1], 4, [7, 0], 16, [5, 0]] \\
C_2 =& \;  [1, [2, 0], 3, [5, 1], 17, [3, 1], 14, [6, 0]] \\
C_3 =& \;  [12, [7, 1], 5, [3, 0], 2, [4, 0], 6, [6, 1]]
\end{align*}
\end{minipage}

\item 2-factor type $[6, 5, 5, 4]$: one starter \vspace*{-3mm} \\
\begin{minipage}{8cm}
\begin{align*}
C_0 =& \;  [19, [19, 0], 18, [1, 0], 0, [1, 1], 1, [4, 1], 5, [8, 1], 13, [19, 1]] \\
C_1 =& \;  [14, [4, 0], 10, [2, 0], 12, [3, 1], 15, [9, 0], 6, [8, 0]] \\
C_2 =& \;  [3, [5, 0], 8, [6, 1], 2, [9, 1], 11, [6, 0], 17, [5, 1]] \\
C_3 =& \;  [9, [2, 1], 7, [3, 0], 4, [7, 1], 16, [7, 0]]
\end{align*}
\end{minipage}

\newpage

\item 2-factor type $[16, 2, 2]$: one starter \vspace*{-3mm} \\
\begin{minipage}{8cm}
\begin{align*}
C_0 =& \;  [3, [1, 1], 2, [5, 0], 16, [2, 1], 14, [9, 0], 5, [2, 0], 7, [6, 0], 13, [7, 1], 1, [8, 1], 12, [6, 1], 18, [9, 1], 9, [3, 1], 6, \\ & \; [8, 0], 17, [7, 0], 10, [1, 0], 11, [3, 0], 8, [5, 1]] \\
C_1 =& \;  [19, [19, 0], 4, [19, 1]] \\
C_2 =& \;  [15, [4, 0], 0, [4, 1]]
\end{align*}
\end{minipage}

\item 2-factor type $[15, 3, 2]$: two starters \vspace*{-3mm} \\
\begin{minipage}{8cm}
\begin{align*}
C_0 =& \;  [15, [2, 2], 17, [6, 0], 11, [2, -1], 9, [9, 1], 18, [5, 1], 4, [1, -1], 3, [5, 2], 8, [8, 1], 0, [5, -1], 5, [5, 0], 10, \\ & \; [6, 1], 16, [3, -1], 13, [7, -1], 6, [1, 2], 7, [8, 0]] \\
C_1 =& \;  [12, [2, 0], 14, [7, 1], 2, [9, 2]] \\
C_2 =& \;  [19, [19, -1], 1, [19, 1]]
\end{align*}
\end{minipage}

\vspace*{-4mm}
\begin{minipage}{8cm}
\begin{align*}
C_0' =& \;  [15, [3, 2], 18, [8, -1], 7, [2, 1], 5, [4, 0], 1, [19, 2], 19, [19, 0], 12, [9, -1], 3, [7, 2], 10, [1, 1], 9, [7, 0], 16, \\ & \; [3, 1], 0, [8, 2], 11, [3, 0], 8, [6, 2], 14, [1, 0]] \\
C_1' =& \;  [13, [9, 0], 4, [6, -1], 17, [4, 2]] \\
C_2' =& \;  [6, [4, -1], 2, [4, 1]]
\end{align*}
\end{minipage}

\item 2-factor type $[14, 4, 2]$: two starters \vspace*{-3mm} \\
\begin{minipage}{8cm}
\begin{align*}
C_0 =& \;  [7, [9, 1], 17, [6, 0], 4, [2, -1], 2, [1, 1], 3, [3, -1], 0, [3, 2], 16, [19, -1], 19, [19, 1], 5, [3, 0], 8, [7, 1], 15, \\ & \; [1, 2], 14, [4, 0], 10, [8, 1], 18, [8, 2]] \\
C_1 =& \;  [1, [8, 0], 9, [4, 2], 13, [2, 1], 11, [9, -1]] \\
C_2 =& \;  [6, [6, -1], 12, [6, 1]]
\end{align*}
\end{minipage}

\vspace*{-4mm}
\begin{minipage}{8cm}
\begin{align*}
C_0' =& \;  [1, [8, -1], 12, [4, -1], 8, [9, 2], 18, [7, -1], 6, [5, 0], 11, [5, 2], 16, [9, 0], 7, [2, 2], 9, [4, 1], 5, [3, 1], 2, \\ & \; [7, 0], 14, [1, -1], 13, [6, 2], 0, [1, 0]] \\
C_1' =& \;  [3, [5, 1], 17, [2, 0], 15, [5, -1], 10, [7, 2]] \\
C_2' =& \;  [19, [19, 0], 4, [19, 2]]
\end{align*}
\end{minipage}

\item 2-factor type $[13, 5, 2]$: two starters \vspace*{-3mm} \\
\begin{minipage}{8cm}
\begin{align*}
C_0 =& \;  [13, [7, -1], 1, [2, -1], 3, [3, 0], 0, [6, 1], 6, [6, 2], 12, [1, 0], 11, [1, 2], 10, [5, 0], 15, [2, 1], 17, [7, 2], 5, \\ & \; [3, 1], 2, [2, 0], 4, [9, 2]] \\
C_1 =& \;  [19, [19, 0], 7, [8, -1], 18, [4, -1], 14, [5, -1], 9, [19, 2]] \\
C_2 =& \;  [8, [8, 0], 16, [8, 2]]
\end{align*}
\end{minipage}

\vspace*{-4mm}
\begin{minipage}{8cm}
\begin{align*}
C_0' =& \;  [18, [7, 0], 11, [2, 2], 13, [8, 1], 5, [9, -1], 14, [4, 1], 10, [7, 1], 3, [9, 0], 12, [3, 2], 9, [9, 1], 0, [6, 0], 6, \\ & \; [1, 1], 7, [3, -1], 4, [5, 2]] \\
C_1' =& \;  [15, [5, 1], 1, [4, 2], 16, [1, -1], 17, [4, 0], 2, [6, -1]] \\
C_2' =& \;  [19, [19, -1], 8, [19, 1]]
\end{align*}
\end{minipage}

\item 2-factor type $[12, 6, 2]$: two starters \vspace*{-3mm} \\
\begin{minipage}{8cm}
\begin{align*}
C_0 =& \;  [18, [9, 1], 8, [3, -1], 5, [4, 2], 1, [9, 0], 10, [2, 2], 12, [9, -1], 2, [2, 0], 0, [9, 2], 9, [3, 0], 6, [2, -1], 4, \\ & \; [7, -1], 16, [2, 1]] \\
C_1 =& \;  [19, [19, 1], 3, [5, 1], 17, [4, 1], 13, [6, 1], 7, [4, -1], 11, [19, -1]] \\
C_2 =& \;  [15, [1, 0], 14, [1, 2]]
\end{align*}
\end{minipage}

\newpage

\begin{minipage}{8cm}
\begin{align*}
C_0' =& \;  [18, [7, 1], 6, [3, 2], 3, [6, 0], 16, [6, -1], 10, [3, 1], 13, [19, 2], 19, [19, 0], 11, [8, 1], 0, [5, -1], 14, [6, 2], \\ & \; 8, [1, 1], 7, [8, 0]] \\
C_1' =& \;  [15, [8, -1], 4, [5, 0], 9, [8, 2], 17, [4, 0], 2, [1, -1], 1, [5, 2]] \\
C_2' =& \;  [12, [7, 0], 5, [7, 2]]
\end{align*}
\end{minipage}

\item 2-factor type $[11, 7, 2]$: two starters \vspace*{-3mm} \\
\begin{minipage}{8cm}
\begin{align*}
C_0 =& \;  [0, [1, 0], 18, [2, 2], 16, [4, 0], 1, [5, 2], 6, [3, -1], 9, [4, -1], 13, [9, -1], 3, [2, -1], 5, [7, -1], 12, [8, 1], \\ & \; 4, [4, 1]] \\
C_1 =& \;  [8, [19, -1], 19, [19, 0], 7, [7, 2], 14, [1, -1], 15, [5, 0], 10, [8, -1], 2, [6, 2]] \\
C_2 =& \;  [11, [6, -1], 17, [6, 1]]
\end{align*}
\end{minipage}

\vspace*{-4mm}
\begin{minipage}{8cm}
\begin{align*}
C_0' =& \;  [2, [8, 2], 13, [8, 0], 5, [5, 1], 10, [1, 1], 11, [19, 2], 19, [19, 1], 14, [6, 0], 8, [1, 2], 9, [2, 1], 7, [7, 1], 0, [2, 0]] \\
C_1' =& \;  [3, [4, 2], 18, [5, -1], 4, [7, 0], 16, [9, 1], 6, [9, 2], 15, [3, 1], 12, [9, 0]] \\
C_2' =& \;  [1, [3, 0], 17, [3, 2]]
\end{align*}
\end{minipage}

\item 2-factor type $[10, 8, 2]$: two starters \vspace*{-3mm} \\
\begin{minipage}{8cm}
\begin{align*}
C_0 =& \;  [9, [6, 1], 15, [5, 1], 1, [9, -1], 11, [7, 1], 18, [19, 1], 19, [19, -1], 2, [8, 2], 10, [3, 0], 7, [1, 2], 8, [1, 0]] \\
C_1 =& \;  [5, [9, 2], 14, [8, 0], 3, [3, -1], 0, [7, 2], 12, [6, 0], 6, [2, -1], 4, [7, -1], 16, [8, 1]] \\
C_2 =& \;  [17, [4, 0], 13, [4, 2]]
\end{align*}
\end{minipage}

\vspace*{-4mm}
\begin{minipage}{8cm}
\begin{align*}
C_0' =& \;  [12, [2, 1], 14, [3, 2], 11, [1, 1], 10, [5, -1], 15, [1, -1], 16, [8, -1], 5, [3, 1], 2, [4, -1], 6, [7, 0], 18, [6, -1]] \\
C_1' =& \;  [17, [6, 2], 4, [9, 0], 13, [5, 2], 8, [5, 0], 3, [4, 1], 7, [2, 2], 9, [9, 1], 0, [2, 0]] \\
C_2' =& \;  [1, [19, 0], 19, [19, 2]]
\end{align*}
\end{minipage}

\item 2-factor type $[9, 9, 2]$: two starters \vspace*{-3mm} \\
\begin{minipage}{8cm}
\begin{align*}
C_0 =& \;  [7, [4, 2], 11, [4, 0], 15, [7, 2], 8, [8, 0], 0, [2, 2], 17, [6, 0], 4, [6, 1], 10, [4, 1], 14, [7, -1]] \\
C_1 =& \;  [16, [19, 2], 19, [19, 0], 5, [1, 2], 6, [5, 1], 1, [1, 0], 2, [3, -1], 18, [6, -1], 12, [3, 2], 9, [7, 0]] \\
C_2 =& \;  [13, [9, 0], 3, [9, 2]]
\end{align*}
\end{minipage}

\vspace*{-4mm}
\begin{minipage}{8cm}
\begin{align*}
C_0' =& \;  [18, [4, -1], 3, [5, 0], 17, [8, -1], 9, [3, 1], 12, [8, 2], 4, [7, 1], 16, [5, -1], 2, [2, 1], 0, [1, 1]] \\
C_1' =& \;  [1, [6, 2], 7, [1, -1], 8, [2, 0], 6, [8, 1], 14, [9, -1], 5, [5, 2], 10, [3, 0], 13, [2, -1], 11, [9, 1]] \\
C_2' =& \;  [19, [19, -1], 15, [19, 1]]
\end{align*}
\end{minipage}

\item 2-factor type $[14, 3, 3]$: one starter \vspace*{-3mm} \\
\begin{minipage}{8cm}
\begin{align*}
C_0 =& \;  [12, [6, 1], 18, [5, 0], 13, [7, 1], 1, [7, 0], 8, [5, 1], 3, [19, 0], 19, [19, 1], 5, [3, 0], 2, [2, 1], 4, [6, 0], 10, [4, 0], \\ & \; 14, [3, 1], 17, [1, 1], 16, [4, 1]] \\
C_1 =& \;  [7, [1, 0], 6, [9, 0], 15, [8, 1]] \\
C_2 =& \;  [9, [2, 0], 11, [8, 0], 0, [9, 1]]
\end{align*}
\end{minipage}

\item 2-factor type $[13, 4, 3]$: one starter \vspace*{-3mm} \\
\begin{minipage}{8cm}
\begin{align*}
C_0 =& \;  [13, [1, 0], 12, [6, 1], 18, [8, 0], 10, [1, 1], 11, [4, 0], 15, [7, 1], 3, [2, 1], 5, [3, 0], 2, [7, 0], 9, [19, 1], 19, \\ & \; [19, 0], 4, [4, 1], 8, [5, 1]] \\
C_1 =& \;  [7, [9, 0], 16, [3, 1], 0, [2, 0], 17, [9, 1]] \\
C_2 =& \;  [1, [6, 0], 14, [8, 1], 6, [5, 0]]
\end{align*}
\end{minipage}

\item 2-factor type $[12, 5, 3]$: one starter \vspace*{-3mm} \\
\begin{minipage}{8cm}
\begin{align*}
C_0 =& \;  [5, [8, 1], 13, [1, 0], 12, [4, 1], 16, [6, 1], 3, [7, 0], 15, [9, 1], 6, [4, 0], 10, [3, 1], 7, [7, 1], 14, [5, 0], 9, [2, 0], \\ & \; 11, [6, 0]] \\
C_1 =& \;  [19, [19, 0], 1, [1, 1], 2, [3, 0], 18, [5, 1], 4, [19, 1]] \\
C_2 =& \;  [17, [2, 1], 0, [8, 0], 8, [9, 0]]
\end{align*}
\end{minipage}

\item 2-factor type $[11, 6, 3]$: one starter \vspace*{-3mm} \\
\begin{minipage}{8cm}
\begin{align*}
C_0 =& \;  [3, [6, 1], 16, [8, 1], 8, [5, 1], 13, [3, 0], 10, [4, 1], 14, [4, 0], 18, [2, 0], 1, [1, 0], 0, [5, 0], 5, [7, 1], 12, [9, 0]] \\
C_1 =& \;  [2, [9, 1], 11, [7, 0], 4, [3, 1], 7, [1, 1], 6, [19, 0], 19, [19, 1]] \\
C_2 =& \;  [9, [8, 0], 17, [2, 1], 15, [6, 0]]
\end{align*}
\end{minipage}

\item 2-factor type $[10, 7, 3]$: one starter \vspace*{-3mm} \\
\begin{minipage}{8cm}
\begin{align*}
C_0 =& \;  [8, [7, 0], 15, [9, 1], 6, [1, 1], 7, [6, 1], 13, [4, 1], 9, [3, 0], 12, [19, 1], 19, [19, 0], 2, [5, 0], 16, [8, 1]] \\
C_1 =& \;  [18, [7, 1], 11, [8, 0], 3, [2, 0], 1, [3, 1], 4, [6, 0], 17, [2, 1], 0, [1, 0]] \\
C_2 =& \;  [14, [4, 0], 10, [5, 1], 5, [9, 0]]
\end{align*}
\end{minipage}

\item 2-factor type $[9, 8, 3]$: one starter \vspace*{-3mm} \\
\begin{minipage}{8cm}
\begin{align*}
C_0 =& \;  [14, [4, 1], 10, [6, 1], 4, [1, 0], 5, [9, 0], 15, [2, 0], 17, [4, 0], 13, [7, 0], 6, [3, 0], 3, [8, 1]] \\
C_1 =& \;  [19, [19, 0], 16, [5, 0], 11, [3, 1], 8, [7, 1], 1, [1, 1], 2, [2, 1], 0, [9, 1], 9, [19, 1]] \\
C_2 =& \;  [7, [5, 1], 12, [6, 0], 18, [8, 0]]
\end{align*}
\end{minipage}

\item 2-factor type $[11, 5, 4]$: one starter \vspace*{-3mm} \\
\begin{minipage}{8cm}
\begin{align*}
C_0 =& \;  [15, [3, 1], 12, [9, 0], 3, [1, 1], 4, [4, 0], 0, [7, 1], 7, [5, 1], 2, [4, 1], 17, [7, 0], 10, [9, 1], 1, [19, 0], 19, [19, 1]] \\
C_1 =& \;  [13, [2, 1], 11, [5, 0], 16, [2, 0], 18, [6, 0], 5, [8, 0]] \\
C_2 =& \;  [6, [3, 0], 9, [1, 0], 8, [6, 1], 14, [8, 1]]
\end{align*}
\end{minipage}

\item 2-factor type $[10, 6, 4]$: one starter \vspace*{-3mm} \\
\begin{minipage}{8cm}
\begin{align*}
C_0 =& \;  [5, [1, 0], 4, [2, 0], 6, [8, 1], 17, [2, 1], 0, [3, 1], 3, [4, 1], 18, [8, 0], 10, [5, 0], 15, [6, 0], 9, [4, 0]] \\
C_1 =& \;  [19, [19, 1], 14, [7, 1], 7, [5, 1], 12, [1, 1], 13, [3, 0], 16, [19, 0]] \\
C_2 =& \;  [1, [9, 1], 11, [9, 0], 2, [6, 1], 8, [7, 0]]
\end{align*}
\end{minipage}

\item 2-factor type $[9, 7, 4]$: one starter \vspace*{-3mm} \\
\begin{minipage}{8cm}
\begin{align*}
C_0 =& \;  [16, [5, 0], 11, [2, 0], 9, [6, 0], 15, [19, 0], 19, [19, 1], 10, [9, 0], 1, [4, 1], 5, [7, 0], 17, [1, 1]] \\
C_1 =& \;  [4, [3, 1], 7, [1, 0], 8, [6, 1], 2, [2, 1], 0, [3, 0], 3, [9, 1], 12, [8, 1]] \\
C_2 =& \;  [18, [5, 1], 13, [7, 1], 6, [8, 0], 14, [4, 0]]
\end{align*}
\end{minipage}

\item 2-factor type $[10, 5, 5]$: one starter \vspace*{-3mm} \\
\begin{minipage}{8cm}
\begin{align*}
C_0 =& \;  [13, [3, 0], 16, [2, 0], 14, [6, 1], 8, [6, 0], 2, [4, 1], 6, [5, 0], 1, [2, 1], 18, [5, 1], 4, [1, 0], 5, [8, 0]] \\
C_1 =& \;  [7, [4, 0], 11, [1, 1], 12, [3, 1], 9, [8, 1], 17, [9, 0]] \\
C_2 =& \;  [19, [19, 1], 15, [7, 1], 3, [7, 0], 10, [9, 1], 0, [19, 0]]
\end{align*}
\end{minipage}

\item 2-factor type $[9, 6, 5]$: one starter \vspace*{-3mm} \\
\begin{minipage}{8cm}
\begin{align*}
C_0 =& \;  [10, [7, 0], 17, [4, 0], 13, [8, 0], 5, [2, 1], 3, [4, 1], 18, [5, 0], 4, [9, 0], 14, [3, 1], 11, [1, 0]] \\
C_1 =& \;  [15, [5, 1], 1, [6, 1], 7, [2, 0], 9, [3, 0], 6, [9, 1], 16, [1, 1]] \\
C_2 =& \;  [12, [19, 0], 19, [19, 1], 2, [6, 0], 8, [8, 1], 0, [7, 1]]
\end{align*}
\end{minipage}

\item 2-factor type $[8, 7, 5]$: one starter \vspace*{-3mm} \\
\begin{minipage}{8cm}
\begin{align*}
C_0 =& \;  [16, [5, 1], 11, [1, 0], 10, [5, 0], 15, [6, 1], 9, [1, 1], 8, [4, 1], 12, [9, 0], 3, [6, 0]] \\
C_1 =& \;  [18, [19, 0], 19, [19, 1], 6, [8, 1], 14, [3, 0], 17, [2, 0], 0, [7, 1], 7, [8, 0]] \\
C_2 =& \;  [2, [3, 1], 5, [4, 0], 1, [7, 0], 13, [9, 1], 4, [2, 1]]
\end{align*}
\end{minipage}

\item 2-factor type $[8, 6, 6]$: one starter \vspace*{-3mm} \\
\begin{minipage}{8cm}
\begin{align*}
C_0 =& \;  [15, [2, 1], 13, [8, 0], 5, [1, 1], 4, [3, 0], 7, [6, 1], 1, [5, 0], 6, [19, 0], 19, [19, 1]] \\
C_1 =& \;  [8, [6, 0], 14, [3, 1], 17, [1, 0], 18, [4, 0], 3, [7, 1], 10, [2, 0]] \\
C_2 =& \;  [12, [4, 1], 16, [5, 1], 11, [8, 1], 0, [9, 1], 9, [7, 0], 2, [9, 0]]
\end{align*}
\end{minipage}

\item 2-factor type $[7, 7, 6]$: one starter \vspace*{-3mm} \\
\begin{minipage}{8cm}
\begin{align*}
C_0 =& \;  [5, [5, 0], 10, [6, 0], 16, [2, 1], 14, [8, 1], 6, [6, 1], 0, [2, 0], 17, [7, 0]] \\
C_1 =& \;  [18, [5, 1], 4, [3, 1], 1, [8, 0], 9, [3, 0], 12, [1, 0], 13, [9, 0], 3, [4, 1]] \\
C_2 =& \;  [19, [19, 0], 7, [1, 1], 8, [7, 1], 15, [4, 0], 11, [9, 1], 2, [19, 1]]
\end{align*}
\end{minipage}

\item 2-factor type $[18, 2]$: two starters \vspace*{-3mm} \\
\begin{minipage}{8cm}
\begin{align*}
C_0 =& \;  [11, [1, 0], 10, [7, 2], 3, [19, -1], 19, [19, 1], 2, [2, -1], 4, [4, 1], 0, [6, 0], 6, [6, 2], 12, [1, -1], 13, [7, 0], 1, \\ & \; [8, 2], 9, [1, 1], 8, [6, -1], 14, [4, -1], 18, [8, -1], 7, [9, 0], 16, [1, 2], 17, [6, 1]] \\
C_1 =& \;  [15, [9, -1], 5, [9, 1]]
\end{align*}
\end{minipage}

\vspace*{-4mm}
\begin{minipage}{8cm}
\begin{align*}
C_0' =& \;  [13, [5, 0], 8, [7, -1], 1, [3, 1], 4, [9, 2], 14, [5, -1], 0, [3, 0], 3, [2, 2], 5, [8, 0], 16, [2, 1], 18, [8, 1], 7, [4, 2], \\ & \; 11, [5, 1], 6, [3, -1], 9, [7, 1], 2, [4, 0], 17, [5, 2], 12, [2, 0], 10, [3, 2]] \\
C_1' =& \;  [19, [19, 0], 15, [19, 2]]
\end{align*}
\end{minipage}

\item 2-factor type $[17, 3]$: one starter \vspace*{-3mm} \\
\begin{minipage}{8cm}
\begin{align*}
C_0 =& \;  [3, [1, 0], 4, [3, 1], 7, [8, 0], 18, [6, 0], 5, [3, 0], 8, [6, 1], 2, [8, 1], 10, [7, 0], 17, [19, 0], 19, [19, 1], 1, [4, 0], \\ & \; 16, [1, 1], 15, [4, 1], 0, [9, 1], 9, [5, 0], 14, [2, 0], 12, [9, 0]] \\
C_1 =& \;  [6, [5, 1], 11, [2, 1], 13, [7, 1]]
\end{align*}
\end{minipage}

\item 2-factor type $[15, 5]$: one starter \vspace*{-3mm} \\
\begin{minipage}{8cm}
\begin{align*}
C_0 =& \;  [2, [5, 1], 16, [4, 1], 12, [19, 1], 19, [19, 0], 3, [3, 1], 0, [4, 0], 4, [1, 1], 5, [5, 0], 10, [1, 0], 11, [2, 0], 13, [7, 0], \\ & \; 6, [2, 1], 8, [9, 1], 18, [3, 0], 15, [6, 0]] \\
C_1 =& \;  [17, [8, 0], 9, [8, 1], 1, [6, 1], 14, [7, 1], 7, [9, 0]]
\end{align*}
\end{minipage}

\item 2-factor type $[14, 6]$: one starter \vspace*{-3mm} \\
\begin{minipage}{8cm}
\begin{align*}
C_0 =& \;  [9, [7, 0], 2, [8, 0], 13, [2, 1], 15, [8, 1], 7, [4, 0], 3, [19, 0], 19, [19, 1], 10, [9, 0], 0, [1, 1], 1, [7, 1], 8, [3, 1], \\ & \; 11, [5, 1], 16, [9, 1], 6, [3, 0]] \\
C_1 =& \;  [14, [4, 1], 18, [6, 1], 5, [1, 0], 4, [6, 0], 17, [5, 0], 12, [2, 0]]
\end{align*}
\end{minipage}

\item 2-factor type $[13, 7]$: one starter \vspace*{-3mm} \\
\begin{minipage}{8cm}
\begin{align*}
C_0 =& \;  [18, [4, 0], 3, [5, 1], 8, [1, 1], 7, [3, 1], 10, [6, 1], 4, [7, 1], 16, [19, 1], 19, [19, 0], 13, [2, 1], 15, [5, 0], 1, [8, 0], \\ & \; 9, [9, 0], 0, [1, 0]] \\
C_1 =& \;  [14, [2, 0], 12, [7, 0], 5, [6, 0], 11, [9, 1], 2, [4, 1], 6, [8, 1], 17, [3, 0]]
\end{align*}
\end{minipage}

\item 2-factor type $[11, 9]$: one starter \vspace*{-3mm} \\
\begin{minipage}{8cm}
\begin{align*}
C_0 =& \;  [11, [7, 0], 4, [19, 0], 19, [19, 1], 5, [3, 0], 2, [4, 0], 6, [6, 1], 0, [9, 1], 9, [9, 0], 18, [8, 0], 10, [3, 1], 13, [2, 1]] \\
C_1 =& \;  [16, [1, 1], 15, [7, 1], 8, [1, 0], 7, [4, 1], 3, [5, 0], 17, [5, 1], 12, [8, 1], 1, [6, 0], 14, [2, 0]]
\end{align*}
\end{minipage}

\end{itemizenew}

\bigskip

\normalsize

\subsection*{Acknowledgements}

The second author gratefully acknowledges support by the Natural Sciences and Engineering Research Council of Canada (NSERC).

\end{document}